\PassOptionsToPackage{dvipsnames,table}{xcolor}
\documentclass[twoside]{article}
\usepackage[a4paper,top=2cm,bottom=2.2cm,left=2.5cm,right=2.5cm]{geometry}

\usepackage[utf8]{inputenc}
\usepackage[T1]{fontenc}
\usepackage{lmodern}

\newif\ifusebiblatex
\usebiblatextrue
\newif\iffullsurvey
\fullsurveyfalse
\newif\iffullappendix
\fullappendixfalse
\newif\ifsolvealg
\solvealgfalse

\usepackage{amsmath,amssymb}

\usepackage{xcolor}
\usepackage{stmaryrd}
\usepackage{mathtools}
\usepackage{tikz}
\usetikzlibrary{calc,matrix,intersections,arrows.meta}
\usepackage{nicefrac}
\usepackage{enumitem}
\usepackage[algo2e, ruled, vlined, linesnumbered]{algorithm2e}
\usepackage{changepage}
\usepackage[normalem]{ulem}
\usepackage[scr=esstix]{mathalpha}
\usepackage{array}
\usepackage{nicematrix}
\usepackage{multirow}
\usepackage{booktabs}
\usepackage{makecell}

\usepackage{inconsolata}
\usepackage[font=small,labelfont={sf,bf},labelsep=space]{caption}
\usepackage{subcaption}
\usepackage{cancel}
\usepackage{setspace}
\usepackage{float}

\usepackage{etoolbox}
\usepackage{hyperref}

\usepackage[english]{babel}
\usepackage{csquotes}
\ifusebiblatex
\usepackage[maxbibnames=999,style=numeric-comp,backend=biber]{biblatex}
\fi

\usepackage{titlesec}
\titleformat{\section}
  {\large\bfseries}
  {\thesection}
  {1em}
  {}
\titleformat{\subsection}
  {\normalfont\bfseries}
  {\thesubsection}
  {1em}
  {}

\usepackage{hyperref} %
\usepackage[capitalise,nameinlink]{cleveref} %
\crefname{algocf}{algorithm}{algorithms}
\Crefname{algocf}{Algorithm}{Algorithms}
\usepackage{thmtools}

\NiceMatrixOptions{
    code-for-first-row = \color{gray},
    code-for-last-row = \color{gray},
    code-for-first-col = \color{gray},
    code-for-last-col = \color{gray}
}

\ifusebiblatex

\DeclareLanguageMapping{english}{UKenglish}
\AtEveryBibitem{%
  \clearlist{language}%
}
\DeclareSourcemap{
  \maps[datatype=bibtex]{
    \map{
      \step[fieldsource=publisher,
            match=\regexp{arXiv},
            final]
      \step[fieldset=url,null]
      \step[fieldset=urldate,null]
      \step[fieldset=note,fieldvalue={(Preprint)}]
    }
  }
}
\AtEveryBibitem{%
  \ifentrytype{misc}
    {}%
    {%
        \ifentrytype{software}
        {}%
        {%
            \iffieldundef{doi}
              {}%
              {%
                  \clearfield{url}%
                  \clearlist{urldate}%
                  \clearfield{urlyear}%
                  \clearfield{urlmonth}%
                  \clearfield{urlday}%
              }%
        }%
    }%
}
\fi

\makeatletter
\patchcmd\algocf@Vline{\vrule}{\vrule \kern-0.4pt}{}{}
\patchcmd\algocf@Vsline{\vrule}{\vrule \kern-0.4pt}{}{}
\makeatother

\makeatletter
\renewcommand{\algocf@captiontext}[2]{%
  {#1\space}%
  {#2}%
}
\makeatother

\newcommand{\tikzmark}[1]{%
    \tikz[remember picture, overlay] \node (#1) at (0,0) {};%
    \ignorespaces%
}
\def\colorbarxshift{4em}
\def\colorbarxshift{9.8pt}
\def\colorbarwidth{2.4pt}
\newcommand{\ColorBar}[3]{%
    \begin{tikzpicture}[overlay, remember picture]
        \newdimen\startx
        \newdimen\starty
        \newdimen\endy
        \pgfextractx{\startx}{\pgfpointanchor{xmark}{center}}
        \pgfextracty{\starty}{\pgfpointanchor{#2}{center}}
        \pgfextracty{\endy}{\pgfpointanchor{#3}{center}}
        \draw[#1, line width=\colorbarwidth, line cap=round]
            ($(\startx-\colorbarxshift, \starty+0.75em)$) -- ($(\startx-\colorbarxshift, \endy-0.28em)$);
    \end{tikzpicture}%
}
\newcommand{\ColorBarTwo}[4]{%
    \begin{tikzpicture}[overlay,remember picture]
        \newdimen\startx
        \newdimen\starty
        \newdimen\endy
        \pgfextractx{\startx}{\pgfpointanchor{xmark}{center}}%
        \pgfextracty{\starty}{\pgfpointanchor{#3}{center}}%
        \pgfextracty{\endy}{\pgfpointanchor{#4}{center}}%
        \coordinate (Top) at ($(\startx-\colorbarxshift,\starty+0.75em)$);
        \coordinate (Bot) at ($(\startx-\colorbarxshift,\endy-0.28em)$);
        \coordinate (Ltop)  at ($(\startx-\colorbarxshift-\colorbarwidth/2,\starty+0.75em+\colorbarwidth)$);
        \coordinate (Rtop)  at ($(\startx-\colorbarxshift+\colorbarwidth/2,\starty+0.75em+\colorbarwidth)$);
        \coordinate (Lbot)  at ($(\startx-\colorbarxshift-\colorbarwidth/2,\endy-0.28em-\colorbarwidth)$);
        \coordinate (Rbot)  at ($(\startx-\colorbarxshift+\colorbarwidth/2,\endy-0.28em-\colorbarwidth)$);
        \coordinate (Mleft) at ($(Ltop)!0.9!(Lbot)$);
        \coordinate (Mright) at ($(Rtop)!0.1!(Rbot)$);
        \begin{scope}
            \draw[#1, line width=\colorbarwidth, line cap=round] (Top) -- (Bot);
        \end{scope}
        \begin{scope}
            \clip (Mleft) -- (Mright) -- (Rbot) -- (Lbot) -- cycle;
            \draw[#2, line width=\colorbarwidth, line cap=round] (Top) -- (Bot);
        \end{scope}
    \end{tikzpicture}%
}

\renewcommand\L[1]{L^{\hspace{-1pt}#1}}
\newcommand\Lt[1]{L^{\hspace{-1pt}#1\,\top}}
\newcommand\Lit[1]{L^{\hspace{-1pt}#1\;-\!\top}}
\newcommand\Li[1]{L^{\hspace{-1pt}#1\;-\!1}}

\newcommand\tilL[1]{\tilde L^{\hspace{-1pt}#1}}
\newcommand\tilLt[1]{\tilde L^{\hspace{-1pt}#1\,\top}}

\newcommand\Acl{A^{\mathrm{cl}}}%
\newcommand\Acltp{A^{\mathrm{cl}\,\top}}%
\newcommand\tilAcl{\tilde A^{\mathrm{cl}}}%
\newcommand\Aprop{{\hat A}}%

\makeatletter
\newcommand{\verbatimfont}[1]{\renewcommand{\verbatim@font}{\ttfamily#1}}
\makeatother

\newcommand\smalltriangleright{\triangleright}

\def\dx{\Delta\hspace{-0.2pt}x}
\def\du{\Delta\hspace{-0.8pt}u}
\def\dxbf{d_{\bf x}}
\def\dubf{d_{\bf u}}

\def\shortminus{\scalebox{0.6}[0.8]{$-$}}

\SetAlFnt{\setstretch{1.12}}
\SetKwProg{Fn}{Function}{}{}
\SetKwFor{While}{while}{}{}%
\SetKwFor{GotoLoop}{}{}{}%
\SetKwFor{For}{for}{}{}%
\SetKwIF{If}{ElseIf}{Else}{if}{}{else if}{else}{end if}%
\newcommand\mykwcomment{\quad$\smalltriangleright$\ }
\SetKwComment{Comment}{\mykwcomment}{}
\SetKwComment{NakedComment}{}{}
\SetKwComment{FullLineComment}{$\smalltriangleright$\ }{}
\newcommand\mycommentstyle[1]{\textnormal{#1}}
\SetCommentSty{mycommentstyle}
\newcommand\mycomment[1]{\hfill\Comment{#1}}
\newcommand\mycommentnoline[1]{\hfill\mycommentstyle{\hfill\mykwcomment #1}}
\newcommand\mycommentnofill[1]{\quad\mycommentstyle{\mykwcomment #1}}

\newcommand\rmrk[1]{\textcolor{gray}{\small(#1)}}
\newcommand\rmrknp[1]{\textcolor{gray}{\small #1}}

\colorlet{algblue}{blue!90!black!90!red!85!green}
\colorlet{alggreen}{green!60!black}
\colorlet{algpurple}{purple}
\colorlet{algorange}{orange}
\colorlet{algred}{red!70!black}
\colorlet{algteal}{teal!80!blue}
\colorlet{algmagenta}{magenta!70!black}
\colorlet{algolive}{olive!70!black}
\colorlet{algpink}{pink!70!red!90!black}

\colorlet{colorA}{red!60!white}
\colorlet{colorB}{blue!60!white}
\colorlet{colorL}{red!40!orange!90!blue!50!white}
\colorlet{colorU}{green!80!white!85!black!60!white}
\colorlet{colorY}{yellow!40!orange!60!white}
\colorlet{colorQ}{red!40!blue!70!white}

\colorlet{colorAd}{red!85!black!60!white}
\colorlet{colorBd}{blue!70!white}
\colorlet{colorLd}{red!40!orange!85!blue!50!white}
\colorlet{colorUd}{green!90!white!80!black!70!white}
\colorlet{colorYd}{yellow!40!orange!90!white}
\colorlet{colorQd}{red!40!blue!80!white}

\SetKwInOut{KwOpCnt}{Operation count}

\newcommand\Mdiag{M}
\newcommand\Ksub{K}

\newcommand{\Mfwd}[1][]{%
    \rlap{\smash{\raisebox{1.25ex}{$\hspace{0.45em}\scriptveryshortrightarrow$}}}M%
    \if\relax\detokenize{#1}\relax\else^{(\mkern-1.5mu#1\mkern-1.5mu)}\fi
}
\newcommand{\Mbwd}[1][]{%
    \rlap{\smash{\raisebox{1.25ex}{$\hspace{0.4em}\scriptveryshortleftarrow$}}}M%
    \if\relax\detokenize{#1}\relax\else^{(\mkern-1.5mu#1\mkern-1.5mu)}\fi
}
\newcommand{\Kfwd}[1][]{%
    \rlap{\smash{\raisebox{1.35ex}{$\hspace{0.4em}\scriptveryshortrightarrow$}}}K%
    \if\relax\detokenize{#1}\relax\else^{(\mkern-1.5mu#1\mkern-1.5mu)}\fi
}
\newcommand{\Kbwd}[1][]{%
    \rlap{\smash{\raisebox{1.35ex}{$\hspace{0.34em}\scriptveryshortleftarrow$}}}K%
    \if\relax\detokenize{#1}\relax\else^{(\mkern-1.5mu#1\mkern-1.5mu)}\fi
}
\newcommand{\KfwdT}[1][]{%
    \rlap{\smash{\raisebox{1.35ex}{$\hspace{0.4em}\scriptveryshortrightarrow$}}}K%
    \if\relax\detokenize{#1}\relax^\top\else^{(\mkern-1.5mu#1\mkern-1.5mu) \top}\fi
}
\newcommand{\KbwdT}[1][]{%
    \rlap{\smash{\raisebox{1.35ex}{$\hspace{0.34em}\scriptveryshortleftarrow$}}}K%
    \if\relax\detokenize{#1}\relax^\top\else^{(\mkern-1.5mu#1\mkern-1.5mu)\top}\fi
}

\input{mynotation.tex}

\usepackage[notext,nomath]{stix}

\ifusebiblatex
\fi

\usepackage{fancyhdr}
\title{{\sc Cyqlone}: A Parallel, High-Performance Linear Solver for Optimal Control}
\date{\today{}}

\begin{document}

\author{Pieter Pas \and Panagiotis Patrinos%
\thanks{\noindent
The authors are with the STADIUS Center for Dynamical Systems, Signal Processing and Data Analytics,
Department of Electrical Engineering (ESAT),
KU Leuven, Belgium. \quad
Email: \texttt{\{pieter.pas,panos.patrinos\}@esat.kuleuven.be}
}
}

\makeatletter
\let\@fnsymbol@orig\@fnsymbol
\def\@fnsymbol#1{}
\makeatother

\maketitle

\makeatletter
\let\@fnsymbol\@fnsymbol@orig
\makeatother

\begin{abstract}
We present \cyqlone{}, a solver for linear systems with a stage-wise optimal control
structure that fully exploits the various levels of parallelism
available in modern hardware. \cyqlone{} unifies algorithms based on the
sequential Riccati recursion, parallel Schur complement methods, and cyclic reduction methods,
thereby minimizing the required number of floating-point operations, while allowing
parallelization across a configurable number of processors.
Given sufficient parallelism, the solver run time scales with the logarithm of
the horizon length (in contrast to the linear scaling of sequential Riccati-based methods),
enabling real-time solution of long-horizon problems.
Beyond multithreading on multi-core processors, implementations of \cyqlone{} can also
leverage vectorization using batched linear algebra routines. Such batched
routines exploit data parallelism using single instruction, multiple data
(SIMD) operations, and expose a higher degree of instruction-level
parallelism than their non-batched counterparts. This enables them to significantly
outperform BLAS and BLASFEO for the small matrices that arise in optimal control.
Building on this high-performance linear solver, we develop \cyqpalm{},
a parallel and optimal-control-specific variant of the QPALM quadratic
programming solver. It combines the parallel and vectorized linear algebra
operations from \cyqlone{} with a parallel line search and parallel
factorization updates, resulting in order-of-magnitude speedups over the
state-of-the-art HPIPM solver.
Open-source \Cpp{} implementations of \cyqlone{} and \cyqpalm{} are available
at {\fontsize{9}{11}\selectfont \url{https://github.com/kul-optec/cyqlone}}.

\end{abstract}

\medskip

\section{Introduction}
The need for the efficient solution of linear-quadratic optimal control problems
(OCPs) arises in popular engineering applications such as linear model
predictive control (MPC), moving horizon estimation (MHE) and data
assimilation (DA). Moreover, numerical solvers that target nonlinear MPC, MHE or DA
problems are often based on iterative procedures that solve linearized
subproblems, which can in turn be reformulated as linear-quadratic OCPs \cite{pas_gaussnewton_2023}.
As these subproblems account for the majority of the computational effort of the
overall solver, improving the performance of linear-quadratic OCP solvers can
greatly reduce the solver run time for a wide class of problems.
Furthermore, applications like MPC and MHE impose hard limits on the solver run time,
while often running on constrained embedded devices.
Although general-purpose QP solvers can be used in such settings,
OCP-specific solvers such as HPIPM \cite{frison_hpipm_2020} and PIQP multi-stage
\cite{schwan_exploiting_2025-1} take advantage of the
particular OCP structure, typically resulting in faster run times,
lower power consumption, and a smaller memory footprint.

With these requirements in mind, the goal of the present work is the development
of a performant numerical linear-quadratic OCP solver that exploits both the
inherent parallel structure of OCPs and the various levels of parallelism
available in modern processors to minimize solver run time.
Specifically, we consider the parallelization and vectorization of the
factorization and solution of linear systems with
optimal control structure originating from the optimality conditions
of linear-quadratic OCPs. This parallel linear solver, which we call \cyqlone{},
is then used to solve Newton systems in a fully parallel implementation of the
augmented Lagrangian-based \qpalm{} method \cite{hermans_qpalm_2022,hermans_qpalm_2019-1,lowenstein_qpalm-ocp_2024}.
The resulting \cyqpalm{} solver, a highly performant quadratic programming solver for optimal control problems with inequality constraints,
achieves speedups of an order of magnitude or higher compared to the state-of-the-art
HPIPM solver, particularly for problems with long horizons.

\medskip

\subsection{Exploiting the parallelism available in modern hardware}

\begin{figure}
    \centering
    \includegraphics[scale=0.8, clip, trim=0 0.45cm 0 0.3cm]{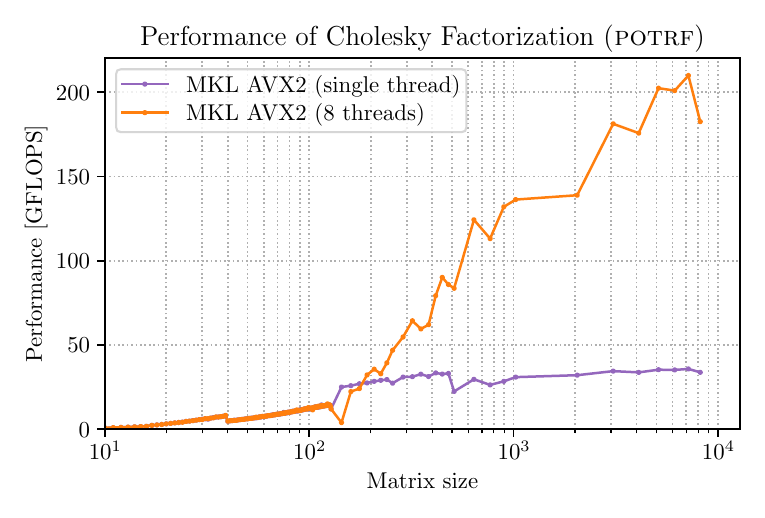}
    \caption{Performance of single-threaded and multithreaded Cholesky factorization (\textsc{potrf}) routines
    from the Intel MKL, on an Intel Core Ultra 7 265. The multithreaded variant
    becomes more performant than the single-threaded variant for sizes $192\times 192$ and larger.
    }
    \label{fig:potrf-threaded-perf}
\end{figure}

From dual-core Raspberry Pi Pico microcontrollers to 192-core
x86-64 server CPUs, multiprocessors dominate the hardware landscape.
As gains in clock speed and instructions per cycle have slowed over
the past decades, parallel processing has become essential for reaching peak
performance on modern systems \cite[Ch.\,5]{hennessy_computer_2011}.
In the specific context of optimal control, existing methods such as the
well-known Riccati recursion used in solvers like HPIPM and \fatrop{} \cite{vanroye_fatrop_2023-1}
are inherently sequential and not directly amenable to parallelization at the algorithmic level
because of the recursive backward dynamic programming approach.
For problems with large numbers of states, individual linear algebra
operations can be parallelized using multithreaded BLAS (basic linear algebra subprograms)
and LAPACK (linear algebra package) implementations \cite{smith_anatomy_2014}. However, for small matrices,
parallelization within a single linear algebra operation comes with significant
overhead compared to single-threaded alternatives, largely canceling out the
performance gained from using multiple processors
\cite[Fig.\,16]{catalan_programming_2020}, as shown in \Cref{fig:potrf-threaded-perf}.

For this reason, the \cyqlone{} method presented in this paper considers
parallelization along the time dimension of the OCPs, resulting in high
speedups compared to state-of-the-art serial methods, even for small
matrices.
Existing solvers that parallelize along the time dimension of OCPs
such as
parallel Schur complement methods \cite{pas_exploiting_2024,frasch_parallel_2015},
partitioned dynamic programming \cite{wright_partitioned_1991},
related parametric subdivision approaches \cite{nielsen_structure_2015,jallet_parallel_2024-1}
and cyclic reduction methods \cite{nicholson_parallel_2018}
suffer from a number of limitations, as discussed in more detail in \cref{sec:compare-ocp-lin-solver}:
some methods contain serial bottlenecks by only parallelizing part of the algorithm,
they may require significantly higher FLOP (floating-point operation)
counts compared to serial alternatives,
or they may be suboptimal when the number of available processors is smaller than
the horizon length of the OCP.
\cyqlone{} unifies the sequential Riccati recursion, parallel Schur complement
methods and cyclic reduction methods to
minimize the FLOP count by exploiting the problem-specific structure, while
allowing the user to select the amount of parallelism that is most suitable
for the hardware at hand.

Efficient linear algebra routines are fundamental to an optimized implementation of the \cyqlone{} solver.
Many existing optimization solvers for optimal control, such as HPIPM, \fatrop{}, and PIQP multi-stage \cite{schwan_exploiting_2025-1} rely on the BLASFEO library \cite{frison_blasfeo_2018}:
BLASFEO stores the matrices in a specialized packed format, referred to as \textit{panel-major} storage order.
This storage format obviates the need for the online packing step performed by general-purpose BLAS implementations \cite{van_zee_blis_2015},
improves cache locality, and enables vectorization.
These properties make BLASFEO an attractive choice for the small to medium-sized
matrices commonly used in solvers for optimal control.

\begin{figure}
    \centering
    \includegraphics[scale=0.8, clip, trim=0 0.45cm 0 0.3cm]{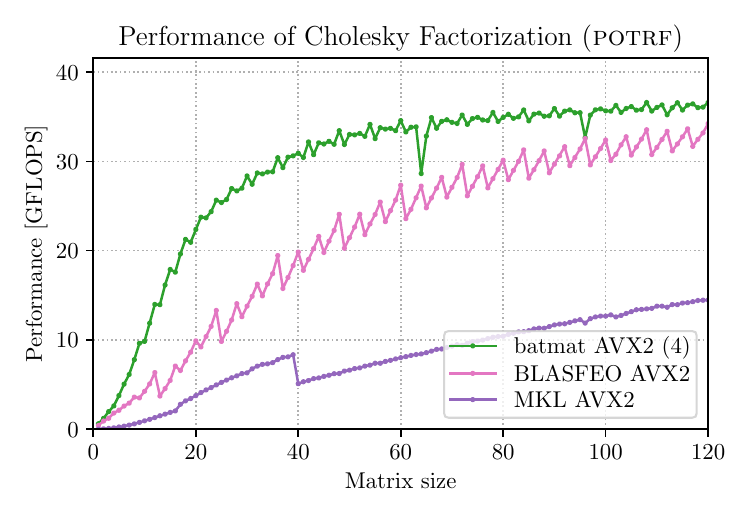}
    \caption{Performance of different Cholesky factorization (\textsc{potrf}) routines
    when applied to batches of eight square matrices.
    The \textsc{Batmat} implementation using batch-wise vectorization \cite{tttapa_batmat_2025}
    significantly outperforms the BLASFEO and Intel MKL BLAS implementations for small matrices.
    Experiments were carried out on an Intel Core Ultra 7 265 with a sustained clock speed of 5.2\,GHz and a theoretical peak AVX2 throughput of 41.6 GFLOPS.
    }
    \label{fig:syrk-potrf-perf}
\end{figure}

One limitation of the vectorization technique used by BLASFEO is that it is
less efficient for small or structured (triangular or symmetric) matrices.
Additionally, some critical steps in e.g. the Cholesky factorization routine cannot be
vectorized effectively. We will see that the parallelism exposed by the
\cyqlone{} solver not only enables parallelization across different processors,
but also vectorization along the time dimension of the OCP, or more generally,
\textit{batch-wise vectorization}.
This vectorization strategy, which was previously inapplicable to solvers for optimal control,
is able to outperform BLAS and BLASFEO for small matrices such as the ones arising in optimal control, as demonstrated for the Cholesky factorization in \Cref{fig:syrk-potrf-perf}.

In the development of \cyqlone{} and \cyqpalm{}, we consider the following levels of parallelism
to fully exploit the performance of modern hardware \cite{hennessy_computer_2011}:
\begin{enumerate}
    \setlength\itemsep{0em}
    \item \textbf{Thread-level parallelism (TLP)}: The proposed \cyqlone{} solver
    for linear systems with optimal control structure
    enables parallelization across different stages of the OCP by partitioning
    the horizon into smaller sub-intervals which are solved by different
    processors in parallel.
    This allows the method to leverage the full
    computational power of modern multi-core processors,
    with run times scaling logarithmically with the horizon length
    if enough processors are available
    (in contrast to sequential solvers such as the Riccati-based solvers
    used by HPIPM and \fatrop{}, whose run times scale linearly with
    horizon length).
    \item \textbf{Data-level parallelism (DLP)}: Batches of four or eight matrices from different stages are
    interleaved into a specialized compact storage format,
    enabling the use of fast batched linear algebra routines.
    By using single instruction, multiple data (SIMD) operations across such a batch, these routines can be implemented
    more efficiently than standard BLAS routines using row-wise vectorization (for small matrices).
    When applied to matrices whose dimensions are not multiples of the hardware's vector length,
    or when applied to triangular or symmetric matrices, batch-wise vectorization along the time dimension of the OCP
    results in improved utilization of the available vector lanes compared to row-wise vectorization.
    \item \textbf{Instruction-level parallelism (ILP)}:
    Our implementation of \cyqlone{} makes use of custom batched
    linear algebra routines built on top of optimized micro-kernels \cite{tttapa_batmat_2025,van_zee_blis_2015}:
    these low-level implementations of fundamental linear algebra operations for small blocks of matrices are hand-coded to exploit pipelining,
    out-of-order execution and superscalar capabilities of modern CPUs. By building blocked algorithms on top of these micro-kernels,
    keeping in mind the properties of the memory hierarchy, high performance can be achieved for larger matrices as well.
    Furthermore, by processing matrices in batches rather than one at a time,
    bottlenecks posed by data dependencies within a matrix (such as dependencies on the pivot in Gaussian elimination) can be avoided, further improving ILP.
\end{enumerate}

\subsection{Contributions}

The main contributions of this work can be summarized as follows:

\begin{enumerate}
    \setlength\itemsep{0em}
    \item We present \cyqlone{}, a new \textbf{parallel linear solver for linear-quadratic control problems}.
    This solver exploits the particular structure of optimal control problems by using a modified
    Riccati recursion on different intervals of the OCP horizon, followed by cyclic reduction of the remaining Schur complement. Furthermore, it is
    able to adapt to various degrees of parallelism, depending on how many
    processors are available, without restrictions on the horizon length,
    enabling excellent scaling to long horizons and large numbers of processors. \rmrk{\Cref{sec:cyqlone}}
    \item A \textbf{parallel solver for quadratic programs with optimal control structure}
    named \cyqpalm{} is developed that leverages \cyqlone{} as its linear solver.
    This solver is based on the recently proposed \qpalmocp{} method \cite{lowenstein_qpalm-ocp_2024}, and
    aims to fully parallelize all steps of the algorithm,
    including parallel matrix factorizations and back substitutions using \cyqlone{},
    parallel evaluation of the necessary gradients and residuals,
    and a parallel exact line search. \rmrk{\Cref{sec:qpalm-ocp}} \\
    A popular optimal control benchmark is used to demonstrate
    that the \cyqlone{} and \cyqpalm{} solvers outperform
    existing state-of-the-art \\solvers by a significant margin. \rmrk{\Cref{sec:results}}
    \item Parallel \textbf{factorization update routines} for \cyqlone{} are derived. We show how they further improve the performance of \cyqpalm{}. \rmrk{\Cref{sec:fact-upd}}
    \item \textbf{Optimized parallel and vectorized \Cpp{} implementations} of the \cyqlone{} and \cyqpalm{} solvers are developed,
    which achieve excellent performance, even for small matrices.
    The source code is made available under an open-source license \cite{kul-optec_cyqlone_2025}. \rmrk{\Cref{sec:vectorization}}
    \item We perform a detailed \textbf{comparison of different solution methods for linear--quadratic control problems},
    discussing the tradeoffs in total computational cost and parallelizability.
    We show how \cyqlone{} relates to and improves upon existing methods. \rmrk{\Cref{sec:compare-ocp-lin-solver}}
\end{enumerate}

\section{Problem statement and notation \done}\label{sec:formulation}
We consider discrete-time linear-quadratic optimal control problems with inequality constraints of the form
\begin{equation} \label{eq:ocp}
     \begin{aligned}
        &\minimize_{\mathbf{u}, \mathbf{x}}&& \sum_{j=0}^{N-1} \ell_j(u^j, x^j) + \ell_N(x^N) \\
        &\subjto && E_0 x^0 = x_\text{init} \\
        &&& E_{j+1} x^{j+1} = A_j x^j + B_j u^j + f^j &&{\color{gray}\scriptstyle(0 \le j \lt N)} \\
        &&& b_l^j \le C_j x^j + D_j u^j \le b_u^j &&{\color{gray}\scriptstyle(0 \le j \lt N)} \\
        &&& b_l^N \le C_N x^N \le b_u^N,
    \end{aligned} 
\end{equation}
where the convex stage-wise and terminal cost functions are given by \\
$
    \ell_j(u, x) = \tfrac12 \begin{pmatrix}
        u \\ x
    \end{pmatrix}^{\!\!\top}\!\!\! \begin{pmatrix}
        R^\ell_j & S^\ell_j \\
        S^{\ell\top}_j\!\!\!\! & Q^\ell_j
    \end{pmatrix}\! \begin{pmatrix}
        u \\ x
    \end{pmatrix} + \begin{pmatrix}
        u \\ x
    \end{pmatrix}^{\!\!\top}\!\!\! \begin{pmatrix}
        r_\ell^j \\ q_\ell^j
    \end{pmatrix}
    \;\text{and}\;
    \ell_N(x) = \tfrac12 \tp x Q^\ell_N x + \tp x q_\ell^N,
$
with $Q^\ell_j \in \possdefset{\R^{n_x}}$, $S^\ell_j \in \R^{n_u\times n_x}$ and $R^\ell_j \in \possdefset{\R^{n_u}}$.
The state variables $\mathbf{x} = [x^0\, \cdots\, x^N] \in \R^{n_x\times(N+1)}$ and controls $\mathbf{u} = [u^0\, \cdots\, u^{N-1}] \in \R^{n_u\times N}$
are governed by the linear dynamics described by $A_j \in \R^{n_x\times n_x}$,
$B_j \in \R^{n_x\times n_u}$ and $f^j \in \R^{n_x}$. Linear stage-wise mixed state--input
constraints are encoded using $C_j\in\R^{n_y \times n_x}$ and $D_j\in\R^{n_y \times n_u}$,
with $C_N \in\R^{n_y \times n_x}$ for the terminal constraints.
The invertible matrices $E_j \in \R^{n_x\times n_x}$ allow for implicit state equations.
Although the matrices $E_j$ are used by some of the methods discussed in \Cref{sec:compare-ocp-lin-solver},
we will usually assume without loss of generality that $E_j=\I$, which can be achieved
by transforming $E_j x^{j+1} = A_j x^j + B_j u^j + f^j$ into
$x^{j+1} = \inv E_j A_j x^j + \inv E_j B_j u^j + \inv E_j f^j$.\,\footnote{In fact,
this transformation also improves performance in practice, since it lowers the
number of operations required when solving linear systems involving the constraints.
The upfront cost of computing an LU factorization of $E_j$
is often negligible compared to the resulting performance gains in the optimization solver.}

In the first part of the paper, and in \Cref{sec:cyqlone} specifically, we focus on linear-quadratic optimal control problems subject to
equality constraints only:
\begin{equation} \label{eq:ocp-eq}
    \begin{aligned}
         & \minimize_{\mathbf{u}, \mathbf{x}} &  & \sum_{j=0}^{N-1} \varphi_j(x^j, u^j) + \varphi_N(x^N)                                    \\
         & \subjto                            &  & E_0 x^0 = x_\text{init}                                                                          \\
         &                                    &  & E_{j+1} x^{j+1} = A_j x^j + B_j u^j + f^j,                     &  & {\color{gray}\scriptstyle(0 \le j \lt N)} \\
    \end{aligned}
\end{equation}
where the strongly convex stage-wise and terminal cost functions are given by \\
$
    \varphi_j(u, x) = \tfrac12 \begin{pmatrix}
        u \\ x
    \end{pmatrix}^{\!\!\top}\!\!\! \begin{pmatrix}
        R_j & S_j \\
        S^{\top}_j\!\!\!\! & Q_j
    \end{pmatrix} \begin{pmatrix}
        u \\ x
    \end{pmatrix} + \begin{pmatrix}
        u \\ x
    \end{pmatrix}^{\!\!\top}\!\!\! \begin{pmatrix}
        r^j \\ q^j
    \end{pmatrix}
    \;\text{and}\;
    \varphi_N(x) = \tfrac12 \tp x Q_N x + \tp x q^N.
$
The connection between the general optimal control problem \eqref{eq:ocp} and the strongly convex equality constrained variant \eqref{eq:ocp-eq}
will be discussed in \Cref{sec:qpalm-ocp}.

After elimination of the fixed initial state $x^0$, the Karush--Kuhn--Tucker (KKT) optimality conditions of \eqref{eq:ocp-eq}
are given by
\begin{equation} \label{eq:opt-cond-ocp-eq}
    \begin{aligned}
        \left\{\begin{aligned}
                    & Q_N x^N - \ttp E_N \lambda^{N-1}                                                   &  & = -q^N \quad &                                \\
                    & R_j u^j + S_j x^j + \tp B_j \lambda^j                                       &  & = -r^j \quad & {\color{gray} \scriptstyle (0 \lt j \lt N)} \\
                    & \ttp{S_j} \hspace{1pt}u^j + Q_j x^j + \tp A_j \lambda^j - \ttp E_j \lambda^{j-1} &  & = -q^j \quad & {\color{gray} \scriptstyle (0 \lt j \lt N)} \\
                    & R_0 u^0 + \tp B_0 \lambda^0                                                   &  & = -\tilde r^0                                        \\
                    & B_j u^j + A_j x^j - E_{j+1} x^{j+1}                                       &  & = -f^j \quad & {\color{gray} \scriptstyle (0 \lt j \lt N)} \\
                    & B_0 u^0 - E_1 x^1                                                           &  & = -\tilde f^0,
               \end{aligned}\right.
    \end{aligned}
\end{equation}
with $\tilde r^0 = r^0 + S_0 \inv E_0 x_\text{init}$ and $\tilde f^0 = f^0 + A_0 \inv E_0 x_\text{init}$.
The Lagrange multipliers $\lambda^j \in \R^{n_x}$ correspond to the dynamics constraints between stage $j$ and stage $j+1$.
We will refer to linear systems of the form \eqref{eq:opt-cond-ocp-eq} as \textit{KKT systems with optimal control structure}.

\section{Preliminaries}
We start by reviewing some facts about the block Cholesky factorization that will
be used throughout the remainder of this paper.
The block Cholesky factorization formula introduced below is key in deriving
the \cyqlone{} algorithm; the cyclic reduction (CR) algorithm is used to solve
block-tridiagonal systems in parallel;
and the concepts of fill-in and the elimination tree are essential in analyzing
the efficiency and parallelizability of linear solvers.

\subsection{Block Cholesky factorizations \done} \label{sec:block-chol}

A $2\times 2$ block Cholesky factorization $LD\ttp L$ of a symmetric matrix $H$
with $L$ block lower triangular and $D$ block diagonal satisfies \vspace{-0.4em}
\begin{equation} \label{eq:blk-chol}
    \begin{aligned}
    H &= \begin{pmatrix}
        H_{11} & \ttp H_{21} \\
        H_{21} & H_{22}
    \end{pmatrix} =
    \begin{pmatrix}
        L_{11} \\ L_{21}
    \end{pmatrix} D_1
    \ttp{\begin{pmatrix}
            L_{11} \\ L_{21}
        \end{pmatrix}} +
    \begin{pmatrix}
        0 & 0 \\ 0 & H / H_{11}
    \end{pmatrix} \\
    &=
    {
        \underbrace{
    \begin{pmatrix}
        L_{11} & \\ L_{21} & L_{22}
    \end{pmatrix}}_{L}
    \underbrace{\begin{pmatrix}
        D_1 & \\ & D_2
    \end{pmatrix}}_{D}
    \underbrace{\begin{pmatrix}
            L_{11} & \\ L_{21} & L_{22}
        \end{pmatrix}^{\mathclap{\!\top\!}}\!}_{\ttp L}}
    \\
    &\text{where }\begin{aligned}[t]
        H_{11} &= L_{11} D_1 \ttp L_{11}, \qquad
        L_{21} = H_{21} \invtp{L}_{11} \inv D_1, \\
        H/H_{11} &\defeq H_{22} - H_{21} \inv H_{11} \ttp H_{21} = H_{22} - L_{21} D_1 \ttp L_{21} = L_{22} D_2 \ttp L_{22}.
    \end{aligned}
    \end{aligned}
\end{equation}
The block $H/H_{11}$ is known as the Schur complement of $H_{11}$ in $H$.
The matrices $L_{ii}$ and $D_{i}$ are not unique: we only require that the inverses of $L_{ii}$ and $D_{i}$ are easy to apply, and that $D_{i}$ is symmetric.
In practice, we will restrict ourselves to lower and upper triangular matrices $L_{ii}$, and signature matrices $D_i$.\,\footnote{Signature matrices are diagonal matrices with elements $\pm1$.}
This leaves some freedom in how we choose $L_{ii}$ and $D_i$: a first design parameter is the size of the block $H_{11}$, a second is the factorization of $H_{11}$,
and a third is the factorization of $H/H_{11}$. The size of $H_{11}$ will usually follow from the block structure of the original matrix $H$.
If an ``obvious'' factorization of $H_{11}$ or $H/H_{11}$ exists, we can use it as our choice of $L_{ii}$ and $D_i$.
Specifically, if $H_{11}$ is a $1\times 1$ matrix, we use the scalar factorization $L_{11} = \sqrt{|H_{11}|}$ and $D_1 = \sgn(H_{11})$, and similarly for the case where $H/H_{11}$ is a $1\times 1$ matrix.
If no such factorization is available,
we instead use another smaller block Cholesky factorization, obtained by recursive application of \eqref{eq:blk-chol} to $H_{11}$ and/or $H/H_{11}$.

The block Cholesky factorization procedure can be summarized as follows: \vspace{-1em}
\begin{enumerate}
    \setlength\itemsep{0em}
    \setcounter{enumi}{-1}
    \item Isolate the top-left block $H_{11}$; \label{it:isolate-H11}
    \item Factorize $H_{11}$: select $L_{11}$ and $D_1$ such that $H_{11} = L_{11} D_1 \ttp L_{11}$, either using specific knowledge about $H_{11}$ or by recursively computing a (block) Cholesky factorization;
    \item Compute the subdiagonal block $L_{21} = H_{21} \invtp L_{11} \inv D_1$;
    \item Subtract the product $L_{21} D_1 \ttp L_{21}$ from the block $H_{22}$ to compute the Schur complement $H/H_{11}$; \label{it:update-schur-compl}
    \item Factorize $H/H_{11}$: select $L_{22}$ and $D_2$ such that $H/H_{11} = L_{22} D_2 \ttp L_{22}$, either using specific knowledge about $H/H_{11}$ or by recursively computing a (block) Cholesky factorization.
\end{enumerate}

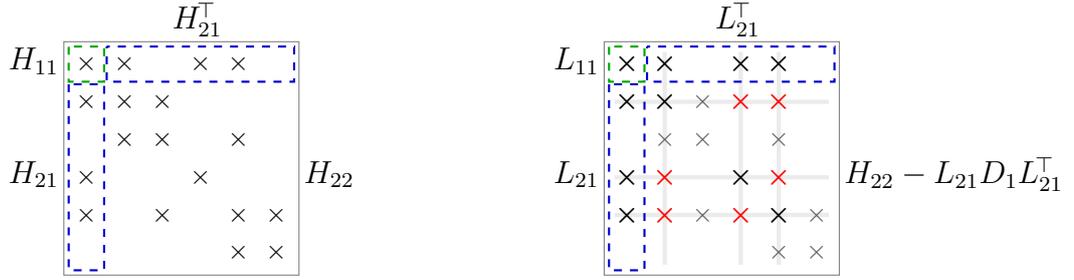
\begin{figure}
\centering
\hfill
\begin{tikzpicture}[scale=0.7]
  \matrix[matrix of nodes,nodes={minimum size=4mm, anchor=center},column sep=0mm, row sep=0mm] (m) {
    ~ & ~ & ~ & ~ & ~ & ~ \\
    ~ & ~ & ~ & ~ & ~ & ~ \\
    ~ & ~ & ~ & ~ & ~ & ~ \\
    ~ & ~ & ~ & ~ & ~ & ~ \\
    ~ & ~ & ~ & ~ & ~ & ~ \\
    ~ & ~ & ~ & ~ & ~ & ~ \\
  };

  \draw[gray] 
    ($(m-1-1.north west)+(-0.05,0.05)$) rectangle ($(m-6-6.south east)+(0.05,-0.05)$);

  \foreach \i/\j in {3/2, 5/3, 6/5} {
    \node at (m-\i-\j) {$\times$};
    \node at (m-\j-\i) {$\times$};
  }
  \foreach \i/\j in {2/1, 4/1, 5/1} {
    \node at (m-\i-\j) {$\times$};
    \node at (m-\j-\i) {$\times$};
  }

  \foreach \i in {1,...,6} {
    \node at (m-\i-\i) {$\times$};
  }

  \draw[blue!80!black,thick,dashed] 
    ($(m-2-1.north west)+(0.03,-0.03)$) rectangle ($(m-6-1.south east)+(-0.03,0.03)$);
  \draw[blue!80!black,thick,dashed,dash phase=2.3] 
    ($(m-1-2.north west)+(0.03,-0.03)$) rectangle ($(m-1-6.south east)+(-0.03,0.03)$);
  \draw[green!66!black,thick,dashed,dash phase=-0.1]
    ($(m-1-1.north west)+(0.03,-0.03)$) rectangle ($(m-1-1.south east)+(-0.03,0.03)$);

  \node at ($(m-1-1)+(-0.3,0.05)$) [left] {$H_{11}$};
  \node at ($(m-4-1)+(-0.3,0.05)$) [left] {$H_{21}$};
  \node at ($(m-1-4)+(-0.05,0.3)$) [above] {$\ttp H_{21}$};
  \node at ($(m-4-6)+(0.3,0.05)$) [right] {$H_{22}$};
\end{tikzpicture}
\hfill
\begin{tikzpicture}[scale=0.7]
  \matrix[matrix of nodes,nodes={minimum size=4mm, anchor=center},column sep=0mm, row sep=0mm] (m) {
    ~ & ~ & ~ & ~ & ~ & ~ \\
    ~ & ~ & ~ & ~ & ~ & ~ \\
    ~ & ~ & ~ & ~ & ~ & ~ \\
    ~ & ~ & ~ & ~ & ~ & ~ \\
    ~ & ~ & ~ & ~ & ~ & ~ \\
    ~ & ~ & ~ & ~ & ~ & ~ \\
  };

  \draw[gray] 
    ($(m-1-1.north west)+(-0.05,0.05)$) rectangle ($(m-6-6.south east)+(0.05,-0.05)$);

  \foreach\i in {2, 4, 5} {
    \draw[lightgray,ultra thick, opacity=0.3]
        ($(m-\i-1)-(0.2,0)$) -- ($(m-\i-6)+(0.2,0)$);
    \draw[lightgray,ultra thick, opacity=0.3]
        ($(m-1-\i)+(0,0.2)$) -- ($(m-6-\i)-(0,0.2)$);
  }

  \foreach \i/\j in {3/2, 5/3, 6/5} {
    \node at (m-\i-\j) {$\color{gray!70!black}\times$};
    \node at (m-\j-\i) {$\color{gray!70!black}\times$};
  }
  \foreach \i/\j in {2/1, 4/1, 5/1} {
    \node at (m-\i-\j) {$\boldsymbol\times$};
    \node at (m-\j-\i) {$\boldsymbol\times$};
  }

  \foreach \i in {3, 6} {
    \node at (m-\i-\i) {$\color{gray!70!black}\times$};
  }
  \foreach \i in {1, 2, 4, 5} {
    \node at (m-\i-\i) {$\boldsymbol\times$};
  }

  \node[text=red] at (m-5-4) {$\boldsymbol\times$};
  \node[text=red] at (m-4-5) {$\boldsymbol\times$};
  \node[text=red] at (m-2-4) {$\boldsymbol\times$};
  \node[text=red] at (m-2-5) {$\boldsymbol\times$};
  \node[text=red] at (m-4-2) {$\boldsymbol\times$};
  \node[text=red] at (m-5-2) {$\boldsymbol\times$};

  \draw[blue!80!black,thick,dashed] 
    ($(m-2-1.north west)+(0.03,-0.03)$) rectangle ($(m-6-1.south east)+(-0.03,0.03)$);
  \draw[blue!80!black,thick,dashed,dash phase=2.3] 
    ($(m-1-2.north west)+(0.03,-0.03)$) rectangle ($(m-1-6.south east)+(-0.03,0.03)$);
  \draw[green!66!black,thick,dashed,dash phase=-0.1]
    ($(m-1-1.north west)+(0.03,-0.03)$) rectangle ($(m-1-1.south east)+(-0.03,0.03)$);

  \node at ($(m-1-1)+(-0.3,0.05)$) [left] {$L_{11}$};
  \node at ($(m-4-1)+(-0.3,0.05)$) [left] {$L_{21}$};
  \node at ($(m-1-4)+(-0.05,0.3)$) [above] {$\ttp L_{21}$};
  \node at ($(m-4-6)+(0.3,0.05)$) [right] {$H_{22} - L_{21}D_1\ttp L_{21}$};
\end{tikzpicture}
\hfill
\caption{Visualization of steps \ref{it:isolate-H11}--\ref{it:update-schur-compl} of the block Cholesky factorization of a sparse symmetric matrix $H$ as in \eqref{eq:blk-chol}. Nonzero blocks of $H$ are indicated by crosses ($\times$).
Blocks that are modified by the elimination of $H_{11}$ and $H_{21}$ are shown in bold, and blocks of fill-in in the Schur complement are shown in red.}
\label{fig:blk-chol}
\end{figure}

When computing $LD\ttp L$ factorizations of sparse matrices, an important consideration is \textit{fill-in}:
the outer product $L_{21} D_1 \ttp L_{21}$ is nonzero at the intersections of the rows where $L_{21}$ is nonzero and
the columns where $\ttp L_{21}$ is nonzero. Consequently, the matrix $H/H_{11}$ may be significantly less sparse than $H_{22}$. Zero entries of $H_{22}$ that become nonzero in $H/H_{11}$ are referred to as fill-in.
\Cref{fig:blk-chol} visualizes the phenomenon for a small example matrix.
Such fill-in may cause more fill-in down the line, and causes the factorization of $H/H_{11}$ to require more floating-point operations.
Practical sparse solvers therefore use (heuristic) fill-in-reducing permutations of the original matrix $H$ to maintain sparsity.

\subsection{Cyclic reduction of block-tridiagonal linear systems \done} \label{sec:block-tri-diag-chol-cr}

Symmetric (block-)tridiagonal linear systems play a crucial role in the solution of optimal control and other dynamic optimization problems,
because they can be used to represent coupling between successive time steps, as is the case in \eqref{eq:opt-cond-ocp-eq}.
In this section, we discuss solution methods for block-tridiagonal systems, and
in particular, the CR algorithm for solving such systems in parallel.
Consider a symmetric positive definite block-tridiagonal system of equations of size $N$
with diagonal blocks $\Mdiag_k \in \posdefset{\R^{n}}$ and subdiagonal blocks $\Ksub_k \in \R^{n\times n}$:
\begin{equation} \label{eq:block-tri}
    \mathscr M\bf x=\scalebox{1.0}{$\setlength{\arraycolsep}{1pt}\begin{pmatrix}
                \Mdiag_{0\phantom{\shortminus N}} & \tp \Ksub_{0\phantom{\shortminus N}}                                                                                              \\
                \Ksub_{0\phantom{\shortminus N}} & \Mdiag_{1\phantom{\shortminus N}}     & \tp \Ksub_{1}                                                                                  \\
                                             & \Ksub_{1\phantom{\shortminus N}}     & \scalebox{0.6}{$\ddots$} & \scalebox{0.6}{$\ddots$}                                        \\
                                             &                                  & \scalebox{0.6}{$\ddots$} & \Mdiag_{N\shortminus 2}                  & \; \tp \Ksub_{N\shortminus 2} \\[2pt]
                                             &                                  &                          & \Ksub_{N\shortminus 2}^{\phantom{\top}} & \; \Mdiag_{N\shortminus 1}     \\
            \end{pmatrix}
            \begin{pmatrix}
                x^0 \\ x^1 \\ \scalebox{0.6}{$\vdots$} \\ x^{N\shortminus 2} \\ x^{N\shortminus 1}
            \end{pmatrix}
            =
            \begin{pmatrix}
                b^0 \\ b^1 \\ \scalebox{0.6}{$\vdots$} \\ b^{N\shortminus 2} \\ b^{N\shortminus 1}
            \end{pmatrix} = \bf b.
        $}
\end{equation}
Such block-tridiagonal matrices can be factorized using the block Cholesky factorization outlined in the previous section. Thanks to the special structure, only the top block of $L_{21}$ will be nonzero, so only a single diagonal block of $H_{22}$ needs to be updated, without any fill-in. If $\Mdiag_0 = L_{11}D_1 \ttp L_{11}$, then
\begin{equation} \label{eq:block-tri-cholesky}
    \begin{aligned}
        L_{21} &= \scalebox{1.0}{$\setlength{\arraycolsep}{1pt}\left(\begin{array}{c}
            \Ksub_0 \invtp L_{11}\! \inv D_1 \\\hline 0 \\ \scalebox{0.6}{$\vdots$} \\ 0
        \end{array}\right)$}\qquad\text{and}\qquad
        \mathscr M/\Mdiag_0 &= \scalebox{1.0}{$\setlength{\arraycolsep}{1pt}\left(\begin{array}{c|ccc}
            \Mdiag_1 - \Ksub_0 \invtp L_{11}\! \inv D_1 \inv L_{11} \tp \Ksub_0 & \quad\tp \Ksub_{1} \\\hline
            \quad \Ksub_{1\phantom{\shortminus N}}     & \scalebox{0.6}{$\ddots$} & \scalebox{0.6}{$\ddots$}                                        \\
                                            & \scalebox{0.6}{$\ddots$} & \Mdiag_{N\shortminus 2}                  & \; \tp \Ksub_{N\shortminus 2} \\[2pt]
                                            &                          & \Ksub_{N\shortminus 2}^{\phantom{\top}} & \; \Mdiag_{N\shortminus 1}     \\
        \end{array}\right).$}
    \end{aligned}
\end{equation}
This block-tridiagonal factorization procedure requires $N$ steps and produces no fill-in.
However, due to the coupling through the off-diagonal blocks, this process is inherently serial: the block $\Mdiag_{k}$ can only be
updated and factorized once $\Mdiag_{k-1}$ and $\Ksub_{k-1}$ have been eliminated, because of the direct coupling between
equations $k-1$ and $k$. A well-known approach that avoids this serial bottleneck
is \textit{cyclic reduction} (CR) \cite{heller_aspects_1976,gander_cyclic_1998,bini_cyclic_2009}.

CR is built on the idea of simultaneous elimination of the odd block rows:
Since there is no direct coupling between equations $2k-1$ and $2k+1$,
their elimination can be performed in parallel.
After eliminating the odd equations of a system of size $N$, the remaining equations form a smaller
block-tridiagonal system of size $\lceil \frac{N}2 \rceil$, and the process can
be repeated recursively, for a total of $\lceil\log_2(N)\rceil + 1$ steps.
Because equations $2k-1$ and $2k+1$ are still coupled (indirectly) through equation $2k$,
this remaining even equation needs to be updated. Consequently, we will find that each
block row elimination in CR causes one block of fill-in. This is a disadvantage compared
to the straightforward block-tridiagonal Cholesky factorization from \eqref{eq:block-tri-cholesky},
but the increased number of operations per row is often an acceptable trade-off for the gained parallelism, especially for large $N$.

\subsubsection{Derivation via simultaneous elimination of odd equations \done} \label{subsec:cr-elim-deriv}

A single application of the CR algorithm involves updating the
even block rows $2k$ by eliminating the odd block rows $2k\pm1$ around them.
Let us consider three equations from \eqref{eq:block-tri} in isolation, for any index $0 \le i \lt \frac {N}2$, with the boundary conditions $\Ksub_{-1} = 0 = \Ksub_{N-1}$:
\begin{equation} \label{eq:tridiagonal-equations-3}
    \scalebox{1.0}{$\left\{
            \begin{aligned}
                 & \Ksub_{2i - 2}\, x^{2i - 2} + \Mdiag_{2i-1}\, x^{2i-1} + \tp \Ksub_{2i-1}\, x^{2i} \hspace{-1em}&& = b^{2i - 1} \\
                 & \Ksub_{2i - 1}\, x^{2i - 1} + \Mdiag_{2i}\, x^{2i} + \tp \Ksub_{2i}\, x^{2i + 1} \hspace{-1em}&& = b^{2i}     \\
                 & \Ksub_{2i}\, x^{2i} + \Mdiag_{2i + 1}\, x^{2i + 1} + \tp \Ksub_{2i + 1}\, x^{2i + 2} \hspace{-1em}&& = b^{2i + 1}. \\
            \end{aligned}
            \right.$}
\end{equation}
We can eliminate $x^{2i-1}$ from equation $2i$ by multiplying equation $2i-1$ by $\Ksub_{2i-1} \inv \Mdiag_{2i-1}$ and subtracting it from equation $2i$.
Similarly, subtracting $\tp \Ksub_{2i} \inv \Mdiag_{2i+1}$ times equation $2i+1$ eliminates $x^{2i+1}$.
\begin{align}
     & \scalebox{1.0}{$\left\{
            \begin{aligned}
                 & \Ksub_{2i-1} \inv \Mdiag_{2i-1} \Ksub_{2i - 2}\, x^{2i - 2} + \cancel{\Ksub_{2i-1} \, x^{2i-1}} + \Ksub_{2i-1} \inv \Mdiag_{2i-1} \tp \Ksub_{2i-1}\, x^{2i} \hspace{-1em}&& = \Ksub_{2i-1} \inv \Mdiag_{2i-1}\, b^{2i - 1} \tikzmark{oddprev}                  \\
                 & \cancel{\Ksub_{2i - 1}\, x^{2i - 1}} + \Mdiag_{2i}\, x^{2i} + \cancel{\tp \Ksub_{2i}\, x^{2i + 1}} \hspace{-1em}&& = b^{2i}\tikzmark{even}                                                   \\
                 & \tp \Ksub_{2i} \inv \Mdiag_{2i+1} \Ksub_{2i}\, x^{2i} + \cancel{\tp \Ksub_{2i}\, x^{2i + 1}} + \tp \Ksub_{2i} \inv \Mdiag_{2i+1} \tp \Ksub_{2i + 1}\, x^{2i + 2} \hspace{-1em}&& = \tp \Ksub_{2i} \inv \Mdiag_{2i+1}\, b^{2i + 1}.\hspace{0.66em}\tikzmark{oddnext}
            \end{aligned}
            \begin{tikzpicture}[overlay, remember picture]
                \path[draw, ->, ]($(oddprev) + (0.1cm,0.10cm)$) to[bend left]  node[midway,right,inner sep=4pt] {$-$} ([xshift=1.9cm, yshift=0.2cm]even);
                \path[draw, ->, ]($(oddnext) + (0.1cm,0.05cm)$) to[bend right] node[midway,right,inner sep=4pt] {$-$}  ([xshift=1.9cm, yshift=0cm]even);
            \end{tikzpicture}%
            \right.$} \notag
    \intertext{By introducing $L_{k} \defeq \operatorname{chol}(\Mdiag_{k})$,\; $Y_{k} \defeq \Ksub_{k} \invtp{L_{k}}$,\; $U_{k} \defeq \tp \Ksub_{k-1} \invtp{L_{k}}$ and $\tilde b^{k} \defeq \inv L_{k} b^{k}$, this simplifies to}
     & \scalebox{1.0}{$\left\{
            \begin{aligned}
                 & Y_{2i-1} \tp U_{2i-1}\, x^{2i - 2} + \cancel{\Ksub_{2i-1} \, x^{2i-1}} + Y_{2i-1}\tp Y_{2i-1}\, x^{2i} \hspace{-1em}&& = Y_{2i-1} \, \tilde b^{2i - 1} \hspace{1.01em}\tikzmark{oddprev} \\
                 & \cancel{\Ksub_{2i - 1}\, x^{2i - 1}} + \Mdiag_{2i}\, x^{2i} + \cancel{\tp \Ksub_{2i}\, x^{2i + 1}} \hspace{-1em}&& = b^{2i}\hspace{0.6em}\tikzmark{even}                             \\
                 & U_{2i+1}\tp U_{2i+1}\, x^{2i} + \cancel{\tp \Ksub_{2i}\, x^{2i + 1}} + U_{2i+1} \tp Y_{2i+1}\, x^{2i + 2} \hspace{-1em}&& = U_{2i+1}\, \tilde b^{2i + 1}.\hspace{0.6em}\tikzmark{oddnext}
            \end{aligned}
            \begin{tikzpicture}[overlay, remember picture]
                \path[draw, ->, ]($(oddprev) + (0.1cm,0.10cm)$) to[bend left]  node[midway,right,inner sep=4pt] {$-$} ([xshift=1.1cm, yshift=0.2cm]even);
                \path[draw, ->, ]($(oddnext) + (0.1cm,0.05cm)$) to[bend right] node[midway,right,inner sep=4pt] {$-$}  ([xshift=1.1cm, yshift=0cm]even);
            \end{tikzpicture}%
            \right.$} \label{eq:cr-updated-odd}
\end{align}
After subtracting equations $2i\pm 1$ from equation $2i$, the resulting equation contains only even variables:
\begin{equation}\label{eq:cr-even}
    \hspace{-4em}%
    \begin{aligned}[b]
        &-Y_{2i-1}\tp U_{2i-1} \, x^{2i-2} + \left( \Mdiag_{2i} \!-\! Y_{2i-1}\tp Y_{2i-1} \!-\! U_{2i+1}\tp U_{2i+1} \right) x^{2i} - U_{2i+1}\tp Y_{2i+1}\, x^{2i+2} \\
        &\quad = b^{2i} - Y_{2i-1}\, \tilde b^{2i-1} - U_{2i+1}\, \tilde b^{2i+1}.
    \end{aligned}%
    \hspace{-4em}%
\end{equation}
This is another tridiagonal system, half the size of the original system, with diagonal blocks $M^{(1)}_{2i} \defeq \Mdiag_{2i} - Y_{2i-1}\tp Y_{2i-1} - U_{2i+1}\tp U_{2i+1}$ and subdiagonal blocks $K^{(1)}_{2i-2} \defeq -Y_{2i-1}\tp U_{2i-1}$.
It can either be solved directly, or by another recursive application of CR.
Once the even variables $x^{2i}$ have been determined,
the odd variables can be recovered from the odd equations in \eqref{eq:tridiagonal-equations-3}. Multiplying by $\inv L_{2i-1}$ and rearranging, clearly $x^{2i-1}$ solves
\begin{equation}\label{eq:cr-odd}
    \tp L_{2i-1}\, x^{2i-1} = \tilde b^{2i-1} - \tp U_{2i-1}\, x^{2i-2} - \tp Y_{2i-1} x^{2i}.
\end{equation}
The full recursive procedure for solving a symmetric block-tridiagonal linear system using CR is described by
\Cref{alg:cr} in \Cref{app:cr}.

\subsubsection{Cyclic reduction as block Cholesky factorization \done} \label{sec:cr-as-chol}

Both (block) CR and (block) Cholesky factorization are instances
of Gaussian elimination, and there exists a natural equivalence between the two methods:
Block CR of a block-tridiagonal matrix can be interpreted as the block Cholesky factorization of
a specific permutation of the original matrix \cite[(3)]{gander_cyclic_1998}, where the block rows
and columns are ordered by increasing 2-adic valuation $\nu_2$,
i.e. the multiplicity of two in the prime factorization of their index $i$ in the original block-tridiagonal matrix:
\begin{equation}
    \nu_2(i) \defeq \begin{cases}
        \infty & \text{if } i = 0, \\
        \max \defset{m \in \N\;}{\;2^m \mathbin{|} i\vphantom{X^X_X}} & \text{otherwise.}
    \end{cases}
\end{equation}
All block rows/columns with odd indices (valuation 0) are ordered first, then the block rows/columns whose indices are divisible by 2 but not 4 (valuation 1),
then indices divisible by 4 but not 8 (valuation 2), etc. The block row/column with index 0 (valuation $\infty$) is ordered last.
This permutation gives rise to block diagonal submatrices of decreasing sizes ($N/2$, $N/4$, $\dots$) that can be factorized in parallel,
as visualized in \Cref{fig:cyclic-reduction-permutation}.
\Cref{fig:cyclic-reduction-permutation-labels-32-vl1} shows the same representation of CR for a larger matrix, clearly highlighting the self-similar, recursive structure of the method.

For practical reasons, it is easier to label the subdiagonal blocks in the permuted block-tridiagonal matrix by the column they belong to:
we define $\Kfwd_{k} \defeq \Ksub_{k}$ and $\Kbwd_{k} \defeq \ttp\Ksub_{k-1}$ (for odd $k$).
For example, as can be seen in \Cref{fig:cyclic-reduction-permutation},
$\Kfwd_{1} \defeq \Ksub_{1}$ is the subdiagonal block that couples column 1 to column 2 (increasing the index, hence the forward arrow),
and $\Kbwd_{1} \defeq \ttp\Ksub_{0}$ is the subdiagonal block that couples column 1 to column 0 (decreasing the index, hence the backward arrow).
This results in more symmetric formulas $Y_{k} = \Kfwd_{k} \invtp{L_{k}}$ and $U_{k} = \Kbwd_{k} \invtp{L_{k}}$.

\subsubsection{Parallelism and the elimination tree \done}

\begin{figure}
    \begin{minipage}{0.77\textwidth}
        \includegraphics[width=\textwidth]{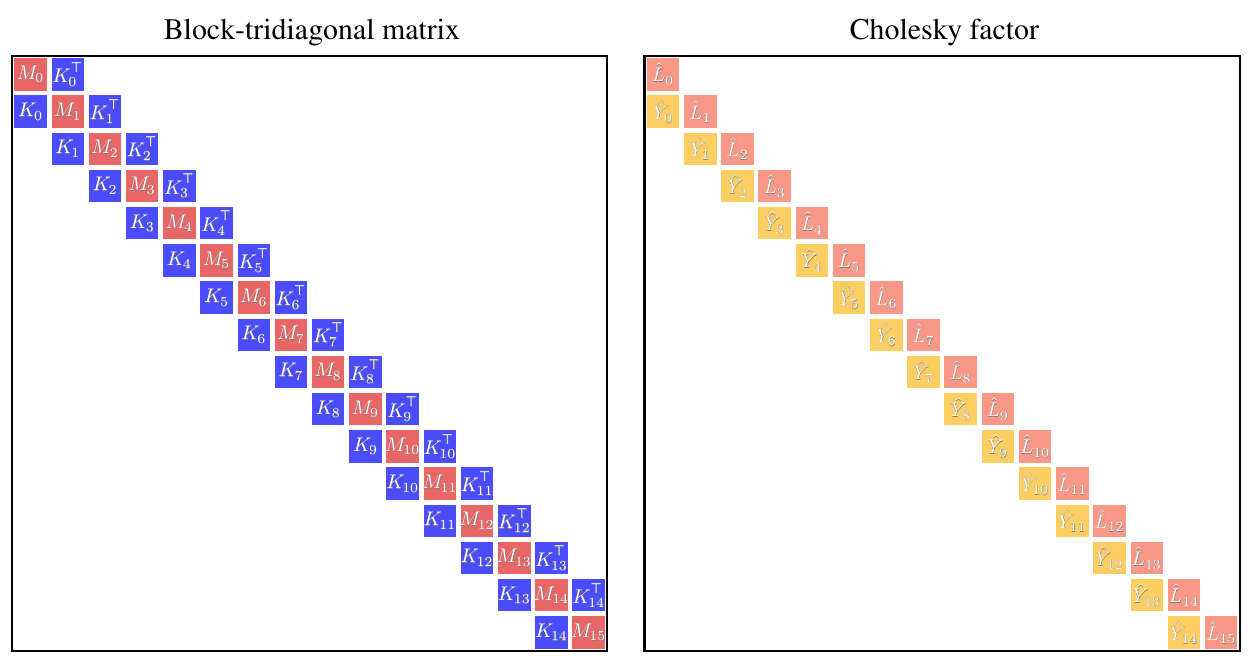}
    \end{minipage}%
    \begin{minipage}{0.23\textwidth}
        \centering
        \includegraphics[scale=0.43]{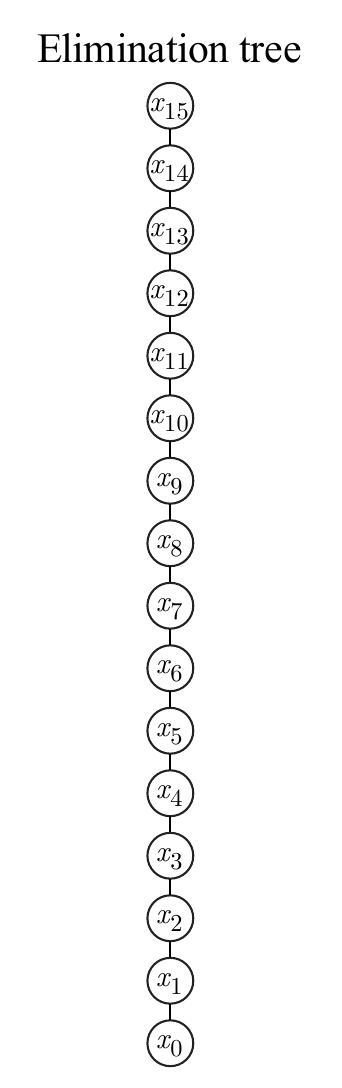}
        \vspace{-3.62em}%
    \end{minipage}
    \begin{minipage}{0.77\textwidth}
        \includegraphics[width=\textwidth]{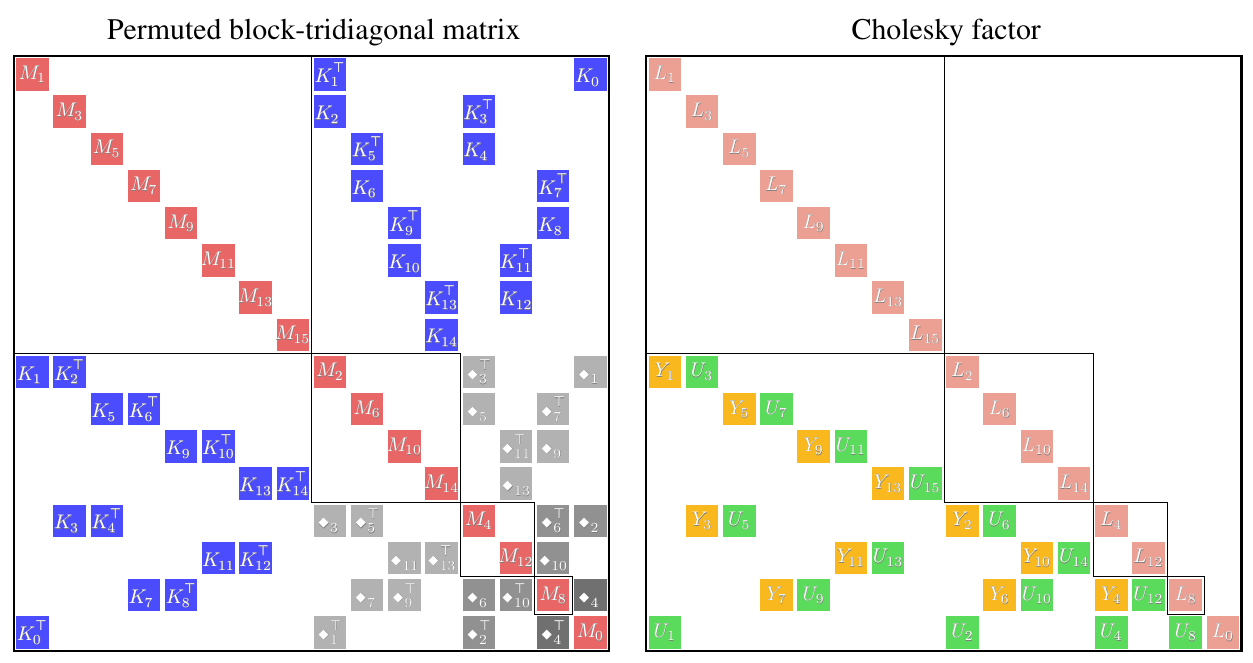}
    \end{minipage}%
    \begin{minipage}{0.23\textwidth}
        \centering
        \includegraphics[scale=0.43]{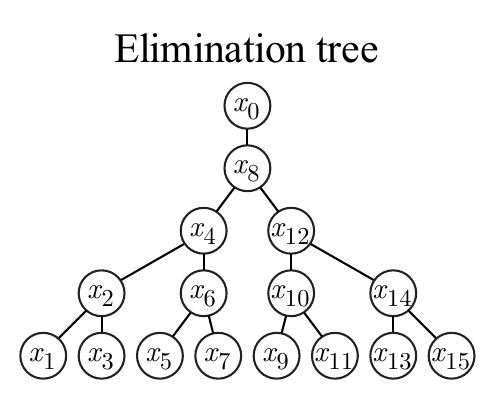}
    \end{minipage}%
    \caption{Top: a $16\times 16$ symmetric block-tridiagonal matrix with its block Cholesky factor and the corresponding elimination tree.
    No fill-in is incurred, but the single linear path in the elimination tree highlights the sequential nature of this factorization procedure.
    Bottom: CR as the factorization of a permutation of the same block-tridiagonal matrix,
    along with its Cholesky factor (with fill-in generated during the factorization shown in gray and labels matching the derivation in \Cref{subsec:cr-elim-deriv}).
    Blocks in the square outlines along the diagonal can be factorized simultaneously.
    The parallel branches in the elimination tree visualize the available parallelism during the Cholesky factorization of the permuted matrix.}
    \label{fig:cyclic-reduction-permutation}
\end{figure}

A tool that helps visualize the available parallelism in the block Cholesky factorization of a given sparse matrix is the
\textit{elimination tree} \cite[Ch.\,4]{davis_direct_2006}. Each node of such a tree represents a
variable in the system (or a block column in the matrix),
and the edges indicate the order in which the variables are eliminated:
a node can only be eliminated once all its children have been eliminated.
The elimination tree of the original block-tridiagonal matrix \eqref{eq:block-tri} is a linear graph: the variables have to be eliminated in order and one at a time.
In contrast, the elimination tree of the permuted matrix in \Cref{fig:cyclic-reduction-permutation}
is a perfect binary tree with many leaf nodes that can be eliminated in parallel. These leaf nodes are exactly the odd variables that are eliminated in the first step of CR.
Thanks to this parallelism, the total height of the elimination tree for the CR method ($\log_2(N) + 1$) is much lower than that of the original block-tridiagonal matrix ($N$).
This height difference indicates that the wall-clock time of an optimal parallel implementation of CR could be
significantly lower than for the block-tridiagonal Cholesky method.
Because of the additional fill-in, the asymptotic reduction in wall time of CR compared to the Cholesky factorization of the block-tridiagonal system for the case of $P=N$ processors
turns out to be $\frac{10 \log_2(N) \;+\; 1}{7 (N-1) \;+\; 1}$\footnote{Based on the theoretical FLOPs in the critical path,
and performing as much of \Cref{alg:cr} in parallel as possible, including the loops, as well as \cref{ln:cr-Y,ln:cr-U} and \cref{ln:cr-syrk,ln:cr-gemm}}.

\medskip
\section{{\normalfont\cyqlone{}}: Parallel factorization and solution of KKT systems with optimal control structure} \label{sec:cyqlone}
\begin{figure}[H]
    \begin{minipage}{0.74\textwidth}
        \includegraphics[width=1\textwidth]{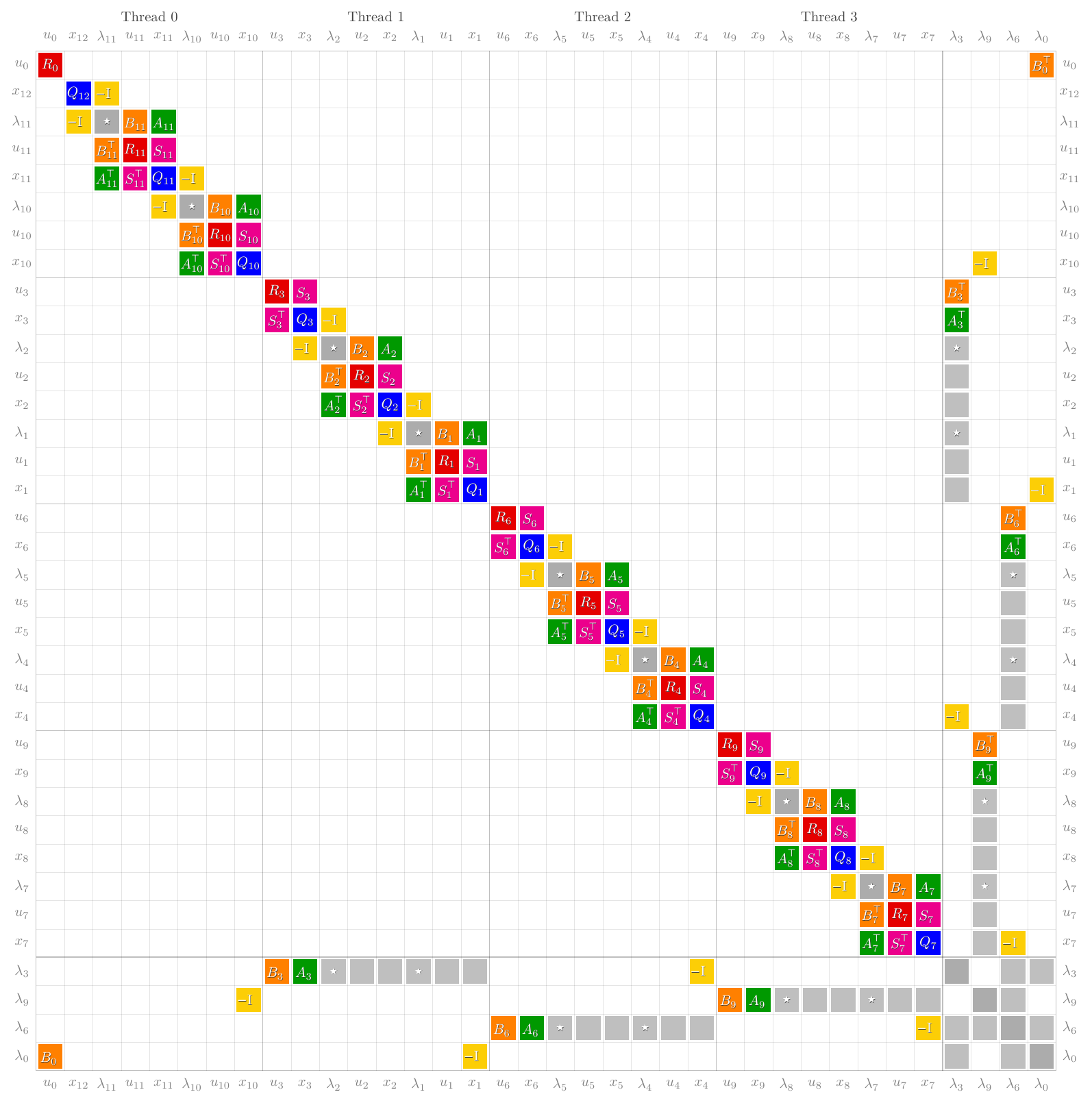}
    \end{minipage}%
    \begin{minipage}{0.26\textwidth}
        \centering
        \includegraphics[scale=0.44]{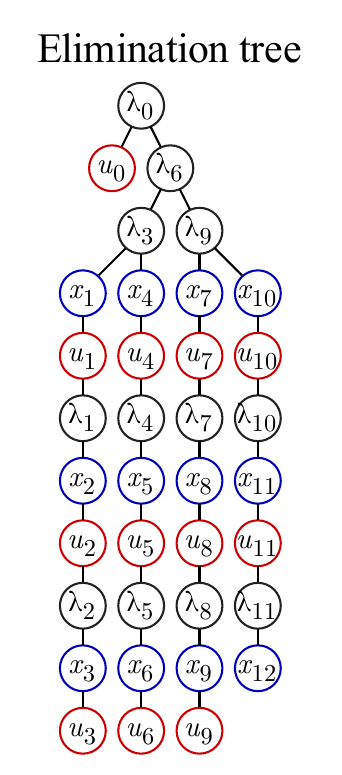}
    \end{minipage}
    \caption{\textsc{Cyqlone} ordering of the coefficient matrix of a KKT system with optimal control structure \eqref{eq:opt-cond-ocp-eq} with horizon $N=12$ for $P=4$ processors.
        Grey blocks indicate fill-in incurred during the factorization of the matrix.
        Stars ($\star$) mark structured or redundant fill-in (see \Cref{sec:mod-ricc}).
        The corresponding elimination tree visualizes the available parallelism during the factorization, with full utilization of all four processors throughout most of the process.}
    \label{mat:pdp-cr-N12-P4}
\end{figure}

This section describes a parallelizable algorithm for the factorization of a particular permutation of the
coefficient matrix of the KKT system with optimal control structure introduced in \eqref{eq:opt-cond-ocp-eq}.
We will refer to this ordering of the coefficient matrix and the corresponding factorization algorithm as the \textsc{Cyqlone} method.
\Cref{mat:pdp-cr-N12-P4} shows an illustrative example for an OCP with horizon $N=12$ and parallelization across $P=4$ processors.
This matrix visualization gives a clear idea of the computational cost (by considering the fill-in) and the parallelizability (using the elimination tree).
The specific permutations and the structure of the matrix for the general case will be discussed below.

\subsection{Matrix structure and high-level factorization procedure \done} \label{sec:mat-struc-hi-lev}
\newcommand\Perm[1]{\mathrm P_{\!#1}}

We partition the matrix visualized in \Cref{mat:pdp-cr-N12-P4} as follows:%
{\setstretch{1.25}
\begin{equation} \label{eq:block-matrix-example-4}
    (\text{Fig.~\ref{mat:pdp-cr-N12-P4}}):\;\; \mathscr K = \begin{pNiceArray}{c||c}[columns-width=1.8em]
        \mathscr R & \tp{\mathscr A} \\\hline\hline
        \mathscr A & 0
    \end{pNiceArray} =
    \begin{pNiceArray}{c|c|c|c||c}[columns-width=2.6em,last-row]
        \mathscr R_0         &                      &                      &                      & \tp{\mathscr A}_0 \tp{\Perm0} \\\hline
        & \mathscr R_1         &                      &                      & \tp{\mathscr A}_1 \tp{\Perm1} \\\hline
        &                      & \mathscr R_2         &                      & \tp{\mathscr A}_2 \tp{\Perm2} \\\hline
        &                      &                      & \mathscr R_3         & \tp{\mathscr A}_3 \tp{\Perm3} \\\hline\hline
        \Perm0 \mathscr A_0 & \Perm1 \mathscr A_1 & \Perm2 \mathscr A_2 & \Perm3 \mathscr A_3 & 0 \\
        c\!=\!0 & c\!=\!1 & c\!=\!2 & c\!=\!3 &
    \end{pNiceArray}.
\end{equation}}%
The structure of this matrix is based on the partitioning of the full horizon of length $N$ into $P$ smaller intervals of length $n\defeq N/P$.
Each block column $c \in \{0, 1, \dotsc, P - 1\}$
contains the states and controls of the stages in the interval $\{n c - n + 1,\; n c - n + 2,\, \dotsc,\; nc\}$ (in reverse order),
as well as the Lagrange multipliers corresponding to the dynamics constraints in the interior of the interval.
The purpose of this partitioning is to allow the different intervals to be solved in parallel on $P$ processors.

Block column $c=0$ is special, because it contains the controls $u^0$ of the first stage, combined with the variables of the last $n$ stages.
It is natural to handle this edge case uniformly by reducing the appropriate indices modulo $N$ and closing the horizon in a circular or periodic fashion. One could imagine
introducing matrices $S_0 = S_{12} = 0$ and $A_0 = A_{12} = 0$, ensuring that all block columns have the same structure.

The matrices $\mathscr R_c$ represent optimality conditions of smaller OCPs with horizon length $n$ on the different intervals of the partitioning.
The dynamics constraints describing the coupling between different intervals are found in the bottom row, described by matrices $\mathscr A_c$,
and they serve as initial and terminal state constraints for the smaller OCPs on each interval.
An odd--even permutation of these coupling constraints, represented by the permutation matrices $\Perm c$, is used to enable cyclic reduction of the bottom-right block, as introduced in \Cref{sec:cr-as-chol}:
the constraints corresponding to Lagrange multipliers
$\lambda^{ni}$ for odd $i$ are ordered first,
then the ones for which $i$ is divisible by 2 but not 4, then $i$ divisible by 4 but not 8, etc., with $\lambda^0$ last.

The \textsc{Cyqlone} algorithm for solving KKT systems with optimal control structure is derived by
systematically performing block Cholesky factorization of different blocks of \eqref{eq:block-matrix-example-4}, %
leading to a factorization of the full matrix $\mathscr K$ that can be used to solve \eqref{eq:opt-cond-ocp-eq} by forward and back substitution.
Specifically, we obtain a factorization of $\mathscr K$ by applying the straightforward block Cholesky procedure outlined in \Cref{sec:block-chol}:

\begin{enumerate}
    \setlength{\itemsep}{0.6em}
    \item Compute block Cholesky factorizations of the diagonal blocks \\$\mathscr R_c = \L{\mathscr R}_c \mathscr D_c \Lt{\mathscr R}_c$ (in parallel). \label{it:bchol-diag}
    \item Use these factorizations to solve $\L{\mathscr A}_c = \mathscr A_c \Lit{\mathscr R}_c \inv{\mathscr D}_c$ (in parallel). \label{it:solve-subdiag}
    \item Compute the Schur complement $\displaystyle \mathscr M = -\mathscr K / \mathscr R = -\smash{\sum_{c=0}^{P-1}} \Perm c \L{\mathscr A}_c \mathscr D_c \Lt{\mathscr A}_c \tp{\Perm c}$ \\(in parallel). \label{it:compute-schur}
    \item Compute a block Cholesky factorization of the Schur complement \\$\mathscr M = \L{\mathscr M} \Lt{\mathscr M}$ (using CR). \label{it:cr-schur}
\end{enumerate}

Each step makes use of the internal block structure of the matrices $\mathscr R_c$ and $\mathscr A_c$: For example, steps \ref{it:bchol-diag} and \ref{it:solve-subdiag} can be combined into a single Riccati-like recursion that exploits the optimal control structure;
many of the terms in step \ref{it:compute-schur} cancel out; and the structure of symmetric or triangular blocks is exploited during the individual block matrix operations.

Following \eqref{eq:blk-chol}, the full Cholesky factorization of $\mathscr K$ can then be written as
\begin{equation} \label{eq:chol-fac-K-cyqlone}
    \setstretch{1.2}
    \mathscr K =
    \begin{pNiceArray}{c||c}[columns-width=1.8em]
        \L{\mathscr R} & \\\hline\hline
        \L{\mathscr A} & \L{\mathscr M}
    \end{pNiceArray}
    \begin{pNiceArray}{c||c}[columns-width=1.8em]
        \mathscr D \\\hline\hline
        & \!\!-\I
    \end{pNiceArray}
    \ttp{
    \begin{pNiceArray}{c||c}[columns-width=1.8em]
        \L{\mathscr R} & \\\hline\hline
        \L{\mathscr A} & \L{\mathscr M}
    \end{pNiceArray}
    },
\end{equation}
where the matrices $\L{\mathscr R}$, $\L{\mathscr A}$ and $\mathscr D$ consist of the blocks $\L{\mathscr R}_c$, $\Perm c\L{\mathscr A}_c$ and $\mathscr D_c$, respectively.

\subsection{Modified Riccati recursion \rmrk{steps \ref{it:bchol-diag} and \ref{it:solve-subdiag}} \done} \label{sec:mod-ricc}

Consider one of the first four block columns of $\mathscr K$ in \eqref{eq:block-matrix-example-4} in isolation, for example the second block column ($c=1$),
corresponding to variables $(u_3, x_3, \lambda_2, \allowbreak u_2, x_2, \lambda_1, u_1, x_1)$.
Together with the corresponding initial and terminal constraint equations associated with Lagrange multipliers $\lambda_0$ and $\lambda_3$, this results in the following smaller block matrix:
\begin{align}
    \mathscr K_1 \defeq
    &\left(
    \begin{array}{c|c}
            \mathscr R_1 & \ttp{\mathscr A_1} \\\hline
            \mathscr A_1 & 0
        \end{array}
    \right) \\ 
    ={}
             & \scalebox{0.87}{\setstretch{1.25}$
            \begin{pNiceArray}{cccccccc|cc}[columns-width=3.4em,first-row]
                u^3      & x^3 & \lambda^2 & u^2      & x^2 & \lambda^1 & u^1      & x^1 & \lambda^3 & \lambda^0 \\
                R_3      & S_3 &           &          &     &           &          &     & \ttp B_3  &           \\
                \ttp S_3 & Q_3 & -\I       &          &     &           &          &     & \ttp A_3  &           \\
                & -\I & 0         & B_2      & A_2 &           &          &     &           & -\I       \\
                &     & \ttp B_2  & R_2      & S_2 &           &          &     &           &           \\
                &     & \ttp A_2  & \ttp S_2 & Q_2 & -\I       &          &     &           &           \\
                &     &           &          & -\I & 0         & B_1      & A_1 &           &           \\
                &     &           &          &     & \ttp B_1  & R_1      & S_1 &           &           \\
                &     &           &          &     & \ttp A_1  & \ttp S_1 & Q_1 &           &           \\\hline
                B_3      & A_3 &           &          &     &           &          &     & 0         &           \\
                &     &           &          &     &           &          & -\I &           & 0
            \end{pNiceArray}%
    $}\notag                                           \\
    \intertext{Performing block Cholesky factorization of $\mathscr R_1$ using \eqref{eq:blk-chol}, we have}
    \mathscr K_1 ={}
             & \left(
    \begin{array}{c}
            \L{\mathscr R}_1 \\\hline
            \L{\mathscr A}_1
        \end{array}
    \right)
    \mathscr D_1
    \ttp{
        \left(
        \begin{array}{c}
            \L{\mathscr R}_1 \\\hline
            \L{\mathscr A}_1
        \end{array}
        \right)
    }
    +
    \left(
    \begin{array}{c|c}
            0 & 0                           \\\hline
            0 & \mathscr K_1 / \mathscr R_1
        \end{array}
    \right),                                           \displaybreak[1] \\
    \L{\mathscr R}_1                                                \defeq{}
             & \scalebox{0.87}{\setstretch{1.25}$\left(
            \begin{NiceArray}{cccccccc}[columns-width=3.4em]
                    \L R_3 &           &                 &        &           &                 &        &        \\
                    \L S_3 & \L Q_3    &                 &        &           &                 &        &        \\
                           & -\Lit Q_3 & -\Lit Q_3       &        &           &                 &        &        \\
                           &           & \ttp B_2 \L Q_3 & \L R_2 &           &                 &        &        \\
                           &           & \ttp A_2 \L Q_3 & \L S_2 & \L Q_2    &                 &        &        \\
                           &           &                 &        & -\Lit Q_2 & -\Lit Q_2       &        &        \\
                           &           &                 &        &           & \ttp B_1 \L Q_2 & \L R_1 &        \\
                           &           &                 &        &           & \ttp A_1 \L Q_2 & \L S_1 & \L Q_1
                \end{NiceArray}
    \right),$} \label{eq:factor-R1}                                          \\
    \L{\mathscr A}_1 \defeq{} &\mathscr A_1 \Lit{\mathscr R}_1 \inv{\mathscr D}_1 \\
    ={}
             & \scalebox{0.87}{\setstretch{1.25}$\left(
            \begin{NiceArray}{cccccccc}[columns-width=3.34em]
                    \L{B}_3 & \!\!\Acl_3 \Lit Q_3\!\! & \!\!\Acl_3 \Lit Q_3\!\! & \L{B}_2 & \!\!\Acl_2 \Lit Q_2\!\! & \!\!\Acl_2 \Lit Q_2\!\! & \L{B}_1 & \L{A}_1   \\
                            &                         &                         &         &                         &                         &         & -\Lit Q_1
                \end{NiceArray}
            \right),$} \label{eq:LA-single-col-example-4}
    \\
    \mathscr D_1 \defeq
    \blkdiag &\mkern1mu \scalebox{0.87}{\setstretch{1.8}$\left(
            \begin{NiceArray}{cccccccc}[columns-width=3.41em]
                    \I & \I & -\I & \I & \I & -\I & \I & \I
                \end{NiceArray}
            \right).$}
\end{align}
The matrices $\L R_j$, $\L S_j$, and $\L Q_j$, which make up $\L{\mathscr R}_1$ in \eqref{eq:factor-R1}, are computed using the block Cholesky factorization procedure described
in \Cref{sec:block-chol}.
Recall that there is some freedom in the choice of the factors of the diagonal blocks.
Thanks to the special structure of $\mathscr R_1$, the following factorization can be used for the blocks corresponding to the states and the Lagrange multipliers:
\begin{equation} \label{eq:redundant-fact}
    \scalebox{1}{$
    \begin{pmatrix}
        P & -\I \\ -\I & \phantom{-}0
    \end{pmatrix}$} =
    \scalebox{1}{$\begin{pmatrix}
        L \\ -\invtp L & -\invtp L
    \end{pmatrix}$}
    \scalebox{1}{$
    \begin{pmatrix}
        \I \\ & -\I
    \end{pmatrix}$}
    \scalebox{1}{$
    \begin{pmatrix}
        \tp L & -\inv L \\ & -\inv L
    \end{pmatrix}$}, \quad \text{where} \; P = L \tp L.
\end{equation}
Note how the bottom blocks $\invtp L$ do not need to be computed or stored:
Storing $L$ is sufficient, because multiplication
by $\invtp L$ (required during the solution of the linear system)
can be implemented as back substitution using $L$;
its inverse does not need to be formed explicitly.
The factorization in \eqref{eq:redundant-fact} is also responsible for the duplicate $\Acl_j \Lit Q_j$ blocks in $\L{\mathscr A}_1$.
Since these redundant blocks do not have to be computed or stored explicitly, we denoted them using a star in \Cref{mat:pdp-cr-N12-P4}.
The optimized block Cholesky factorization of $\mathscr R_1$
is described by \Cref{alg:fact-mod-riccati} (invoked with $c=1$, $N=12$ and $P=4$ for the example considered here).
By using \eqref{eq:redundant-fact}, this procedure fully exploits the optimal control structure of $\mathscr R_1$, and it is similar%
\,\footnote{In the case where $P=1$, we have $A_0=A_N = 0 = S_N=S_0$, and \cref{ln:riccati-Acl,ln:riccati-mul-Acl} in \Cref{alg:fact-mod-riccati} can be skipped since $\Acl_j$ and $\L S_j$ are zero as well. The resulting algorithm corresponds to \cite[Alg.\,3]{frison_efficient_2013} almost exactly.}
to the factorized Riccati recursion described in \cite[Alg.~3]{frison_efficient_2013}.
It is also closely related to backward dynamic programming (see \Cref{sec:riccati-permutation}).
Additionally, \Cref{alg:fact-mod-riccati} returns the blocks $\L {B}_j$, $\Acl_j$ and $\L {A}_1$ of $\L {\mathscr A}_1$. Merging the evaluation of $\L{\mathscr R}_1$ and $\L{\mathscr A}_1$ results in improved locality and cache performance.
The products ${\Acl_j \Lit Q_j}$ are not carried out explicitly, except for the very first stage $j_1$ in each interval.

To obtain an alternative block Cholesky factorization of $\mathscr R_1$,
one could use the trivial block Cholesky factorization of the entire $2\times2$
block of states and Lagrange multipliers, leading to an equivalent algorithm that is closer to the asymmetric implementation of the Riccati recursion from \cite[Alg.\,1]{frison_efficient_2013}:
\begin{equation}
    \scalebox{1}{$\begin{pmatrix}
        P & -\I \\ -\I & \phantom{-}0
    \end{pmatrix}$} =
    \underbrace{
        \scalebox{1}{$\begin{pmatrix}
            \I & \\ & \phantom{-}\I
        \end{pmatrix}$}}_{L}
    \underbrace{
        \scalebox{1}{$\begin{pmatrix}
            P & -\I \\ -\I & \phantom{-}0
        \end{pmatrix}$}}_{D}
    \underbrace{
        \scalebox{1}{$\begin{pmatrix}
            \I & \\ & \phantom{-}\I
        \end{pmatrix}$}}_{\ttp L},
    \;\;\text{with}\;\;
    \inv{\scalebox{1}{$\begin{pmatrix}
            P & -\I \\ -\I & \phantom{-}0
        \end{pmatrix}$}}\!\!
    = \scalebox{1}{$\begin{pmatrix}
        0 & -\I \\ -\I & -P
    \end{pmatrix}$}.
\end{equation}

\begin{algorithm2e}[htbp]
    \def\bbwd{b_\mathrm{bwd}}
    \def\bfwd{b_\mathrm{fwd}}
    \caption{Factorization of a single modified Riccati block column}
    \label{alg:fact-mod-riccati}
    \DontPrintSemicolon
    \KwIn{$A, B, Q, R, S$: OCP data matrices}
    \KwIn{$N$: horizon length}
    \KwIn{$P$: number of intervals to partition the horizon into}
    \KwIn{$c \in \N_{[0, P)}$: index of the block column to factor (selects the interval)}
    \KwOut{$\L{\mathscr R}_c$ ($\L{R}, \L{S}, \L{Q}$), $\L{\mathscr A}_c$ ($\L{B}, \L{A}, \Acl$): blocks of the Cholesky factor of $\mathscr K_c$}
    \vspace{0.3em}
    \Fn{\normalfont\textsc{factor--block--column--riccati}$(c)$}{
    $n = N / P$\mycommentnofill{\rmrknp{Number of stages per interval}}\;
    $j_{1} = n(c-1)+1$,\quad $j_{n} = nc$\mycommentnofill{\rmrknp{First/last stage indices in block column $c$}}\;
    $\begin{pmatrix}
            \hat B_{j_{n}} & \Aprop_{j_{n}}
        \end{pmatrix}
        \assignpc \begin{pmatrix}
            B_{j_{n}} & A_{j_{n}}
        \end{pmatrix}$\;
    $\begin{pmatrix}
            \hat R_{j_{n}}       & \hat S_{j_{n}} \\
            \tp {\hat S_{j_{n}}} & \hat Q_{j_{n}}
        \end{pmatrix} \assignpc \begin{pmatrix}
            R_{j_{n}}       & S_{j_{n}} \\
            \tp {S_{j_{n}}} & Q_{j_{n}}
        \end{pmatrix}$\;

    \For{$j = j_{n}, j_{n} - 1, j_{n} - 2, \dots, j_{1}$}{
        $\begin{pmatrix}
                L^{\!R}_{j}               \\
                L^{\!S}_{j} & L^{\!Q}_{j} \\
            \end{pmatrix} \assignpc \operatorname{chol} \begin{pmatrix}
                \hat R_{j}       & \hat S_{j} \\
                \tp {\hat S_{j}} & \hat Q_{j}
            \end{pmatrix}$
        \hspace{-5pt}
        \mycomment{\texttt{\rmrk{potrf}} \quad $\hspace{-0.45em}\quad\mathclap{\tfrac16 n_{ux}^3}\quad$}
        $L^{\!B}_{j} \assignpc \hat B_{j} L^{\!R\; -\!\top}_{j}$
        \mycomment{\texttt{\rmrk{trsm}} \quad $\quad\mathclap{\tfrac12 n_{x} n_{u}^2}\quad$}
        $\Acl_{j} \assignpc \Aprop_{j} - L^{\!B}_{j} L^{\!S\,\top}_{j}$ \label{ln:riccati-Acl}
        \mycomment{\texttt{\rmrk{gemm}} \quad $\quad\mathclap{n_x^2 n_u}\quad$}
        \vspace{0.3em}
        \If{$j > j_{1}$\mycommentnofill{\rmrknp{Propagate dynamics and cost-to-go to the previous stage}}}{
            $\begin{pmatrix}
                    \hat B_{j-1} & \Aprop_{j-1}
                \end{pmatrix}
                \assignpc \Acl_{j} \begin{pmatrix}
                    B_{j-1} & A_{j-1}
                \end{pmatrix}$ \label{ln:riccati-mul-Acl}
            \mycomment{\texttt{\rmrk{gemm}} \quad $\quad\mathclap{n_{ux} n_{x}^2}\quad$}
            $V_{j-1} \assignpc \begin{pmatrix}
                    \ttp B_{j-1} \\ \ttp A_{j-1}
                \end{pmatrix} L^{\!Q}_{j}$
            \mycomment{\texttt{\rmrk{trmm}} \quad $\quad\mathclap{\tfrac12 n_{ux} n_{x}^2}\quad$}
            $\begin{pmatrix}
                    \hat R_{j-1}       & \hat S_{j-1} \\
                    \tp {\hat S_{j-1}} & \hat Q_{j-1}
                \end{pmatrix} \assignpc \begin{pmatrix}
                    R_{j-1}       & S_{j-1} \\
                    \tp {S_{j-1}} & Q_{j-1}
                \end{pmatrix} + V_{j-1}\tp V_{j-1}$\mycomment{\texttt{\rmrk{syrk}} \quad $\quad\mathclap{\tfrac12 n_{ux}^2 n_{x}}\quad$}
        }
    }
    $L^{\!A}_{j_{1}} \assignpc \Acl_{j_{1}} L^{\!Q\;-\!\top}_{j_{1}} $\mycommentnofill{\rmrknp{Only the first stage needs $\L{A}$}} \mycomment{\texttt{\rmrk{trsm}} \quad $\quad\mathclap{\tfrac12 n_{x}^3}\quad$}
    }
    \vspace{0.3em}
    \setstretch{0.8}
    To avoid unintelligible index manipulations in the pseudocode, we use the\;
    convention that stage indices are reduced modulo $N$, e.g. $R_{j}$
    represents $R_{j\,\mathbin{\%}N}$.\;
    The indices of $Q_j$ and $\smash{\L Q_j}$ are
    one-based because of the elimination of $x^0$,\;
    hence $Q_{j}$ represents $Q_{((j-1)\mathbin{\%}N)+1}$.
    Furthermore, $A_0 = 0$ and $S_0 = 0$.\;
\end{algorithm2e}

\subsection{Computation of the Schur complement \rmrk{step \ref{it:compute-schur}} \done} \label{sec:compute-schur}

Once $\L{\mathscr A}_1$ has been computed by \Cref{alg:fact-mod-riccati}, evaluating the contribution of each block column to the Schur complement is done using a straightforward structured symmetric matrix product:
\vspace{0pt}
\begin{align}
    -\mathscr K_1 / \mathscr R_1 ={} &\L{\mathscr A}_1 \mathscr D_1 \Lt{\mathscr A}_1    \notag\\
                                   ={} & \scalebox{0.87}{\setstretch{1.25}$\left(
            \begin{NiceArray}{cc}[columns-width=3.4em]
                    \L{B}_3 \Lt{B}_3 + \L{B}_2 \Lt{B}_2 + \L{B}_1 \Lt{B}_1 + \L{A}_1 \Lt{A}_1 & -\L{A}_1 \Li Q_1 \\
                    -\Lit Q_1 \Lt{A}_1                                                        & \Lit Q_1 \Li Q_1
                \end{NiceArray}
            \right)$}
    \defeq \scalebox{0.87}{\setstretch{1.25}$\left(
            \begin{NiceArray}{cc}[columns-width=2.5em]
                    \Mfwd_1 & \;\KbwdT_1 \\
                    \Kbwd_1  & \Mbwd_0
                \end{NiceArray}
            \right)$}
    \label{eq:schur-compl-single-col-example-4-c1}
    \intertext{Notice how the blocks involving $\Acl_j \Lit{Q}_j$ from \eqref{eq:LA-single-col-example-4} cancel out for all but the first stage.
    Arrows are used to indicate the direction of the coupling to other block columns.
    For example, $\Mfwd_1$ is the contribution to the Schur complement from the constraints representing the coupling forward in time, from the second to the third block columns ($c=1$ to $c=2$),
    and $\Mbwd_0$ is the contribution of the coupling back in time, from the second to the first block column ($c=1$ to $c=0$). The updated forward coupling matrices $\Mfwd_i$
    are derived from the matrices $A_j$ and $B_j$ in the eliminated dynamics constraints between intervals,
    whereas the backward coupling matrices $\Mbwd_i$ are derived from the corresponding matrix $E_{j+1}=\I$.
    Because of the odd--even structure of CR, the order of the block rows in $\mathscr A_c$ is reversed for even $c$. For the third block column ($c=2$):}
    -\mathscr K_2 / \mathscr R_2= {} & \scalebox{0.87}{\setstretch{1.25}$\left(
            \begin{NiceArray}{cc}[columns-width=3.4em]
                    \Lit Q_4 \Li Q_4 & -\Lit Q_4 \Lt{A}_4                                                        \\
                    -\L{A}_4 \Li Q_4 & \L{B}_6 \Lt{B}_6 + \L{B}_5 \Lt{B}_5 + \L{B}_4 \Lt{B}_4 + \L{A}_4 \Lt{A}_4
                \end{NiceArray}
            \right)$}
    \defeq \scalebox{0.87}{\setstretch{1.25}$\left(
            \begin{NiceArray}{cc}[columns-width=2.5em]
                    \Mbwd_1  & \;\KfwdT_1 \\
                    \Kfwd_1 & \Mfwd_2
                \end{NiceArray}
            \right)$}
    \label{eq:schur-compl-single-col-example-4-c2}
    \intertext{The indices of the subdiagonal blocks $\Kbwd_i$ and $\Kfwd_i$ are chosen to match the ones of the diagonal blocks in the same column, as this simplifies the notation in the upcoming CR step.
    For completeness, we also list the Schur complement matrices for the first and last block columns. The matrix $\Kfwd_3$ is zero because the original
    OCP has no terminal constraints ($A_{12} = 0$). Blocks that are zero are shown in gray and have been included to demonstrate the regular structure of the factorization.
    A straightforward generalization to problems with coupling between the first and last stage is possible, but is not considered here for the sake of brevity.}
    -\mathscr K_0 / \mathscr R_0= {} & \scalebox{0.87}{\setstretch{1.25}$\left(
            \begin{NiceArray}{cc}[columns-width=3.4em]
                    \Lit Q_{10} \Li Q_{10} & {\color{gray} -\Lit Q_{10} \L{A}_{10}}                                                                                                                            \\
                    {\color{gray} -\Lt{A}_{10} \Li{Q}_{10}}                      & \L{B}_{0} \Lt{B}_{0} + {\color{gray}\L{B}_{11} \Lt{B}_{11} + \L{B}_{10} \Lt{B}_{10} + \L{A}_{10} \Lt{A}_{10}}
                \end{NiceArray}
            \right)$}
    \defeq \scalebox{0.87}{\setstretch{1.25}$\left(
            \begin{NiceArray}{cc}[columns-width=2.5em]
                    \Mbwd_3                & {\color{gray}\;\KfwdT_3} \\
                    {\color{gray}\Kfwd_3} & \Mfwd_0
                \end{NiceArray}
            \right)$}
    \label{eq:schur-compl-single-col-example-4-c0}                    \\
    -\mathscr K_3 / \mathscr R_3= {} & \scalebox{0.87}{\setstretch{1.25}$\left(
            \begin{NiceArray}{cc}[columns-width=3.4em]
                    \L{B}_9 \Lt{B}_9 + \L{B}_8 \Lt{B}_8 + \L{B}_7 \Lt{B}_7 + \L{A}_7 \Lt{A}_7\;\;\;\,{} & -\L{A}_7 \Li Q_7 \\
                    -\Lit Q_7 \Lt{A}_7                                                        & \Lit Q_7 \Li Q_7
                \end{NiceArray}
            \right)$}
    \defeq \scalebox{0.87}{\setstretch{1.25}$\left(
            \begin{NiceArray}{cc}[columns-width=2.5em]
                    \Mfwd_3 & \;\KbwdT_3 \\
                    \Kbwd_3  & \Mbwd_2
                \end{NiceArray}
            \right)$}
    \label{eq:schur-compl-single-col-example-4-c3}
\end{align}
Combining the contributions from all four block columns in \eqref{eq:schur-compl-single-col-example-4-c1}--\eqref{eq:schur-compl-single-col-example-4-c3}, the full Schur complement reads:
\begin{equation}
    \mathscr M = -\mathscr K / \mathscr R ={} \scalebox{0.87}{\setstretch{1.5}$
            \begin{pNiceArray}{cccc}[columns-width=3.4em,first-row,first-col]
                &\shortstack{$\;\lambda^3$\\$\scriptsize i=1$}& \shortstack{$\;\lambda^9$\\$\scriptsize i=3$} & \shortstack{$\;\lambda^6$\\$\scriptsize i=2$} & \shortstack{$\;\lambda^0$\\$\scriptsize i=0$} \\
                \lambda^3 \;\; & \Mbwd_1 + \Mfwd_1 &                                                                    & \KfwdT_1                              & \KbwdT_1                               \\
                \lambda^9 \;\; & & \Mbwd_3 + \Mfwd_3 & \KbwdT_3                               & {\color{gray}\KfwdT_3}                                                                   \\
                \lambda^6 \;\; & \Kfwd_1                                   & \Kbwd_3                                    & \Mbwd_2 + \Mfwd_2 &                                                                    \\
                \lambda^0\;\; & \Kbwd_1                                    & {\color{gray}\Kfwd_3}                                                                   &                                                                    & \Mbwd_0 + \Mfwd_0  \\
            \end{pNiceArray}$}.
    \label{eq:schur-compl-example-4}
\end{equation}
In the following sections and in the pseudocode, we use indices $j$ to refer to the original OCP stages, and indices $i$ to refer to the block columns of the (permuted) Schur complement \eqref{eq:schur-compl-example-4}.
Column $i$ corresponds to the Lagrange multipliers $\lambda^{ni}$.
The rows and columns of the Schur complement are permuted to enable CR: they are ordered by increasing 2-adic valuation $\nu_2(i)$ (see \Cref{sec:cr-as-chol}).
Consequently, the off-diagonal blocks $\Kbwd_i$ (representing coupling back in time)
are found in column $i$ and row $i - 2^{\nu_2(i)}$;
and blocks $\Kfwd_i$ (coupling forward in time) are found in column $i$ and row $i + 2^{\nu_2(i)}$.

The construction of $\mathscr M$ can be fully parallelized, assuming that the final additions of the two terms in the diagonal blocks are synchronized correctly.
A possible parallel procedure for constructing $\mathscr M$ is given in the \textsc{compute--schur} function of \Cref{alg:fact-parallel-riccati-cr}.
Since the diagonal blocks will be updated during the CR algorithm in the following step, we use an explicit superscript
to indicate the level $l$ in the CR recursion (as in \Cref{alg:cr}), with initial values
$\Mdiag^{(0)}_i \defeq \Mbwd_i + \Mfwd_i$.

\subsection{Factorization of the Schur complement \rmrk{step \ref{it:cr-schur}} \done}

Thanks to the odd--even ordering of the Lagrange multipliers in \Cref{mat:pdp-cr-N12-P4}, cyclic reduction of the Schur complement $\mathscr M$ corresponds directly to its block Cholesky factorization.
The resulting block Cholesky factor $\L{\mathscr M}$ is a sparse matrix with at most two subdiagonal blocks per block column, with a particular structure as visualized in \Cref{fig:cyclic-reduction-permutation}.
For the specific case of \Cref{mat:pdp-cr-N12-P4} with only $P=4$ processors,
we have:
\vspace{-1.4em}
\begin{align} \label{eq:chol-fac-schur-compl-M}
    \L{\mathscr M} ={}
     & \scalebox{0.87}{\setstretch{1.5}$
            \begin{pNiceArray}{cccc}[columns-width=3.4em,first-row,first-col]
                &\lambda^3 & \lambda^9 & \lambda^6 & \lambda^0 \\
                \lambda^3 \;\; & L_1 & & & \\
                \lambda^9 \;\; & & L_3 & & \\
                \lambda^6 \;\; & Y_1 & U_3 & L_2 & \\
                \lambda^0\;\; & U_1 & {\color{gray} Y_3} & U_2 & L_0 \\
            \end{pNiceArray}.$}
\end{align}
A procedure for the factorization of $\mathscr M$ is given in \Cref{alg:fact-parallel-riccati-cr}, which is a specialization of \Cref{alg:cr}.
The main difficulty here is keeping track of the appropriate indices. Studying the indices in \Cref{fig:cyclic-reduction-permutation} may be helpful to better understand the pseudocode.
To clarify the structure of the algorithm, we note that at level $l$ in the CR recursion, the columns $i$ for which $\nu_2(i)=l$ are eliminated.
For these columns, the diagonal blocks $\Mdiag^{(l)}_i = L_i \ttp L_i$ are factorized first, and then the two subdiagonal blocks
$U_i=\Kbwd_i \invtp L_i$ and $Y_i=\Kfwd_i \invtp L_i$ are computed.
Finally, the Schur complement of the eliminated columns is computed by subtracting products of $U_i$ and $Y_i$ from the bottom-right corner:
The symmetric products $U_i \ttp U_i$ and $Y_i \ttp Y_i$ are subtracted from the diagonal blocks $\Mdiag_{i-2^l}^{(l)}$ and $\Mdiag_{i+2^l}^{(l)}$ respectively,
labeling the updated matrices $\smash{\Mdiag_{i-2^l}^{(l+1)}}$ and $\smash{\Mdiag_{i+2^l}^{(l+1)}}$; and
the asymmetric products $U_i \ttp Y_i$ and $Y_i \ttp U_i$ result in off-diagonal fill-in, labeled $\Kbwd_{i+2^l}$ and $\Kfwd_{i-2^l}$.
Since the matrix is symmetric, only the diagonal and subdiagonal blocks are computed.
The stride between the indices of successive rows/columns of $\mathscr M$ and $\L{\mathscr M}$ at level $l$ is $2^l$,
so index offsets of $\pm2^l$ can be interpreted as the next/previous row/column of the reduced system.

\begin{figure}
    \centering
    \includegraphics[width=\textwidth]{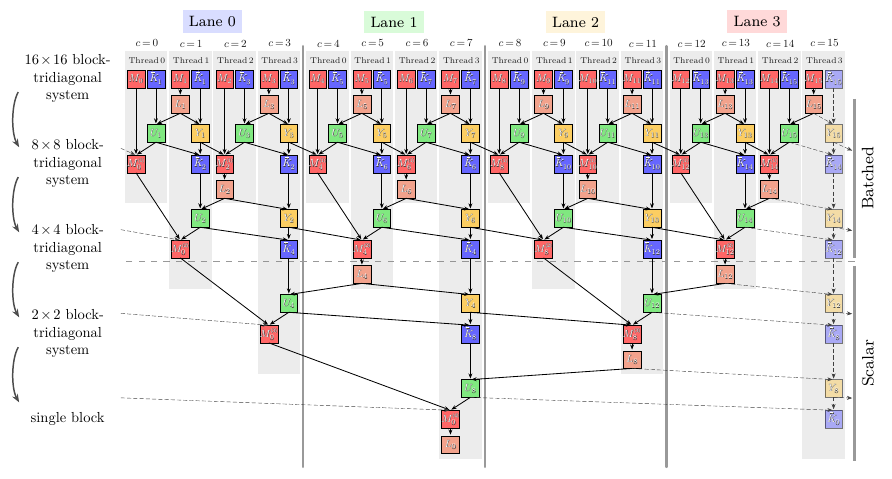}
    \caption{Graph representation of the different steps required for the cyclic reduction of a symmetric $16\times 16$
        block-tridiagonal matrix.
        The 16 block columns of the original matrix are distributed across 4 different threads and 4 vector lanes ($P=16$).
        Each square node represents one of the blocks of the block-tridiagonal matrix or its Cholesky factor;
        nodes including a superscript are blocks of the intermediate block-tridiagonal systems of reduced size from \eqref{eq:cr-even}.
        Node colors match the colors in \Cref{alg:fact-parallel-riccati-cr}.
        Edges represent data dependencies between steps of the factorization. For instance, following \Cref{subsec:cr-elim-deriv}, $U_1 \defeq \Kbwd_0 \invtp L_1$, so there are edges from $\Kbwd_0$ and $L_1$ to $U_1$.
        The vertical axis represents the wall time.
        The horizontal position of each block indicates the thread and vector lane where it is computed (variable $c$ in \Cref{alg:fact-parallel-riccati-cr}).
        Notice the recursive, self-similar structure of the method: the operations performed at each level are the same, just spread out horizontally and executed by half as many threads than the level before it.
        Dashed arrows complete the structure, and are required when solving a system with circular coupling $\Kfwd_{15} \neq 0$ or when extending the figure to larger matrix sizes.
        Vector lane assignment, the distinction between scalar and batched levels, and other vectorization aspects
        will be discussed in \Cref{sec:vectorization}.
    }
    \label{fig:cyclic-reduction-graph-16}
\end{figure}

A visual representation of the steps of the CR algorithm is shown in
\Cref{fig:cyclic-reduction-graph-16} (with colors of the different matrices matching the colors in the margin of \Cref{alg:fact-parallel-riccati-cr}).
Further details about the thread scheduling and vector lane assignment will be discussed in \Cref{sec:vectorization}.

\begin{algorithm2e}[htbp]
    \def\bbwd{i_{\scriptveryshortleftarrow}}
    \def\bfwd{i_{\scriptveryshortrightarrow}}
    \def\spaceif{\hspace{0.4em}}
    \def\spaceflop{\hspace{0.1em}}
    \def\spaceflopright{\hspace{-1em}}
    \caption{\cyqlone{} factorization}
    \label{alg:fact-parallel-riccati-cr}
    \DontPrintSemicolon
    \KwIn{$A, B, Q, R, S$ and $N$: OCP data matrices and horizon length}
    \KwIn{$P$: number of intervals to partition the horizon into}
    \KwOut{$\L{\mathscr R}$ ($\L{R}, \L{S}, \L{Q}$), $\L{\mathscr A}$ ($\L{B}, \L{A}, \Acl$), $\L{\mathscr M}$ ($L, Y, U$): blocks of the Cholesky factor of $\mathscr K$}
    \vspace{0.3em}
        \For{$c=0,...,P-1$ {\normalfont (in parallel)}}{\label{ln:loop-processors}
        \textsc{factor--block--column--riccati}$(c)$\mycomment{\rmrk{steps \ref{it:bchol-diag} and \ref{it:solve-subdiag}}}
        \textsc{compute--schur}$(c)$\mycomment{\rmrk{step \ref{it:compute-schur}}}
        \textsc{factor--schur}$(c)$\mycomment{\rmrk{step \ref{it:cr-schur}}}
    }
    \vspace{0.3em}
    \Fn{\normalfont\textsc{compute--schur}$(c)$\mycommentnoline{Compute part of the Schur complement $\mathscr M$}}{
        \tikzmark{xmark}
        $n = N / P$\mycommentnofill{\rmrknp{Number of stages per interval}}\;
        $j_{1} = n(c-1)+1$,\quad $j_{n} = nc$\mycommentnofill{\rmrknp{First/last stage indices in block column $c$}}\;
        $\bfwd = c$,\quad %
        $\bbwd = c - 1$\mycommentnofill{\rmrknp{Row indices of forward/backward coupling in $\Perm{c} \mathscr A_c$}}\; %
        $T_c \assignpc L^{\!Q\;-\!\top}_{j_{1}}$\mycomment{\texttt{\rmrk{trtri}\hspace{2pt}} \spaceflop $\tfrac16 n_{x}^3$\spaceflopright}
        \vspace{0.1em}
        \tikzmark{updateB0}
        \textbf{if} $\nu^P_2(\bbwd) > \nu^P_2(\bfwd)$: \spaceif
            $\Kbwd_{\bfwd} \assignpc -T_c L^{\!A\; \top}_{j_{1}}$
            \spaceif \textbf{else}: \spaceif
            $\Kfwd_{\bbwd} \assignpc -L^{\!A}_{j_{1}} \tp T_c$\!\!
            \tikzmark{updateB1}
            \ColorBar{colorBd}{updateB0}{updateB1}
            \mycomment{\texttt{\rmrk{trmm}\hspace{2pt}\phantom{i}} \spaceflop $\tfrac12 n_{x}^3$\spaceflopright}
        --- sync --- \mycommentnofill{\rmrknp{Wait for $T_{c+1}$}}\;
        \tikzmark{updateA0}
        $\Mbwd_{c} \assignpc T_{c+1}\tp T_{c+1}$
        \tikzmark{updateA1}
        \ColorBar{colorAd}{updateA0}{updateA1}
        \mycomment{\texttt{\rmrk{lauum}\hspace{2pt}} \spaceflop $\tfrac16 n_{x}^3$\spaceflopright}
        \vspace{0.15em}
        \tikzmark{updateA0}
        $\Mfwd_{c} \assignpc W \ttp W$\qquad where $W = \begin{pmatrix}
            L^{\!B}_{j_{n}} \;\cdots\;\; L^{\!B}_{j_{1}} \;\; L^{\!A}_{j_{1}}
        \end{pmatrix}$
        \tikzmark{updateA1}
        \ColorBar{colorAd}{updateA0}{updateA1}
        \mycomment{\texttt{\rmrk{syrk}} \spaceflop $\tfrac n2 n_{u} n_{x}^2 + \tfrac12 n_{x}^3$\spaceflopright}
        $\Mdiag^{(0)}_{c} \assignpc \Mbwd_{c} + \Mfwd_{c}$\;
    }
    \vspace{0.3em}
    \Fn{\normalfont\textsc{factor--schur}$(c)$\mycommentnoline{Cyclic reduction of the Schur complement $\mathscr M$}}{
        \tikzmark{xmark}
        \tikzmark{factorL0}
        \textbf{if} $\nu^P_2(c) = 0$: \spaceif
            $L_{c} \assignpc \operatorname{chol} \big(\Mdiag^{(0)}_{c}\big)$
            \tikzmark{factorL1}
            \ColorBar{colorLd}{factorL0}{factorL1}
            \mycomment{\texttt{\rmrk{potrf}} \spaceflop $\tfrac16 n_{x}^3$\spaceflopright}
        \For{$l = 0,...,\log_2(P) - 1$\mycommentnofill{\rmrknp{Recursion level of CR}}} {
            $i_U = c + 1$, \quad
            $i_Y = c + 1 - 2^l$
            \;
            --- sync --- \mycommentnofill{\rmrknp{Wait for $L$}}\;
            \tikzmark{solveU0}
            \textbf{if}\phantom{\textbf{el}} $\nu^P_2(i_U) = l$: \spaceif
                $U_{i_U} \assignpc \Kbwd_{i_U} \invtp L_{i_U}$
                \tikzmark{solveU1}
                \ColorBar{colorUd}{solveU0}{solveU1}
                \mycomment{\texttt{\rmrk{trsm}} \spaceflop $\tfrac12 n_{x}^3$\spaceflopright}
            \vspace{0.1em}
            \tikzmark{solveY0}
            \textbf{elif} $\nu^P_2(i_Y) = l$: \spaceif
                $Y_{i_Y} \assignpc \Kfwd_{i_Y} \invtp L_{i_Y}$
                \tikzmark{solveY1}
                \ColorBar{colorYd}{solveY0}{solveY1}
                \mycomment{\texttt{\rmrk{trsm}} \spaceflop $\tfrac12 n_{x}^3$\spaceflopright}
            --- sync --- \mycommentnofill{\rmrknp{Wait for $U, Y$}}\;
            \textbf{if}\phantom{\textbf{el}} $\nu^P_2(i_U) = l$: \spaceif
                \textsc{factor--L}$(l, i_Y)$\;
            \textbf{elif} $\nu^P_2(i_Y) = l$: \spaceif
                \textsc{update--K}$(l, i_Y)$\;
        }
    }
    \vspace{0.3em}
    \Fn{\normalfont\textsc{factor--L}$(l, i)$\mycommentnoline{Update \& factorize $\Mdiag^{(l+1)}_{i}$ after elimination of level $l$}}{
        \tikzmark{xmark}
        $i_U = i + 2^l$, \quad $i_Y = i - 2^l$\mycommentnofill{\rmrknp{Column indices of $U_{i_U}$/$Y_{i_Y}$ in row $i$ and level $l$}}\;
        \tikzmark{updateA0}
        $\Mdiag^{(l+1)}_{i} \assignpc \Mdiag^{(l)}_{i} - U_{i_U} \ttp U_{i_U} - Y_{i_Y} \ttp Y_{i_Y}$
        \tikzmark{updateA1}
        \ColorBar{colorAd}{updateA0}{updateA1}
        \mycomment{\texttt{\rmrk{$2\times$syrk}} \spaceflop $n_{x}^3$\spaceflopright}
        \tikzmark{factorL0}
        \textbf{if} $\nu^P_2(i) = l + 1$: \spaceif
            $L_{i} \assignpc \operatorname{chol} \big(\Mdiag^{(l+1)}_{i}\big)$
            \tikzmark{factorL1}
            \ColorBar{colorLd}{factorL0}{factorL1}
            \mycomment{\texttt{\rmrk{potrf}} \spaceflop $\tfrac16 n_{x}^3$\spaceflopright}
            \vspace{0.2em}
    }
    \vspace{0.3em}
    \Fn{\normalfont\textsc{update--K}$(l, i)$\mycommentnoline{Compute fill-in from elimination of column $i$ of $\mathscr M$}}{
        \tikzmark{xmark}
        $\bbwd = i - 2^l$, \quad
        $\bfwd = i + 2^l$\mycommentnofill{\rmrknp{Row indices of $U_i$/$Y_i$ in column $i$}}\;
        \tikzmark{updateB0}
        \textbf{if} $\nu^P_2(\bbwd) > \nu^P_2(\bfwd)$: \spaceif
            $\Kbwd_{\bfwd} \assignpc -U_{i} \ttp Y_{i}$ \spaceif
            \textbf{else}: \spaceif
            $\Kfwd_{\bbwd} \assignpc -Y_{i} \ttp U_{i}$
            \tikzmark{updateB1}
            \ColorBar{colorBd}{updateB0}{updateB1}
            \hspace{-0.5em}
            \mycomment{\texttt{\rmrk{gemm}} \spaceflop $n_{x}^3$\spaceflopright}
        \vspace{0.2em}
    }
    \vspace{1em}
    We again use implicit reduction of the indices $i$ modulo $P$ to simplify the notation,\;
    and we define $\nu^P_2(0) \defeq \nu_2(P)$ and $\nu^P_2(i) \defeq \nu_2(i)$ for $0 \lt i \lt P$. \quad Color coding\;
    matches \Cref{fig:cyclic-reduction-permutation,fig:cyclic-reduction-graph-16}.\;
    \vspace{0.3em}
\end{algorithm2e}

\subsection{Operation counts \done}

\Cref{alg:fact-mod-riccati,alg:fact-parallel-riccati-cr} list the names of the different BLAS and LAPACK routines that can be used to implement the different steps, together with the cubic terms of the FLOP count%
\footnote{A fused multiply--add (FMA) is counted as a single operation, reflecting the properties of modern hardware where FMA instructions have the same latency and throughput as a standalone multiplication. Using this convention, $n\times n$ matrix--matrix multiplication requires $n^3$ FLOPs, and Cholesky factorization of an $n\times n$ matrix requires $\frac16n^3$ operations.}
for each routine (with $n_{ux} = n_u + n_x$).
Ignoring lower-order terms, the total number of floating-point operations in the critical path of \Cref{alg:fact-parallel-riccati-cr} is given by
\begin{subequations}
    \label{eq:flop}
    \begin{align}
        \tfrac NP \left(
        \tfrac16 n_{ux}^3 %
        + \tfrac12 n_{x} n_{u}^2 %
        + n_{x}^2 n_{u} %
        \right)
        + \left(\tfrac NP - 1\right) \left(
        \tfrac32 n_{ux} n_{x}^2 %
        + \tfrac12 n_{ux}^2 n_{x} %
        \right)
        + \tfrac12 n_{x}^3           \label{eq:flop-riccati}            \\ %
        {} + \tfrac43 n_{x}^3 %
        + \tfrac NP \tfrac12\, n_{u} n_{x}^2 \label{eq:flop-eval-schur} \\ %
        {} + \tfrac16 n_{x}^3 + \log_2(P)\, \tfrac53 n_{x}^3\rlap{,} \label{eq:flop-cr-schur} %
    \end{align}
\end{subequations}
where \eqref{eq:flop-riccati} is the cost of the modified Riccati recursion from \Cref{alg:fact-mod-riccati},
\eqref{eq:flop-eval-schur} is the cost of computing the Schur complement $\mathscr M$ in the \textsc{compute--schur} function of \Cref{alg:fact-parallel-riccati-cr},
and \eqref{eq:flop-cr-schur} is the cost of the cyclic reduction of $\mathscr M$ in the \textsc{factor--schur} function of \Cref{alg:fact-parallel-riccati-cr}.

The Riccati-like elimination of the different intervals in steps \ref{it:bchol-diag} and \ref{it:solve-subdiag}, and the evaluation of the Schur complement in step \ref{it:compute-schur} are embarrassingly parallel:
Each interval can be processed fully independently, without communication (which can be seen from the parallel branches at the bottom of the elimination tree of \Cref{mat:pdp-cr-N12-P4}). For \eqref{eq:flop-riccati}--\eqref{eq:flop-eval-schur}, the number of operations in the critical path scales linearly with $\landauO(N/P)$.
CR in step \ref{it:cr-schur} requires $\log_2(P)+1$ steps, with communication in each step (each node near the root of the elimination tree receives updated matrices from its two children).
If the number of stages is significantly larger than the number of processors ($N \gg P$), the run time of the first three steps dominates
and the expected overall speedup is close to linear.
If the number of stages is closer to the number of processors, ($N \approx P$), the fourth step dominates the overall run time (scaling logarithmically with $P$),
resulting in diminishing returns when doubling $P$ from $N/2$ to $N$, for example.

\begin{figure}
    \centering
    \includegraphics[scale=0.71, clip, trim=0 0.4cm 0 0.4cm]{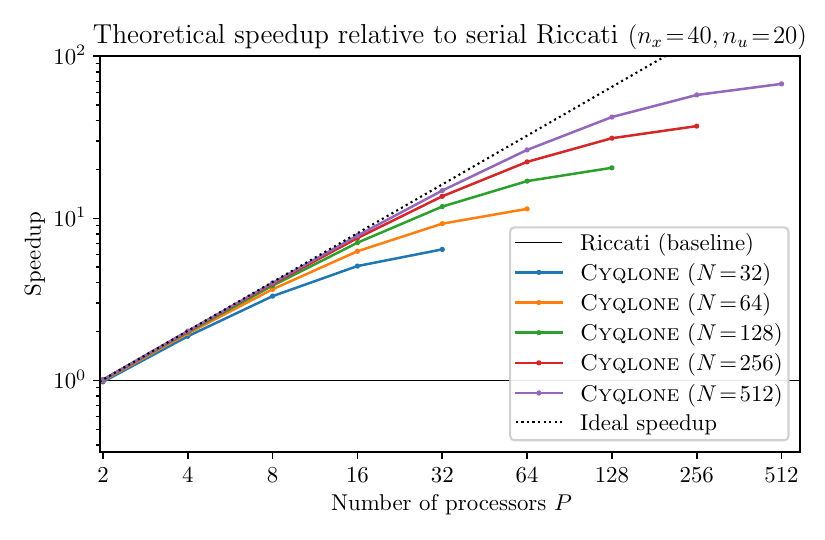}
    \caption{Speedup of the parallel \textsc{Cyqlone} factorization compared to the serial Riccati recursion for various horizon lengths, considering the number of floating-point operations in the critical path.}
    \label{fig:theoretical-speedup}
\end{figure}

\Cref{fig:theoretical-speedup} shows the theoretical limits of the speedup of \textsc{Cyqlone} over a serial implementation based on the Riccati recursion \cite[Alg.\,3]{frison_efficient_2013},
which requires
$\tfrac16 n_x^3 + N \left( \tfrac12 n_{ux} n_x^2 + \tfrac12 n_{ux}^2 n_x + \tfrac16 n_{ux}^3 \right) + \landauO(N n_{ux}^2)$ operations.
The speedup is computed by considering the number of floating-point operations in the critical path of both methods (which is what determines the wall-clock time).
The partitioning into parallel sub-intervals in \cyqlone{} results in some additional fill-in compared to the serial Riccati recursion,
which explains why \cyqlone{} does not see any speedup in the case with two threads.
When the number of processors $P$ is greater than two, the parallelism gained from the partitioning already makes up for this extra work, and doubling the number of processors approximately halves the run time for sufficiently large $N$.
For a single processor, no partitioning of the horizon is performed, and \cyqlone{}
reduces to an algorithm that is almost equivalent to the factorized Riccati recursion from \cite[Alg.\,3]{frison_efficient_2013},
as discussed in the footnote in \Cref{sec:mod-ricc}.
The only difference is the early elimination of $u^0\!$, with negligible impact on the run time for large $N$.

\subsection{Remainders and padding} \label{sec:cyqlone-padding}

Since \cyqlone{} relies on CR,
performance is optimal when the number of processors $P$ is a power of two that divides the horizon length $N$, but this is not a strict requirement:
\cyqlone{} can be used for any $N$ by rounding it up to the next multiple of $P$, and
introducing padding blocks $Q_{j+1}=\I$, $R_j=\I$, $S_{j+1}=0$ and $E_{j+1}=\I$, $A_j=0$, $B_j=0$ for $N \le j \lt P\lceil N/P \rceil$.

\section{Solving problems with inequality constraints \done} \label{sec:qpalm-ocp}
To handle optimal control problems with inequality constraints of the form \eqref{eq:ocp},
we make use of \qpalm{}, a proximal augmented Lagrangian method
for quadratic programs with inequality constraints \cite{hermans_qpalm_2022,hermans_qpalm_2019-1}.
This section briefly summarizes the structure and properties of the inner
optimization problems that need to be solved iteratively by \qpalm{}'s
semismooth Newton solver, which is responsible for the
main computational cost of the method.
When preserving the states and the dynamics constraints in the inner problems
as in the \qpalmocp{} solver described in \cite[Sec.~III]{lowenstein_qpalm-ocp_2024},
the resulting Newton systems are KKT systems with optimal control structure \eqref{eq:opt-cond-ocp-eq},
which we solve using the \textsc{Cyqlone} method described in \Cref{sec:cyqlone}.
Other operations required in the \qpalm{} algorithm, such as matrix--vector products, dot products, vector addition, etc.,
can be parallelized as well. Furthermore, we discuss how \qpalm{}'s exact line search
procedure can be optimized and parallelized.
We will refer to the parallel \qpalm{}-based solver for linear--quadratic optimal control 
problems of the form \eqref{eq:ocp} that uses \cyqlone{} as its linear solver as \cyqpalm{}.

\subsection{The augmented Lagrangian inner problem \done}

\def\qpalmiter{\langle\hspace{-0.8pt}\nu\hspace{-0.8pt}\rangle}

\qpalm{} and \qpalmocp{} solve problem \eqref{eq:ocp} by relaxing the inequality
constraints $b_l^j \le C_j x^j + D_j u^j \le b_u^j$ and $b_l^N \le C_N x^N \le b_u^N$
using an augmented Lagrangian method (ALM).
The inner problem to be solved at each outer iteration of this method amounts to the minimization
of a piecewise quadratic augmented Lagrangian function, subject to
equality constraints:
\begin{equation} \label{eq:ocp-al}
    \begin{aligned}
         & \minimize_{\mathbf{u}, \mathbf{x}} &  & \sum_{j=0}^{N-1} \ell^\mathrm\Sigma_j(x^j, u^j) + \ell^\mathrm\Sigma_N(x^N)                                    \\
         & \subjto                            &  & E_0 x^0 = x_\text{init}                                                                          \\
         &                                    &  & E_{j+1} x^{j+1} = A_j x^j + B_j u^j + f^j                     &  & {\scriptstyle(0 \le j \lt N)} \\
    \end{aligned}
\end{equation}
where the strongly convex stage-wise cost functions are given by the $\mathrm\Gamma$-proximal $\mathrm\Sigma_j$-augmented Lagrangians
\begin{equation} \label{eq:ocp-al-funcs}
    \begin{aligned}
        \ell^\mathrm\Sigma_j(x, u)               & =
        \ell_j(x, u) + \tfrac12 \dist^2_{\mathrm\Sigma_j}\!\!\big(
        C_j x + D_j u + \inv {\mathrm\Sigma}_j y^j; \, [b_l^j, b_u^j]\big) \\
        &\quad
        + \tfrac12 \tnormsq{x - \bar x^j}_{\inv{\mathrm\Gamma_x}}
        + \tfrac12 \tnormsq{u - \bar u^j}_{\inv{\mathrm\Gamma_u}} \\
        \text{and }\quad \ell^\mathrm\Sigma_N(x) & =
        \ell_N(x) + \tfrac12 \dist^2_{\mathrm\Sigma_N}\!\!\big(
        C_N x + \inv {\mathrm\Sigma}_N y^N; \, [b_l^N, b_u^N]\big) \\
        &\quad
        + \tfrac12 \tnormsq{x - \bar x^N}_{\inv{\mathrm\Gamma_x}}.
    \end{aligned}
\end{equation}
In these expressions, the vectors $y^j \in \R^{n_y}$ represent the current estimate of the
Lagrange multipliers corresponding to the inequality constraints in stage $j$, and
$\mathrm\Sigma_j \in \posdefset{\R^{n_y}}$ are diagonal matrices containing the
ALM penalty weights. The matrices $\mathrm\Gamma_x \in \posdefset{\R^{n_x}}$ and $\mathrm\Gamma_u \in \posdefset{\R^{n_u}}$
add primal regularization through a proximal term, relative to the fixed vectors
$\bar x^j \in \R^{n_x}$ and $\bar u^j \in \R^{n_u}$.
The function $\dist^2_\mathrm\Sigma(w;\,\mathrm\Omega)$ denotes the squared distance in $\mathrm\Sigma$-norm between a point $w$ and a set $\mathrm\Omega$.
The intervals $[b_l^j, b_u^j]$
should be interpreted component-wise, describing $n_y$-dimensional boxes.
Because of the squared distance, the augmented Lagrangians are not twice
continuously differentiable. Therefore, QPALM--OCP employs a semismooth Newton
method for solving \eqref{eq:ocp-al}, using the generalized Hessian matrices $H_j(u, x)$ \cite[Eq.~15]{lowenstein_qpalm-ocp_2024}:
\def\setJ{\mathcal{J}}
\begin{align}
        H_j(u, x) & \defeq \begin{pmatrix}
                               R_j(u, x) & S_j(u, x) \\ \tp {S_j(u, x)}\!\!\!\! & Q_j(u, x)
                           \end{pmatrix} \notag \\
                  & = \begin{pmatrix}
                          R^\ell_j & S^\ell_j \\ S^{\ell\top}_j\!\!\!\! & Q^\ell_j
                      \end{pmatrix}
        + \begin{pmatrix}
              \tp D_j \\ \tp C_j
          \end{pmatrix} \mathrm\Sigma_j^\setJ\!(u, x)
        \begin{pmatrix}
            D_j & C_j
        \end{pmatrix}
        +
        \begin{pmatrix}
            \inv{\mathrm\Gamma_u} \\ & \inv{\mathrm\Gamma_x}
        \end{pmatrix} \label{eq:qpalm-ocp-hessians} \\
        &\in \partial_C \big(\nabla \ell^\mathrm\Sigma_j(u, x)\big), \notag
\end{align}
where we defined $\mathrm\Sigma_j^\setJ\!(u, x)$ as $\mathrm\Sigma_j$ with the entries corresponding to inactive constraints set to zero:
\begin{equation}
    \Big(\mathrm\Sigma_j^\setJ\!(u, x)\Big)_{ii} \defeq \begin{cases}
        \big(\mathrm\Sigma_j\big)_{ii} \hspace{-0.2em}                                       & \text{if } i \in \setJ_j(u, x), \\
        0 & \text{otherwise.}
    \end{cases}
\end{equation}
The index set
$\setJ_j(u, x) \defeq \defset{i \in \N_{[1, m]}}{\big(C_j x + D_j u + \inv {\mathrm\Sigma}_j y^j\big)_i \not\in \big[(b_l^j)_i, (b_u^j)_i\big]}$
contains the indices of active constraints at the given point $(u, x)$.
The quantities $H_N(x)$, $\mathrm\Sigma_N^\setJ(x)$ and $\setJ_N(x)$ are defined analogously.

\subsection{Structured Newton systems \done}

The Newton step
$(\dxbf, \dubf)$ with
$\dxbf \defeq [\dx^0 \cdots \dx^N] \in \R^{n_x\times (N+1)}$,
$\dubf \defeq [\du^0 \cdots \du^{N-1}] \in \R^{n_u\times N}$
and the Lagrange multipliers corresponding to the equality constraints
$\boldsymbol \lambda \defeq [\lambda^{-1} \cdots \lambda^{N-1}]$\footnote{The multiplier $\lambda^{-1}$ corresponds to the initial state constraint $E_0 x^0 = x_\text{init}$, which will be eliminated later. The remaining $N$ multiplier vectors are then numbered starting from 0.} $\in \R^{n_x\times (N+1)}$
are found by solving a Newton system of the linearized KKT conditions of \eqref{eq:ocp-al} at a point $(\mathbf u, \mathbf x)$:
\begin{equation} \label{eq:newton-sys-with-x0}
    \hspace{-1em}%
    \begin{aligned}
        \left\{\begin{aligned}
                    & Q_N(x^N)\, \dx^N - E_N \lambda^{N-1}                                                                   \hspace{-1em}&  & = -q^N(x^N) \\[0.8em]
                    & R_j(u^j, x^j)\, \du^j + S_j(u^j, x^j)\, \dx^j + \tp B_j \lambda^j                                      \hspace{-1em}&  & = -r^j(u^j, x^j)    \\[-0.6em] &&&\quad\qquad {\color{gray}\scriptstyle (0 \le j \lt N)} \\[-0.25em]
                    & \ttp{S_j(u^j, x^j)} \hspace{-1pt}\du^j + Q_j(u^j, x^j)\, \dx^j + \tp A_j \lambda^j - E_j \lambda^{j-1} \hspace{-1em}&  & = -q^j(u^j, x^j)    \\[-0.6em] &&&\quad\qquad {\color{gray}\scriptstyle (0 \le j \lt N)} \\[-0.25em]
                    & B_j \du^j + A_j \dx^j - E_{j+1} \dx^{j+1}                                                              \hspace{-1em}&  & = -c^j(u^j, x^j, x^{j+1})    \\[-0.6em] &&&\quad\qquad {\color{gray}\scriptstyle (0 \le j \lt N)} \\[-0.25em]
                    & \mathord-E_{0} \dx^{0}                                                                                 \hspace{-1em}&  & = -c^{-1}(x^0).                    &                                \\
               \end{aligned}\right.
    \end{aligned}%
    \hspace{-1em}%
\end{equation}
The right-hand side is given by the negative gradients of the augmented Lagrangians and the negative residuals of the equality constraints:
{\setstretch{1.2}
\begin{equation*}
    \begin{aligned}
        q^N(x^N)
        & \defeq \grad \ell^\mathrm\Sigma_N(x^N) =
        Q^\ell_N x^N
        +
        q_\ell^N
        +
        \ttp C_N
        \, \hat y^N(x^N) +
        \inv{\mathrm\Gamma_x} (x^N - \bar x^N) \\
        \hat y^N(x^N) & \defeq y^N + \mathrm\Sigma_N \Big( C_N x^N - \proj{}\!\big(C_N x^N + \inv {\mathrm\Sigma}_N y^N; \, [b_l^N, b_u^N]\big) \Big) \\
        \begin{pmatrix}
            r^j(u^j, x^j) \\ q^j(u^j, x^j)
        \end{pmatrix} & \defeq \grad \ell^\mathrm\Sigma_j(u^j, x^j) \\[-1em]
        &=
        \begin{pmatrix}
            R^\ell_j & S^\ell_j \\ S^{\ell\top}_j\!\!\!\! & Q^\ell_j
        \end{pmatrix}\!\!
        \begin{pmatrix}
            u^j \\ x^j
        \end{pmatrix} \!+\!
        \begin{pmatrix}
            r_\ell^j \\ q_\ell^j
        \end{pmatrix} \!+\!
        \begin{pmatrix}
            \ttp D_j \\ \ttp C_j
        \end{pmatrix} \hat y^j(u^j, x^j) \!+\!
        \begin{pmatrix}
            \inv{\mathrm\Gamma_u} (u^j - \bar u^j) \\ \inv{\mathrm\Gamma_x} (x^j - \bar x^j)
        \end{pmatrix}                                                                                                                                         \\
        \hat y^j(u^j, x^j) & \defeq y^j + \mathrm\Sigma_j \Big( D_j u^j + C_j x^j - \proj{}\!\big(D_j u^j + C_j x^j + \inv {\mathrm\Sigma}_j y^j; \, [b_l^j, b_u^j]\big) \Big) \\
        c^j(u^j, x^j, x^{j+1})         & \defeq B_j u^j + A_j x^j + f^j - E_{j+1} x^{j+1}                                                                                                                     \\
        c^{-1}(x^0)                    & \defeq x_\text{init} - E_0 x^0.
    \end{aligned}
\end{equation*}}%
Computation of the quantities $\hat y^j$, $r^j$, $q^j$ can be carried out independently
and in parallel for all stages $j$. The equality constraint residuals $c^j$ can also be
computed in parallel, but require communication between stages.

The initial state $x^0$ and the corresponding steps $\dx^0$ and $\lambda^{-1}$ are often eliminated for practical reasons.
The stationarity conditions and equality constraints for the first stage then become
\begin{equation} \label{eq:elim-x0}
    \begin{aligned}
        R_j(u^0, x^0)\, \du^0 + \tp B_0 \lambda^0 & = -r^0(u^0, \inv E_0 x_\text{init}) \\
        B_0 \du^0 - E_1 \dx^1                     & = -c^0(u^0, \inv E_0 x_\text{init}, x^1). %
    \end{aligned}
\end{equation}
The result is a KKT system with optimal control structure that can be written in the form of \eqref{eq:opt-cond-ocp-eq}.

\subsection{Parallel exact line search}

A key component in the semismooth Newton solver used in \qpalm{} is the line
search procedure that minimizes the augmented Lagrangian along the direction
of the Newton step. Thanks to the specific
piecewise quadratic structure of the objective, the minimizer
can be computed exactly \cite[\S3.2]{hermans_qpalm_2022}.

Below, we derive how a bisection algorithm can be used to find the optimal step size
in linear time. This is an improvement compared to the linearithmic complexity of the sorting-based algorithm in \cite[Alg.\,2]{hermans_qpalm_2022}.
Finally, we discuss how the line search algorithm can be parallelized and vectorized.

For the sake of simplicity, consider the abstract QP formulation of \eqref{eq:ocp}:
\begin{equation} \label{eq:ocp-qp-abstract}
\begin{aligned}
    &\minimize_z && \tfrac12 \ttp z Q z + \ttp q z + c \;\defeq\; \ell(z) \\
    &\subjto && M z = b \\
    &&& b_l \le G z \le b_u.
\end{aligned}
\end{equation}
We denote the number of inequality constraints (i.e. the number of rows of $G$) by $m$.
Given the current iterate $z$ and a descent direction $d$, the line search procedure aims to find the step size $\tau$
that minimizes the merit function $\psi(\tau) \defeq \ell^\mathrm\Sigma(z + \tau d)$,
where $\ell^\mathrm\Sigma(z) \defeq \ell(z) + \tfrac12 \dist^2_{\mathrm\Sigma}\!\big(
        G z + \inv {\mathrm\Sigma} y; \, [b_l, b_u]\big)
        + \tfrac12 \tnormsq{z - \bar z}_{\inv{\mathrm\Gamma}}$ is the $\mathrm\Gamma$-proximal $\mathrm\Sigma$-augmented
Lagrangian of \eqref{eq:ocp-qp-abstract} with respect to the inequality constraints, analogous to \eqref{eq:ocp-al-funcs}.
By convexity, the optimal step size $\tau$ is the root of $\psi^\prime(\tau)=0$.
If $z$ is feasible w.r.t. the equality constraints, the Newton direction $d$
lies in the null space of $M$, maintaining feasibility of $z+\tau d$ for any $\tau$,
so the equality constraints need not be included in the merit function \cite{lowenstein_qpalm-ocp_2024}.

\subsubsection{Bisection of the step size}
It can be shown that the derivative of $\psi$ with respect to the
step size is given by \cite[Eq.\,3.5]{hermans_qpalm_2022}
\begin{equation} \label{eq:psi-prime}
    \psi^\prime(\tau) = \tinprod{\grad\ell^\mathrm\Sigma(z + \tau d)}{d} = \eta \tau + \beta + \tinprod{\delta}{\left[\delta \tau - \alpha\right]_+},
\end{equation}
where $\eta \defeq \tp d (Q + \inv{\mathrm\Gamma}) d$, $\beta \defeq \tp d (Q z + q + \inv{\mathrm\Gamma}(z - \bar z))$,
$\smash{\delta \defeq \scalebox{0.9}{$\begin{pmatrix}
            - \sqrt{\mathrm\Sigma} Gd \vphantom{\inv{\sqrt{\mathrm\Sigma}}} \\
            \phantom- \sqrt{\mathrm\Sigma} Gd \vphantom{\inv{\sqrt{\mathrm\Sigma}}}
\end{pmatrix}$}}$, \\ and
$\alpha \defeq \scalebox{0.9}{$\begin{pmatrix}
            \phantom- \inv{\sqrt{\mathrm\Sigma}} \big(y + \mathrm\Sigma (Gz - b_l)\big) \\
            - \inv{\sqrt{\mathrm\Sigma}} \big(y + \mathrm\Sigma (Gz - b_u)\big) \\
        \end{pmatrix}$}$. \\
The values of $\tau$ where one or more of the constraints change activity are given by the
breakpoints $t_i = \alpha_i/\delta_i$.
The derivative of the merit function $\psi$ can be evaluated efficiently for a given breakpoint $t_j$:
\def\setI{\mathcal{I}}
\begin{equation} \label{eq:psi-prime-at-breakpoint}
        \psi^\prime(t_j) %
        \stackrel{\eqref{eq:psi-prime}}= t_j \eta + \beta + \!\sum_{i \in \setI_j} \delta_i \left(\delta_i \tfrac{\alpha_j}{\delta_j} \!-\! \alpha_i\right)
        = {t_j \underbrace{\!\bigg(\eta + \!{\sum_{i\in\setI_j} \delta_i^2}\bigg)\!}_{\eta(t_j)}\, + \underbrace{\beta -\! \sum_{i\in\setI_j} \delta_i\alpha_i}_{\beta(t_j)}},
\end{equation}
\begin{equation*}
    \begin{aligned}
        \text{where } \setI_j^{\phantom+} \!&\defeq \setI_j^- \union \setI_j^+ \\
        \setI_j^- \!&\defeq \defset{i \in \N_{[1,2m]}}{\tfrac{\alpha_j}{\delta_j}\delta_i \!-\! \alpha_i \ge 0 \wedge \delta_i \lt 0}
        = \defset{i \in \N_{[1,2m]}}{\scalebox{0.86}{$\begin{aligned} &t_i \ge t_j \\ {}\mathbin{\wedge}{} &\delta_i \lt 0 \end{aligned}$}} \\
        \setI_j^+ \!&\defeq \defset{i \in \N_{[1,2m]}}{\tfrac{\alpha_j}{\delta_j}\delta_i \!-\! \alpha_i \gt 0 \wedge \delta_i \gt 0}
        = \defset{i \in \N_{[1,2m]}}{\scalebox{0.86}{$\begin{aligned} &t_i \lt t_j \\ {}\mathbin{\wedge}{} &\delta_i \gt 0 \end{aligned}$}}\!.
    \end{aligned}
\end{equation*}
A simple bisection technique can be used to find the root $\tau$ of the
monotonically increasing and piecewise linear function $\psi^\prime(\tau)$.
\begin{enumerate}
    \setlength\itemsep{0em}
    \item Initialize the search interval $\underline \tau = 0$, $\overline \tau = +\infty$
    \item Select some $j$ for which $\underline \tau \lt t_j \lt \overline \tau$  \label{item:label-loop-bisect}
    \item If no such $j$ exists, find $\tau$ by interpolation between $\underline \tau$ and $\overline \tau$, and return  \label{item:step-interpolation-tau}
    \item If $\psi^\prime(t_j) = 0$, return $\tau = t_j$
    \item If $\psi^\prime(t_j) < 0$, set $\underline \tau = t_j$
    \item If $\psi^\prime(t_j) > 0$, set $\overline \tau = t_j$
    \item Go to \ref{item:label-loop-bisect}.
\end{enumerate}
As an invariant, we have $\psi^\prime(\underline \tau) < 0$ and $\psi^\prime(\overline \tau) > 0$.\,\footnote{This holds initially because $d$ is a descent direction and $\psi$ is strongly convex.}
If no breakpoint $t_j$ can be found in the interval $(\underline \tau,\; \overline \tau)$ in step \ref{item:step-interpolation-tau},
then $\psi^\prime$ is linear on that interval, and finding its root $\tau$ is trivial.
Thanks to \eqref{eq:psi-prime-at-breakpoint}, the value of $\psi^\prime(t_j)$ does not have to be recomputed from scratch at each iteration:
if $t_k > t_j$, we have
\begin{equation} \label{eq:ls-update-eta}
    \eta(t_k) = \eta(t_j) -\!\! \sum_{i\in\setI_{j,k}^-} \delta_i^2 +\!\! \sum_{i\in\setI_{j,k}^+} \delta_i^2,
    \vspace{-0.5em}
\end{equation}
where $
    \setI_{j,k}^- \!\defeq \defset{i \in \N_{[1,2m]}}{\scalebox{0.86}{$\begin{aligned} & t_k \gt t_i \ge t_j \\ & {}\mathbin{\wedge}{} \delta_i < 0 \end{aligned}$}}
$ and $
    \setI_{j,k}^+ \!\defeq \defset{i \in \N_{[1,2m]}}{\scalebox{0.86}{$\begin{aligned} & t_j \le t_i \lt t_k \\ & {}\mathbin{\wedge}{} \delta_i > 0 \end{aligned}$}}\!,
$
significantly reducing the number of terms that need to be summed compared to a naive implementation of \eqref{eq:psi-prime-at-breakpoint}. Similar update formulas
can be obtained for $\beta(t_k)$ and for the case where $t_k \lt t_j$.

Unlike the original implementation in \cite[Alg.\,2]{hermans_qpalm_2022}, we do not sort the breakpoints $t_i$, as this would require
$\landauO(m \log m)$ operations.
An algorithm like quickselect or introselect \cite{musser_introspective_1997_manual} can be used to select a breakpoint in
$(\underline \tau,\; \overline \tau)$. As a side effect, it partitions the interval into breakpoints smaller and greater than
the selected pivot $t_j$, and the same partitioning is applied to $\alpha$ and $\delta$.
This partitioning helps avoid scanning the full vectors $\alpha$ and $\delta$ for the
evaluation of \eqref{eq:ls-update-eta} at each trial breakpoint $t_j$.
The cost of performing introselect is linear in the number of breakpoints in $(\underline \tau,\; \overline \tau)$.
If the median element is selected, this number is halved at each iteration, with linear overall complexity, $\landauO(m)$.

Rather than selecting the median breakpoint as the pivot,
it is also possible to select a specific value $t_j$ directly,
and partition the interval based on that choice. While the average cost of this
approach may be lower than a version using introselect, its worst case
complexity is quadratic in the number of breakpoints. Despite the asymptotic
concerns, such a strategy could nevertheless be useful to prune the search space
early on: in practice, the optimal step size is rarely greater than one, so
selecting the largest breakpoint below one as $t_j$ in the first iteration
often eliminates a potentially large number of uninteresting breakpoints.

Practical implementations can be further optimized by storing just $t_i$ and $\delta_i |\delta_i|$
instead of $t_i$, $\alpha_i$ and $\delta_i$. This choice reduces memory usage,
memory bandwidth and the number of branches, and avoids computing square roots.
The update formulas from \eqref{eq:ls-update-eta} can then be written as
\begin{equation*}
    \eta(t_k) = \eta(t_j) +\!\! \sum_{i\in\setI_{j,k}} \delta_i|\delta_i| \;\; \text{and} \;\;
    \beta(t_k) = \beta(t_j) -\!\! \sum_{i\in\setI_{j,k}} t_i \delta_i|\delta_i|, %
    \vspace{-0.5em}
\end{equation*}
where $\setI_{j,k} \defeq \defset{i \in \N_{[1,2m]}}{t_j \le t_i \lt t_k}$.

\subsubsection{Parallelization of the line search algorithm}

The evaluation of the vectors $\alpha$, $\delta$ and $t$ requires only
stage-wise matrix-vector products and element-wise operations, and can thus
be fully parallelized and vectorized. The breakpoints are then partitioned in
parallel: non-finite values of $t_i$, which do not affect $\psi^\prime(\tau)$, are removed
(e.g. due to one-sided constraints or constraints with gradients orthogonal to $d$),
and negative breakpoints are isolated,
as they are not considered during the bisection (but they do affect $\psi^\prime(\tau)$).
Each thread produces its local partitions in parallel.
After selecting a trial breakpoint $t_j$,
all threads update their partitions using $t_j$ as a pivot
and compute their local sums to evaluate $\eta(t_j)$ and $\beta(t_j)$,
which are then combined to compute $\psi^\prime(t_j)$.
We can leverage a vast literature on optimal and parallel selection
algorithms \cite{ribizel_parallel_2020,siebert_scalable_2014,martinez_sesquickselect_2019};
although these operations are unlikely to be a bottleneck in the full implementation
of \cyqpalm{}.
Once the remaining
interval in the bisection algorithm is sufficiently small, the local partitions
of all threads can be merged to complete the procedure on a single thread,
thereby reducing synchronization overhead.
When the size of the interval drops below a second threshold,
the remaining breakpoints are sorted, and a linear search is used to find the
optimal step size $\tau$ (this is faster than continuing the bisection all the
way down to a single element).

\section{Factorization updates \done} \label{sec:fact-upd}
\def\Ssign{\mathscr S}

At every iteration of \qpalm{}'s inner semismooth Newton method,
the generalized Hessian matrices $H_j(u,x)$ defined in \eqref{eq:qpalm-ocp-hessians}
may vary because of changes in the active set
at the current iterate, $\setJ_j(u,x)$.
If the activity of only a limited number of constraints
changes, the factorization of the new Newton system can be expressed as a low-rank
modification of the factorization of the previous system. We will refer to such a
modification of triangular factors as a \textit{factorization update}.
If the rank of the modification is sufficiently small, factorization updates
require fewer operations than a full factorization of the new matrix
(without reusing the previous factors) \cite{gill_methods_1974}
thereby enabling significant performance gains in the resulting \qpalm{} solver.

To derive factorization update formulas for the full block Cholesky factor of
the KKT system with optimal control structure from \Cref{sec:cyqlone}, we 
make use of the hyperbolic Householder transformations described
by \cite{pas_blocked_2025-1} as the main building block.

\newcommand\pmJ{\setJ_{\!\pm}}
\def\Sigmaud{\mathrm\Sigma_{\pmJ}}
\newcommand\Sigmaudk[1]{\mathrm\Sigma_{\pmJ\!, #1}}
\newcommand\Sigmabarudk[1]{\bar\mathrm\Sigma_{\pmJ\!, #1}}

Consider an iterate $(\mathbf u, \mathbf x)$ for which an existing factorization of
the Newton system \eqref{eq:opt-cond-ocp-eq} is available. We now wish to obtain
a factorization of the corresponding Newton system at a different point
$(\tilde {\mathbf u}, \tilde {\mathbf x})$. By \eqref{eq:qpalm-ocp-hessians},
we have that
\begin{equation} \label{eq:low-rank-update-Hj}
    H_j(\tilde u^j, \tilde x^j) = H_j(u^j, x^j) + \begin{pmatrix}
        \ttp D_j \\ \ttp C_j
    \end{pmatrix}
    \Big(\mathrm\Sigma_j^\setJ\!(\tilde u^j, \tilde x^j) - \mathrm\Sigma_j^\setJ\!(u^j, x^j)\Big)
    \begin{pmatrix}
        D_j & C_j
    \end{pmatrix}.
\end{equation}
\def\DSigma{\Delta\hspace{-1pt}\mathrm\Sigma}%
\def\tilDSigma{\Delta\hspace{-1pt}\tilde\mathrm\Sigma}%
When the active sets at $(u^j, x^j)$ and $(\tilde u^j, \tilde x^j)$ differ
for a small number of constraints only, the diagonal matrix
$\DSigma_j \defeq \mathrm\Sigma_j^\setJ\!(\tilde u^j, \tilde x^j) - \mathrm\Sigma_j^\setJ\!(u^j, x^j)$ has a small number of nonzero entries, and
$H_j(\tilde u^j, \tilde x^j)$ and $H_j(u^j, x^j)$ differ by a low-rank term.

\subsection{Parallel factorization updates in {\normalfont\cyqlone} \done}

To see how the low-rank updates in \eqref{eq:low-rank-update-Hj} affect the
factorization computed by the \cyqlone{} method, we will derive a
factorization of the updated KKT matrix $\widetilde{\mathscr K}$, which we define as the matrix $\mathscr K$ from 
\eqref{eq:block-matrix-example-4} with a term $\mathit\Upsilon \Ssign \hspace{-1pt}\ttp{\mathit\Upsilon}$
added to the top-left block:
{
\setstretch{1.2}
\begin{align}
    \label{eq:cyqlone-K-tilde}
    \widetilde {\mathscr K} \defeq{}&
    \begin{pNiceArray}{c||c}[columns-width=1.8em]
        \widetilde{\mathscr R} & \tp{\mathscr A} \\\hline\hline
        \mathscr A & 0
    \end{pNiceArray} \defeq
    \begin{pNiceArray}{c||c}[columns-width=1.8em]
        \mathscr R & \tp{\mathscr A} \\\hline\hline
        \mathscr A & 0
    \end{pNiceArray}
    +
    \begin{pNiceArray}{c}[columns-width=1.8em]
        \mathit\Upsilon \\\hline\hline
        0
    \end{pNiceArray}
    \Ssign
    \ttp{
    \begin{pNiceArray}{c}[columns-width=1.8em]
        \mathit\Upsilon \\\hline\hline
        0
    \end{pNiceArray}
    } \\
    \stackrel{\eqref{eq:chol-fac-K-cyqlone}}={}&
    \begin{pNiceArray}{c||c}[columns-width=1.8em]
        \L{\mathscr R} & \\\hline\hline
        \L{\mathscr A} & \L{\mathscr M}
    \end{pNiceArray}
    \begin{pNiceArray}{c||c}[columns-width=1.8em]
        \mathscr D \\\hline\hline
        & \!\!-\I
    \end{pNiceArray}
    \ttp{
    \begin{pNiceArray}{c||c}[columns-width=1.8em]
        \L{\mathscr R} & \\\hline\hline
        \L{\mathscr A} & \L{\mathscr M}
    \end{pNiceArray}
    }
    +
    \begin{pNiceArray}{c}[columns-width=1.8em]
        \mathit\Upsilon \\\hline\hline
        0
    \end{pNiceArray}
    \Ssign
    \ttp{
    \begin{pNiceArray}{c}[columns-width=1.8em]
        \mathit\Upsilon \\\hline\hline
        0
    \end{pNiceArray}
    } \label{eq:chol-fac-K-cyqlone-K-tilde}. %
\end{align}}%
Following the structure of the algorithm in \Cref{sec:cyqlone},
we apply a factorization update to $\mathscr R$ first, and then update the Schur
complement. To this end, we use a hyperbolic Householder transformation
$\breve{\mathscr Q}^{\mathscr R}$ such that
\begin{align}
    \label{eq:hyp-house-R}
    \begin{pmatrix}
        \tilL{\mathscr R} & 0 & 0 \\
        \tilL{\mathscr A} & \L{\mathscr M} & \mathit\Xi
    \end{pmatrix} &=
    \begin{pmatrix}
        \L{\mathscr R} & 0 & \mathit\Upsilon \\
        \L{\mathscr A} & \L{\mathscr M} & 0
    \end{pmatrix} \breve{\mathscr Q}^{\mathscr R}, &&\text{with }
    \breve{\mathscr Q}^{\mathscr R} \scalebox{0.87}{$\begin{pmatrix}
        \mathscr D \\ & \!\!-\I \\ && \Ssign
    \end{pmatrix}$} \breve{\mathscr Q}^{\mathscr R\top} \hspace{-0.7em}&= \scalebox{0.87}{$\begin{pmatrix}
        \mathscr D \\ & \!\!-\I \\ && \Ssign
    \end{pmatrix}$},
\intertext{followed by another hyperbolic Householder transformation $\breve{\mathscr Q}^{\mathscr M}$
satisfying}
    \label{eq:hyp-house-M}
    \begin{pmatrix}
        \tilL{\mathscr R} & 0 & 0 \\
        \tilL{\mathscr A} & \tilL{\mathscr M} & 0
    \end{pmatrix} &=
    \begin{pmatrix}
        \tilL{\mathscr R} & 0 & 0 \\
        \tilL{\mathscr A} & \L{\mathscr M} & \mathit\Xi
    \end{pmatrix} \breve{\mathscr Q}^{\mathscr M}, &&\text{with }
    \breve{\mathscr Q}^{\mathscr M} \scalebox{0.87}{$\begin{pmatrix}
        \mathscr D \\ & \!\!-\I \\ && \Ssign
    \end{pmatrix}$} \breve{\mathscr Q}^{\mathscr M\,\top} \hspace{-0.7em}&= \scalebox{0.87}{$\begin{pmatrix}
        \mathscr D \\ & \!\!-\I \\ && \Ssign
    \end{pmatrix}$}.
\end{align}
Such transformations $ \breve{\mathscr Q}^{\mathscr R}$ and $ \breve{\mathscr Q}^{\mathscr M}$
can be constructed using the algorithms outlined in \cite[\S III]{pas_blocked_2025-1}.
From \eqref{eq:hyp-house-R}--\eqref{eq:hyp-house-M}, we conclude that $\tilL{\mathscr R}$, $\tilL{\mathscr A}$ and $\tilL{\mathscr M}$ indeed form the desired factorization of $\widetilde{\mathscr K}$:
\begin{align}
    & \begin{pmatrix}
        \tilL{\mathscr R} & 0 \\
        \tilL{\mathscr A} & \tilL{\mathscr M}
    \end{pmatrix} \scalebox{1}{$\begin{pmatrix}
        \mathscr D \\ & \!\!-\I
    \end{pmatrix}$} \ttp{
    \begin{pmatrix}
        \tilL{\mathscr R} & 0 \\
        \tilL{\mathscr A} & \tilL{\mathscr M}
    \end{pmatrix}
    } \notag \\
    ={}& \begin{pmatrix}
        \tilL{\mathscr R} & 0 & 0 \\
        \tilL{\mathscr A} & \tilL{\mathscr M} & 0
    \end{pmatrix} \scalebox{0.87}{$\begin{pmatrix}
        \mathscr D \\ & \!\!-\I \\ && \Ssign
    \end{pmatrix}$} \ttp{
    \begin{pmatrix}
        \tilL{\mathscr R} & 0 & 0 \\
        \tilL{\mathscr A} & \tilL{\mathscr M} & 0
    \end{pmatrix}
    } \\
    ={}& \begin{pmatrix}
        \L{\mathscr R} & 0 & \mathit\Upsilon \\
        \L{\mathscr A} & \L{\mathscr M} & 0
    \end{pmatrix}
    \breve{\mathscr Q}^{\mathscr R} \breve{\mathscr Q}^{\mathscr M}
    \scalebox{0.87}{$\begin{pmatrix}
        \mathscr D \\ & \!\!-\I \\ && \Ssign
    \end{pmatrix}$} 
    \breve{\mathscr Q}^{\mathscr M\,\top} \breve{\mathscr Q}^{\mathscr R\top}
    \ttp{
    \begin{pmatrix}
        \L{\mathscr R} & 0 & \mathit\Upsilon \\
        \L{\mathscr A} & \L{\mathscr M} & 0
    \end{pmatrix}
    } \\
    ={}& \begin{pmatrix}
        \L{\mathscr R} & 0 & \mathit\Upsilon \\
        \L{\mathscr A} & \L{\mathscr M} & 0
    \end{pmatrix}
    \scalebox{0.87}{$\begin{pmatrix}
        \mathscr D \\ & \!\!-\I \\ && \Ssign
    \end{pmatrix}$}
    \ttp{
    \begin{pmatrix}
        \L{\mathscr R} & 0 & \mathit\Upsilon \\
        \L{\mathscr A} & \L{\mathscr M} & 0
    \end{pmatrix}
    } \stackrel{\eqref{eq:chol-fac-K-cyqlone-K-tilde}}= \widetilde{\mathscr K}.
\end{align}

\subsection{Modified Riccati recursion update \rmrk{steps \ref{it:bchol-diag} and \ref{it:solve-subdiag}} \done}

Computing the updated factors $\tilL{\mathscr R}$ and $\tilL{\mathscr A}$
as in \eqref{eq:hyp-house-R} can be done using a modified version of
\cite[Alg.\,4]{pas_blocked_2025-1}, listed in \Cref{alg:fact-mod-riccati-update}.
This algorithm mirrors the factorization in \Cref{alg:fact-mod-riccati}.
In addition to the updated factors, it also returns blocks of the matrix $\mathit\Xi$
that can be used to perform the update of the Schur complement.
Following the specific example from \Cref{sec:cyqlone}, the second block column of $\mathscr K$ is updated as follows:
\begin{align}
    \label{eq:riccati-update-tilde-K1}
    \widetilde{\mathscr K}_1 \defeq{}& \left(\begin{NiceArray}{c|c}[columns-width=2em]
        \widetilde{\mathscr R}_1 & \ttp {\mathscr A}_1 \\\hline
        \mathscr A_1 & 0
    \end{NiceArray}\right) = 
    \mathscr K_1 + \left(\begin{NiceArray}{c}[columns-width=1.2em] \mathit\Upsilon_1 \mkern8mu \\\hline 0 \end{NiceArray}\right)
    \Ssign_1\hspace{-1pt} \ttp {\left(\begin{NiceArray}{c}[columns-width=1.2em] \mathit\Upsilon_1 \mkern8mu \\\hline 0 \end{NiceArray}\right)}, \\
    \text{where}\quad {\mathit\Upsilon}_1                                                \defeq{}
             & \scalebox{0.75}{\setstretch{1.25}$\left(
            \begin{NiceArray}{c|c|c}[columns-width=3em]
                    \ttp D_3 &          &          \\
                    \ttp C_3 &          &          \\
                    0        &          &          \\\hline
                             & \ttp D_2 &          \\
                             & \ttp C_2 &          \\
                             &        0 &          \\\hline
                             &          & \ttp D_1 \\
                             &          & \ttp C_1 %
            \end{NiceArray}
    \right)$} \quad \text{and} \quad
    {\Ssign}_1                                                \defeq{}
             \scalebox{0.75}{\setstretch{1.25}$\left(
            \begin{NiceArray}{c|c|c}[columns-width=2em]
                    \DSigma_3 \\\hline
                    & \DSigma_2 \\\hline
                    && \DSigma_1
            \end{NiceArray}
    \right)$}. \label{eq:ups-sigma-upd-riccati}
\end{align}

\newcommand{\Ups}[1][]{%
    \mathrel{\mathrm\Upsilon}%
    \if\relax\detokenize{#1}\relax\else^{(\mkern-1.5mu#1\mkern-1.5mu)}\fi
}
\newcommand{\Xifwd}[1][]{%
    \mathrel{\rlap{\smash{\raisebox{1.55ex}{$\hspace{0.25em}\scriptveryshortrightarrow$}}}\mathrm\Upsilon}%
    \if\relax\detokenize{#1}\relax\else^{(\mkern-1.5mu#1\mkern-1.5mu)}\fi
}
\newcommand{\Xibwd}[1][]{%
    \mathrel{\rlap{\smash{\raisebox{1.55ex}{$\hspace{0.09em}\scriptveryshortleftarrow$}}}\mathrm\Upsilon}%
    \if\relax\detokenize{#1}\relax\else^{(\mkern-1.5mu#1\mkern-1.5mu)}\fi
}
\newcommand{\Ssignric}[1]{%
    {\mathcal S}_{\mkern-1.5mu#1\mkern-1.5mu}
}
\newcommand{\Ssigncr}[1][]{%
    \Ssign%
    \if\relax\detokenize{#1}\relax\else^{(\mkern-1.5mu#1\mkern-1.5mu)}\fi
}

\begin{algorithm2e}[htbp]
    \def\bbwd{b_\mathrm{bwd}}
    \def\bfwd{b_\mathrm{fwd}}
    \caption{Factorization update of a single modified Riccati block column}
    \label{alg:fact-mod-riccati-update}
    \DontPrintSemicolon
    \KwIn{$A, B, C, D$ and $N$: OCP data matrices and horizon length}
    \KwIn{$P$: number of intervals to partition the horizon into}
    \KwIn{$c \in \N_{[0, P)}$: index of the block column to update (determines the interval)}
    \KwIn{$\DSigma$: changes in the penalty factors}
    \KwIn{$\L{\mathscr R}_c$ ($\L{R}, \L{S}, \L{Q}$) and $\L{\mathscr A}_c$ ($\L{B}, \L{A}, \Acl$): existing factorization of $\mathscr K_c$}
    \KwOut{$\tilL{\mathscr R}_c$ ($\tilL{R}, \tilL{S}, \tilL{Q}$) and $\tilL{\mathscr A}_c$ ($\tilL{B}, \tilL{A}, \tilAcl$): factorization of $\widetilde{\mathscr K}_c\; \eqref{eq:riccati-update-tilde-K1}$}
    \KwOut{$\Xifwd_{c}, \Xibwd_{c-1}, \Ssign_c$: contribution to the update of the Schur complement}
    \vspace{0.3em}
    \Fn{\normalfont\textsc{update--block--column--riccati}$(c)$}{
    $n = N / P$\mycommentnofill{\rmrknp{Number of stages per interval}}\;
    $j_{1} = n(c-1)+1$,\quad $j_{n} = nc$\mycommentnofill{\rmrknp{First/last stage indices in block column $c$}}\;

    $\setstretch{1.2} \begin{pmatrix}
        \mathrm\Upsilon^{u}_{j_n} \\ \mathrm\Upsilon^{x}_{j_n} \\ \mathrm\Upsilon^{\lambda}_{j_n}
    \end{pmatrix} \assignpc \begin{pmatrix}
            \ttp D_{j_n} \\
            \ttp C_{j_n} \\
            0
    \end{pmatrix},\quad
    \Ssignric{j_n} \assignpc \DSigma_{j_n}%
    $\;

    \For{$j = j_{n}, j_{n} - 1, j_{n} - 2, \dots, j_{1}$}{
        $m_j = \operatorname{rank}{\Ssignric{j}}$\;
        $\setstretch{1.2} \begin{pmatrix}
            \tilde L^{\!R\vphantom Q}_{j} & 0 \\
            \tilde L^{\!S}_{j} & \mathrm\Phi^{x}_{j} \\
            \tilde L^{\!B}_j & \mathrm\Phi^{\lambda}_{j}
        \end{pmatrix} \assignpc \begin{pmatrix}
            L^{\!R \vphantom Q}_{j} & \mathrm\Upsilon^{u}_{j} \\
            L^{\!S}_{j} & \mathrm\Upsilon^{x}_{j} \\
            L^{\!B}_{j} & \mathrm\Upsilon^{\lambda}_{j}
        \end{pmatrix} \breve Q^{u}_{j}$, where $\breve Q^{u}_{j}$ is $\left( \begin{smallmatrix}
            \I \\ &\Ssignric{j}
        \end{smallmatrix} \right)$-orth. \cite{pas_blocked_2025-1}%
        \mycomment{\texttt{\rmrk{hyh}} $\hspace{1.5em}\mathclap{\substack{m_jn_u^2 \\ + \\ 4m_jn_un_x}}\hspace{0.5em}$}

        \vspace{0.3em}
        \If{$j > j_{1}$}{
            $\tilAcl_j \assignpc \Acl_j + \mathrm\Phi^{\lambda}_{j}\, \Ssignric{j}\, \mathrm\Phi^{x\,\top}_{j}$
            \mycomment{\texttt{\rmrk{gemm}} \quad $\hspace{0.55em}\mathclap{m_jn_x^2}\hspace{0.5em}$}

            $\setstretch{1.2} \begin{pmatrix}
                \mathrm\Upsilon^{u}_{j-1} \\ \mathrm\Upsilon^{x}_{j-1} \\ \mathrm\Upsilon^{\lambda}_{j-1}
            \end{pmatrix} \assignpc \begin{pmatrix}
                    \ttp B_{j-1} \,\mathrm\Phi^{x}_{j} & \phantom{.}\ttp D_{j-1} \\
                    \ttp A_{j-1} \,\mathrm\Phi^{x}_{j} & \phantom{.}\ttp C_{j-1} \\
                \mathrm\Phi^{\lambda}_{j} & 0
            \end{pmatrix}$
            \mycomment{\texttt{\rmrk{gemm}} \quad $\hspace{2.6em}\mathclap{m_jn_x(n_x\!+\!n_u)}\hspace{1.5em}$}
            $\Ssignric{j-1} \assignpc \left(\begin{smallmatrix}
                \Ssignric{j} \\ & \DSigma_{j-1}
            \end{smallmatrix}\right)%
            $\;
        }
        \vspace{0.3em}

        $\begin{pmatrix}
            \tilde L^{Q}_{j} & 0
        \end{pmatrix} \assignpc \begin{pmatrix}
            L^{Q}_{j} & \mathrm\Phi^{x}_{j}
        \end{pmatrix} \breve Q^{x}_{j}$ \quad where $\breve Q^{x}_{j}$ is $\left( \begin{smallmatrix}
            \I \\ &\Ssignric{j}
        \end{smallmatrix} \right)$-orth.%
        \mycomment{\texttt{\rmrk{hyh}} \quad $\quad\mathclap{m_jn_x^2}\hspace{0.5em}$}
    }
    $\Ssign_c \assignpc \Ssignric{j_1}$\;
    $\setstretch{1.2} \begin{pmatrix}
        \phantom-\tilde L^{A\phantom{\,-1}}_{j_1} & \Xifwd_{c\phantom{-1}} \\
        -\tilde L^{Q\,-\!\top}_{j_1} & \Xibwd_{c-1}
    \end{pmatrix} \assignpc
    \begin{pmatrix}
        \phantom-L^{A\phantom{\,-1}}_{j_1} & \mathrm\Phi^{\lambda}_{j_1} \\
        -L^{Q\,-\!\top}_{j_1} & 0
    \end{pmatrix}
    \breve Q^{x}_{j_1}$
    \mycomment{\texttt{\rmrk{hyh apply}} \quad $\quad\mathclap{4m_{j_1}n_x^2}\hspace{0.5em}$}
}
\end{algorithm2e}

\subsection{Low-rank update of the Schur complement \rmrk{step \ref{it:compute-schur}} \done} \label{sec:update-schur-complement}

Combining the contributions of \Cref{alg:fact-mod-riccati-update} from all block columns, the matrix $\mathit\Xi$ from \eqref{eq:hyp-house-R}
that updates the Schur complement $\tilde{\mathscr{M}} = \mathscr M + \mathit\Xi \Ssign \ttp{\mathit\Xi} = \tilL{\mathscr M}\tilLt{\mathscr M}$
in \eqref{eq:hyp-house-M} and the updated factor $\tilL{\mathscr M}$ are given by
\begin{equation} \label{eq:def-xi-update-example}
    \hspace{-1em}%
    \mathit\Xi = \scalebox{0.75}{\setstretch{1.5}$
            \begin{pNiceArray}{cccc}[columns-width=3em,first-row,first-col]
                & c=0 & c=1 & c=2 & c=3 \\
                \lambda^3 \;\; & & \Xifwd_1 & \Xibwd_1 & \\
                \lambda^9 \;\; & \Xibwd_3 & & & \Xifwd_3 \\
                \lambda^6 \;\; & & & \Xifwd_2 & \Xibwd_2 \\
                \lambda^0\;\; & {\color{gray}\Xifwd_0} & \Xibwd_0 & & \\
            \end{pNiceArray}$}
            \text{ and }\tilL{\mathscr M} = 
            \scalebox{0.75}{\setstretch{1.5}$
            \begin{pNiceArray}{cccc}[columns-width=3em,first-row,first-col]
                &\shortstack{$\;\lambda^3$\\$\scriptsize i=1$}& \shortstack{$\;\lambda^9$\\$\scriptsize i=3$} & \shortstack{$\;\lambda^6$\\$\scriptsize i=2$} & \shortstack{$\;\lambda^0$\\$\scriptsize i=0$} \\
                \lambda^3 \;\; & \tilde L_1 & & & \\
                \lambda^9 \;\; & & \tilde L_3 & & \\
                \lambda^6 \;\; & \tilde Y_1 & \tilde U_3 & \tilde L_2 & \\
                \lambda^0\;\; & \tilde U_1 & {\color{gray}\tilde Y_3} & \tilde U_2 & \tilde L_0 \\
            \end{pNiceArray}$}.
    \hspace{-0.5em}%
\end{equation}
In the absence of coupling between the first and final stages, the block $\tilde Y_3$ is zero,
and the blocks $\Xifwd_0$ and $\Xibwd_3$ have complementary sparsity patterns.
Consequently, the updates of $\L{R}_0$ and $\L{B}_0$ can be performed in isolation
(cf. $u^0$ in the elimination tree in \Cref{mat:pdp-cr-N12-P4}), as long as the contribution to the block $L_0$ of $\L{\mathscr M}$ is taken into account in the very last step
of the update of the Schur complement. This simplifies the following section, allowing us to assume that
$\Xifwd_0 = 0$.

\subsection{Cyclic reduction factorization updates of the Schur complement \rmrk{step \ref{it:cr-schur}} \done} \label{sec:cr-update-schur}

To update the CR factorization, we systematically apply the following
operation to the appropriate blocks of the factors computed by \Cref{alg:fact-parallel-riccati-cr}:
\vspace{0pt}
\begin{equation} \label{eq:update-cr-building-block}
    \scalebox{0.87}{$
    \left( \begin{NiceArray}{c|cc}[columns-width=2.4em]
        \tilde L_{11} & 0 & 0 \\
        \tilde L_{21} & \smash{\widetilde{\mathrm\Xi}_{21}} & \smash{\widetilde{\mathrm\Xi}_{22}} \\
        \tilde L_{31} & \smash{\widetilde{\mathrm\Xi}_{31}} & \smash{\widetilde{\mathrm\Xi}_{32}}
    \end{NiceArray} \right)$} =
    \scalebox{0.87}{$\left( \begin{NiceArray}{c|cc}[columns-width=2.4em]
        L_{11} & {\mathrm\Xi}_{11} & {\mathrm\Xi}_{12} \\
        L_{21} & {\mathrm\Xi}_{21} & 0 \\
        L_{31} & 0 & {\mathrm\Xi}_{32}
    \end{NiceArray} \right)$} \breve Q.
\end{equation}
This operation can be implemented using the algorithms in \cite[\S III]{pas_blocked_2025-1}.

\begin{figure}[p]
    \centering
    \includegraphics[width=\textwidth, clip, trim=0 0.15cm 0 0]{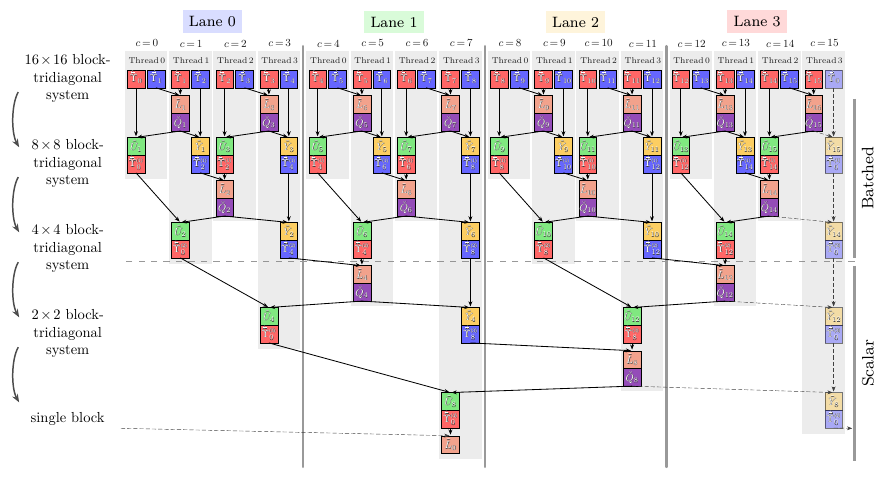}
    \caption{Graph representation of the different steps required for the CR factorization updates of a symmetric $16\times 16$
        block-tridiagonal matrix.
        Node colors match the colors in \Cref{alg:fact-parallel-riccati-cr-update}.
       Note the structural similarity to \Cref{fig:cyclic-reduction-graph-16}.
       Vector lane assignment, the distinction between scalar and batched levels, and other vectorization aspects
       will be discussed in \Cref{sec:vectorization}.
    }
    \label{fig:cyclic-reduction-update-graph-16}
\end{figure}

For each column $i$ of $\L{\mathscr M}$, consider the nonzero blocks, and update them
using the columns of $\mathit\Xi$ that have nonzero elements in row $i$, as described
by \eqref{eq:update-cr-building-block}. The order in which the columns are updated
is the same as the order used during the factorization in \Cref{alg:fact-parallel-riccati-cr}:
in the first level, all odd columns are updated, then all multiples of two that are not multiples of four in the second level, etc.
Within each level, different columns can be updated in parallel.
Note that applying \eqref{eq:update-cr-building-block} introduces the additional nonzero blocks
$\smash{\widetilde{\mathrm\Xi}_{22}}$ and $\smash{\widetilde{\mathrm\Xi}_{31}}$, which have to be taken into account in the next level of the algorithm.
For the specific matrix $\L{\mathscr M}$ defined in \eqref{eq:chol-fac-schur-compl-M}
and with $\mathit\Xi$ and $\tilL{\mathscr M}$ from \eqref{eq:def-xi-update-example},
the updating procedure amounts to
\begin{align*}
    \rotatebox[origin=c]{90}{\underline{Level 0}}\quad && \scalebox{0.87}{$\left( \begin{NiceArray}{c|c}[columns-width=2.0em]
        \tilde L_{1} & 0 \\
        \tilde U_{1} & \Xibwd[1]_0 \\
        \tilde Y_{1} & \Xifwd[1]_2
    \end{NiceArray} \right)$} &=
    \scalebox{0.87}{$\left( \begin{NiceArray}{c|cc}[columns-width=2.0em]
        L_{1} & \Xifwd_1 & \Xibwd_1 \\
        U_{1} & \Xibwd_0 & 0 \\
        Y_{1} & 0 & \Xifwd_2
    \end{NiceArray} \right)$} \breve Q_1, \quad
    \scalebox{0.87}{$\left( \begin{NiceArray}{c|c}[columns-width=2.0em]
        \tilde L_{3} & 0 \\
        \tilde U_{3} & \Xibwd[1]_2 \\
        {\color{gray}\tilde Y_{3}} & {\color{gray}\Xifwd[1]_0}
    \end{NiceArray} \right)$} =
    \scalebox{0.87}{$\left( \begin{NiceArray}{c|cc}[columns-width=2.0em]
        L_{3} & \Xifwd_3 & \Xibwd_3 \\
        U_{3} & \Xibwd_2 & 0 \\
        {\color{gray}Y_{3}} & {\color{gray}0} & {\color{gray}\Xifwd_0}
    \end{NiceArray} \right)$} \breve Q_3, \\[0.6em]
    \rotatebox[origin=c]{90}{\underline{Level 1}}\quad && \scalebox{0.87}{$\left( \begin{NiceArray}{c|c}[columns-width=2.14em]
        \tilde L_{2} & 0 \\
        \tilde U_{2} & \Ups[2]_0
    \end{NiceArray} \right)$} &=
    \scalebox{0.87}{$\left( \begin{NiceArray}{c|cc}[columns-width=2.14em]
        L_{2} & \Xifwd[1]_2 & \Xibwd[1]_2 \\
        U_{2} & \Xibwd[1]_0 & {\color{gray}\Xifwd[1]_0}
    \end{NiceArray} \right)$} \breve Q_2, \\[0.5em]
    \rotatebox[origin=c]{90}{\underline{Level 2}}\quad && \scalebox{0.87}{$\left( \begin{NiceArray}{c|c}[columns-width=2.42em]
        \tilde L_{0} & 0
    \end{NiceArray} \right)$} &=
    \scalebox{0.87}{$\left( \begin{NiceArray}{c|c}[columns-width=2.34em]
        L_{0} & \Ups[2]_0
    \end{NiceArray} \right)$} \breve Q_0.
\end{align*}
The orthogonality metrics of the hyperbolic Householder transformations $\breve Q_i$
are determined by $\Ssign_c$ \eqref{eq:ups-sigma-upd-riccati}:
\begin{align*}
    &\breve Q_1 \text{ and } \breve Q_3 \text{ are } \scalebox{0.87}{$\begin{pmatrix}
        -\I \\ &\Ssign_1 \\ &&\Ssign_2
    \end{pmatrix}$}\text{- and }
    \scalebox{0.87}{$\begin{pmatrix}
        -\I \\ &\Ssign_3 \\ &&\Ssign_0
    \end{pmatrix}$}\text{-orthogonal}, \\
    &\breve Q_2 \text{ and } \breve Q_0 \text{ are } \scalebox{0.87}{$\begin{pmatrix}
        -\I \\ &\Ssign_3 \\ &&\Ssign_0 \\
        &&&\Ssign_1 \\ &&&&\Ssign_2
    \end{pmatrix}$}\text{-orthogonal}.
\end{align*}
As discussed above, the blocks $\Xifwd[l]_0$ shown in gray are zero in the absence of coupling between the first and final stages.
They are included here because the same structure generalizes to larger problems.
Levels $l=\log_2\! P$ and $l=\log_2\! P - 1$ are special cases (at these levels, the system has been reduced to a block $2\times 2$ matrix), while all other
levels $l < \log_2\! P - 1$ follow the same pattern as shown in \eqref{eq:update-cr-building-block}.
The resulting procedure is summarized in \Cref{alg:fact-parallel-riccati-cr-update},
which follows the same structure as the factorization in \Cref{alg:fact-parallel-riccati-cr}.
A graph representation of the algorithm is shown in \Cref{fig:cyclic-reduction-update-graph-16}.
Further implementation details such as the data structures used to store the update matrices are discussed in \Cref{app:data-structures}.
For brevity, the factorization and update algorithms listed here assume that
there is no coupling between the first and final stages;\,\footnote{%
This allows levels $\log_2\! P - 1$ and $\log_2\! P$ to be handled in the same way as earlier levels, simplifying the pseudocode.
In vectorized implementations of \cyqlone{}, the last levels are handled using PCR or PCG instead of CR, avoiding the block $2\times 2$ case entirely.}
they can easily be extended to include such coupling by introducing a special case for the last two levels of the algorithm, as described in \Cref{app:gen-cr-periodic}.

\begin{algorithm2e}[p]
    \caption{\cyqlone{} factorization updates}
    \label{alg:fact-parallel-riccati-cr-update}
    \DontPrintSemicolon
    \setstretch{1}
    \KwIn{$A, B, C, D$ and $N$: OCP data matrices and horizon length}
    \KwIn{$P$: number of intervals to partition the horizon into}
    \KwIn{$\L{\mathscr R}$ ($\L{R}, \L{S}, \L{Q}$), $\L{\mathscr A}$ ($\L{B}, \L{A}, \Acl$), $\L{\mathscr M}$ ($L, Y, U$): existing factorization of $\mathscr K\!$}
    \KwOut{$\tilL{\mathscr R}$ ($\tilL{R}, \tilL{S}, \tilL{Q}$), $\tilL{\mathscr A}$ ($\tilL{B}, \tilL{A}, \tilAcl$) and $\tilL{\mathscr M}$ ($\tilde L, \tilde Y, \tilde U$): factorization of $\widetilde{\mathscr K}\; \mathrlap{\eqref{eq:cyqlone-K-tilde}}$}
    \vspace{0.5em}
    \For{$c=0,...,P-1$ {\normalfont (in parallel)}}{\label{ln:loop-processors-update}
        \vspace{0.1em}
        \tikzmark{xmark}
        \tikzmark{updateXi0}
        $\Xifwd[0]_c$, $\Xibwd[0]_{c-1}$, $\Ssigncr[0]_c = {}$\textsc{update--block--column--riccati}$(c)$
        \tikzmark{updateXi1}
        \ColorBarTwo{colorBd}{colorAd}{updateXi0}{updateXi1}
        \mycomment{\rmrk{steps \ref{it:bchol-diag} and \ref{it:solve-subdiag}}}
        \vspace{0.1em}
        \textsc{update--schur}$(c)$
        \mycomment{\rmrk{step \ref{it:cr-schur}}}
    }
    \vspace{0.5em}
    \Fn{\normalfont\textsc{update--schur}$(c)$}{
        --- sync --- \mycommentnofill{\rmrknp{Wait for $\Xifwd[0]$, $\Xibwd[0]$}}\;
        \textbf{if} $\nu^P_2(c) = 0$: \quad
            \textsc{update--L}$(0, c)$\;
        \For{$l = 0,...,\log_2(P)-1$\mycommentnofill{\rmrknp{Recursion level of CR}}} {
            $i_U = c + 1$, \quad
            $i_Y = c + 1 - 2^l$
            \;
            --- sync --- \mycommentnofill{\rmrknp{Wait for $\breve Q$}}\;
            \textbf{if}\phantom{\textbf{el}} $\nu^P_2(i_U) = l$: \quad
                \textsc{update--U}$(l, i_U)$\;
            \textbf{elif} $\nu^P_2(i_Y) = l$: \quad
                \textsc{update--Y}$(l, i_Y)$\;
            --- sync --- \mycommentnofill{\rmrknp{Wait for $\Xifwd[l]$, $\Xibwd[l]$}}\;
            \textbf{if} $\nu^P_2(i_Y) = l+1$: \quad
                \textsc{update--L}$(l+1, i_Y)$\;
        }
    }
    \vspace{0.2em}
    \Fn{\normalfont\textsc{update--L}$(l, i)$}{
        \tikzmark{xmark}
        $\Ssigncr[l+1]_i \assignpc \blkdiag\big(\Ssigncr[l]_i, \;\Ssigncr[l]_{i + 2^l}
        \big), \quad m_i = \operatorname{rank} \Ssigncr[l+1]_i$\;
        \tikzmark{updateL0}
        $\begin{pNiceArray}{c|c}
            \tilde L_i & 0
        \end{pNiceArray}
        \assignpc
        \begin{pNiceArray}{c|cc}
            L_i & \Xifwd[l]_i & \Xibwd[l]_i
        \end{pNiceArray} \breve Q_i$, where $\breve Q_i$ is $\left(\!\begin{smallmatrix}
            -\I \\ & \phantom-\Ssigncr[l+1]_i
        \end{smallmatrix}\!\right)$-orth. \cite{pas_blocked_2025-1}
        \tikzmark{updateL1}
        \ColorBarTwo{colorLd}{colorQd}{updateL0}{updateL1}
        \mycomment{\texttt{\rmrk{hyh}} $\hspace{0.2em} m_in_x^2\hspace{-1em}$}
    }
    \vspace{0.2em}
    \Fn{\normalfont\textsc{update--U}$(l, i)$}{
        \tikzmark{xmark}
        \tikzmark{updateU0}
        $\begin{pNiceArray}{c|c}
            \tilde U_i & \Xibwd[l+1]_{i-2^l}
        \end{pNiceArray} \assignpc
        \begin{pNiceArray}{c|cc}
            U_i & \Xibwd[l]_{i-2^l} & 0
        \end{pNiceArray} \breve Q_i$
        \tikzmark{updateU1}
        \ColorBarTwo{colorUd}{colorAd}{updateU0}{updateU1}
        \mycomment{\texttt{\rmrk{hyh apply}} \quad $\hspace{0.35em} 2m_in_x^2 \hspace{-1em}$}
    }
    \vspace{0.2em}
    \Fn{\normalfont\textsc{update--Y}$(l, i)$}{
        \tikzmark{xmark}
        \tikzmark{updateY0}
        $\begin{pNiceArray}{c|c}
            \tilde Y_i & \Xifwd[l+1]_{i+2^l}
        \end{pNiceArray} \assignpc
        \begin{pNiceArray}{c|cc}
            Y_i & 0 & \Xifwd[l]_{i+2^l}
        \end{pNiceArray} \breve Q_i$
        \tikzmark{updateY1}
        \ColorBarTwo{colorYd}{colorBd}{updateY0}{updateY1}
        \mycomment{\texttt{\rmrk{hyh apply}} \quad $\hspace{0.35em} 2m_in_x^2 \hspace{-1em}$}
    }
    \vspace{1em}
    This algorithm parallels \Cref{alg:fact-parallel-riccati-cr}. \quad
    Color coding matches \Cref{fig:cyclic-reduction-update-graph-16}.\;
    \vspace{0.2em}
\end{algorithm2e}

Thanks to the similar structure between the factorization (\Cref{fig:cyclic-reduction-graph-16})
and factorization update procedures (\Cref{fig:cyclic-reduction-update-graph-16}),
it is possible to switch from factorization updates to re-factorizations at any level in the algorithm.
This is useful when the sizes $m_i$ of the update matrices become large relative to the block size $n_x$ at a certain level in the recursion,
making re-factorization -- with cost $\landauO(n_x^3)$ -- more economical than factorization updates -- with cost $\landauO(m_i n_x^2)$.

\section{Vectorization} \label{sec:vectorization}
Thanks to the parallelism exposed by the \cyqlone{} permutation from
\Cref{sec:cyqlone}, operations applied to different stages of the OCP can not
only be parallelized across different processors, but they can also be \textit{vectorized}.
Vectorization refers to the process of transforming an algorithm into a form
where operations are carried out on arrays or vectors rather than on individual
values. Consider the example of matrix-vector multiplication $Ax$: a scalar
implementation of this algorithm involves multiplying the individual elements of
each row of the matrix $A$ by the corresponding elements of the vector $x$.
One possible vectorized approach could be to compute the sum of the columns of
$A$, weighted by the elements of $x$:
\begin{align}
    Ax &=
    \underbrace{
    \scalebox{0.85}{$
    \begin{pmatrix}
        a_{11}x_1 + a_{12}x_2 + a_{13}x_3 + a_{14}x_4 \\
        a_{21}x_1 + a_{22}x_2 + a_{23}x_3 + a_{24}x_4 \\
        a_{31}x_1 + a_{32}x_2 + a_{33}x_3 + a_{34}x_4 \\
        a_{41}x_1 + a_{42}x_2 + a_{43}x_3 + a_{44}x_4
    \end{pmatrix}$}
    }_{\text{Scalar}}
    \label{eq:matvec-scalar}
    \displaybreak[1]\\
    &= \scalebox{0.85}{$\begin{pmatrix}
        | & | & | & | \\
        a_1 & a_2 & a_3 & a_4 \\
        | & | & | & | \\
    \end{pmatrix}\!\!
    \begin{pmatrix}
        x_1 \\
        x_2 \\
        x_3 \\
        x_4
    \end{pmatrix}$} =\,
    \underbrace{
    \scalebox{0.85}{$a_1 x_1 + a_2 x_2 + a_3 x_3 + a_4 x_4$}
    }_{\text{Vectorized}}\!.
    \label{eq:matvec-vectorize}
\end{align}
The advantages of vectorization become apparent when considering that modern
processors implement instructions that operate on an entire vector at once,
referred to as \textit{single instruction, multiple data} (SIMD). Examples
of SIMD instruction set architecture extensions include AVX2, AVX-512, AVX10,
NEON and ARM SVE. Modern SIMD instructions often have the same or a similar
throughput as their scalar counterparts, and operate on vector registers whose size
is a small power of two, referred to as the \textit{vector length} (e.g. four
double-precision elements per vector for AVX2, or eight for AVX-512).
The \textit{lane} of a scalar value refers to its position within a vector register.
The first lane is lane zero, and the index of the highest lane is one less than
the vector length. Many SIMD instructions operate on all lanes in parallel
(e.g. an element-wise sum of two vector registers), but
instructions that perform horizontal reductions across lanes or that permute
the values between lanes are also available.
In ideal cases, a vectorized implementation using SIMD can achieve
speedups by a factor of the vector length: the scalar variant of the
matrix-vector example in \eqref{eq:matvec-scalar} requires 16 multiplication instructions, whereas the
vectorized implementation requires just four vector multiplication instructions
(assuming a vector length of four).
Furthermore, the number of elements that can be stored in fast registers inside
of the processor generally also increases when using SIMD instructions, possibly
resulting in lower memory traffic and improved arithmetic intensity.
Finally, thanks to the higher throughput, the lower number of instructions, and
the reduced and coalesced memory traffic, implementations using SIMD may consume less
energy than equivalent scalar variants \cite{moldovanova_automatic_2017,inoue_how_2016}.

An important step in the vectorization process is so-called
\textit{strip-mining} \cite{weiss_strip_1991}, which involves transforming
specific loops in the algorithm into loops with fewer iterations
that operate on chunks with the same size as the platform's SIMD vector length.
In general-purpose implementations of linear algebra algorithms, a common
strategy is to apply strip-mining to the loops that iterate over the rows within
a given column of a matrix, processing multiple elements in the same column
at once using SIMD instructions (similar to \eqref{eq:matvec-vectorize}).
In the following subsection, we discuss the limitations of this row-wise vectorization strategy,
and propose vectorization across different matrices in a batch as
a more performant alternative.

\subsection{Limitations of classical row-wise vectorization} \label{subsec:stripmining}

Most vectorized general-purpose linear algebra libraries such as BLAS, LAPACK,
Eigen and BLASFEO make use of row-wise vectorization where multiple adjacent
elements in the same column are loaded in a single SIMD register.
If the number of rows of the matrix is not an integer multiple of the vector
length, the remainders at the bottom of each column do not fill an entire SIMD
register, leading to poor vector lane utilization.\,%
\footnote{The remainder could be processed using scalar operations,
or using operations with successively smaller vector lengths, e.g. a
remainder of three could be implemented using a SIMD instruction with vector
length two and a scalar instruction rather than three scalar instructions.
Alternatively, padding could be used to round up the number of rows to the next
multiple of the vector length. However, this places the burden of properly padding the input
data on the user of the software. If the target architecture supports it
(e.g. AVX-512 or ARM SVE), masked or predicated SIMD instructions could be used
to avoid padding or scalar loops.}
When dealing with large matrices, the fraction of time taken up by the
processing of this remainder is negligible. However, for small matrices, this
part of the algorithm may take up a significant portion of the overall run time.
Furthermore, triangular matrices --- which are important workhorses in numerical
linear algebra, used in Cholesky factorizations, QR factorizations,
block-Householder reflectors, etc. --- have the property that the number of
nonzero elements differs in each column. As a result, there will always be a
remainder in most of the columns, even if the number of rows of the matrix is a
multiple of the vector length. For example, if the matrix $A$ in \eqref{eq:matvec-vectorize}
is lower triangular, 6 of the 16 multiplications are redundant multiplications by zero,
reducing the efficiency of a vectorized implementation.

In addition to the loss of efficiency because of nonuniform remainders,
vectorization can also be prevented by inherent dependencies encountered in
certain linear algebra operations. For example, Cholesky factorization requires
computing the inverse square root of the pivot element, multiplying the subdiagonal
elements by this value, and then subtracting the symmetric outer product of the
result from the trailing submatrix to obtain the Schur complement.
The top left element of the Schur
complement is then used as the next pivot (cf. \S\ref{sec:block-chol}). Because of the dependency of the
second pivot on the first pivot, it is not possible to compute the inverse
square roots of multiple pivots at once, preventing the vectorization of this
operation. Vectorization can still be used to multiply multiple subdiagonal
elements by the inverse square root of the pivot at once, but the
scalar divisions and square roots dominate the run time for small matrices.\,\footnote{The
latency and reciprocal throughput of division and (inverse) square root
instructions are much higher than for multiplication and addition instructions.}
For the Cholesky factorization specifically, vector lane utilization is further
reduced because of the triangular structure of the resulting Cholesky factor and
the symmetry of the Schur complement, as discussed in the previous paragraph.

\subsection{Batch-wise vectorization}

When the same operation is applied to multiple matrices of the same size
(as is the case in \Cref{alg:fact-parallel-riccati-cr}),
the limitations discussed in the previous section can be avoided by vectorizing
across different matrices rather than along the rows of a single matrix.
We will refer to the strip-mining of the loop that iterates over the different matrices in the batch
as \textit{batch-wise} vectorization.
SIMD instructions are then applied to the elements of just one row of the
matrices at a time,
and there are no dependencies between elements of different matrices in the batch.
For example, when applying batch-wise vectorization to the Cholesky factorization,
the inverse square roots of the first pivots of all matrices in the batch can
be computed at once, instead of one at a time (as was the case for the row-wise vectorization case discussed earlier).

Furthermore, the batch-wise vectorization of batched linear algebra routines is also
easier to implement than row-wise vectorization: the scalar versions of
these routines can simply be applied element-wise to tuples containing the
values for the corresponding elements of all matrices in the batch. In contrast,
row-wise vectorization requires some operations to be carried out using scalar
arithmetic, scalars may need to be broadcast to a full vector register, and
masked operations are needed for the remainders.

A downside of batch-wise vectorization is that the size of the working set is
larger than for the case where a single matrix is processed at a time.
Since all matrices in the batch are processed simultaneously, the cache
utilization is higher, which affects performance for large matrices. %
In optimal control applications, the sizes of the matrices are determined
by the number of states and the number of controls, which are generally
relatively small. Consequently, cache utilization is less of a concern than optimal vector
lane utilization, especially on modern hardware, and we find that batch-wise vectorization is highly
effective for this application.

\subsubsection{Compact batched matrix storage format} \label{subsec:compact}

\begin{figure}[h]
    \centering
    \includegraphics[width=0.75\textwidth]{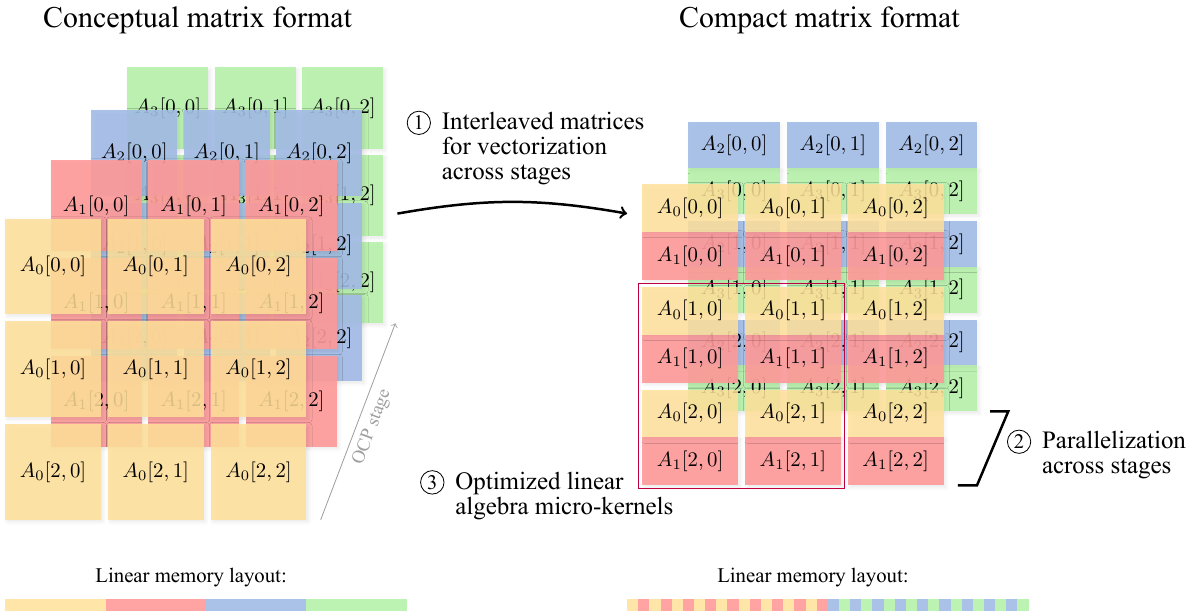}
    \caption{Visual representation of the conversion to compact storage format for four 3-by-3 matrices. Each batch of two matrices is interleaved to match a vector length of two, enabling DLP using SIMD (1). Multiple batches can be processed in parallel by different processors, exploiting TLP using multi-threading (2). Finally, the compactly stored matrices can be processed using specialized batched linear algebra routines with micro-kernels that maximize ILP (3).}
    \label{fig:compact-storage}
\end{figure}

SIMD load and store operations are most efficient when the elements of the SIMD
vector are stored contiguously in memory. When using row-wise vectorization,
this is achieved by storing the matrix in \textit{column-major order}, which is
a commonly used storage format for dense matrices.
When using batch-wise vectorization, e.g. by strip-mining along the time dimension of an OCP,
contiguous storage can be achieved by storing the matrices in memory in such a
way that elements of matrices of different stages are interleaved.
We refer to this interleaved matrix storage format as \textit{compact storage}.
Step 1 of \Cref{fig:compact-storage} visualizes the conversion to a compact
storage format for four 3-by-3 matrices $A_0$ through $A_3$ and a vector length
of two.

Practically, the most common and intuitive way to store $N$ different $m$-by-$n$ matrices
is to store the matrices $A_k$ back-to-back, each in column-major order. We will
refer to this format as the \textit{conceptual storage format}.
Using zero-based indexing, the element at row $r$ and column $c$ of the $j$-th
matrix would then be stored in memory at offset $r + mc + mnj$.
In contrast, the offset of the same element in compact storage with vector length
$v$ would be $j\rem v + vr + vmc + vmn \lfloor j / v \rfloor$.
The conceptual storage format can be written as $(m,n,N){:}(1,m,mn)$ in CuTe notation \cite{cris_cecka_cute_nodate},
while the CuTe layout of the compact storage format can be written as
$\big(m, n, (v, \lceil N/v\rceil)\big){:}\big(v, vm, (1, vmn)\big)$.

\subsection{Vectorized modified Riccati recursion \rmrk{steps \ref{it:bchol-diag} and \ref{it:solve-subdiag}} \done} \label{subsec:vec-mod-ric}

\begin{figure}
\centering
\begin{minipage}{0.495\textwidth}
    \includegraphics[width=1\textwidth]{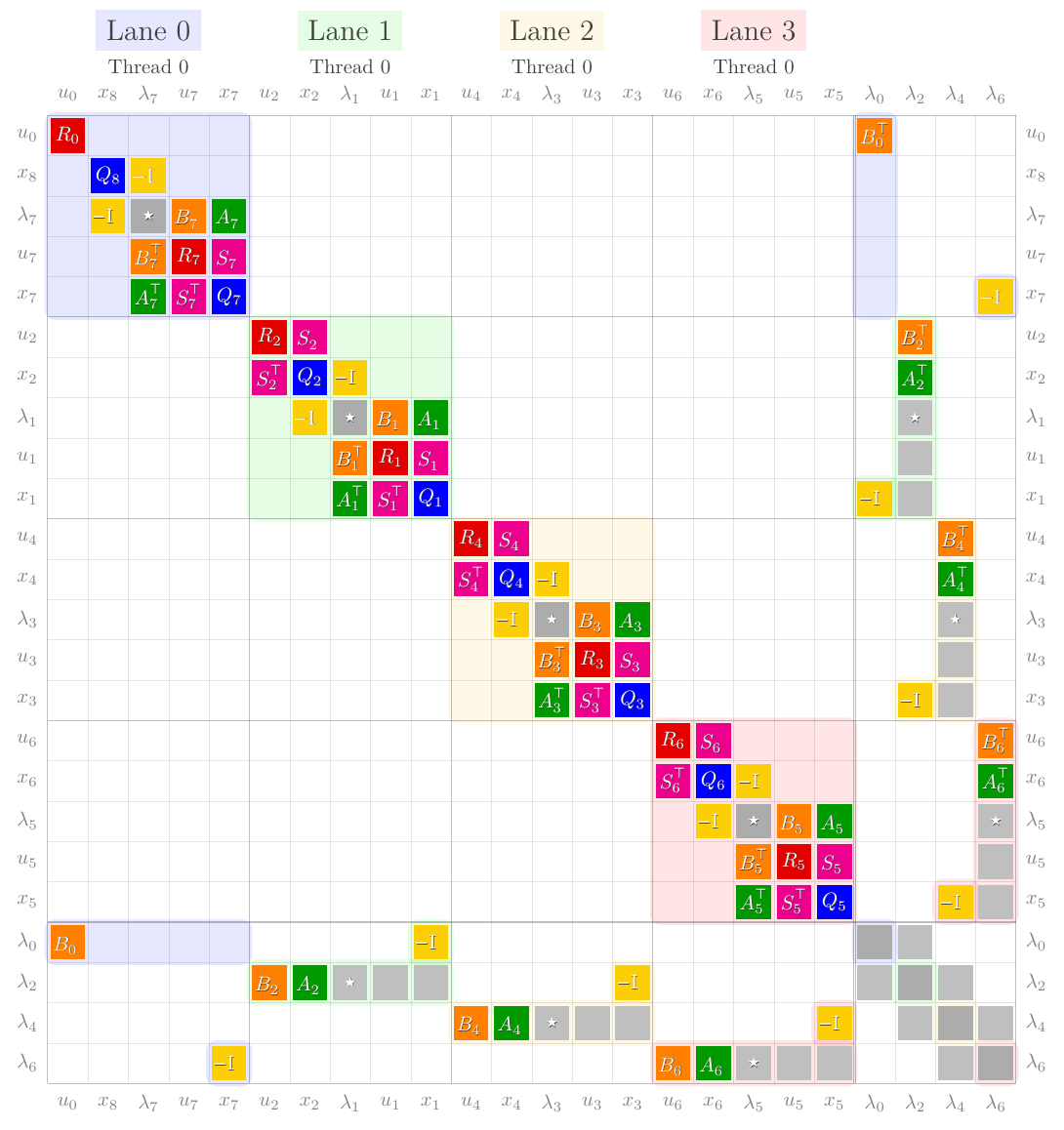}
\end{minipage}
\begin{minipage}{0.495\textwidth}
    \includegraphics[width=1\textwidth]{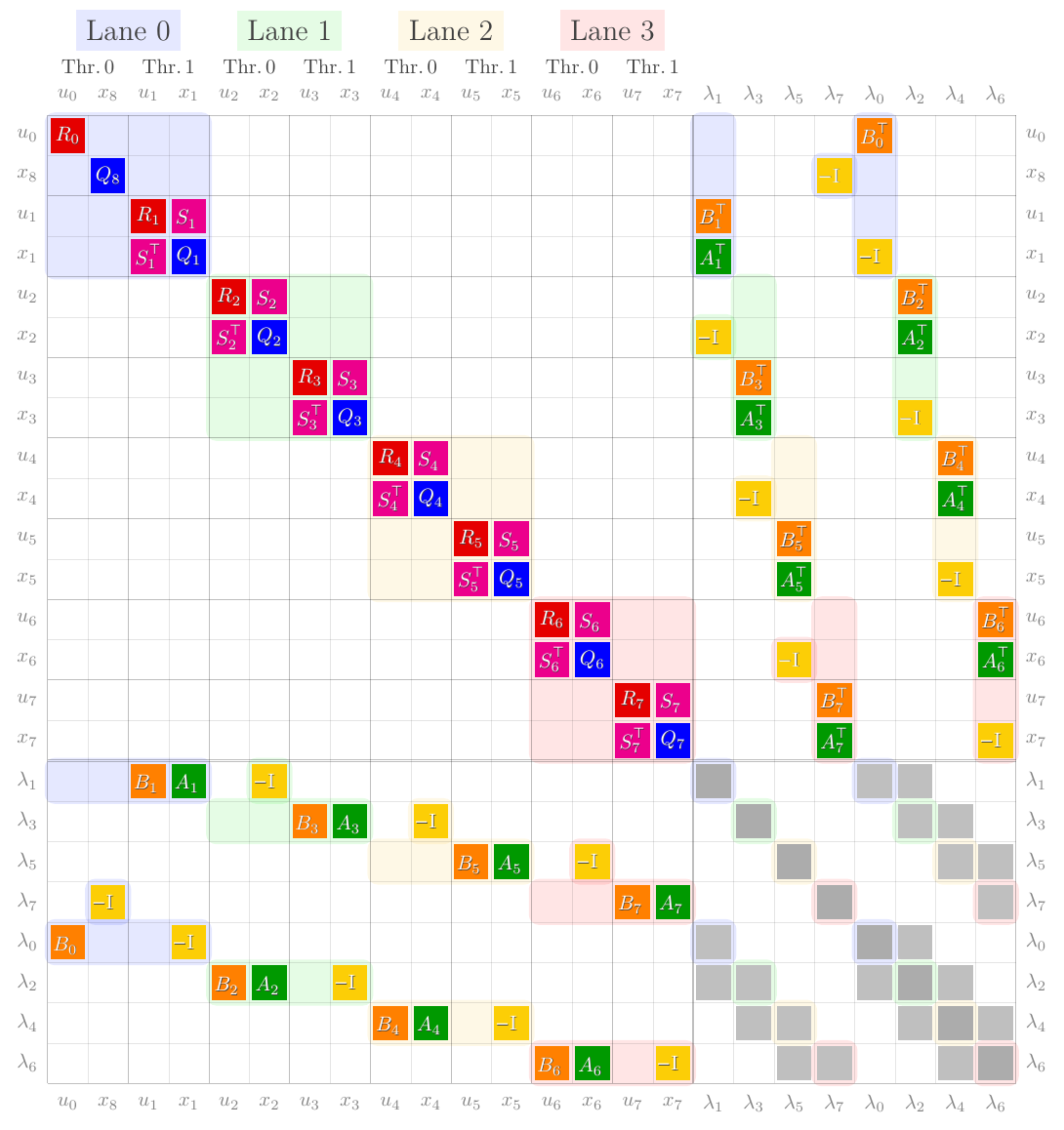}
\end{minipage}
\vspace{-0.8em}
\caption{Vectorization of the modified Riccati recursion (\Cref{alg:fact-mod-riccati}) with vector length $v=4$ on $p=1$ (left) and $p=2$ (right) processors.
The distribution of the data across the four vector lanes is indicated by a blue, green, yellow or red background.}
\label{mat:cyqlone-vec-ric}
\vspace{-0.5em}
\end{figure}

The batch-wise vectorization of \Cref{alg:fact-mod-riccati} is relatively straightforward,
because it closely resembles the existing parallel implementation:
First, the horizon is partitioned into $v$ partitions, where $v$ is the available
vector length. All stages within each partition are assigned to the same vector
lane, corresponding to strip mining across the loop of \cref{ln:loop-processors}
in \Cref{alg:fact-parallel-riccati-cr}.
\Cref{mat:cyqlone-vec-ric} illustrates how this corresponds to distributing the
block columns of the matrix $\mathscr K$ across different vector lanes: All blocks
within a given block column $c$ are assigned to the same lane, highlighted by the four background colors.
Since the control flow in \Cref{alg:fact-mod-riccati} does not depend on the
selected interval $c$, and since there are no dependencies between different intervals,
all operations can be effectively vectorized.
Second, after partitioning the horizon according to the vector length $v$,
further subdivide these intervals to distribute them across $p$ available
processors. The total parallelism is then $P=v\cdot p$.
We will refer to these two levels of partitions as the $v$-partition and 
$p$-partition, respectively.
An example of this type of partitioning with $v=4$ and $p=2$ is shown
in the right matrix of \Cref{mat:cyqlone-vec-ric}:
stage 0/8 is processed by lane 0 of processor 0, stage 1 is processed by lane 0
of processor 1; stage 2 is processed by lane 1 of processor 0, stage 3 is
processed by lane 1 of processor 1, etc.

\subsection{Vectorized Schur complement computation \rmrk{step \ref{it:compute-schur}}} \label{subsec:vec-schur-comp}

The evaluation of the blocks of the Schur complement (see \Cref{sec:compute-schur}) may involve data from
stages in different vector lanes. For reasons that will become apparent in the following section, we do not yet
use an odd--even permutation of $\lambda^0,\lambda^2,\lambda^4,\lambda^6$,
resulting in a $4\times 4$ block-tridiagonal Schur complement matrix as shown in \Cref{mat:cyqlone-vec-ric}:
\begin{equation}
    -\mathscr K / \mathscr R = \scalebox{0.87}{\setstretch{1.5}$\begin{pNiceArray}{cccc}[first-col,first-row,columns-width=2em]
    & \lambda^0 & \lambda^2 & \lambda^4 & \lambda^6 \\
    \lambda^0 & \Mdiag_0 & \tp \Ksub_0 & & \\
    \lambda^2 & \Ksub_0 & \Mdiag_1 & \tp \Ksub_1 & \\
    \lambda^4 & & \Ksub_1 & \Mdiag_2 & \tp \Ksub_2 \\
    \lambda^6 & & & \Ksub_2 & \Mdiag_3
    \end{pNiceArray}$}.
\end{equation}

The blocks $\Mdiag_i$ and $\Ksub_i$ are stored in vector lane $i$, although the blocks used to compute them may be stored in different lanes.
For example, consider the blocks in column
$x^1$ and rows $\lambda^0,\lambda^2$ in
the left matrix of \Cref{mat:cyqlone-vec-ric}.
The result of the product $\Ksub_0 \defeq \KbwdT_1 = -\L{A}_1 \Li Q_1$ (where both of the operands
are stored in lane 1) ends up in the first subdiagonal block of the Schur
complement (row $\lambda^2$, column $\lambda^0$), which is stored in lane 0. Similarly, the operands of
$\Mbwd_0 = \Lit Q_1 \Li Q_1$ are in lane 1, but the result
is added to the first diagonal block $\Mdiag_0 \defeq \Mbwd_0 + \Mfwd_0$ of the Schur complement (row $\lambda^0$, column $\lambda^0$),
which is stored in lane 0.
On the other hand, the sum of the symmetric products
$\L{B}_3 \Lt{B}_3 + \L{B}_2 \Lt{B}_2 + \L{B}_1 \Lt{B}_1 + \L{A}_1 \Lt{A}_1$
is stored in the same lane as the operands, so no special handling is required for this operation.

In general, we see that the results of operations involving the blocks $\Lit Q_{k_1}$ of the first $p$-partition
within each $v$-partition end up in a different vector lane. Practically, this is achieved by performing a
lane-wise rotation before storing the results of these operations to memory.
Such rotations can be implemented without significant overhead using one of the various lane permutation
instructions available in modern SIMD ISA extensions.

\subsection{Vectorized cyclic reduction \rmrk{step \ref{it:cr-schur}} \done} \label{subsec:vec-cr}

Historically, CR has been used on vector machines \cite{lambiotte_solution_1975,kershaw_solution_1982}, shared \cite{hirshman_bcyclic_2010} and distributed memory multicore systems \cite{krechel_parallelization_1990,reale_tridiagonal_1990}, and more recently on GPUs \cite{zhang_fast_2010}.
In this section, we combine both vectorization and shared-memory multithreading to increase the available parallelism on modern consumer CPUs.
To the best of the authors' knowledge, this approach has not yet been described in the literature and is not available in any open-source CR software packages.

The basic idea is simple: At each level $l$ in the CR algorithm,
$N/2^{l+1}$ columns are factorized and eliminated. The operations required are
the same for each column, so vectorization and parallelization across these columns is straightforward (similar to \Cref{subsec:vec-mod-ric}).
However, updating the remaining even equations at each level is more involved,
and requires vector lane crossings (as in \Cref{subsec:vec-schur-comp}).
Furthermore, because of the halving of the problem size at each level, strides and indices change in nontrivial ways.
In the final levels of CR, the size of the reduced system drops below the vector length, leading to poor resource utilization if not handled as a special case. 

We will now consider the CR method in isolation, independent of any optimal control
structure of the original system.
After all, once the Schur complement has been evaluated as described by \Cref{subsec:vec-schur-comp},
we are left with a general block-tridiagonal matrix.
The graphical representation of the CR algorithm applied to a $16\times 16$
block-tridiagonal matrix in \Cref{fig:cyclic-reduction-graph-16} will be used
as the main visual tool to motivate the approach. It visualizes the blocks of
the original system and the intermediate reduced systems
($\Mdiag^{(l)}_i$ and $\Ksub^{(l)}_i$), and the blocks of the resulting block Cholesky
factor ($L_i$, $U_i$ and $Y_i$). The colors of these blocks match the ones used in \Cref{alg:fact-parallel-riccati-cr}.
Additionally, \Cref{fig:cyclic-reduction-graph-16} shows which of the four vector lanes
on which of the four threads is used to compute each block, with arrows to
indicate data dependencies between the blocks.
We will refer to matrices that belong to the four vector lanes within the same
thread as a \textit{batch}. For example, matrices $L_1, L_5, L_9$ and $L_{13}$
form a batch (mapped to thread 1).
The structure of the scheduling in \Cref{fig:cyclic-reduction-graph-16}
is chosen in such a way that all matrices in the same batch
always undergo the same operation. The following sections explain the derivation
of this scheduling in more detail.

\subsubsection{Assignment of matrices to batches for effective vectorization \done}

In order to derive a practical scheduling and vector lane assignment scheme that allows
efficient application of vectorized operations, let us first consider the level in the
recursion where the number of block columns being eliminated is equal to the
vector length.
For a $16\times 16$ block matrix and a vector length of four, this
happens at the second level, as shown in \Cref{fig:cyclic-reduction-graph-16}
(at the dividing line between \textit{Batched} and \textit{Scalar} operations).
This level is the last level at which the same operations are applied to at least
four blocks, and we impose that all four blocks being computed (e.g.
$\Mdiag^{(2)}_0$, $\Mdiag^{(2)}_4$, $\Mdiag^{(2)}_8$ and $\Mdiag^{(2)}_{12}$ in the figure)
belong to a single batch, so that they are processed by different vector lanes.
This choice minimizes unused lanes and thus
maximizes the number of levels that can be processed using fully vectorized operations.
From there on, the lane assignment of the remaining blocks is performed by
working upwards to the earlier levels,
introducing additional batches in such a way that blocks with dependencies
between them are stored in the corresponding lanes of their respective batches.
This is not possible for every operation in the CR algorithm: In particular, the subtraction of
$Y_k \tp Y_k$ from the diagonal blocks necessarily crosses over into different
lanes for the last batch of each level. In \Cref{fig:cyclic-reduction-graph-16},
the arrows from blocks $Y_3, Y_7$ and $Y_{11}$ to the subsequent diagonal blocks
cross the thick gray lines that separate the vector lanes.
These cross-lane operations are limited to a single operation on a single
batch per level, and will be handled in the same way as the cross-lane products
in \Cref{subsec:vec-schur-comp}.
To implement vectorized CR in practice, it is convenient to use periodic CR as described in
\Cref{app:cr-periodic-vectorization}.

\subsubsection{Handling of the final scalar levels \done} \label{subsec:vec-cr-last}

Once the size of the reduced system becomes less than or equal to the vector length,
it is no longer possible to fill complete batches, as there are simply not
enough matrices that undergo the same operation. At this point, one could
process the final levels in the recursion using non-batched
linear algebra routines on individual matrices, using vectorization along the
row dimension.

\begin{figure}
    \centering
    \includegraphics[width=0.85\textwidth]{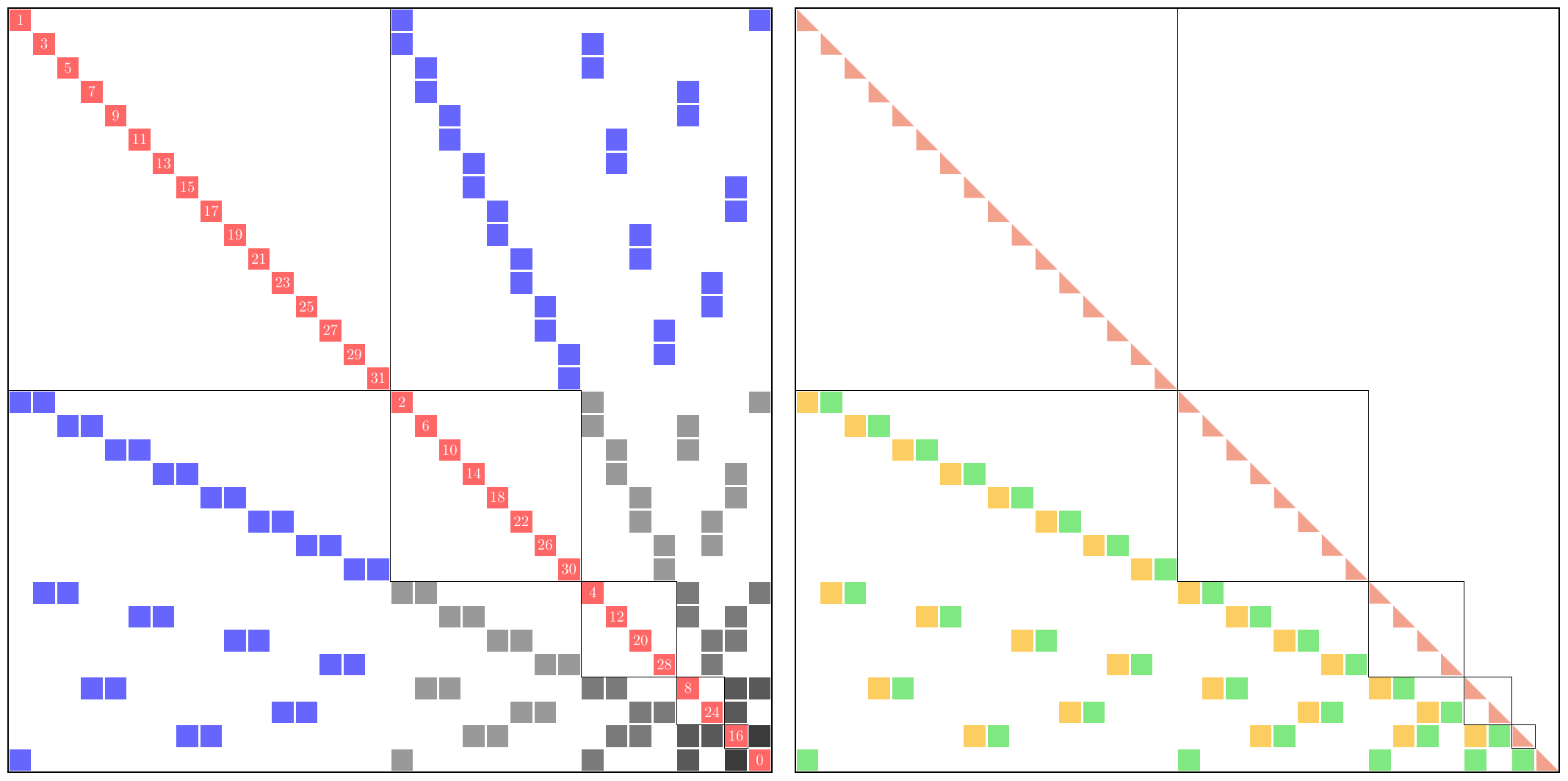}
    \caption{Permutation of a $32\times 32$ block-tridiagonal matrix that enables CR down to a single block (left) and its Cholesky factor (right). Color coding matches \Cref{fig:cyclic-reduction-permutation}.}
    \label{fig:cyclic-reduction-permutation-labels-32-vl1}
\end{figure}
\begin{figure}
    \centering
    \includegraphics[width=0.85\textwidth]{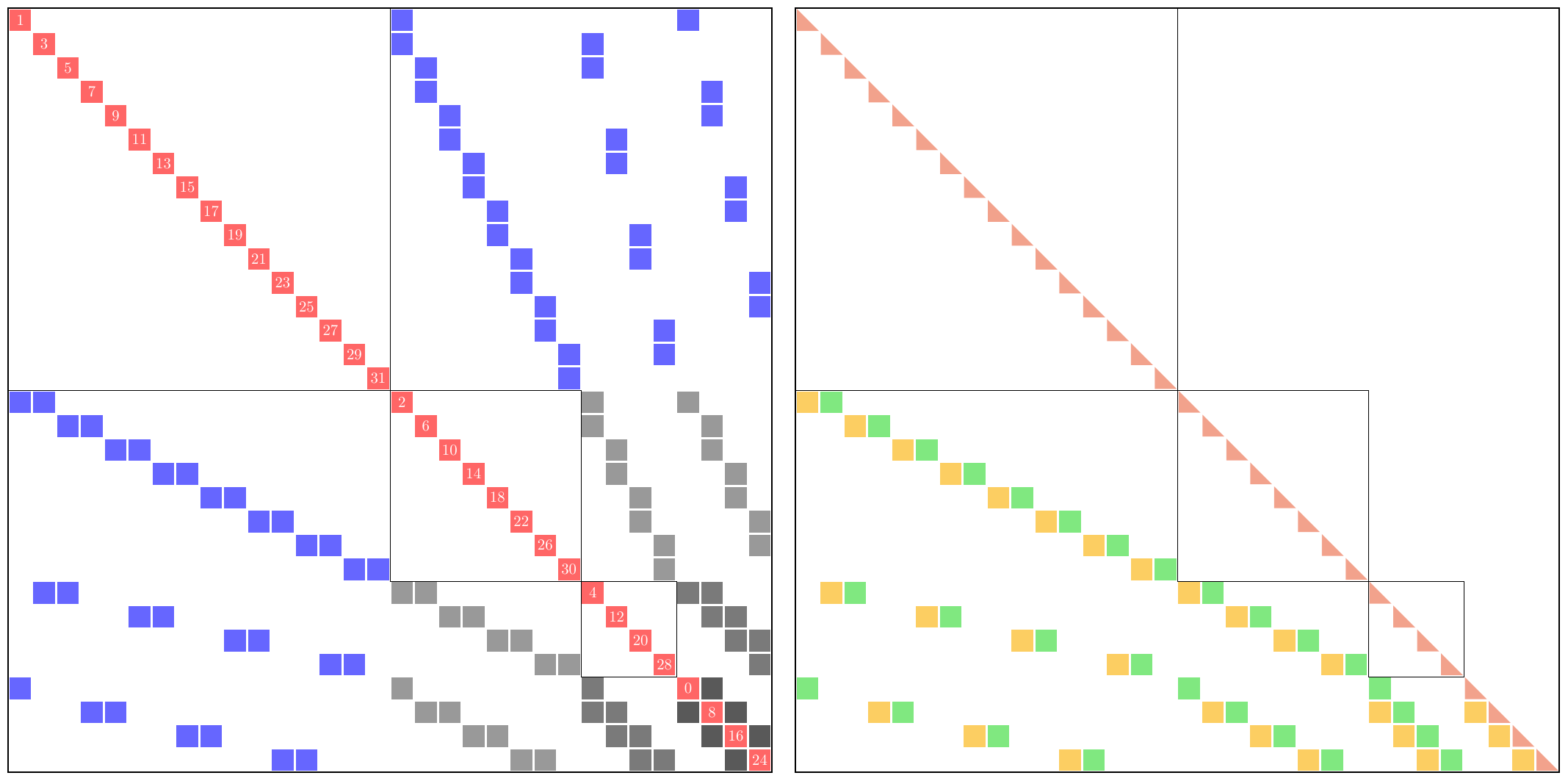}
    \caption{Permutation of a $32\times 32$ block-tridiagonal matrix that enables CR down to an unpermuted $4\times 4$ block-tridiagonal matrix (left) and its Cholesky factor (right). Color coding matches \Cref{fig:cyclic-reduction-permutation}. Unlike \Cref{fig:cyclic-reduction-permutation-labels-32-vl1}, this permutation allows the use of batched operations to compute the blocks of fill-in the last four block columns. Additionally, it leaves the final block matrix in the format of \eqref{eq:cr-pcg}.}
    \label{fig:cyclic-reduction-permutation-labels-32-vl4}
\end{figure}

However, switching to
a scalar implementation of CR may sometimes be suboptimal:
converting from the compact storage
format to back to the conceptual one takes time, and linear algebra operations on individual
matrices are less performant than batched operations (for small matrices).
Performing scalar operations on matrices in compact storage format is also an option,
but this makes it impossible to use coalesced memory accesses, and requires specialized micro-kernels. To address these drawbacks,
we consider two alternatives to scalar CR for the final levels:
a variant of CR known as \textit{parallel cyclic reduction} (PCR)
and a preconditioned conjugate gradient method (PCG).

PCR \cite{sweet_parallel_1988-1,zhang_fast_2010} %
addresses the poor utilization in the final levels of CR as follows:
instead of reducing an $N\times N$ block-tridiagonal matrix to an $N/2\times N/2$
block-tridiagonal matrix by eliminating odd block rows as in CR,
PCR eliminates both the odd and the even block rows, resulting in not one, but
two $N/2\times N/2$ block-tridiagonal matrices. In $\log_2 N$ steps, PCR reduces
an $N\times N$ block-tridiagonal matrix down to $N$ matrices of a
single block each. 
In essence, PCR computes the Cholesky factors of $N$ different permutations of
the original block-tridiagonal matrix (every block row of the original matrix
appears last in exactly one of these permutations), whereas CR computes just a
single Cholesky factor. The PCR algorithm is listed in \Cref{alg:pcr} of \Cref{app:cr}.
Eliminating both the odd and even block rows at each level obviously increases the total number of required operations
compared to CR,
but PCR maintains full vector lane utilization throughout, and completely avoids
the back substitution step needed to recover the solution in standard CR.
Because of the higher throughput of matrix--matrix operations compared to
matrix--vector operations, this is often a worthwhile trade-off, especially
if the matrix dimensions are small or if multiple systems need to be solved using the same factorization (cf. \Cref{sec:fact-upd}).

A second solution strategy is to use an iterative method such as
PCG to solve the small block-tridiagonal
system at the output of the last batched level. This works well because one of
the properties of CR is that the off-diagonal blocks decay faster than the
diagonal ones with increasing levels of recursion \cite{bini_decay_2017,heller_aspects_1976}.
In combination with a suitable preconditioner such as the recently proposed symmetric stair preconditioner
\cite{bu_symmetric_2024}, we can expect PCG to converge in a few iterations.

To end up with a block-tridiagonal system of the same size as the vector length
(which is the desired format for both PCR and PCG),
we permute the final rows of the Schur complement matrix: the permutation
for full CR of a $32\times 32$ matrix is shown in \Cref{fig:cyclic-reduction-permutation-labels-32-vl1} (cf. \Cref{fig:cyclic-reduction-permutation}),
and a variant that uses CR only down to a reduced $4\times4$ block-tridiagonal matrix is shown in \Cref{fig:cyclic-reduction-permutation-labels-32-vl4}.
When using the vectorized approach from \Cref{fig:cyclic-reduction-graph-16} for the initial levels of CR,
this reduced $4\times 4$ matrix will be stored across the four vector lanes of a single batch.
This means that we can readily apply PCR to this batch, and
in the case of PCG, we will make use of this fact to compute a preconditioner for this matrix
using batched operations.

Consider the scenario for $v=4$ where the original matrix has been
cyclically reduced to a $4\times 4$ block-tridiagonal matrix
\begin{equation} \label{eq:cr-pcg}
        \mathscr{M}_{(4)} \defeq \scalebox{0.87}{$\begin{pmatrix}
            \Mdiag_0 & \tp \Ksub_0 \\
            \Ksub_0 & \Mdiag_1 & \tp \Ksub_1 \\
            & \Ksub_1 & \Mdiag_2 & \tp \Ksub_2 \\
            && \Ksub_2 & \Mdiag_3 \\
        \end{pmatrix}$} = \mathscr{L}_{(4)} \tp{\mathscr{L}_{(4)}} + \scalebox{0.87}{$\begin{pmatrix}
            0 & \tp \Ksub_0 \\
            \Ksub_0 & 0 & \tp \Ksub_1 \\
            & \Ksub_1 & 0 & \tp \Ksub_2 \\
            && \Ksub_2 & 0 \\
        \end{pmatrix}$} \defeq \mathscr{L}_{(4)} \tp{\mathscr{L}_{(4)}} + \mathscr{K}_{(4)}, \\
\end{equation}
where $\mathscr{L}_{(4)} \defeq \blkdiag\big( \operatorname{chol}(\Mdiag_0),\, \operatorname{chol}(\Mdiag_1),\, \operatorname{chol}(\Mdiag_2),\, \operatorname{chol}(\Mdiag_3) \big)$. \\
The inverse of the symmetric stair preconditioner is then given by the following expression: \cite[Eq.~38]{bu_symmetric_2024}
\begin{align}
    \setlength{\arraycolsep}{2pt}
    \inv\Phi_{(4)} &= \scalebox{0.87}{$\begin{pmatrix}
        \inv \Mdiag_0 & -\inv \Mdiag_0 \tp \Ksub_0 \inv \Mdiag_1 \\
        -\inv \Mdiag_1 \Ksub_0 \inv \Mdiag_0 & \inv \Mdiag_1 & -\inv \Mdiag_1 \tp \Ksub_1 \inv \Mdiag_2 \\
        & -\inv \Mdiag_2 \Ksub_1 \inv \Mdiag_1 & \inv \Mdiag_2 & -\inv \Mdiag_2 \tp \Ksub_2 \inv \Mdiag_3 \\
        && -\inv \Mdiag_3 \Ksub_2 \inv \Mdiag_2 & \inv \Mdiag_3 \\
    \end{pmatrix}$} \\
    &= \invtp {\mathscr{L}_{(4)}} \inv {\mathscr{L}_{(4)}} \left( \I - \mathscr{K}_{(4)} \invtp {\mathscr{L}_{(4)}} \inv {\mathscr{L}_{(4)}} \right).
\end{align}
The application of $\mathscr{M}_{(4)}$ to a vector requires
multiplication by $\mathscr{K}_{(4)}$ (as defined in \eqref{eq:cr-pcg}) and batched triangular matrix--vector multiplication by $\mathscr{L}_{(4)}$,
whereas application of the preconditioner $\inv\Phi_{(4)}$ requires
multiplication by $\mathscr{K}_{(4)}$ and batched triangular matrix--vector forward and
back substitutions. Thanks to the block-diagonal structure, factorization,
multiplication and substitution of $\mathscr{L}_{(4)}$ can be fully vectorized using the
same batched linear algebra routines as above. Multiplication by $\mathscr{K}_{(4)}$ involves
some lane-wise rotations, but can also be vectorized effectively.

An alternative choice is the block-Jacobi preconditioner $\inv{\tilde \Phi}_{(4)} = \invtp {\mathscr{L}}_{(4)} \inv {\mathscr{L}}_{(4)}$.
It only takes into account the diagonal blocks of $\mathscr{M}_{(4)}$, and is therefore cheaper to evaluate.
It is particularly effective if the entries of $\mathscr{K}_{(4)}$ are sufficiently small (when $\mathscr{M}_{(4)}$ is highly diagonally dominant).

Note that it is possible to combine CR, PCR and PCG in a hybrid solver, switching between algorithms at different levels in the recursion.
For the sake of simplicity, we restrict ourselves to CR+PCR and CR+PCG with
switchover at the level where the number of remaining blocks equals the vector length.

\subsubsection{Constraints on the number of processors}

In contrast to \Cref{sec:cyqlone-padding}, where we did not impose any constraints on the number of processors used,
the vectorized implementation of CR described in this section does require that the number of processors $p$ is a power of two.

\section{Benchmark results \done} \label{sec:results}
To evaluate the performance of the optimized implementation of \cyqlone{} and \cyqpalm{} from
the accompanying open-source software library \cite{kul-optec_cyqlone_2025}, we apply \cyqpalm{} to the commonly used linear mass--spring benchmark described in \cite{wang_fast_2008,wang_performance_2009,domahidi_efficient_2012},
and compare the solver run times to the state-of-the-art HPIPM solver for optimal control problems \cite{frison_hpipm_2020}.

\subsection{Problem description} \label{sec:mass-springs-problem}
The state vector $x \in \R^{2M}$ consists of the positions and velocities of
$M$ masses (with mass $m$) connected by springs (with spring constant $k$).
The leftmost and rightmost masses are connected to fixed walls (one at the origin, and one at distance $w$).
Following \cite{wang_fast_2008}, actuators are attached between pairs of masses. For $M > 6$, we simply repeat this configuration for $n_u = M/2$ actuators (as described by the matrix $W$ below).
The displacements of the masses from their equilibrium positions are bounded in magnitude by 4\,m, and the controls cannot exceed 0.5\,N.
\vspace{-0.5em}
\begin{equation}
    \begin{aligned}
        A_c &= \scalebox{0.9}{$\begin{pmatrix}
            0 & \I_M \\
            V & -\mu \I_M
        \end{pmatrix}$}, \quad \mathrlap{B_c = \scalebox{0.9}{$\begin{pmatrix}
            \Ozero_{M\times M/2} \\
            \operatorname{blkdiag}(W,\, \dotsc,\, W)
        \end{pmatrix}$}, \quad b_c = \frac1m \scalebox{0.9}{$\begin{pmatrix}
            0 \\ \vdots \\ kw
        \end{pmatrix}$}} \\
        V &= \frac1m \scalebox{0.8}{$\begin{pmatrix}
            -2 k & \phantom{-2} k \\
            \phantom{-2} k & -2 k & \phantom{-2} k \\
            & \phantom{-2} k & -2 k & \phantom{-2} k \\
            & & \phantom{-2} \ddots & \phantom{-2} \ddots & \phantom{-2} \ddots \\
            & & & \phantom{-2} k & -2 k
        \end{pmatrix}$}, & \qquad W &= \frac1m \scalebox{0.8}{$\begin{pmatrix}
            \phantom-1 \\
            -1 \\
            & \phantom-1 \\
            & & \phantom-1 \\
            & -1 \\
            & & -1
        \end{pmatrix}$} \\
        C_c &= \scalebox{0.9}{$\begin{pmatrix}
            \I_{M} & \Ozero_{M\times M} \\
            \Ozero_{n_u\times M} & \Ozero_{n_u\times M} \\
        \end{pmatrix}$}, & D_c &= \scalebox{0.9}{$\begin{pmatrix}
            \Ozero_{M\times n_u} \\
            \I_{n_u}
        \end{pmatrix}.$}
    \end{aligned}
\end{equation}

We discretize the continuous-time model using zero-order hold (ZOH) to obtain the discrete-time matrices $A,B,b,C,D$.
The sampling interval is set to $T_s = 15\,\mathrm{s}/N$ (so $T_s=0.5\,\mathrm{s}$ for the reference horizon $N=30$ from \cite{wang_fast_2008}).
For the experiments below, $T_s$ is reduced proportionally as $N$ increases to preserve the original time scale.
The value of $w$ in \cite{wang_fast_2008} and \cite{domahidi_efficient_2012} is
set to zero, so all masses are pulled towards the origin.
The constants $k$ and $m$ are set to one, with no friction ($\mu=0$).
The quadratic MPC cost drives the system to steady state, with $Q=Q_N=\I$ and $R=\I$, that is $\ell(x, u) = \tfrac12 \tp{(x - x_\mathrm{\infty})} Q (x - x_\mathrm{\infty}) + \tfrac12 \tp u R u$.
The purpose of this problem is to gauge the performance of the QP solver and the underlying linear algebra; we do not consider recursive feasibility.

To generate a set of benchmark problems, we vary the number of masses $M\in\{6,12,18,24,30\}$ and the horizon length $N\in\{32, 64, 96, 128, 192, 256\}$. For each combination of $M$ and $N$,
we randomly generate 150 initial states with velocity equal to zero and positions whose deviations from the steady state are uniformly distributed in $[-3, 3]$.
Source code for the benchmarks is available in the \cyqlone{} GitHub repository \cite{kul-optec_cyqlone_2025}.

\subsection{Benchmark results -- {\normalfont\cyqpalm{}} \done} \label{sec:benchmark-results-cyqpalm}

\begin{figure}[htbp]
    \centering
    \includegraphics[scale=0.8]{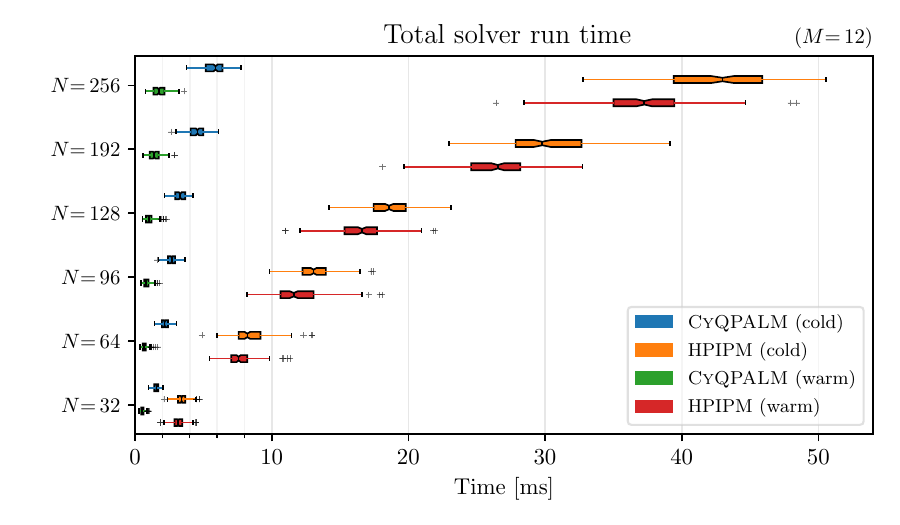}
    \caption{Box plots of the total solver run times for $6\times 150$ instances of the mass--springs benchmark described in
    \Cref{sec:mass-springs-problem}, with $M=12$ masses and varying horizon lengths $N$.
    Timings for \cyqpalm{}
    (with a vector length of four across eight cores) and for HPIPM (using the \textit{speed} profile).
    We compare the solvers without warm starting (cold) and using the solution of the previous MPC iteration shifted by one time step as initialization (warm).}
    \label{fig:box-time-spring-mass-wang-boyd-2008}
\end{figure}

\begin{figure}[htbp]
    \centering
    \includegraphics[scale=0.8]{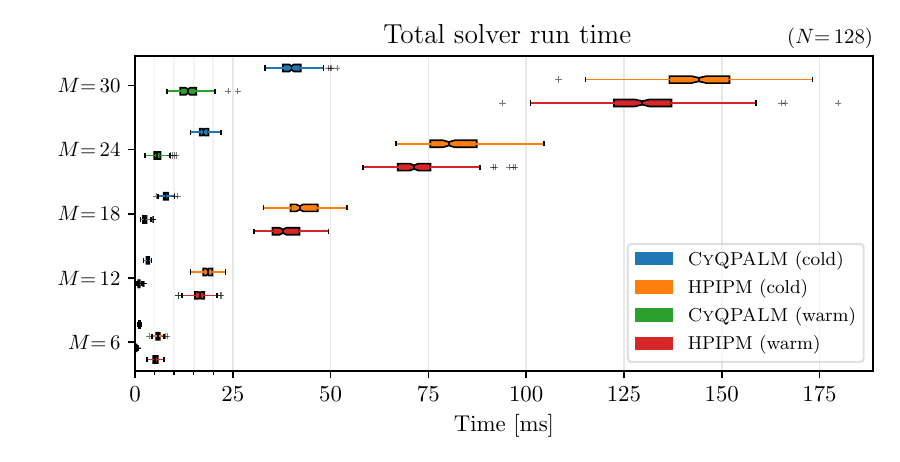}
    \caption{Box plots of the total solver run times for $5\times 150$ instances of the mass--springs benchmark described in
    \Cref{sec:mass-springs-problem}, horizon length $N=128$ and varying numbers of masses $M$.
    Timings for \cyqpalm{}
    (with a vector length of four across eight cores) and for HPIPM (using the \textit{speed} profile).
    We compare the solvers without warm starting (cold) and using the solution of the previous MPC iteration shifted by one time step as initialization (warm).}
    \label{fig:box-time-spring-mass-wang-boyd-2008-M}
\end{figure}

The solver run times for \cyqpalm{} (\qpalm{} using \cyqlone{} as an inner solver) are compared to those of HPIPM with the \textit{speed} profile in \Cref{fig:box-time-spring-mass-wang-boyd-2008,fig:box-time-spring-mass-wang-boyd-2008-M}.
One set of experiments measures the run times when the solvers are initialized using all zeros,
and another measures the run times when the solvers are initialized with the shifted solution of a previous problem.
Specifically, we solve a single optimal control problem (excluded from the measurements) and
apply the first controls from the solution to the system, as in MPC.
The resulting state is used to construct a second OCP, and we measure the solver run time for this second problem,
using the solution to the first problem as an initial guess, after shifting it over by one time step (repeating the variables for the last stage).

Experiments are carried out on an Intel Core i7-11700 at 2.5\,GHz (without frequency scaling to reduce noise in the timing measurements),
using all eight processors and vector length four.\,\footnote{Even though the processor supports a maximum vector length of eight, it lacks a dedicated functional unit for 512-bit floating-point operations.
Instead, two 256-bit units are combined to carry out 512-bit operations \cite[\S 12.9]{fog_microarchitecture_2024}, resulting in a similar throughput for vector lengths four and eight.
Since the working set for vector length eight is larger, using a vector length of four is preferable.}

\Cref{fig:box-time-spring-mass-wang-boyd-2008} shows how the solver run times scale for various horizon lengths $N$,
and \Cref{fig:box-time-spring-mass-wang-boyd-2008-M} shows how the solver run times scale for various numbers of masses $M$.
For the longest horizon length $N=256$ and with $M=12$ masses, the worst-case run time of \cyqpalm{} is over six times faster
than the worst-case run time of HPIPM. This allows the sampling time to be reduced considerably
in real-time applications such as MPC. When warm-starting the solvers, \cyqpalm{} is over thirteen times faster in terms of worst-case run time,
and over twenty times faster on average.
The speedup factors for other combinations of $M$ and $N$ are listed in \Cref{tab:speedup-time-cold,tab:speedup-time-warm}
(without and with warm starting, respectively).
The first value in each cell is the average of the quotients of the
total solver run times for all problem instances,
and the second value is the quotient of the maximum (worst-case) run time of all problem instances for the two solvers.\,\footnote{That is, the quotient of the maxima,
in contrast to the maximum of the quotients. We are interested in the former because it determines the minimum sample time that can be used in real-time control applications.}

\begin{table}[htbp]
    \centering
    \caption{Average speedup of the run time and (speedup of the worst-case run time) for \cyqpalm{} (cold start) compared to \hpipm{} (cold start).}
    \label{tab:speedup-time-cold}
    \fontsize{7}{9}\selectfont
    \begin{tabular}{r|w{c}{4.08em}w{c}{4.08em}w{c}{4.08em}w{c}{4.08em}w{c}{4.08em}w{c}{4.08em}w{c}{4.08em}w{c}{4.08em}}
    \toprule
    $M$ & $N$=32 & $N$=64 & $N$=96 & $N$=128 & $N$=160 & $N$=192 & $N$=224 & $N$=256 \\
    \midrule
    6 & \cellcolor[rgb]{0.913,0.967,0.896}\textcolor{black}{1.7\,(\kern-0.75pt1.8\kern-0.75pt)} & \cellcolor[rgb]{0.709,0.884,0.684}\textcolor{black}{3.2\,(\kern-0.75pt3.2\kern-0.75pt)} & \cellcolor[rgb]{0.507,0.793,0.506}\textcolor{black}{4.3\,(\kern-0.75pt3.9\kern-0.75pt)} & \cellcolor[rgb]{0.270,0.678,0.372}\textcolor{black}{5.4\,(\kern-0.75pt5.1\kern-0.75pt)} & \cellcolor[rgb]{0.142,0.550,0.274}\textcolor{white}{6.3\,(\kern-0.75pt5.5\kern-0.75pt)} & \cellcolor[rgb]{0.048,0.469,0.207}\textcolor{white}{7.0\,(\kern-0.75pt6.1\kern-0.75pt)} & \cellcolor[rgb]{0.000,0.332,0.133}\textcolor{white}{7.8\,(\kern-0.75pt6.8\kern-0.75pt)} & \cellcolor[rgb]{0.000,0.267,0.106}\textcolor{white}{8.2\,(\kern-0.75pt6.8\kern-0.75pt)} \\
    12 & \cellcolor[rgb]{0.857,0.944,0.835}\textcolor{black}{2.2\,(\kern-0.75pt2.3\kern-0.75pt)} & \cellcolor[rgb]{0.591,0.832,0.574}\textcolor{black}{3.9\,(\kern-0.75pt4.3\kern-0.75pt)} & \cellcolor[rgb]{0.351,0.718,0.412}\textcolor{black}{5.1\,(\kern-0.75pt4.8\kern-0.75pt)} & \cellcolor[rgb]{0.216,0.629,0.333}\textcolor{black}{5.8\,(\kern-0.75pt5.5\kern-0.75pt)} & \cellcolor[rgb]{0.134,0.542,0.268}\textcolor{white}{6.4\,(\kern-0.75pt6.2\kern-0.75pt)} & \cellcolor[rgb]{0.074,0.491,0.225}\textcolor{white}{6.8\,(\kern-0.75pt6.4\kern-0.75pt)} & \cellcolor[rgb]{0.005,0.432,0.176}\textcolor{white}{7.2\,(\kern-0.75pt6.9\kern-0.75pt)} & \cellcolor[rgb]{0.000,0.413,0.167}\textcolor{white}{7.4\,(\kern-0.75pt6.5\kern-0.75pt)} \\
    18 & \cellcolor[rgb]{0.824,0.931,0.799}\textcolor{black}{2.5\,(\kern-0.75pt2.7\kern-0.75pt)} & \cellcolor[rgb]{0.585,0.829,0.570}\textcolor{black}{3.9\,(\kern-0.75pt4.0\kern-0.75pt)} & \cellcolor[rgb]{0.389,0.736,0.430}\textcolor{black}{4.9\,(\kern-0.75pt4.6\kern-0.75pt)} & \cellcolor[rgb]{0.257,0.672,0.366}\textcolor{black}{5.5\,(\kern-0.75pt5.0\kern-0.75pt)} & \cellcolor[rgb]{0.194,0.605,0.316}\textcolor{black}{6.0\,(\kern-0.75pt5.5\kern-0.75pt)} & \cellcolor[rgb]{0.168,0.578,0.295}\textcolor{white}{6.2\,(\kern-0.75pt6.3\kern-0.75pt)} & \cellcolor[rgb]{0.164,0.574,0.292}\textcolor{white}{6.2\,(\kern-0.75pt5.7\kern-0.75pt)} & \cellcolor[rgb]{0.175,0.585,0.301}\textcolor{white}{6.1\,(\kern-0.75pt5.8\kern-0.75pt)} \\
    24 & \cellcolor[rgb]{0.835,0.936,0.811}\textcolor{black}{2.4\,(\kern-0.75pt2.5\kern-0.75pt)} & \cellcolor[rgb]{0.618,0.845,0.597}\textcolor{black}{3.8\,(\kern-0.75pt3.9\kern-0.75pt)} & \cellcolor[rgb]{0.496,0.788,0.497}\textcolor{black}{4.4\,(\kern-0.75pt4.1\kern-0.75pt)} & \cellcolor[rgb]{0.445,0.764,0.458}\textcolor{black}{4.6\,(\kern-0.75pt4.8\kern-0.75pt)} & \cellcolor[rgb]{0.480,0.780,0.483}\textcolor{black}{4.5\,(\kern-0.75pt4.1\kern-0.75pt)} & \cellcolor[rgb]{0.507,0.793,0.506}\textcolor{black}{4.3\,(\kern-0.75pt4.1\kern-0.75pt)} & \cellcolor[rgb]{0.524,0.801,0.520}\textcolor{black}{4.2\,(\kern-0.75pt4.2\kern-0.75pt)} & \cellcolor[rgb]{0.507,0.793,0.506}\textcolor{black}{4.3\,(\kern-0.75pt4.8\kern-0.75pt)} \\
    30 & \cellcolor[rgb]{0.842,0.938,0.819}\textcolor{black}{2.3\,(\kern-0.75pt2.5\kern-0.75pt)} & \cellcolor[rgb]{0.681,0.872,0.656}\textcolor{black}{3.4\,(\kern-0.75pt3.1\kern-0.75pt)} & \cellcolor[rgb]{0.648,0.858,0.624}\textcolor{black}{3.6\,(\kern-0.75pt3.4\kern-0.75pt)} & \cellcolor[rgb]{0.653,0.860,0.629}\textcolor{black}{3.6\,(\kern-0.75pt3.4\kern-0.75pt)} & \cellcolor[rgb]{0.624,0.847,0.602}\textcolor{black}{3.8\,(\kern-0.75pt2.6\kern-0.75pt)} & \cellcolor[rgb]{0.591,0.832,0.574}\textcolor{black}{3.9\,(\kern-0.75pt2.9\kern-0.75pt)} & \cellcolor[rgb]{0.579,0.827,0.565}\textcolor{black}{4.0\,(\kern-0.75pt2.8\kern-0.75pt)} & \cellcolor[rgb]{0.552,0.814,0.542}\textcolor{black}{4.1\,(\kern-0.75pt2.7\kern-0.75pt)} \\
    \bottomrule
\end{tabular}

\end{table}
\begin{table}[htbp]
    \centering
    \caption{Average speedup of the run time and (speedup of the worst-case run time) for \cyqpalm{} (warm start) compared to \hpipm{} (warm start).}
    \fontsize{7}{9}\selectfont
    \label{tab:speedup-time-warm}
    \fontsize{7}{9}\selectfont
    \begin{tabular}{r|w{c}{4.08em}w{c}{4.08em}w{c}{4.08em}w{c}{4.08em}w{c}{4.08em}w{c}{4.08em}w{c}{4.08em}w{c}{4.08em}}
    \toprule
    $M$ & $N$=32 & $N$=64 & $N$=96 & $N$=128 & $N$=160 & $N$=192 & $N$=224 & $N$=256 \\
    \midrule
    6 & \cellcolor[rgb]{0.875,0.952,0.854}\textcolor{black}{4.7\,(\kern-0.75pt2.9\kern-0.75pt)} & \cellcolor[rgb]{0.672,0.868,0.647}\textcolor{black}{\phantom{0}9.6\,(\kern-0.75pt6.2\kern-0.75pt)} & \cellcolor[rgb]{0.463,0.773,0.470}\textcolor{black}{13.4\,\phantom{0}(\kern-0.75pt7.8\kern-0.75pt)} & \cellcolor[rgb]{0.227,0.641,0.342}\textcolor{black}{17.4\,\phantom{0}(\kern-0.75pt9.0\kern-0.75pt)} & \cellcolor[rgb]{0.130,0.539,0.265}\textcolor{white}{19.9\,(\kern-0.75pt10.2\kern-0.75pt)} & \cellcolor[rgb]{0.044,0.465,0.204}\textcolor{white}{21.9\,(\kern-0.75pt11.0\kern-0.75pt)} & \cellcolor[rgb]{0.000,0.327,0.131}\textcolor{white}{24.8\,(\kern-0.75pt14.1\kern-0.75pt)} & \cellcolor[rgb]{0.000,0.267,0.106}\textcolor{white}{26.1\,(\kern-0.75pt13.4\kern-0.75pt)} \\
    12 & \cellcolor[rgb]{0.824,0.931,0.799}\textcolor{black}{6.1\,(\kern-0.75pt4.3\kern-0.75pt)} & \cellcolor[rgb]{0.574,0.824,0.561}\textcolor{black}{11.5\,(\kern-0.75pt6.9\kern-0.75pt)} & \cellcolor[rgb]{0.345,0.715,0.409}\textcolor{black}{15.2\,(\kern-0.75pt10.1\kern-0.75pt)} & \cellcolor[rgb]{0.245,0.660,0.357}\textcolor{black}{16.9\,\phantom{0}(\kern-0.75pt9.5\kern-0.75pt)} & \cellcolor[rgb]{0.171,0.582,0.298}\textcolor{white}{18.9\,(\kern-0.75pt11.9\kern-0.75pt)} & \cellcolor[rgb]{0.121,0.531,0.259}\textcolor{white}{20.2\,(\kern-0.75pt11.3\kern-0.75pt)} & \cellcolor[rgb]{0.074,0.491,0.225}\textcolor{white}{21.2\,(\kern-0.75pt12.0\kern-0.75pt)} & \cellcolor[rgb]{0.018,0.443,0.185}\textcolor{white}{22.5\,(\kern-0.75pt13.5\kern-0.75pt)} \\
    18 & \cellcolor[rgb]{0.805,0.924,0.780}\textcolor{black}{6.6\,(\kern-0.75pt4.5\kern-0.75pt)} & \cellcolor[rgb]{0.541,0.809,0.533}\textcolor{black}{12.0\,(\kern-0.75pt7.5\kern-0.75pt)} & \cellcolor[rgb]{0.427,0.755,0.449}\textcolor{black}{13.9\,(\kern-0.75pt10.6\kern-0.75pt)} & \cellcolor[rgb]{0.314,0.699,0.394}\textcolor{black}{15.7\,(\kern-0.75pt10.7\kern-0.75pt)} & \cellcolor[rgb]{0.264,0.675,0.369}\textcolor{black}{16.5\,(\kern-0.75pt12.7\kern-0.75pt)} & \cellcolor[rgb]{0.234,0.648,0.348}\textcolor{black}{17.2\,(\kern-0.75pt11.8\kern-0.75pt)} & \cellcolor[rgb]{0.234,0.648,0.348}\textcolor{black}{17.2\,(\kern-0.75pt12.6\kern-0.75pt)} & \cellcolor[rgb]{0.249,0.664,0.360}\textcolor{black}{16.8\,(\kern-0.75pt11.3\kern-0.75pt)} \\
    24 & \cellcolor[rgb]{0.813,0.927,0.787}\textcolor{black}{6.4\,(\kern-0.75pt4.7\kern-0.75pt)} & \cellcolor[rgb]{0.585,0.829,0.570}\textcolor{black}{11.3\,(\kern-0.75pt7.7\kern-0.75pt)} & \cellcolor[rgb]{0.502,0.791,0.501}\textcolor{black}{12.7\,\phantom{0}(\kern-0.75pt9.4\kern-0.75pt)} & \cellcolor[rgb]{0.474,0.778,0.479}\textcolor{black}{13.2\,\phantom{0}(\kern-0.75pt9.3\kern-0.75pt)} & \cellcolor[rgb]{0.535,0.806,0.529}\textcolor{black}{12.1\,\phantom{0}(\kern-0.75pt9.0\kern-0.75pt)} & \cellcolor[rgb]{0.552,0.814,0.542}\textcolor{black}{11.8\,\phantom{0}(\kern-0.75pt8.9\kern-0.75pt)} & \cellcolor[rgb]{0.563,0.819,0.552}\textcolor{black}{11.6\,\phantom{0}(\kern-0.75pt7.8\kern-0.75pt)} & \cellcolor[rgb]{0.585,0.829,0.570}\textcolor{black}{11.2\,\phantom{0}(\kern-0.75pt7.9\kern-0.75pt)} \\
    30 & \cellcolor[rgb]{0.813,0.927,0.787}\textcolor{black}{6.4\,(\kern-0.75pt3.9\kern-0.75pt)} & \cellcolor[rgb]{0.672,0.868,0.647}\textcolor{black}{\phantom{0}9.6\,(\kern-0.75pt7.7\kern-0.75pt)} & \cellcolor[rgb]{0.634,0.852,0.611}\textcolor{black}{10.3\,\phantom{0}(\kern-0.75pt7.7\kern-0.75pt)} & \cellcolor[rgb]{0.658,0.862,0.633}\textcolor{black}{\phantom{0}9.9\,\phantom{0}(\kern-0.75pt6.8\kern-0.75pt)} & \cellcolor[rgb]{0.653,0.860,0.629}\textcolor{black}{\phantom{0}9.9\,\phantom{0}(\kern-0.75pt7.6\kern-0.75pt)} & \cellcolor[rgb]{0.634,0.852,0.611}\textcolor{black}{10.4\,\phantom{0}(\kern-0.75pt7.5\kern-0.75pt)} & \cellcolor[rgb]{0.644,0.856,0.620}\textcolor{black}{10.2\,\phantom{0}(\kern-0.75pt7.4\kern-0.75pt)} & \cellcolor[rgb]{0.634,0.852,0.611}\textcolor{black}{10.3\,\phantom{0}(\kern-0.75pt7.5\kern-0.75pt)} \\
    \bottomrule
\end{tabular}

\end{table}

When the number of masses is small (e.g. $M=6$ with $n_x=12$, $n_u=3$), the performance
of BLASFEO is relatively low, while \cyqpalm{}'s vectorization is highly
performant, resulting in an significant speedup.
For $M\le 18$, performance of \cyqlone{} generally increases with the horizon length,
as predicted by the theoretical operation counts (cf. \Cref{fig:theoretical-speedup}).
The effect is magnified by the relatively higher synchronization overhead for short horizons:
for small $N$, the embarrassingly parallel modified Riccati recursion accounts for
a smaller fraction of the total run time compared to CR, the latter of which
requires more synchronization between threads.
Such synchronization steps introduce some latency,
and the communication of matrices between threads results in L1 and L2 cache misses,
limiting the performance of subsequent operations on these matrices.
This effect is exacerbated for the smallest problems (e.g. $M=6$ and $N=32$),
where the run time of individual linear algebra operations is
very short (thanks to \cyqlone{}'s batched routines),
so the synchronization overhead accounts for a larger fraction of the total run time.

For larger numbers of masses, as in \Cref{fig:box-time-spring-mass-wang-boyd-2008-M}, \cyqpalm{} also significantly outperforms HPIPM,
although the speedup is lower than for the cases $M=6$ or $M=12$ discussed above.
There are multiple contributing factors:
\begin{enumerate}
    \setlength\itemsep{0em}
    \item The cost of the \cyqlone{} factorization scales
    cubically with increasing numbers of states (cf. \eqref{eq:flop}),
    so the factorization step dominates the overall run time for large $M$.
    In contrast, if $M$ is small,
    quadratic operations such as matrix--vector products also make up a
    considerable fraction of the run time. While matrix--vector products can be
    fully parallelized (with a speedup of close to $8\times$ for $P=8$ processors),
    the speedup for the \cyqlone{} factorization is less than half of the number of
    processors (cf. \Cref{fig:theoretical-speedup}).
    \item Performance of the BLASFEO routines used by HPIPM increases for increasing matrix sizes, whereas
    the batched linear algebra routines used in \cyqlone{} already reach close to
    peak performance for much smaller matrices (cf. \Cref{fig:syrk-potrf-perf}). BLASFEO therefore
    starts to catch up for larger $M$.
    \item For $M\ge20$, batches of matrices of size $n_x$$\times$$n_x$ no longer fit the L1 cache.
    For large problems (e.g. $M=24$ and $N\ge160$), the speedup no longer
    increases with horizon length because the working set exceeds the L2 cache.
    Results for processors with larger caches are included in \Cref{app:results}.
    \item \qpalm{} appears to require slightly more iterations than HPIPM for larger problems (although some iterations are cheaper, thanks to the factorization updates, cf. \Cref{fig:scaling-abs-time_per_iter-spring-mass-wang-boyd-2008-scaling}).
\end{enumerate}

In conclusion, \cyqpalm{} is significantly faster than \hpipm{} across the board,
achieving the largest speedups for a small number of masses
(where the performance of \cyqlone{}'s batched linear algebra is much higher
than that of BLASFEO), and for long horizons
(where the cost and the overhead in the CR phase of
\cyqlone{} are small compared to the duration of the Riccati phase).

\subsection{Benchmark results -- {\normalfont\cyqlone{}}  \done}

\begin{figure}[t]
    \centering
    \includegraphics[scale=0.8]{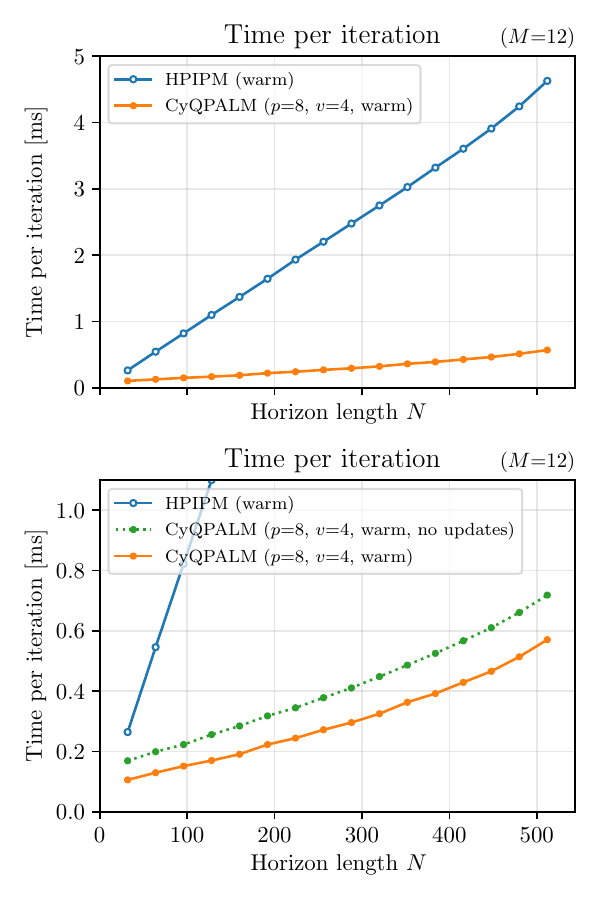}
    \caption{Top: Average run times per iteration of the solvers, for the mass--springs benchmark described in
    \Cref{sec:mass-springs-problem}, with $M=12$ masses and varying horizon lengths $N$.
    Timings for \cyqpalm{}
    (with a vector length of four across eight cores) and for HPIPM (using the \textit{speed} profile), both with warm starting.
    Bottom: The same graph, zoomed in and including \cyqpalm{} without factorization updates.}
    \label{fig:scaling-abs-time_per_iter-spring-mass-wang-boyd-2008-scaling}
\end{figure}

\begin{table}[htbp]
    \caption{Average speedup of the run time per iteration for \cyqpalm{} (cold start) compared to \hpipm{} (cold start).}
    \label{tab:speedup-time-per-iter-cold}
    \centering
    \fontsize{7}{9}\selectfont
    \begin{tabular}{r|w{c}{3.08em}w{c}{3.08em}w{c}{3.08em}w{c}{3.08em}w{c}{3.08em}w{c}{3.08em}w{c}{3.08em}w{c}{3.08em}}
    \toprule
    $M$ & $N$=32 & $N$=64 & $N$=96 & $N$=128 & $N$=160 & $N$=192 & $N$=224 & $N$=256 \\
    \midrule
    6 & \cellcolor[rgb]{0.872,0.950,0.850}\textcolor{black}{2.1} & \cellcolor[rgb]{0.639,0.854,0.615}\textcolor{black}{3.7} & \cellcolor[rgb]{0.414,0.749,0.443}\textcolor{black}{4.8} & \cellcolor[rgb]{0.219,0.633,0.336}\textcolor{black}{5.8} & \cellcolor[rgb]{0.117,0.528,0.256}\textcolor{white}{6.6} & \cellcolor[rgb]{0.013,0.439,0.182}\textcolor{white}{7.3} & \cellcolor[rgb]{0.000,0.342,0.137}\textcolor{white}{7.8} & \cellcolor[rgb]{0.000,0.267,0.106}\textcolor{white}{8.3} \\
    12 & \cellcolor[rgb]{0.802,0.922,0.776}\textcolor{black}{2.7} & \cellcolor[rgb]{0.496,0.788,0.497}\textcolor{black}{4.4} & \cellcolor[rgb]{0.257,0.672,0.366}\textcolor{black}{5.5} & \cellcolor[rgb]{0.153,0.562,0.283}\textcolor{white}{6.3} & \cellcolor[rgb]{0.078,0.494,0.228}\textcolor{white}{6.8} & \cellcolor[rgb]{0.018,0.443,0.185}\textcolor{white}{7.2} & \cellcolor[rgb]{0.000,0.403,0.162}\textcolor{white}{7.5} & \cellcolor[rgb]{0.000,0.363,0.146}\textcolor{white}{7.7} \\
    18 & \cellcolor[rgb]{0.756,0.903,0.729}\textcolor{black}{3.0} & \cellcolor[rgb]{0.474,0.778,0.479}\textcolor{black}{4.5} & \cellcolor[rgb]{0.276,0.681,0.375}\textcolor{black}{5.4} & \cellcolor[rgb]{0.190,0.601,0.313}\textcolor{black}{6.0} & \cellcolor[rgb]{0.146,0.554,0.277}\textcolor{white}{6.4} & \cellcolor[rgb]{0.121,0.531,0.259}\textcolor{white}{6.5} & \cellcolor[rgb]{0.125,0.535,0.262}\textcolor{white}{6.5} & \cellcolor[rgb]{0.146,0.554,0.277}\textcolor{white}{6.4} \\
    24 & \cellcolor[rgb]{0.770,0.909,0.743}\textcolor{black}{2.9} & \cellcolor[rgb]{0.524,0.801,0.520}\textcolor{black}{4.3} & \cellcolor[rgb]{0.395,0.739,0.434}\textcolor{black}{4.9} & \cellcolor[rgb]{0.364,0.724,0.418}\textcolor{black}{5.0} & \cellcolor[rgb]{0.402,0.742,0.437}\textcolor{black}{4.9} & \cellcolor[rgb]{0.463,0.773,0.470}\textcolor{black}{4.6} & \cellcolor[rgb]{0.474,0.778,0.479}\textcolor{black}{4.5} & \cellcolor[rgb]{0.474,0.778,0.479}\textcolor{black}{4.5} \\
    30 & \cellcolor[rgb]{0.775,0.911,0.747}\textcolor{black}{2.8} & \cellcolor[rgb]{0.585,0.829,0.570}\textcolor{black}{4.0} & \cellcolor[rgb]{0.574,0.824,0.561}\textcolor{black}{4.0} & \cellcolor[rgb]{0.579,0.827,0.565}\textcolor{black}{4.0} & \cellcolor[rgb]{0.552,0.814,0.542}\textcolor{black}{4.1} & \cellcolor[rgb]{0.541,0.809,0.533}\textcolor{black}{4.2} & \cellcolor[rgb]{0.530,0.804,0.524}\textcolor{black}{4.2} & \cellcolor[rgb]{0.513,0.796,0.511}\textcolor{black}{4.3} \\
    \bottomrule
\end{tabular}

\end{table}

We also compare the solver run time per iteration between \cyqpalm{} and HPIPM
to gauge the performance of the \cyqlone{} solver, other parallel linear algebra routines, and the parallel line search,
independently of the number of iterations of the optimization solvers.
The top graph of \Cref{fig:scaling-abs-time_per_iter-spring-mass-wang-boyd-2008-scaling} shows that
\cyqpalm{} is over seven times faster per iteration than HPIPM.
When factorization updates are disabled (as shown in the bottom graph of \Cref{fig:scaling-abs-time_per_iter-spring-mass-wang-boyd-2008-scaling}),
performance of \cyqpalm{} reduces by as much as 30\%.
The speedup factors of the solver run times per iteration are reported in \Cref{tab:speedup-time-per-iter-cold} (no warm starting).

\begin{figure}
    \centering
    \includegraphics[width=\textwidth,clip,trim=0.55cm 1.1cm 2.2cm 1.7cm]{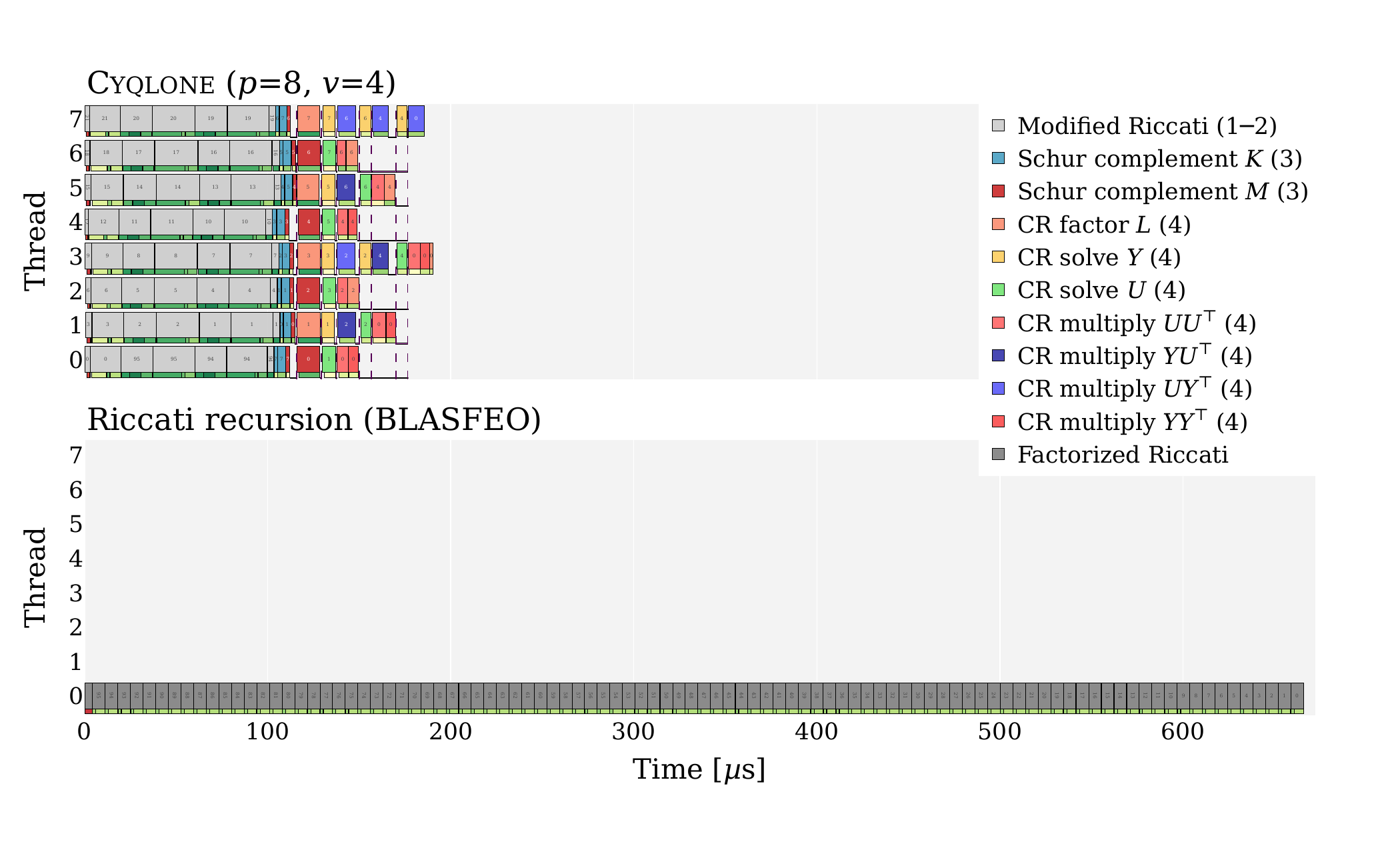}
    \caption{Thread-level execution traces comparing \cyqlone{} and the serial factorized Riccati recursion implemented using BLASFEO,
    applied to the KKT system \eqref{eq:opt-cond-ocp-eq} of an OCP with $n_x$=30, $n_u$=20 and $N$=96.
    Each large rectangle represents a specific step in the algorithm, with the wall time along the horizontal axis, and the thread on which the step is executed along the vertical axis.
    The index in the center of each rectangle refers to the OCP stage $j$ (the stage corresponding to the first vector lane of the batch in the case of \cyqlone{}) or the index $i$ in the Schur complement $\mathscr M$.
    Legend labels include the steps of the \cyqlone{} algorithm: (\ref{it:bchol-diag}--\ref{it:solve-subdiag})
    for the modified Riccati recursion in \Cref{alg:fact-mod-riccati},
    (\ref{it:compute-schur}) for the evaluation of the Schur complement, and (\ref{it:cr-schur}) for its cyclic reduction in \Cref{alg:fact-parallel-riccati-cr}.
    Colored bars at the bottom of each block indicate the performance of individual linear algebra operations:
    red for low performance ($\lt 2$ GFLOPS), yellow for medium performance ($\lt 10$ GFLOPS), and green for high performance (with the darkest green reaching up to 17.7 GFLOPS).
    For example, each gray block in the bottom graph corresponds to the application of the Riccati recursion \eqref{eq:riccati-recursion-P} to a single stage, which is done using two BLASFEO operations:
    \textsc{trmm} to evaluate $\ttp{(B_j\;\;A_j)} L_{j+1}$ ($\approx 13$ GFLOPS), and fused \textsc{syrk+potrf} to update and factorize the cost Hessian ($\approx 12$ GFLOPS).
    For \cyqlone{}, the colors in the CR phase of the algorithm (roughly) match the ones in \Cref{fig:cyclic-reduction-graph-16}, although some \textsc{syrk+potrf} operations are fused.
    Vertical dashed lines represent barriers for synchronization between threads, and horizontal purple lines indicate wait times for threads that already arrived at the barrier.}
    \label{fig:trace}
\end{figure}

Finally, we consider the factorization step in isolation: \Cref{fig:trace} shows detailed
execution traces for the \cyqlone{} factorization (\Cref{alg:fact-parallel-riccati-cr}),
and for the factorized Riccati recursion from \cite[Alg.\,3]{frison_efficient_2013}, implemented using BLASFEO.
Thanks to the parallelization and vectorization,
\cyqlone{} is around 3.5 times faster than BLASFEO for this benchmark.
Even though the matrices are small, the modified Riccati recursion phase
(steps \ref{it:bchol-diag} and \ref{it:solve-subdiag}) and the Schur complement
construction (step \ref{it:compute-schur})
achieve close to the system's peak performance because of the batch-wise vectorization
(the performance of the first stages to be processed is slightly lower because
the necessary data is not yet in the cache).
The linear algebra operations in the CR phase (step \ref{it:cr-schur})
achieve lower performance because of the cache misses caused by inter-processor
communication (which goes through the slower shared L3 cache or triggers a direct cache-to-cache transfer through the ring interconnect).

\section{Comparison of solvers for KKT systems with OCP structure \done}\label{sec:compare-ocp-lin-solver}
This section explores several possible solution methods for linear systems of the form \eqref{eq:opt-cond-ocp-eq},
using the matrix visualization introduced in \Cref{mat:pdp-cr-N12-P4} and the corresponding elimination tree as graphical tools
for comparing the computational cost (by counting fill-in)
and parallelizability (using the elimination tree) of the various methods.
\Cref{tab:compare-ocp-solvers} lists the methods discussed in the remainder of this section,
with some important properties that affect their performance.

\begin{table}
    \centering
    \caption{Comparison of linear solvers for KKT systems with optimal control structure.
    Optimal values for each column are shown in bold.}
    \label{tab:compare-ocp-solvers}
    \renewcommand{\arraystretch}{1}
    \small
    \begin{tabular}{rccccc}
        Method & Figure & Parallelizable & \shortstack{Structure\\exploitation} & \shortstack{Fill-in\\$\smash{(\lambda^j,\lambda^j)}$ blocks} & \shortstack{CR\\block size} \\
        \toprule
        \makecell[r]{Pas et al. \cite{pas_exploiting_2024}} & \ref{mat:schur-compl-1} & Partial & No & $N$ & -- \\
        \cmidrule(lr){1-6}
        \makecell[r]{Frison et al.\\\cite[Alg.\,1]{frison_efficient_2013}} & \ref{mat:riccati} & No & \makecell[c]{Asymmetric\\Riccati} & 0 & -- \\
        \cmidrule(lr){1-6}
        \makecell[r]{Frison et al.\\\cite[Alg.\,3]{frison_efficient_2013}} & \ref{mat:riccati} & No & \makecell[c]{\textbf{Factorized}\\\textbf{Riccati}} & 0 & -- \\
        \cmidrule(lr){1-6}
        \makecell[r]{Wright \cite{wright_partitioned_1991}} & \ref{mat:pdp-pri} & \textbf{Yes} & \makecell[c]{Asymmetric\\Riccati} & $\boldsymbol{P}$ & $2 n_x$ \\
        \cmidrule(lr){1-6}
        \makecell[r]{Nicholson et al.\\\cite[\S 4.2]{nicholson_parallel_2018}} & \ref{mat:direct-cr-1} & \textbf{Yes} & No & $N$ & $2n_x + n_u$ \\
        \cmidrule(lr){1-6}
        \makecell[r]{Jallet et al. \cite{jallet_parallel_2024-1}} & \ref{mat:parallel-xu-elim} &  Partial & \makecell[c]{Generalized\\Riccati} & $\boldsymbol{P}$ & -- \\
        \cmidrule(lr){1-6}
        \makecell[r]{Frasch et al.\\\cite[\S 6.1]{frasch_parallel_2015}} & \ref{mat:schur-cr} & \textbf{Yes} & No & $N$ & $\boldsymbol{n_x}$ \\
        \cmidrule(lr){1-6}
        \cyqlone{} & \ref{mat:cyqlone-N8} & \textbf{Yes} & \makecell[c]{\textbf{Factorized}\\\textbf{Riccati}} & $\boldsymbol{P}$ & $\boldsymbol{n_x}$ \\
        \bottomrule
    \end{tabular}
\end{table}

\begin{figure}
\begin{minipage}{0.725\textwidth}
    \includegraphics[width=1\textwidth]{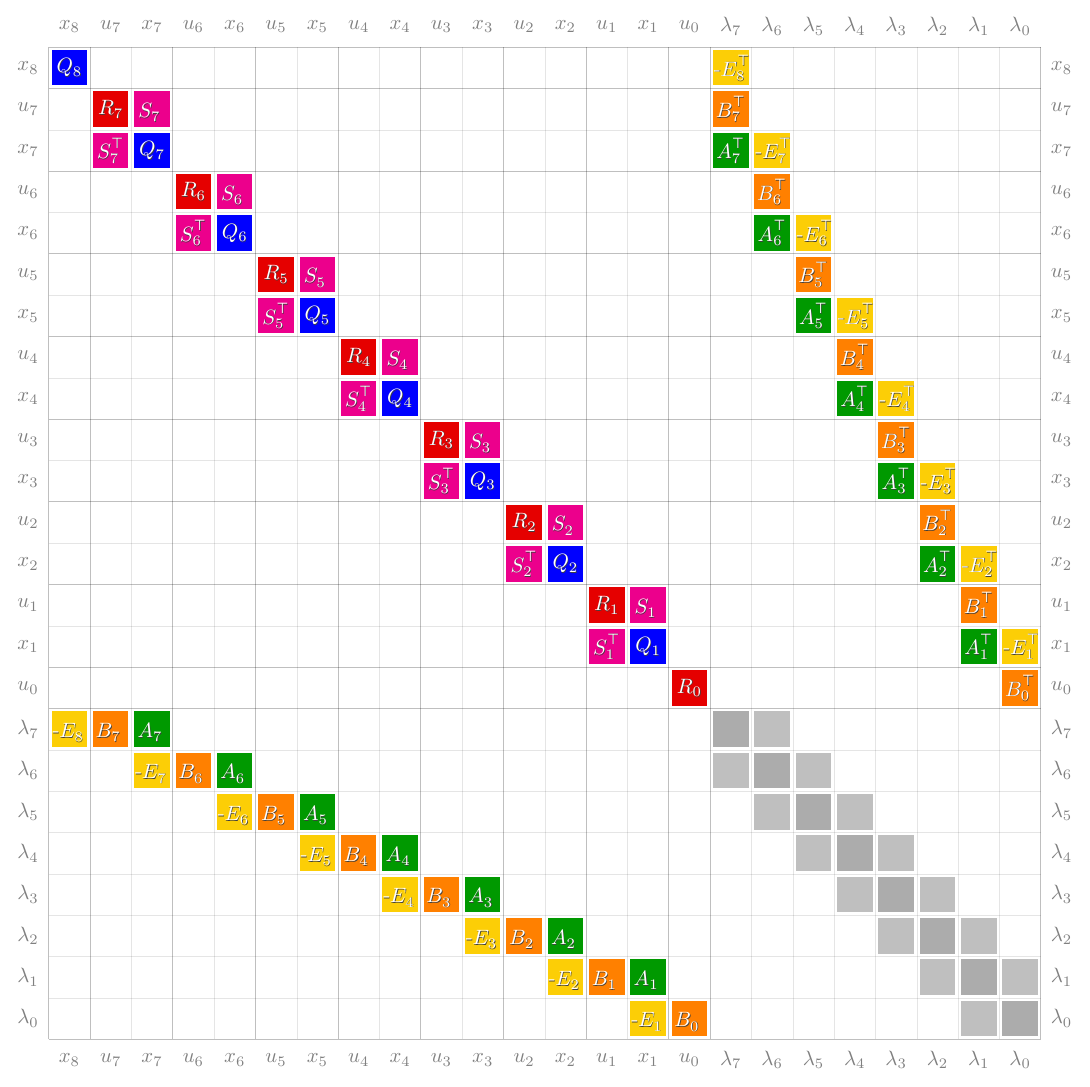}
\end{minipage}%
\begin{minipage}{0.26\textwidth}
    \centering
    \includegraphics[scale=0.44]{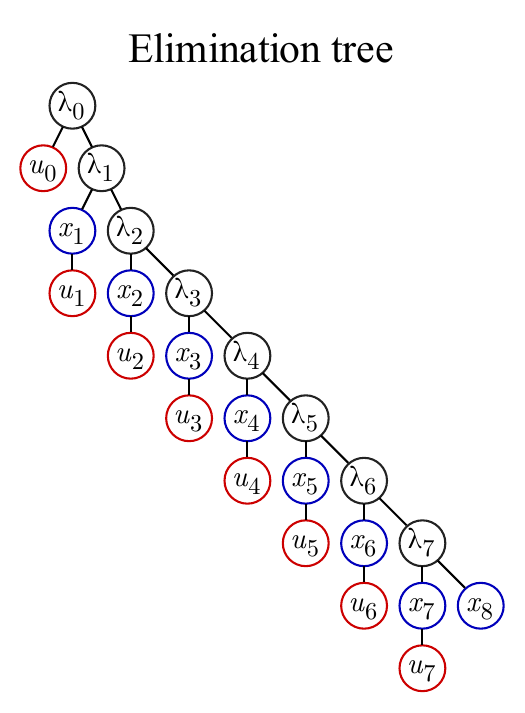}
\end{minipage}%
\vspace{-0.8em}
\caption{Traditional coefficient matrix corresponding to the optimality conditions \eqref{eq:opt-cond-ocp-eq}
of a linear-quadratic OCP \eqref{eq:ocp-eq} with horizon $N=8$. Fill-in incurred during the block Cholesky factorization is shown in gray.}
\label{mat:schur-compl-1}
\vspace{-0.5em}
\end{figure}

The coefficient matrix corresponding to \eqref{eq:opt-cond-ocp-eq} is often given in the standard KKT form as in \Cref{mat:schur-compl-1},
with the cost Hessian in the top left, and constraint Jacobians in the bottom left and top right corners.
The following subsections relate existing methods for solving \eqref{eq:opt-cond-ocp-eq} and \eqref{eq:ocp-eq} to the
block Cholesky factorization of different permutations of this matrix, based on the order in which the different variables $x^j, u^j$ and $\lambda^j$ are eliminated.
Different orderings lead to various degrees of parallelism and different amounts of fill-in, which may preserve or destroy the optimal control structure.
In the last subsection, we discuss how the proposed \cyqlone{} method improves upon these existing methods.

\subsection{Schur complement method \done} \label{sec:schur-compl-method}

The Schur complement method used in \cite{lowenstein_qpalm-ocp_2024,pas_exploiting_2024} corresponds to performing the block Cholesky factorization of the matrix in \Cref{mat:schur-compl-1}.
In the process, some fill-in is generated: elements that were structurally zero in the original sparse matrix are nonzero in its Cholesky factor.
As discussed in \Cref{sec:block-chol}, this loss of sparsity results in a higher number of operations that need to be performed during the factorization and subsequent forward and back substitutions, lowering the efficiency of the overall method.

The block diagonal structure of the cost Hessian allows the different blocks to be eliminated in parallel,
which can also be seen in the corresponding elimination tree. Such parallelization techniques are used in \cite{pas_exploiting_2024}.
Despite the parallelism available in the factorization of the cost Hessian, the generated fill-in is a block-tridiagonal matrix with $N$ blocks of size $n_x\times n_x$, 
limiting parallel scalability (as can be seen from the elimination tree, cf. \Cref{sec:block-tri-diag-chol-cr}).

\subsection{Riccati recursion \done} \label{sec:riccati-permutation}

\begin{figure}[hb]
\begin{minipage}{0.725\textwidth}
    \includegraphics[width=1\textwidth]{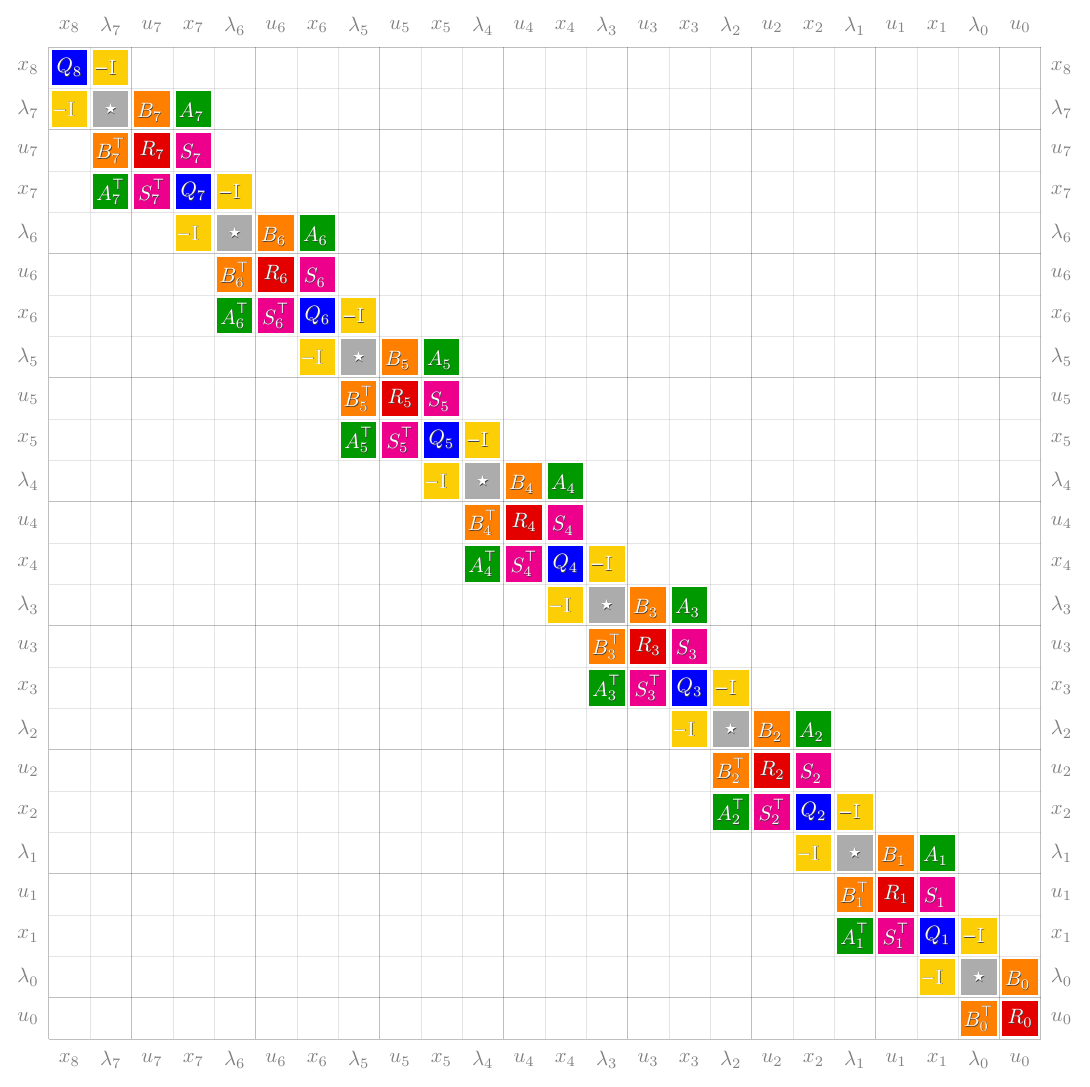}
\end{minipage}%
\begin{minipage}{0.26\textwidth}
    \centering
    \includegraphics[scale=0.44]{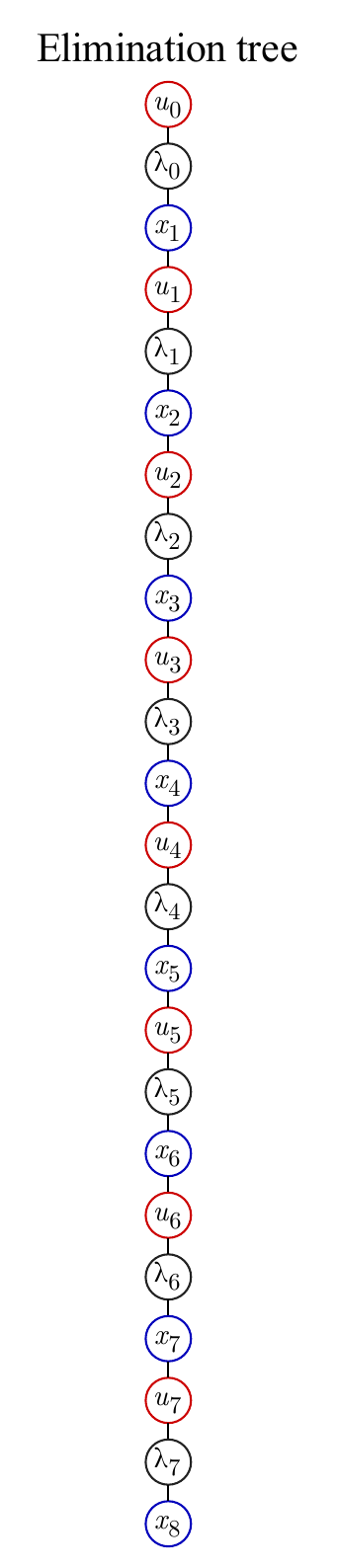}
\end{minipage}%
\vspace{-0.8em}
\caption{Permutation of \Cref{mat:schur-compl-1} corresponding to the Riccati recursion. 
Gray blocks marked by a star ($\star$) represent special redundant fill-in as described in \Cref{sec:mod-ricc}.}
\label{mat:riccati}
\vspace{-0.5em}
\end{figure}

In the case where $E_j = \I$, the Riccati recursion can be used to
solve the Newton system without any fill-in.
This method is usually derived using backward dynamic programming, resulting in the following recursive expression
of the so-called cost-to-go matrix $P_k$
\cite[\S 8.8.3]{rawlings_model_2017}:
\begin{equation} \label{eq:riccati-recursion-P}
    \scalebox{0.8}{$
    \begin{aligned}
        P_N & = Q_N                       \\[-0.5em]
        P_k & = Q_k + \tp A_k P_{k+1} A_k
        - \left( \tp S_k  + \tp A_k P_{k+1} B_k \right)
        \inv{\left( R_k  + \tp B_k P_{k+1} B_k\right)}
        \left( S_k + \tp B_k P_{k+1} A_k \right).
    \end{aligned}
    $}
\end{equation}
A naive implementation of
\eqref{eq:riccati-recursion-P} as in \cite[Alg.\,1]{frison_efficient_2013}
involves asymmetric products (e.g. $P_{k+1} A_k$), which is suboptimal.
Instead, practical implementations often use a factorized variant that exploits the
symmetry of $P_k=L_k\ttp L_k$ \cite[Alg.\,3]{frison_efficient_2013}, requiring fewer
floating-point operations.
Such factorized implementations are
specific instances of the block Cholesky factorization
of the block-tridiagonal permutation of the coefficient matrix shown in \Cref{mat:riccati}.
The Riccati recursion can easily be related to this block Cholesky factorization
by observing that $P_k$ in \eqref{eq:riccati-recursion-P} is the Schur complement of the
top-left $3\times 3$ block in
\begin{equation}
    \scalebox{0.8}{\setstretch{1.25}$\begin{pNiceArray}{ccc|c}
        P_{k+1} & -\I & & \\
        -\I & 0 & B_k & A_k \\
        & \ttp B_k & R_k & S_k \\\hline
        & \ttp A_k & \ttp S_k & Q_k
    \end{pNiceArray}.$}
\end{equation}

\subsection{Permutations that expose parallelism \done}

The downside of the Riccati recursion is that it is an inherently serial method: Equations are eliminated sequentially, with no room for parallelization.
This is evidenced by the shape of the elimination tree from \Cref{sec:riccati-permutation}, which is a single linear path.
Since the Riccati recursion is a specific instance of the block Cholesky factorization of a block-tridiagonal system (cf. \Cref{sec:block-tri-diag-chol-cr}),
a natural question that arises is whether CR-like techniques can also be used to build parallel variants of the Riccati recursion.
The answer turns out to be positive. However, part of the structure is lost in the process.
In the following sections, we summarize some different ways of applying recursive CR-like approaches to the Riccati recursion.

\subsubsection{Partitioned dynamic programming \done} \label{sec:pdp}

In the case of the Riccati recursion, \textit{backward dynamic programming} (DP) is used to summarize all
contributions of later stages using a single parametric problem described by the cost-to-go matrix $P_j$, receding recursively using \eqref{eq:riccati-recursion-P}, one stage at a time. After $P_1$ has been determined using such a backward sweep,
$\du^0$ can be computed, which is then propagated forward in time using the system dynamics.
Both operations are fundamentally sequential.
In an attempt to introduce some parallelism into the problem, and inspired by CR, Wright \cite{wright_partitioned_1991} proposed the method of \textit{partitioned dynamic programming} (PDP),
which considers \textit{multiple} parametric subproblems on different sub-intervals of the horizon rather than a
single parametric subproblem that models the entire tail of the horizon as in backward DP.
The advantage of this method is that the different sub-intervals can be eliminated in parallel.
The downside is that these parametric subproblems introduce an additional term involving the Lagrange multipliers into
the dynamics equations of the reduced problem (i.e. fill-in in the $(\lambda_j, \lambda_j)$ blocks of the matrix),
thus destroying the KKT structure (the bottom-right block of \Cref{mat:schur-compl-1} becomes nonzero) \cite[(37)]{wright_partitioned_1991}.
From the second level onwards, Wright therefore uses a more general method to handle this additional term \cite[\S 3.3]{wright_partitioned_1991},
whereas the first level is solved more efficiently using a method based on the Riccati recursion (using the asymmetric variant).

\begin{figure}
\begin{minipage}{0.725\textwidth}
    \includegraphics[width=1\textwidth]{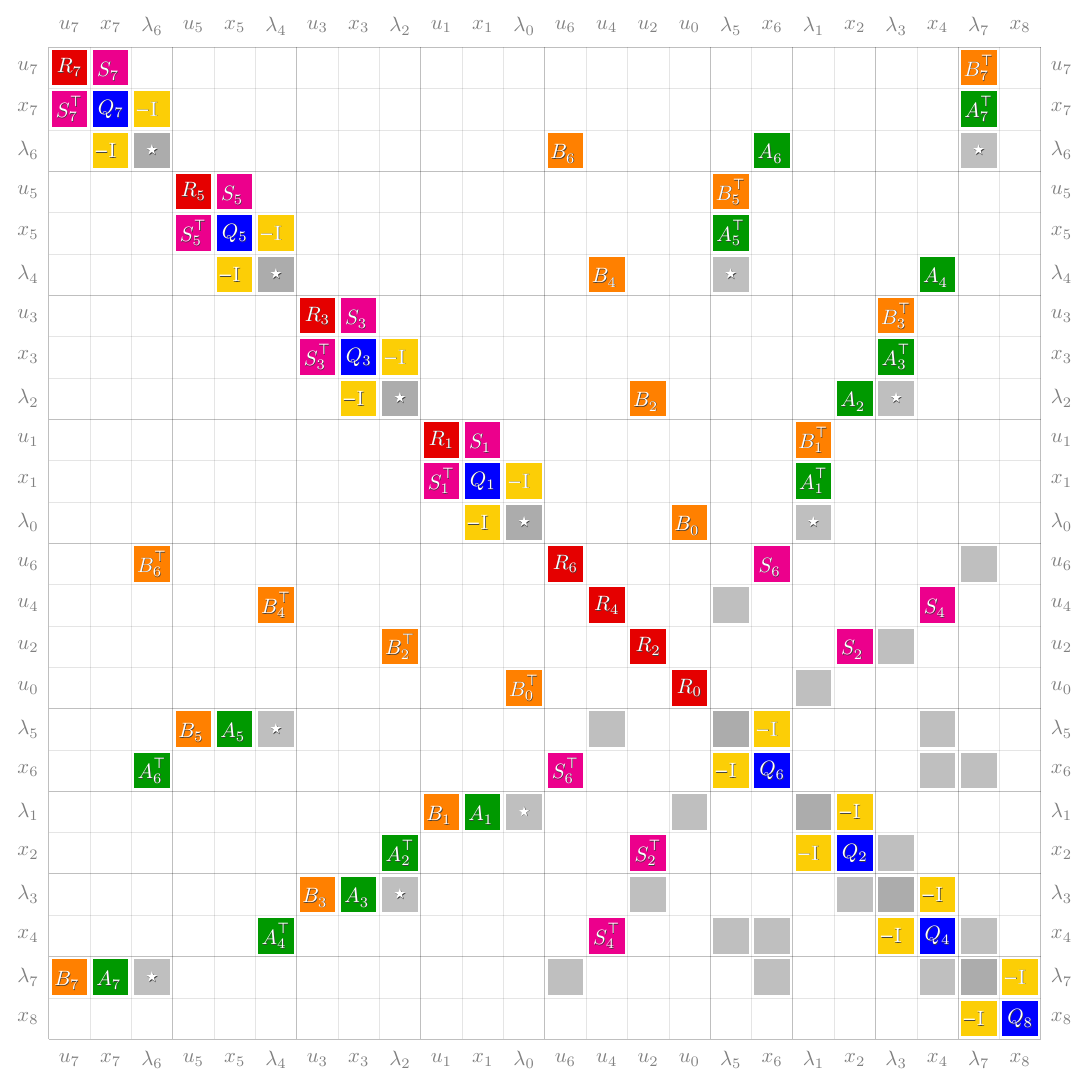}
\end{minipage}%
\begin{minipage}{0.26\textwidth}
    \centering
    \includegraphics[scale=0.44]{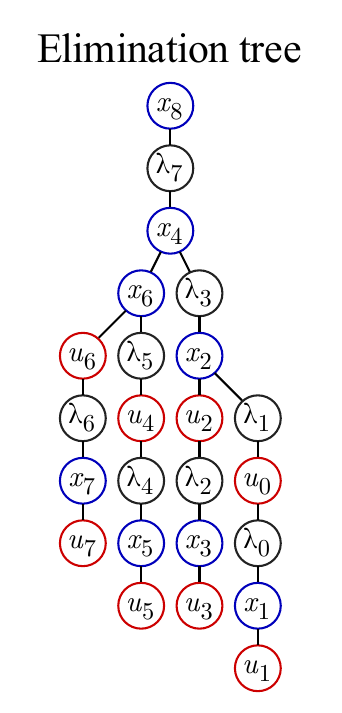}
\end{minipage}%
\vspace{-0.8em}
\caption{Permutation of \Cref{mat:schur-compl-1} corresponding to partitioned dynamic programming followed by elimination of the separator controls and cyclic reduction of the resulting system, as described by \cite[\S 3.3]{wright_partitioned_1991}.
The three steps of the method are clearly visible.}
\label{mat:pdp-pri}
\vspace{-0.5em}
\end{figure}

The full method consists of three steps:
\quad 1. Use PDP to eliminate all internal variables of the stages in the selected sub-intervals (e.g. all odd stages), keeping only the variables of the so-called separator stages (e.g. all even stages) \cite[\S 3.1]{wright_partitioned_1991};
\quad 2. Eliminate all controls in the separator stages \cite[(15)--(16)]{wright_partitioned_1991};
\quad 3. Solve the remaining block-tridiagonal system of the form \cite[(37)]{wright_partitioned_1991} using CR \cite[\S 3.2]{wright_partitioned_1991}.
This corresponds to the block Cholesky factorization of the matrix in \Cref{mat:pdp-pri}.\,\footnote{It should be noted that the blocks $\hat J_i$ in \cite[Alg.\,PRI]{wright_partitioned_1991} are not necessarily invertible, even when the full blocks $\left(\begin{smallmatrix}
    \hat J_i & -\I \\ -I & \hat Q_i
\end{smallmatrix}\right)$ are, so this needs to be considered as a single $2\times 2$ block to prevent the Cholesky factorization from breaking down.}

A related parallel method that also uses parametric subproblems on different sub-intervals of the horizon
is described by Nielsen and Axehill in \cite{nielsen_parallel_2015} and \cite[\S 6.1]{nielsen_structure_2015}.
Their approach is based on a recursive variant of partial condensing, combined with
the elimination of redundant controls in the condensed problems:
The rank of the condensed transfer matrix $(
A_n \cdots A_1 B_0 \;\;\allowbreak
A_n \cdots A_2 B_1 \;\;\allowbreak
\dotsb \;\;\allowbreak
A_n B_{n-1} \;\;\allowbreak
B_n
)$ is at most the number of states $n_x$, and therefore,
the contributions of all controls of a condensed subproblem can be summarized using a
single reduced control vector of dimension at most $n_x$. A \textit{master problem} is then
formed that couples the initial states of all intervals using these reduced controls.
The derivation in \cite[\S 6.1.2]{nielsen_structure_2015} makes use of the
singular value decomposition of the condensed transfer matrix to eliminate the components of the controls
that lie in the null space of the condensed transfer matrix, while \cite[Alg.\,13]{nielsen_structure_2015}
uses a Schur complement to eliminate redundant controls in the subproblems, which can again be viewed as a block Cholesky factorization.

An important difference between the proposed \cyqlone{} method from \Cref{sec:cyqlone} and the methods from \cite{wright_partitioned_1991} and \cite{nielsen_structure_2015}
is that \cyqlone{} eliminates all states and controls before applying CR,
whereas \cite[(42)]{wright_partitioned_1991} keeps the separator states in the reduced problems, and \cite[(6.37)]{nielsen_structure_2015}
keeps both separator states and controls in the master problem. The reason for keeping some states and controls is to obtain reduced problems with a structure similar to the original problem.
However, since the CR phase of these algorithms is the main bottleneck if the number of processors is large,
it is advantageous to reduce the size of the blocks that need to be processed using CR by eliminating the states and controls early on, as in \cyqlone{}.
The effect of the block sizes during the CR phase can be observed in the elimination trees of \Cref{mat:pdp-pri,mat:cyqlone-N8}.
\Cref{tab:compare-ocp-solvers} compares the block sizes used during the CR phase of different methods.

\subsubsection{Direct cyclic reduction \done} \label{sec:direct-cr}

\begin{figure}
\begin{minipage}{0.725\textwidth}
    \includegraphics[width=1\textwidth]{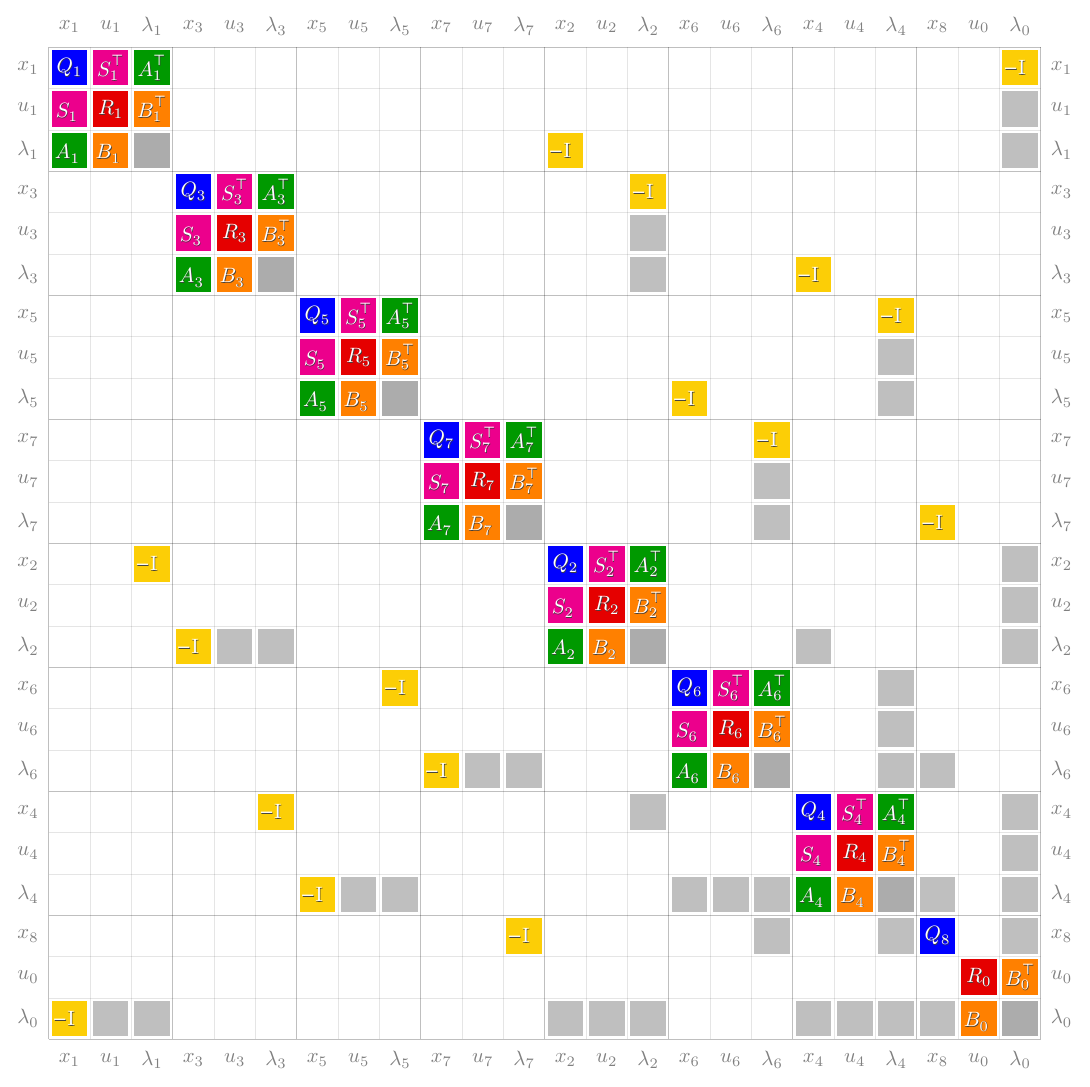}
\end{minipage}%
\begin{minipage}{0.26\textwidth}
    \centering
    \includegraphics[scale=0.44]{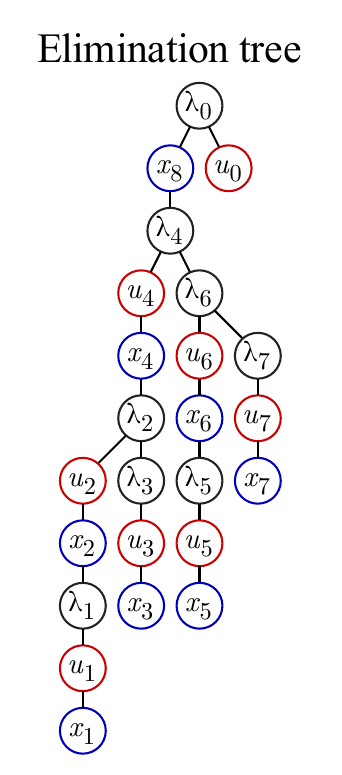}
\end{minipage}%
\vspace{-0.8em}
\caption{Permutation of \Cref{mat:schur-compl-1} corresponding to CR of the block-tridiagonal matrix from \cite[(20)--(21)]{nicholson_parallel_2018}.}
\label{mat:direct-cr-1}
\vspace{-0.5em}
\end{figure}

In \cite[\S 4.2]{nicholson_parallel_2018}, the authors write \eqref{eq:opt-cond-ocp-eq} as a block-tridiagonal system with diagonal blocks of the form
\scalebox{0.87}{$\begin{pmatrix}
    Q_j & \ttp S_j & \ttp A_j \\
    S_j & R_j & \ttp B_j \\
    A_j & B_j & 0
\end{pmatrix}$}
and subdiagonal blocks
\scalebox{0.87}{$\begin{pmatrix}
    0 & \phantom{-}0 & -\I \\
    0 & \phantom{-}0 & \phantom{-}0 \\
    0 & \phantom{-}0 & \phantom{-}0
\end{pmatrix}$}. They then apply CR to this system, which is equivalent
to block Cholesky factorization of the matrix in \Cref{mat:direct-cr-1}.
This
approach has three main disadvantages: \quad 1. It does not exploit the optimal control structure at the first level (general fill-in is incurred in all $(\lambda^j,\lambda^j)$ blocks, unlike in \Cref{mat:riccati,mat:pdp-pri}); \quad 2. CR is applied to large block matrices of size $2n_x + n_u$; and \quad 3. The number of nonzeros in the subdiagonal blocks during CR is not uniform, leading to workload imbalance if distributed across different processors.

\subsubsection{Parallel Riccati-based elimination of all states and controls \done}

In \cite{jallet_parallel_2024-1}, the authors partition the horizon into $P$ intervals which they first eliminate in parallel.
While the approach is inspired by the partitioning techniques from \cite{wright_partitioned_1991,nielsen_structure_2015} described in \Cref{sec:pdp},
the method in \cite{jallet_parallel_2024-1} differs from these methods in that it fully eliminates all state and control variables.
The elimination is done using a parametric optimization problem with optimal control structure on each interval.
This parametric optimization-based derivation is different from the block Cholesky approach used in this paper, and \cite[Alg.\,2 \& 4]{jallet_parallel_2024-1} differ from \Cref{alg:fact-mod-riccati},
but the resulting Schur complement matrices are equivalent up to a permutation, since both methods eliminate the same variables during the first phase of the algorithm.
The \textit{generalized Riccati} method in \cite[Alg.\,1]{jallet_parallel_2024-1} is a generalization of an asymmetric
Riccati-based method similar to \cite[Alg.\,1]{frison_efficient_2013},
whereas \cyqlone{}'s \Cref{alg:fact-mod-riccati} can be related to the more efficient factorized variant \cite[Alg.\,3]{frison_efficient_2013}.
After eliminating all primal variables and the Lagrange multipliers for the coupling between intervals, \cite[\S V.C]{jallet_parallel_2024-1}
uses a block variant of the Thomas algorithm. This is a sequential algorithm,
and the cost of factorizing the Schur complement (consisting of $P$ blocks of size $n_x\times n_x$) scales linearly with the number of processors $P$,
making their method less suitable for highly parallel machines.
In contrast, the cost of the CR-based approach used in \cyqlone{} scales logarithmically with the number of processors, and is therefore preferable for $P>2$.

\begin{figure}
\begin{minipage}{0.725\textwidth}
    \includegraphics[width=1\textwidth]{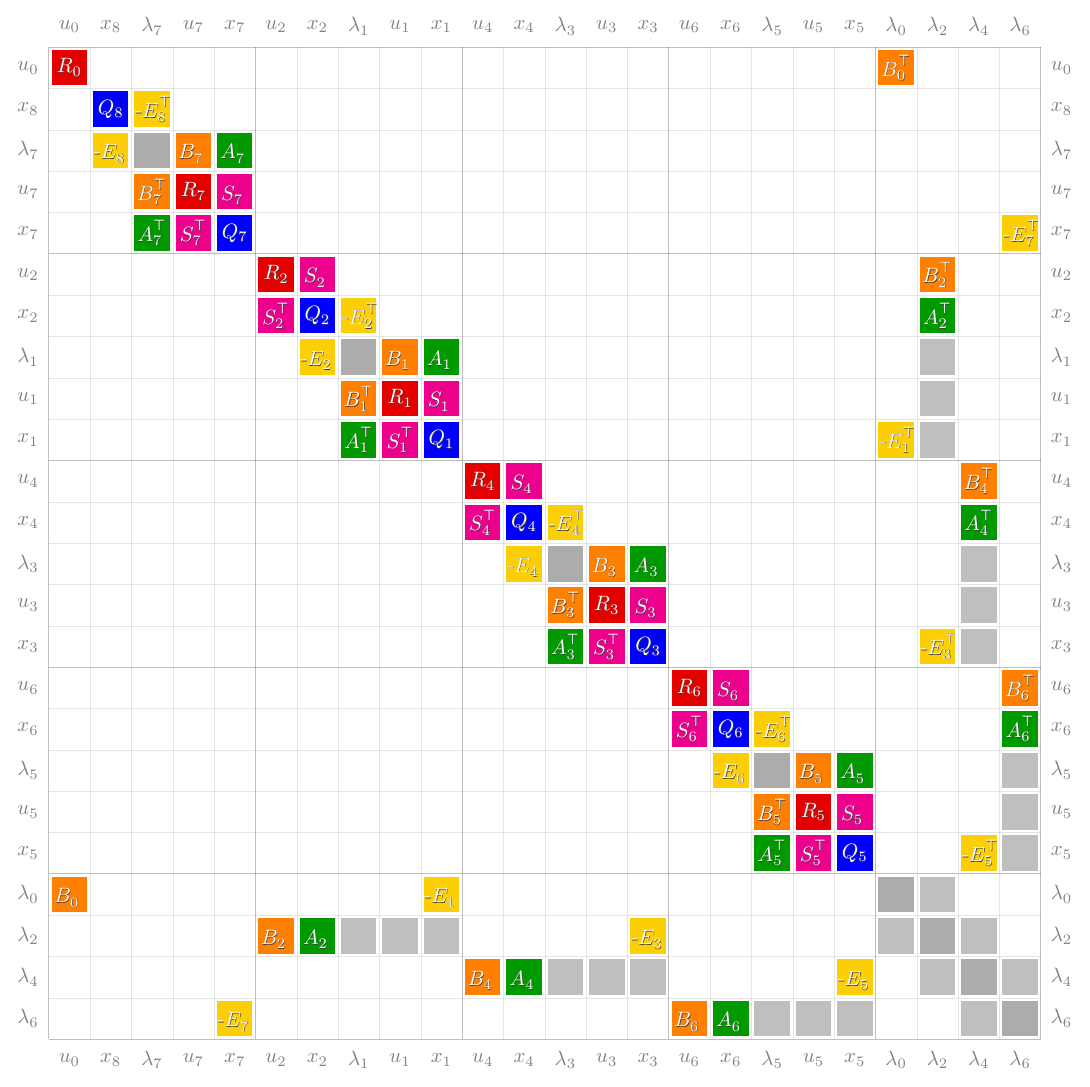}
\end{minipage}%
\begin{minipage}{0.26\textwidth}
    \centering
    \includegraphics[scale=0.44]{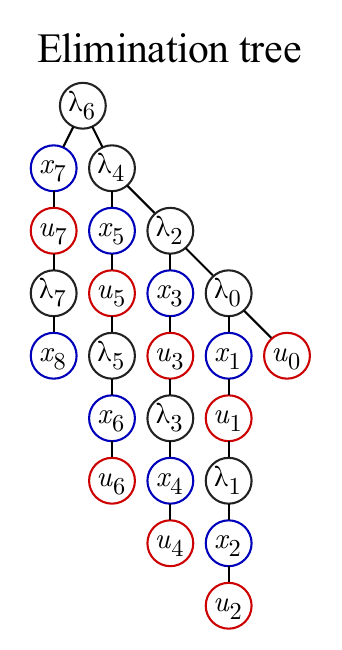}
\end{minipage}%
\vspace{-0.8em}
\caption{Variant of \Cref{mat:schur-compl-1} that eliminates four intervals of two stages in parallel,
and then solves the block-tridiagonal Schur complement, similar to what the method from \cite{jallet_parallel_2024-1} would give rise to if it made use of the Cholesky factorization.}
\label{mat:parallel-xu-elim}
\vspace{-0.5em}
\end{figure}

\Cref{mat:parallel-xu-elim} relates the parallel Riccati-based elimination of all states and controls from \cite{jallet_parallel_2024-1}
to the block Cholesky factorization of a permuted matrix for the case $N=2P=8$ (based on the elimination order of the variables;
the implementation in \cite{jallet_parallel_2024-1} does not use the block Cholesky interpretation directly).
For $P=N$, the method reduces to the standard Schur complement method from \Cref{sec:schur-compl-method}.

\subsubsection{Cyclic reduction of the Schur complement \done} \label{sec-cr-schur}

\begin{figure}
\begin{minipage}{0.725\textwidth}
    \includegraphics[width=1\textwidth]{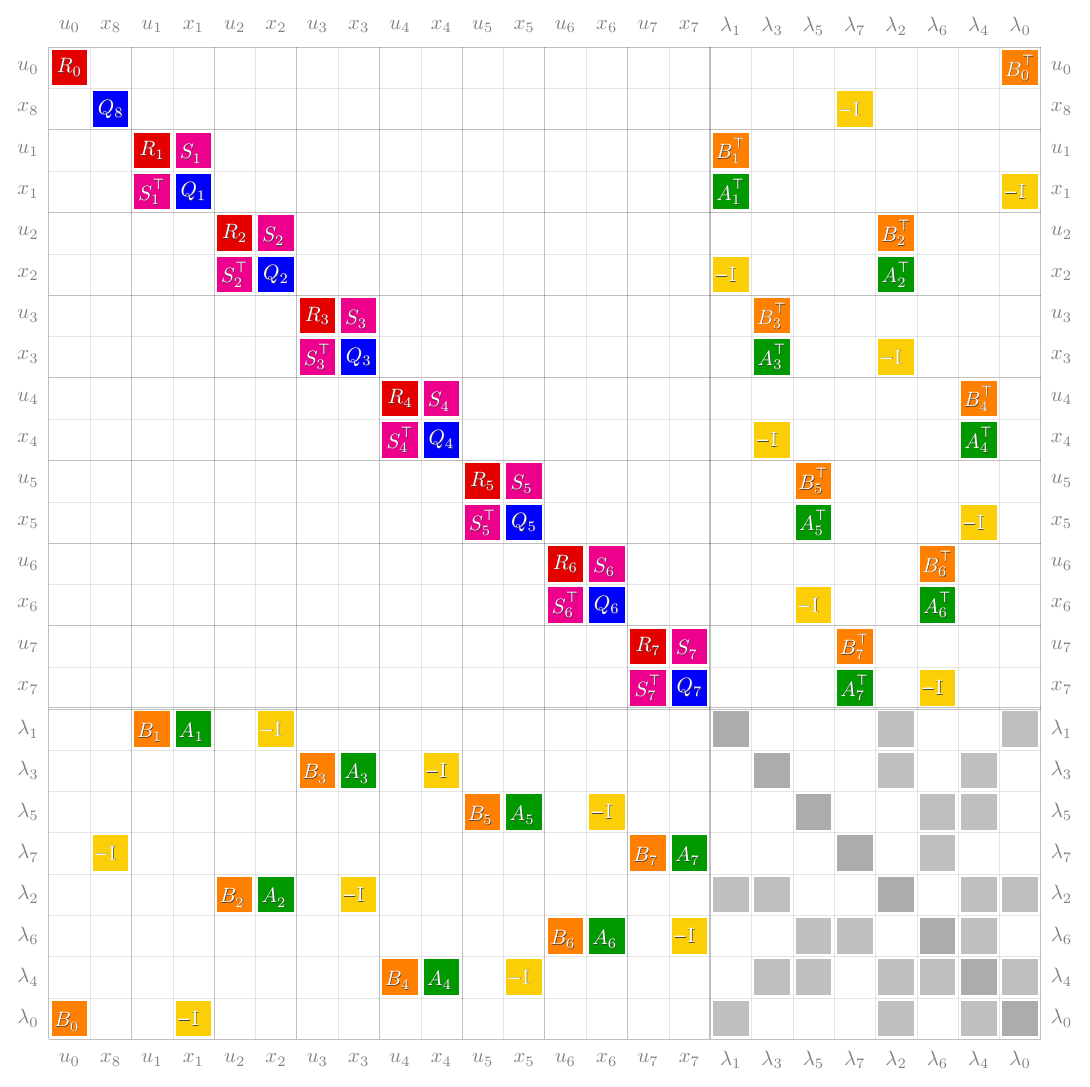}
\end{minipage}%
\begin{minipage}{0.26\textwidth}
    \centering
    \includegraphics[scale=0.44]{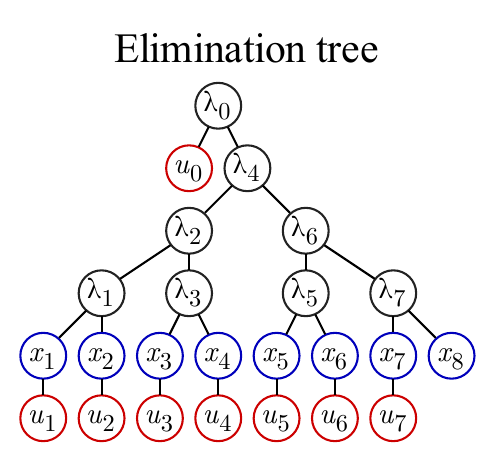}
\end{minipage}%
\vspace{-0.8em}
\caption{Variant of \Cref{mat:schur-compl-1} that uses CR for the factorization of the Schur complement.
Equivalent to \cyqlone{} for $P=8$.}
\label{mat:schur-cr}
\vspace{-0.5em}
\end{figure}

The Schur complement method from \Cref{sec:schur-compl-method} offers great parallelizability during the factorization of the cost Hessian. However, its block-tridiagonal Schur complement presents a serial bottleneck. \Cref{mat:schur-cr} shows how this can be remedied by applying CR to the Schur complement.
This approach reveals an excellent level of parallelism,
with in the shortest elimination tree so far,
thanks to the early elimination of all states and controls before applying CR (cf. \Cref{sec:pdp,sec:direct-cr}).

The qpDUNES solver \cite[\S 6.1]{frasch_parallel_2015} uses cyclic reduction of the block-tridiagonal semismooth Newton system
that arises from the dual of the original OCP, which is equivalent to applying CR to the Schur complement.

In the special case where the number of processors equals the horizon length ($N=P$),
the proposed \cyqlone{} method becomes equivalent to cyclic reduction of the Schur complement.
In this case, each interval in the partitioning from \Cref{sec:mat-struc-hi-lev} has length $n=N/P=1$, so the propagation of the cost-to-go and the dynamics in \Cref{alg:fact-mod-riccati} is skipped.

If the number of available processors is smaller than the horizon length ($P < N$),
using the Schur complement method is suboptimal. For example, in the extreme case where $P = 1$,
it is well known that the Riccati recursion requires fewer FLOPs than the Schur complement method.
One reason is that the Schur complement method does not fully exploit the structure of the constraints:
By eliminating the variables $x^j$ and $x^{j+1}$ simultaneously, the method suffers from fill-in in the
$(\lambda^j,\lambda^j)$ block. In contrast, using the Riccati recursion to eliminate $x^{j+1}$ and $\lambda^{j}$ before eliminating $x^j$
results in redundant fill-in as discussed in \Cref{sec:mod-ricc}, obviating the need for explicit factorization of the $(\lambda^j,\lambda^j)$ block.

This redundant fill-in property of the Riccati recursion is quite fragile: permuting any two rows and columns in \Cref{mat:riccati} either results in additional fill-in in one of the $(\lambda^j,\lambda^j)$ blocks,
or causes the Cholesky factorization to break down due to zero pivots.
As a general rule, permuting the system to expose parallelism introduces additional fill-in.
Therefore, by exposing more parallelism than the hardware can exploit, we incur the cost of fill-in
that cannot be hidden by further parallelizing any computations.
Using these insights, the key to an efficient parallel method for solving KKT systems with optimal control structure lies in:
\begin{enumerate}
    \setlength\itemsep{0em}
    \item Eliminating the states and controls of at least $P$ independent stages simultaneously to maximize utilization of the parallel hardware;
    \item Eliminating no more than $P$ independent stages simultaneously to minimize additional fill-in;
    \item Using the Riccati recursion where possible to exploit redundant fill-in;
    \item Using a parallelizable method like CR, PCR or PCG to factorize the Schur complement after eliminating all states and controls.
\end{enumerate}
If the number of $(\lambda^j,\lambda^j)$ blocks that receives fill-in is less than $P$, the first item is violated: this means that fewer than $P$ states are eliminated simultaneously, with some processors being idle.
If it is greater than $P$, the second item is violated: there exist permutations with the same degree of parallelism that incur less fill-in.
\Cref{tab:compare-ocp-solvers} contains the values corresponding to different methods: no fill-in for the serial Riccati recursion (suboptimal utilization of parallel hardware),
and fill-in in all $N$ of the $(\lambda^j,\lambda^j)$ blocks for Schur complement methods (higher computational cost).
Three methods partition the horizon into exactly $P$ intervals which are then eliminated using some variant of the Riccati recursion. This results in exactly $P$ blocks of fill-in in the $(\lambda^j,\lambda^j)$ blocks, which is optimal.

\clearpage
\subsubsection{\;{\normalfont\textsc{Cyqlone}} \done} \label{sec:compare-cyqlone}

\begin{figure}[hb]
\begin{minipage}{0.725\textwidth}
    \includegraphics[width=1\textwidth]{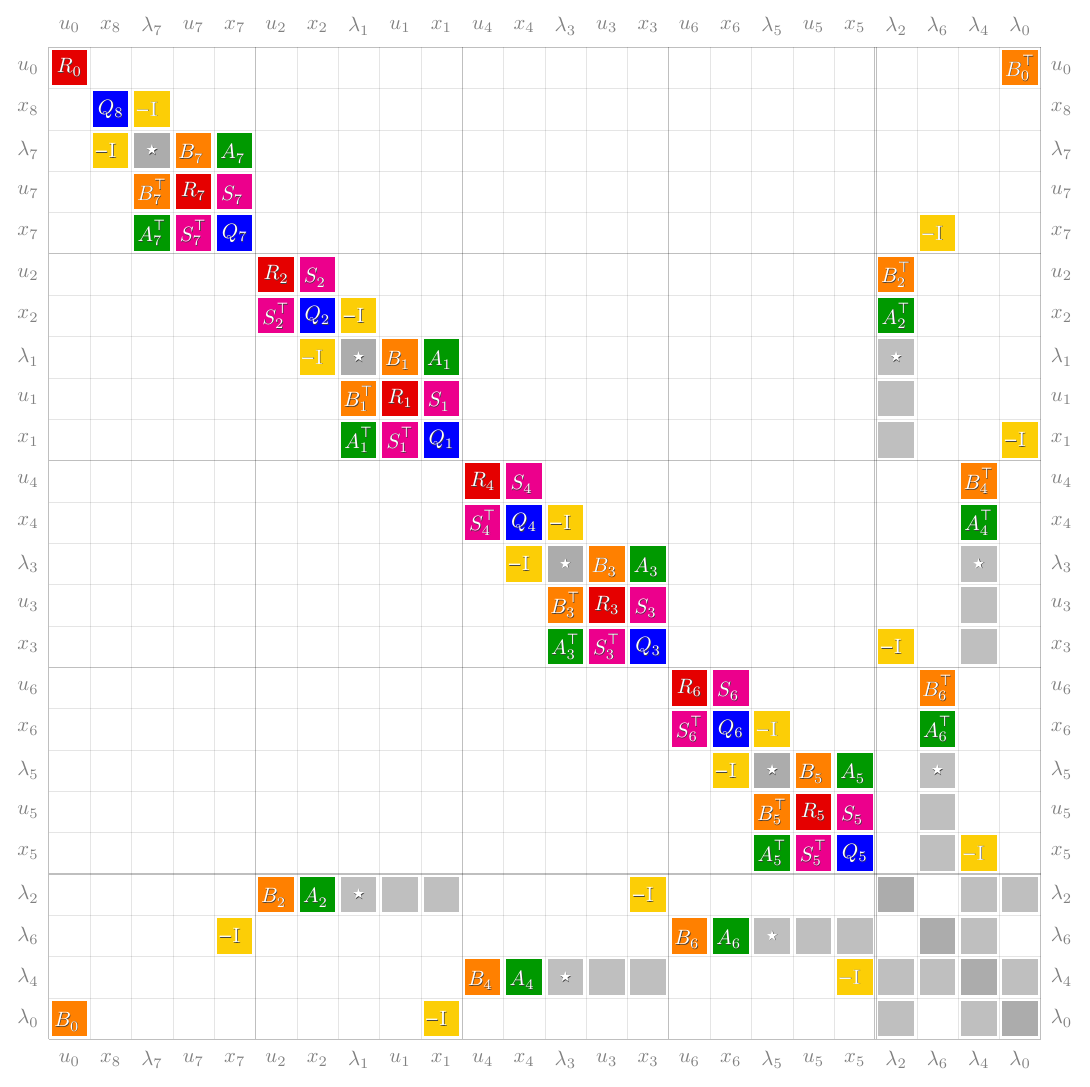}
\end{minipage}%
\begin{minipage}{0.26\textwidth}
    \centering
    \includegraphics[scale=0.44]{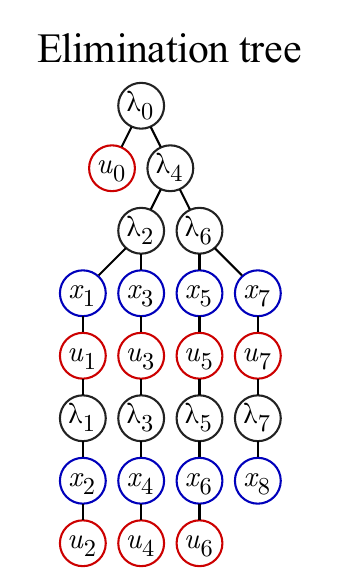}
\end{minipage}%
\vspace{-0.8em}
\caption{Variant of \Cref{mat:schur-compl-1} corresponding to the \cyqlone{} method (\Cref{sec:cyqlone}) with $P=4$.}
\label{mat:cyqlone-N8}
\vspace{-0.5em}
\end{figure}

We conclude by highlighting the advantages of the proposed \cyqlone{} method (shown in \Cref{mat:cyqlone-N8}) compared to the existing methods described above:

\begin{itemize}
    \setlength\itemsep{0em}
    \item If $P=N$, \cyqlone{} is equivalent to a parallel Schur complement method where the block-tridiagonal Schur complement is solved using CR, with a total run time proportional to $\log_2 N$.
    This method has the shortest elimination tree of all methods considered here, making it the fastest (if sufficient parallelism is available).
    \item If $P < N$, sub-intervals of the original horizon are eliminated in parallel using a modified Riccati recursion that fully exploits the optimal control structure (profiting from redundant fill-in).
    \item This elimination procedure makes use of the factorization of the cost-to-go matrix, which requires fewer floating-point operations than the traditional asymmetric variant that uses the recursion \eqref{eq:riccati-recursion-P} directly.
    \item CR is applied to the block-tridiagonal Schur complement that remains after the parallel sub-interval elimination, with a run time proportional to $\log_2 P$ (in contrast to sequential Thomas-like algorithms).
    \item The dimension of the blocks used during CR is $n_x$, as opposed to $2n_x$ or $2n_x + n_u$ in other CR-based methods. This
    results in fewer floating-point operations in the critical path of the method.
    \item If $P \ll N$, the embarrassingly parallel sub-interval elimination dominates the total run time, resulting in near-ideal scaling of the run time of approximately $\landauO(1/P)$ (cf. \Cref{fig:theoretical-speedup}).
\end{itemize}

\par\noindent\rule{\textwidth}{0.5pt}

\clearpage
\appendix

\section{Cyclic reduction and parallel cyclic reduction} \label{app:cr}
In the interest of self-containment, we list the complete CR algorithm \cite{gander_cyclic_1998,heller_aspects_1976} for the specific case of \textit{symmetric}
block-tridiagonal systems in \Cref{alg:cr}.
Compare this straightforward implementation to the \cyqlone{} implementation in the \textsc{factor-schur} function of \Cref{alg:fact-parallel-riccati-cr},
which performs the same operations with explicit multi-threading.

The closely related PCR algorithm is listed in \Cref{alg:pcr}.
There are two important differences between CR and PCR: 1. Each level of PCR processes the same
number of matrices (\scalebox{0.9}{$0\le k\lt N$}), whereas the number of active processors in CR halves at each level;
and 2. PCR is tail recursive, whereas CR requires an additional back substitution step at all levels.
As discussed in \Cref{subsec:vec-cr-last}, this enables vectorization of PCR with full vector lane utilization,
and PCR speeds up the solution of the same linear system with different right-hand sides by avoiding the back substitution step.

\begin{algorithm2e}[htbp]
    \caption{CR: Solution of a symmetric block-tridiagonal system using cyclic \rlap{reduction}}
    \label{alg:cr}
    \DontPrintSemicolon
    \KwIn{$N$: number of diagonal blocks (power of two)}
    \KwIn{$\Mdiag^{(0)}=\{\Mdiag_k\}_{k=0}^{N-1}$: diagonal blocks of the original matrix}
    \KwIn{$\Ksub^{(0)}=\{\Ksub_k\}_{k=0}^{N-2}$: subdiagonal blocks of the original matrix}
    \KwIn{$b_{(0)}=\{b^k\}_{k=0}^{N-1}$: right-hand side of the linear system}
    \KwOut{$\{x^k\}_{k=0}^{N-1}$: solution of the linear system}
    \vspace{0.5em}
    \textsc{cyclic--reduction}$(\Mdiag^{(0)},\; \Ksub^{(0)},\; b_{(0)},\;0)$\;
    \vspace{0.3em}
    \Fn{\normalfont\textsc{cyclic--reduction}$(\Mdiag^{(l)},\; \Ksub^{(l)},\; b_{(l)},\;l)$}{
        \If{$2^l = N$ \mycommentnoline{Base case (single block remaining)}}{
            $L_0 = \operatorname{chol}(\Mdiag_{0}^{(l)})$\;
            $x^0 = \invtp L_0\! \inv L_0 b^0_{(l)}$\;
            \Return\;
        }
        $p_l = N/2^{l+1}$\mycommentnofill{Number of active processors}\;
        \vspace{0.5em}
        $\Ksub_{N-2^l}^{(l)}=0$\;
        \For{$i=0, 1, \dotsc, p_l-1$ \mycommentnoline{Factor odd equations \eqref{eq:cr-updated-odd}}}{
            $k = 2^{l+1} i + 2^l$\;
            $L_{k} = \operatorname{chol}(\Mdiag_{k}^{(l)})$\;
            $Y_{k} = \Ksub_{k}^{(l)} \invtp L_{k}$\label{ln:cr-Y} \;
            $U_{k} = \Ksub_{k-2^l}^{(l) \top} \invtp L_{k}$\label{ln:cr-U} \;
            $\tilde b^{k} = \inv L_{k} b^{k}_{(l)}$\;
        }
        \vspace{0.5em}
        $\Mdiag_0^{(l+1)} = \Mdiag_0^{(l)} - U_{2^l} \ttp U_{2^l}$\mycomment{Update even equations \eqref{eq:cr-even}}
        $b^0_{(l+1)} = b^0_{(l)} - U_{2^l} \tilde b_{2^l}$\;
        \For{$i=1, \dotsc, p_l-1$}{
            $k = 2^{l+1} i$\;
            $\Mdiag_{k}^{(l+1)} = \Mdiag_{k}^{(l)} - Y_{k-2^l} \ttp Y_{k-2^l} - U_{k+2^l} \ttp U_{k+2^l}$\label{ln:cr-syrk} \;
            $\Ksub_{k-2^{l+1}}^{(l+1)} = -Y_{k-2^l} \ttp U_{k-2^l}$\label{ln:cr-gemm} \;
            $b^{k}_{(l+1)} = b^{k}_{(l)} - Y_{k-2^l} \tilde b^{k-2^l} - U_{k+2^l} \tilde b^{k+2^l}$\;
        }
        \vspace{0.4em}
        \textsc{cyclic--reduction}$(\Mdiag^{(l+1)},\; \Ksub^{(l+1)},\; b_{(l+1)},\; l\!+\!1)$\mycomment{\hspace{0.54em}\parbox[c]{9.8em}{Factor and solve even\\[-0.2em]equations (recursively)}}
        \vspace{0.4em}
        \For{$i=0, 1, \dotsc, p_l-1$ \mycommentnoline{Solve odd equations \eqref{eq:cr-odd}}}{
            $k = 2^{l+1} i + 2^l$\;
            $x^{k} = \invtp L_{k}\big(\tilde b^k - \ttp Y_{k} x^{k+2^l} - \ttp U_{k} x^{k-2^l}\big)$\;
        }
    }
\end{algorithm2e}

\begin{algorithm2e}[htbp]
    \caption{PCR: Solution of a symmetric block-tridiagonal system using parallel cyclic \rlap{reduction}}
    \label{alg:pcr}
    \DontPrintSemicolon
    \KwIn{$N$: number of diagonal blocks (power of two)}
    \KwIn{$\Mdiag^{(0)}=\{\Mdiag_k\}_{k=0}^{N-1}$: diagonal blocks of the original matrix}
    \KwIn{$\Ksub^{(0)}=\{\Ksub_k\}_{k=0}^{N-2}$: subdiagonal blocks of the original matrix, $\Ksub_{N-1}^{(0)}=0$}
    \KwIn{$b_{(0)}=\{b^k\}_{k=0}^{N-1}$: right-hand side of the linear system}
    \KwOut{$\{x^k\}_{k=0}^{N-1}$: solution of the linear system}
    \def\forallk{\hfill\NakedComment{\scalebox{0.75}{\color{gray}$(0\le k \lt N)$}\hspace{0.01em}}}
    \vspace{0.7em}

    \textsc{parallel--cyclic--reduction}$(\Mdiag^{(0)},\, \Ksub^{(0)},\, b_{(0)},\, 0)$\;

    \vspace{0.4em}
    \Fn{\normalfont\textsc{parallel--cyclic--reduction}$(\Mdiag^{(l)},\, \Ksub^{(l)},\, b_{(l)}, l)$}{
        $L_k^{(l)} = \operatorname{chol}\big(\Mdiag_k^{(l)}\big)$\forallk
        \If{$2^l = N$ \mycommentnoline{Base case ($N$ individual blocks remaining)}}{
            \Return $x^k = L_k^{(l)\,-\!\top} L_k^{(l)\,-\!1} b^k_{(l)}$\forallk
        }

        \vspace{0.5em}
        $Y_k^{(l)} = \Ksub_k^{(l)} L_k^{(l)\, -\!\top}$\forallk
        $U_k^{(l)} = \Ksub_{k-2^l}^{(l)\top} L_k^{(l)\, -\!\top}$\forallk
        $\tilde b^k = L_k^{(l)\,-\!1} b^k_{(l)}$\forallk

        \vspace{0.5em}
        $\Mdiag_k^{(l+1)} =
            \Mdiag_k^{(l)}
            - Y_{k-2^l}^{(l)} Y_{k-2^l}^{(l) \top}
            - U_{k+2^l}^{(l)} U_{\,k+2^l}^{(l) \top}$\forallk
        $\Ksub_k^{(l+1)} = -Y_{k+2^l}^{(l)} U_{k+2^l}^{(l) \top}$\forallk
        $b^k_{(l+1)} =
              b^k_{(l)}
              - Y_{k-2^l}^{(l)}\,\tilde b^{k-2^l}
              - U_{k+2^l}^{(l)}\,\tilde b^{k+2^l}$\forallk

        \vspace{0.6em}
        \mycommentstyle{$\smalltriangleright$\ {Reduce all equations recursively}} \\
        \Return \textsc{parallel--cyclic--reduction}$(\Mdiag^{(l+1)},\, \Ksub^{(l+1)},\, b_{(l+1)},\, l\!+\!1)$
    }
    \vspace{0.6em}
    For the sake of readability, subscripts $k\pm2^l$ are implicitly reduced modulo $N$.
\end{algorithm2e}

\subsection{Generalization of CR and PCR for periodic block-tridiagonal systems} \label{app:gen-cr-periodic}

\begin{figure}[p]
\begin{minipage}[t][11cm][t]{\linewidth}
\begin{algorithm2e}[H]
    \caption{Periodic PCR factorization of a block-tridiagonal matrix}
    \label{alg:pcr-periodic}
    \DontPrintSemicolon
    \KwIn{$N$: number of diagonal blocks (power of two)}
    \KwIn{$\Mdiag^{(0)}=\{\Mdiag_k\}_{k=0}^{N-1}$: diagonal blocks of the original matrix}
    \KwIn{$\Ksub^{(0)}=\{\Ksub_k\}_{k=0}^{N-1}$: subdiagonal blocks of the original matrix}
    \KwOut{$L_k^{(l)}, U_k^{(l)}, Y_k^{(l)}$: blocks of the Cholesky factors}
    \def\forallk{\hfill\NakedComment{\scalebox{0.75}{\color{gray}$(0\le k \lt N)$}\hspace{0.01em}}}
    \vspace{0.7em}

    \For{$l=0, 1, \dotsc, \log_2 N$}{
        $L_k^{(l)} = \operatorname{chol}\big(\Mdiag_k^{(l)}\big)$\forallk
        \If{$l + 1 < \log_2 N$}{
            $U_k^{(l)} = \Ksub_{k-2^l}^{(l)\top} L_k^{(l)\, -\!\top}$\forallk
            $Y_k^{(l)} = \Ksub_k^{(l)} L_k^{(l)\, -\!\top}$\forallk
            $\Mdiag_k^{(l+1)} =
            \Mdiag_k^{(l)}
                - Y_{k-2^l}^{(l)} Y_{k-2^l}^{(l) \top}
                - U_{k+2^l}^{(l)} U_{\,k+2^l}^{(l) \top}$ \forallk
            $\Ksub_k^{(l+1)} = -Y_{k+2^l}^{(l)} U_{k+2^l}^{(l) \top}$ \forallk
        }\Else(\mycommentnoline{$2\times 2$ block matrices remaining}){
            $\mathring\Ksub_{k}^{(l)} \defeq \Ksub_{k}^{(l)} + \Ksub_{k+2^l}^{(l)\top}$\mycomment{\rmrk{Overlapping forward and backward coupling}}
            $U_k^{(l)} = \mathring\Ksub_{k-2^l}^{(l)\top} L_k^{(l)\, -\!\top}$\forallk
            $\Mdiag_k^{(l+1)} =
            \Mdiag_k^{(l)}
                - U_{k+2^l}^{(l)} U_{\,k+2^l}^{(l) \top}$ \forallk
            $L_k^{(l+1)} = \operatorname{chol}\big(\Mdiag_k^{(l+1)}\big)$\forallk
        }
    }
\end{algorithm2e}
\end{minipage}
\begin{minipage}[t][14.5cm][t]{\linewidth}
\begin{figure}[H]
    \centering
    \includegraphics[width=1\textwidth,trim=0.4cm 0 0.4cm 0,clip]{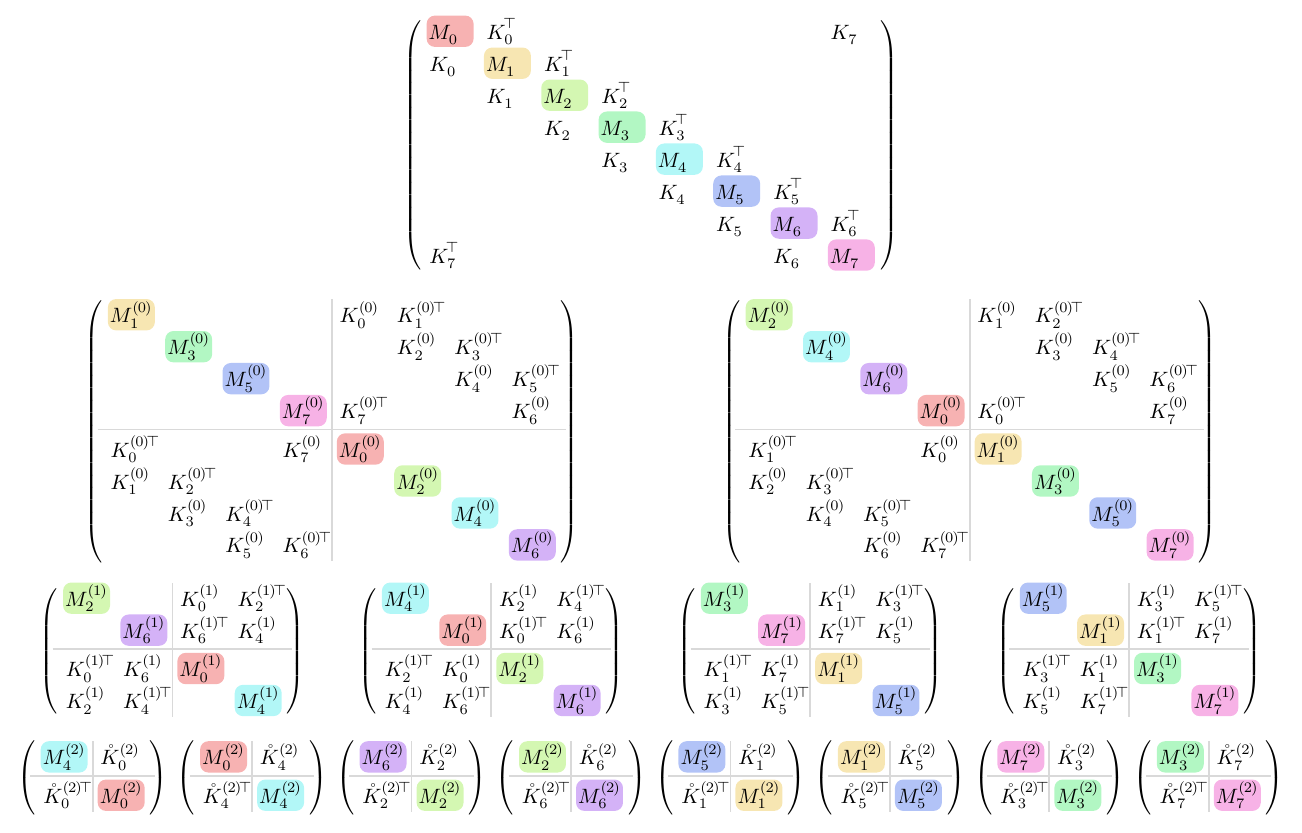}
    \caption{Visualization of PCR (\Cref{alg:pcr-periodic}) for a block-tridiagonal matrix with $N=8$ diagonal blocks of size $n\times n$.
    At each level in the tree, each branch eliminates either the odd (left) or even (right) blocks of the Schur complement of its parent.
    This reduces the coupling between them until $N$ independent dense matrices of size $2n\times 2n$ remain.
    Exactly $N$ diagonal blocks are eliminated in parallel at each level, and all matrices in the same level have the same structure. In the final level, fill-in from the forward and backward coupling blocks overlaps, resulting in the off-diagonal blocks \scalebox{0.9}{$\mathring K_k^{(l)} \defeq K_k^{(l)} + K^{(l)\top}_{k+2^l}$}.}
    \label{fig:pcr}
\end{figure}
\end{minipage}
\end{figure}

\begin{figure}[p]
\begin{minipage}[t][11cm][t]{\linewidth}
\begin{algorithm2e}[H]
    \caption{Periodic PCR factorization updates by a block-bidiagonal matrix}
    \label{alg:pcr-periodic-update}
    \DontPrintSemicolon
    \KwIn{$N$: number of diagonal blocks (power of two)}
    \KwIn{$L_k^{(l)}, U_k^{(l)}, Y_k^{(l)}$: blocks of the original Cholesky factors \rmrk{\Cref{alg:pcr-periodic}}}
    \KwIn{$\Xibwd[0]\,=\{\Xibwd_k\}_{k=0}^{N-1}$: diagonal blocks of the block-bidiagonal update matrix}
    \KwIn{$\Xifwd[0]\,=\{\Xifwd_k\}_{k=0}^{N-1}$: subdiagonal blocks of the block-bidiagonal update matrix\!\!\!\!}
    \KwIn{$\Ssigncr[0]\,=\{\mathscr S_k\}_{k=0}^{N-1}$: diagonal blocks of scaling matrices corresponding to $\Xifwd[0]$\!\!\!\!}
    \KwOut{$\tilde L_k^{(l)}, \tilde U_k^{(l)}, \tilde Y_k^{(l)}$: blocks of the updated Cholesky factors}
    \def\forallk{\hfill\NakedComment{\scalebox{0.75}{\color{gray}$(0\le k \lt N)$}\hspace{0.01em}}}
    \vspace{0.7em}

    \For{$l=0, 1, \dotsc, \log_2 N$}{
        \setstretch{1.4}
        $\mathscr S_{k}^{(l+1)} = \blkdiag\big( \mathscr S_{k}^{(l)}, \mathscr S_{k+2^l}^{(l)} \big)$\forallk
        $\begin{pNiceArray}{c|c}
            \tilde L_k^{(l)} & \;\;0\;\;
        \end{pNiceArray} = \begin{pNiceArray}{c|cc}
            L_k^{(l)} & \Xifwd[l]_k & \Xibwd[l]_k
        \end{pNiceArray} \breve Q_k^{(l)}$, where $\breve Q_k^{(l)}$ is $\left(\!\begin{smallmatrix}
            \I \\ & \phantom x\mathscr S_{k}^{(l+1)}
        \end{smallmatrix}\!\!\right)$-orth. \cite{pas_blocked_2025-1}\forallk
        \If{$l + 1 < \log_2 N$}{
            $\begin{pNiceArray}{c|c}
                \tilde U_k^{(l)} & \Xibwd[l+1]_{k-2^l}
            \end{pNiceArray} = \begin{pNiceArray}{c|cc}
                U_k^{(l)} & \Xibwd[l]_{k-2^l} & \;\;0\;\;
            \end{pNiceArray} \breve Q_{k}^{(l)}$\forallk
            $\begin{pNiceArray}{c|c}
                \tilde Y_k^{(l)} & \Xifwd[l+1]_{k+2^l}
            \end{pNiceArray} = \begin{pNiceArray}{c|cc}
                Y_k^{(l)} & \;\;0\;\; & \Xifwd[l]_{k+2^l}
            \end{pNiceArray} \breve Q_{k}^{(l)}$\forallk
        }\Else(\mycommentnoline{$2\times 2$ block matrices remaining}){
            $\begin{pNiceArray}{c|c}
                \tilde U_k^{(l)} & \Ups[l+1]_{k\pm2^l}
            \end{pNiceArray} = \begin{pNiceArray}{c|cc}
                U_k^{(l)} & \Xibwd[l]_{k-2^l} & \Xifwd[l]_{k+2^l}
            \end{pNiceArray} \breve Q_{k}^{(l)}$\forallk
            $\begin{pNiceArray}{c|c}
                \tilde L_k^{(l+1)} & \;\;0\;\;
            \end{pNiceArray} = \begin{pNiceArray}{c|c}
                L_k^{(l+1)} & \Ups[l+1]_k
            \end{pNiceArray} \breve Q_k^{(l+1)}$\forallk
        }
    }
\end{algorithm2e}
\end{minipage}
\begin{minipage}[t][14.5cm][t]{\linewidth}
\begin{figure}[H]
    \centering
    \includegraphics[width=1\textwidth,trim=0.4cm 0 0.4cm 0,clip]{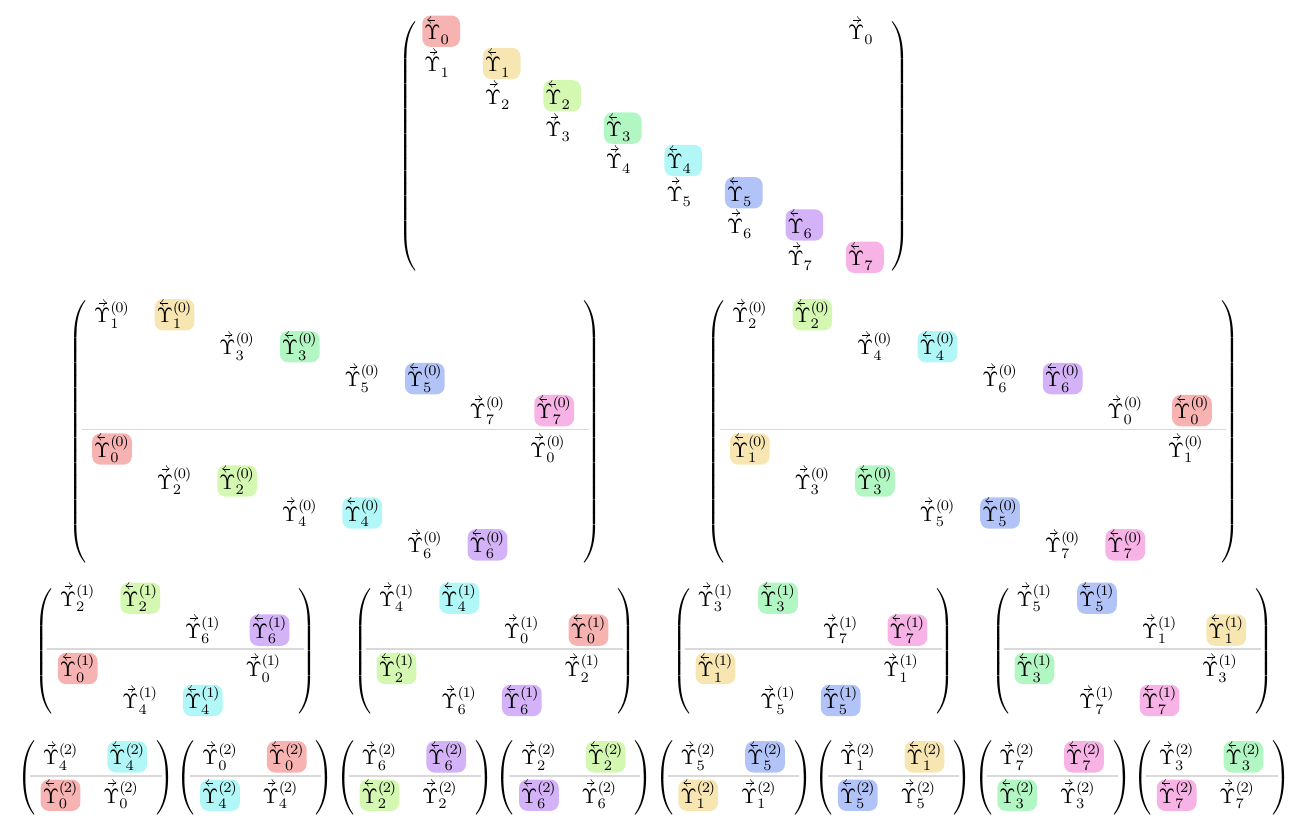}
    \caption{Visualization of factorization updates for PCR (\Cref{alg:pcr-periodic-update}) for a block-tridiagonal matrix with $N=8$ diagonal blocks of size $n\times n$
    by a block-bidiagonal update matrix. This figure shows the structure of the update matrices only; the structure of the updated factors is the same as in \Cref{fig:pcr}.}
    \label{fig:pcr-update}
\end{figure}
\end{minipage}
\end{figure}

Both CR and PCR can easily be modified to solve periodic block-tridiagonal systems, where the first and last variables are coupled by an off-diagonal block $\Ksub_{N-1}$.
The key difference is the overlap of the forward and backward coupling
blocks once the system has been reduced to a block $2\times 2$ system. The generalized algorithm is listed in \Cref{alg:pcr-periodic},
with a visualization of the structure in \Cref{fig:pcr}.
Since PCR contains CR (it is exactly the leftmost branch in \Cref{fig:pcr}, where odd equations are eliminated at each level),
we omit the periodic generalization of CR for brevity.
The factorization update procedure for updates by a block-bidiagonal matrix has a very similar structure,
and is listed in \Cref{alg:pcr-periodic-update}, with a visualization in \Cref{fig:pcr-update}.
Specifically, given the PCR factorization of a block-tridiagonal matrix $\mathscr{M}$, \Cref{alg:pcr-periodic-update} computes the PCR factorization of $\tilde{\mathscr{M}} \defeq \mathscr{M} + \mathit{\Xi} \mathscr{S} \ttp{\mathit{\Xi}}$, with
scaling matrix $\mathscr S = \blkdiag(\mathscr S_1, \dotsc, \mathscr S_{N-1}, \mathscr S_0)$,
\begin{equation}
    \label{eq:appendix-upd-mat-bidiag}
    \mathscr{M} = \begin{pmatrix}
        M_0 & \ttp K_0 & 0 & \cdots & 0 & K_{N\shortminus1} \\
        K_0 & M_1 & \ttp K_1 & \cdots & 0 & 0 \\
        0 & K_1 & M_2 & \cdots & 0 & 0 \\
        \vdots & \vdots & \vdots & \ddots & \vdots & \vdots \\
        0 & 0 & 0 & \cdots & M_{N\shortminus2} & \ttp K_{N\shortminus2} \\
        \ttp K_{\mathrlap{N\shortminus1}\phantom{11}} & 0 & 0 & \cdots & K_{N\shortminus2} & M_{N\shortminus1}
    \end{pmatrix}\!,\;\;\text{and}\;\;\;
    \mathit{\Xi} = \begin{pmatrix}
        \Xibwd_0 & 0 & 0 & \cdots & 0 & \Xifwd_0 \\
        \Xifwd_1 & \Xibwd_1 & 0 & \cdots & 0 & 0 \\
        0 & \Xifwd_2 & \Xibwd_2 & \cdots & 0 & 0 \\
        \vdots & \vdots & \vdots & \ddots & \vdots & \vdots \\
        0 & 0 & 0 & \cdots & \Xibwd_{N\shortminus2} & 0 \\
        0 & 0 & 0 & \cdots & \Xifwd_{N\shortminus1} & \Xibwd_{N\shortminus1}
    \end{pmatrix}.
\end{equation}

\subsection{Implementation of vectorized non-periodic CR using periodic CR} \label{app:cr-periodic-vectorization}

\begin{figure}[htbp]
    \centering
    \includegraphics[width=0.85\textwidth]{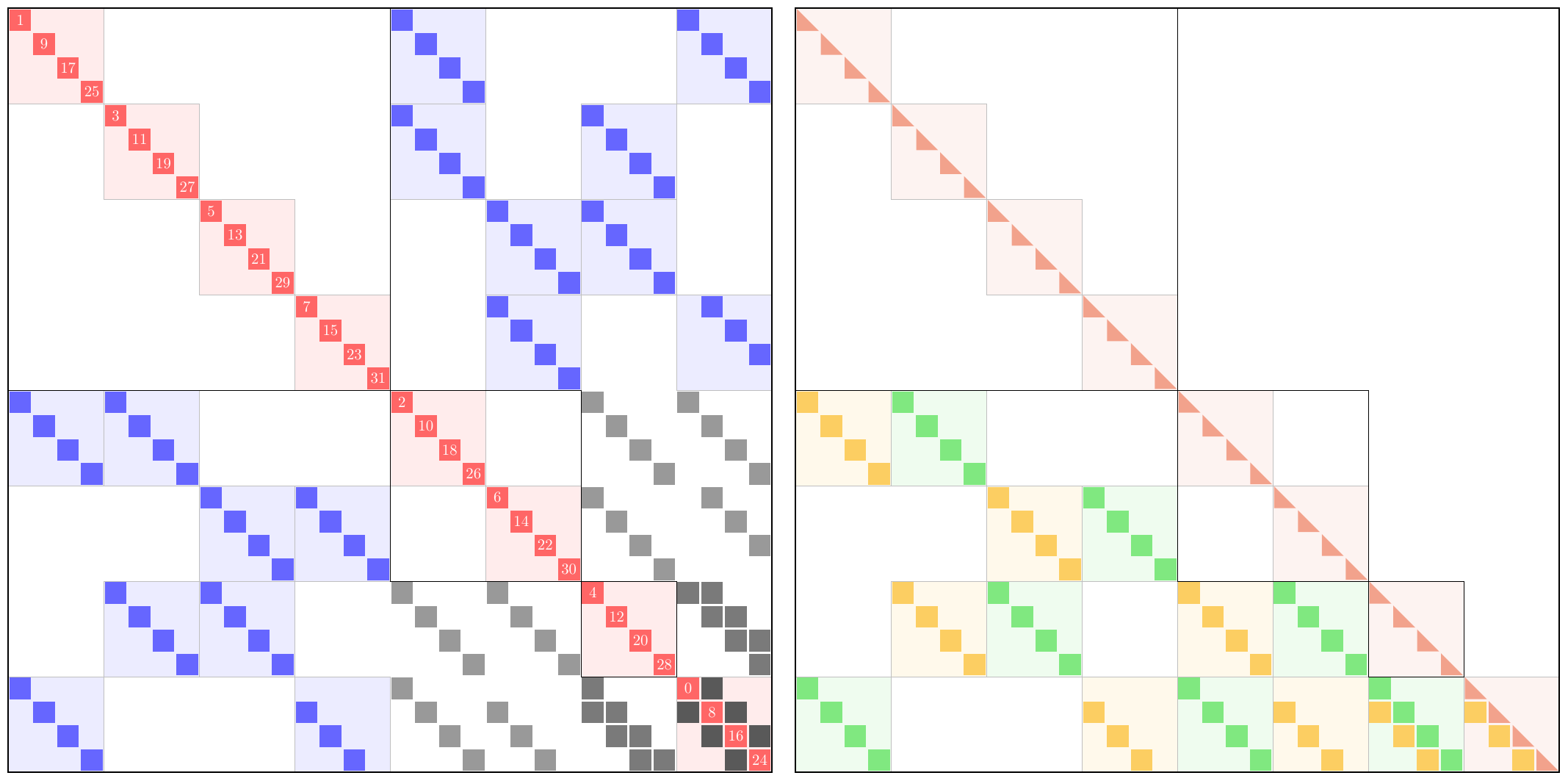}
    \caption{Alternative permutation of a $32\times 32$ block-tridiagonal matrix (left) and its Cholesky factor (right). Compared to \Cref{fig:cyclic-reduction-permutation-labels-32-vl4},
    only block rows/columns that are eliminated within the same level of CR have been swapped.
    By comparing this figure to \Cref{fig:cyclic-reduction-permutation}, we observe that the Cholesky factorization of this permutation
    can be viewed as performing CR on a $8\times 8$ block-tridiagonal matrix, with blocks of size $4n\times 4n$ instead of $n\times n$ (indicated using subtly colored squares).
    Each of these blocks is itself a $4\times 4$ block-tridiagonal matrix,
    except for the blocks coupling the last block column within each CR level with the bottom right block, where the first block-subdiagonal is nonzero.}
    \label{fig:cyclic-reduction-permutation-labels-32-vl4-vectorized}
\end{figure}

The usefulness of the periodic generalization from the previous section becomes
apparent when we consider the implementation of vectorized CR:
\Cref{subsec:vec-cr} derived the algorithm by starting from the standard
CR algorithm, and allocating operations on different blocks to different threads
and vector lanes. Alternatively, the same algorithm can also be derived by permuting
the original block-tridiagonal matrix into a periodic block-tridiagonal matrix where most
blocks are themselves block-diagonal, as shown in \Cref{fig:cyclic-reduction-permutation-labels-32-vl4-vectorized}.
Thanks to this block-diagonal sub-structure,
operations are trivially vectorized.
Implementations using this periodic view of vectorized CR allow the
multithreaded CR routines to be almost fully oblivious to the vectorization,
on the condition that they support periodic block-tridiagonal matrices and that they use batched operations instead of scalar ones.
There is one place where this abstraction leaks: the coupling between the last
and first blocks is not block-diagonal, but rather has a nonzero first block-subdiagonal.
This translates into a lane-wise rotation or shift of the batched operations,
as discussed in \Cref{subsec:vec-schur-comp,subsec:vec-cr}.
Regardless, this approach leads to a clean implementation and enables
straightforward vectorization of \Cref{alg:fact-parallel-riccati-cr}, which we leverage in the \Cpp{} implementation in \cite{kul-optec_cyqlone_2025}.

\section{Data structure details} \label{app:data-structures}
Efficient memory utilization is crucial for the performance of the CR factorization and update procedures described by \Cref{alg:fact-parallel-riccati-cr,alg:fact-parallel-riccati-cr-update}. During the CR factorization phase, the Cholesky factorization of the diagonal blocks can be performed in place, overwriting the input blocks $\Mdiag_i^{(l)}$ by $L_i$. Similarly, the subdiagonal blocks $\Kfwd[l]_i$ and $\Kbwd[l]_i$ can be overwritten by $Y_i$ and $U_i$. By avoiding additional workspaces for these matrices, the memory footprint and the total working set are reduced, improving cache efficiency and lowering the memory bandwidth requirements.

\begin{figure}
    \centering
    \includegraphics[scale=1.55,trim=0 0.08cm 0 0,clip]{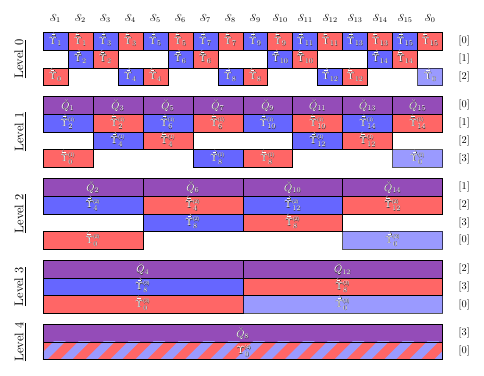}
    \caption{Schematic representation of the contents of the data structure used for storing the update matrices $\Xifwd[l]$ and $\Xibwd[l]$ at different levels of the CR factorization update procedure for the Schur complement described in \Cref{sec:cr-update-schur}, for $P=16$ processors. Indices in square brackets on the right indicate which of the four workspaces is used to store the blocks in each row.
    Each colored block represents either an update matrix $\Xifwd[l]_i$, $\Xibwd[l]_i$ or the matrix of hyperbolic Householder reflectors representing $\breve Q_i$.
    The light blue $\Xifwd[l]_0$ blocks are zero in the absence of coupling between the first and last stages (see \Cref{sec:update-schur-complement}).
    The horizontal position of each block indicates the column offset within its workspace. Note the contiguous storage of $\Xifwd[l]_i$ and $\Xibwd[l]_i$.
    Levels 3 and 4 are special, as they correspond to the update of the final $2\times 2$ block matrix
    $\left(\begin{smallmatrix} L_8 \\ U_8 & L_0 \end{smallmatrix}\right)$.
    Color coding matches that of \Cref{fig:cyclic-reduction-update-graph-16}.
    The update matrices used here are those in the leftmost branch of \Cref{fig:pcr-update} (albeit for a larger matrix).}
    \label{fig:cr-update-storage-16}
\end{figure}

In the factorization update algorithm, we can utilize in-place implementations of the routines from \cite{pas_blocked_2025-1},
where the hyperbolic Householder reflector vectors describing $\breve Q_i$ overwrite the update matrices (which are reduced to zero in the process, and therefore no longer need to be stored),
and where the updated matrices $\tilde L_i$, $\tilde U_i$ and $\tilde Y_i$ overwrite the original matrices $L_i$, $U_i$ and $Y_i$.
Unlike $\tilde L_i$, $\tilde Y_i$ and $\tilde U_i$,
the update matrices $\Xifwd[l]_i$ and $\Xibwd[l]_i$ and the transformations $\breve Q_i$
are not returned to the user, and can therefore be stored in temporary workspaces.
In fact, the update matrices used at a given level $l$ of the CR update procedure can
be discarded in the next level, allowing for efficient reuse of the workspace.
Furthermore, clever allocation of the update matrices to different workspaces
enables in-place operations with minimal data movement. Such an allocation
strategy is illustrated in \Cref{fig:cr-update-storage-16}, and
requires only four workspaces of size $n_x\times m$ each to store all update
matrices and reflector vectors, where $m$ is the total update rank
(i.e. the rank of $\mathscr S$ in \Cref{sec:update-schur-complement}).

At level 0, three of the four workspaces in \Cref{fig:cr-update-storage-16} are initialized to a compressed representation of the
full update matrix $\mathit\Xi$ from \eqref{eq:def-xi-update-example} or \eqref{eq:appendix-upd-mat-bidiag}.
In general, at level $l$ of the CR update procedure from \Cref{alg:fact-parallel-riccati-cr-update},
the update matrices $\Xifwd[l]_i$ and $\Xibwd[l]_i$ for which $\nu_2(i) = l$ are stored
contiguously in workspace $l \bmod 4$,
those for which $\nu_2(i) = l+1$ in workspace $l+1 \bmod 4$,
and those for which $\nu_2(i) > l+1$ in workspace $l+2 \bmod 4$, creating a staggered layout.
The fourth workspace, $l+3 \bmod 4$, is available to store update matrices for the next level.
As a first step in level $l$, the update matrices $\Xifwd[l]_i$ and $\Xibwd[l]_i$ for which $\nu_2(i) = l$
are reduced to zero by constructing and applying an
$\Ssigncr[l]_i$-orthogonal transformation $\breve Q_i$,
updating $L_i$ to $\tilde L_i$.
The update matrices (which are now zero)
are overwritten by the reflector vectors representing $\breve Q_i$ (which is possible thanks to their contiguous storage).
Next, the transformation $\breve Q_i$ is applied to the remaining update matrices
$\Xibwd[l]_{\vphantom i\smash{i-2^l}}$ and $\Xifwd[l]_{\vphantom i\smash{i+2^l}}$
for which $\nu_2(i\pm2^l) > l$, to obtain $\tilde U_i$ and $\tilde Y_i$,
resulting in the update matrices $\Xibwd[l+1]_{\vphantom X\smash{i-2^l}}$ and $\Xifwd[l+1]_{\vphantom X\smash{i+2^l}}$. 
These new update matrices are again stored in a staggered fashion across three workspaces, such that the matrices that are reduced in the next two levels end up in contiguous memory, as described above.
The gaps in workspaces $l+1 \bmod 4$ and $l+2 \bmod 4$ are intentional:
they allow parallel application of the transformation $\breve Q_i$ to $\Xibwd[l]_{\vphantom i\smash{i-2^l}}$ and $\Xifwd[l]_{\vphantom i\smash{i+2^l}}$.
At the end of level $l$, workspace $l \bmod 4$ still contains the reflector vectors,
which are no longer needed and can be overwritten by the next level.
The last two levels are special, as this is where the original system has been reduced to a $2\times2$ block matrix
(see \Cref{sec:cr-update-schur,app:gen-cr-periodic}).

\section{Additional benchmark results} \label{app:results}
To demonstrate the performance of \cyqlone{} across a variety of hardware platforms, we provide
additional benchmark results for the latest generation Intel Core Ultra 7 265 consumer-level CPU,
a high-end Intel Xeon Platinum 8360Y 
server CPU,
and a Raspberry Pi 5 single-board computer
in \Cref{tab:app:speedup-time-cold,tab:app:speedup-time-warm,tab:app:speedup-time-per-iter-cold}.
For larger problems with $M\ge 24$, the Ultra 7 265 and Xeon 8468 maintain excellent performance,
whereas the older Core i7-11700 shows diminishing speedups because of its relatively small L2 and L3 caches,
as compared in \Cref{tab:app:cache-sizes}.
Due to its small caches and limited memory bandwidth, the performance drop
for larger problems is more pronounced on the Raspberry Pi 5. Despite its
limited parallelism (four cores and two vector lanes), the Raspberry Pi 5
achieves impressive speedups for $M = 6$ (up to $16\times$ when warm starting).

As a general rule, we conclude that performance increases for longer horizons, where the modified Riccati recursion dominates
the overall run time. Performance for small problems (in the upper left corner of the tables)
is generally lower because the synchronization overhead is significant compared to the number of floating point operations.
In the opposite corner and near the bottom of the table, performance decreases due to the larger
working set sizes required by these larger problems.
This performance drop is more pronounced on hardware with smaller caches such as the Core i7-11700 and Raspberry Pi 5.

\def\smoll{\small}
\begin{table}[htbp]
    \centering
    \caption{Average speedup of the run time and (speedup of the worst-case run time) for \cyqpalm{} (cold start) compared to \hpipm{} (cold start).}
    \label{tab:app:speedup-time-cold}
    \small{Intel Core i7-11700 \quad ($p$=8, $v$=4) -- Clang 20.1.8} \\[0.2em]
    \smoll
     \\[0.6em]
    \small{Intel Core Ultra 7 265 \quad ($p$=8, $v$=4) -- Clang 21.1.8} \\[0.2em]
    \smoll
    \begin{tabular}{r|w{c}{4.08em}w{c}{4.08em}w{c}{4.08em}w{c}{4.08em}w{c}{4.08em}w{c}{4.08em}w{c}{4.08em}w{c}{4.08em}}
    \toprule
    $M$ & $N$=32 & $N$=64 & $N$=96 & $N$=128 & $N$=160 & $N$=192 & $N$=224 & $N$=256 \\
    \midrule
    6 & \cellcolor[rgb]{0.962,0.986,0.953}\textcolor{black}{1.1\,(\kern-0.75pt1.1\kern-0.75pt)} & \cellcolor[rgb]{0.868,0.949,0.846}\textcolor{black}{2.1\,(\kern-0.75pt2.0\kern-0.75pt)} & \cellcolor[rgb]{0.775,0.911,0.747}\textcolor{black}{2.7\,(\kern-0.75pt2.5\kern-0.75pt)} & \cellcolor[rgb]{0.653,0.860,0.629}\textcolor{black}{3.4\,(\kern-0.75pt3.2\kern-0.75pt)} & \cellcolor[rgb]{0.519,0.798,0.515}\textcolor{black}{4.1\,(\kern-0.75pt3.5\kern-0.75pt)} & \cellcolor[rgb]{0.402,0.742,0.437}\textcolor{black}{4.6\,(\kern-0.75pt4.1\kern-0.75pt)} & \cellcolor[rgb]{0.289,0.687,0.381}\textcolor{black}{5.0\,(\kern-0.75pt4.4\kern-0.75pt)} & \cellcolor[rgb]{0.208,0.621,0.327}\textcolor{black}{5.5\,(\kern-0.75pt4.4\kern-0.75pt)} \\
    12 & \cellcolor[rgb]{0.935,0.975,0.922}\textcolor{black}{1.4\,(\kern-0.75pt1.5\kern-0.75pt)} & \cellcolor[rgb]{0.809,0.925,0.783}\textcolor{black}{2.5\,(\kern-0.75pt2.5\kern-0.75pt)} & \cellcolor[rgb]{0.629,0.850,0.606}\textcolor{black}{3.5\,(\kern-0.75pt3.5\kern-0.75pt)} & \cellcolor[rgb]{0.469,0.775,0.474}\textcolor{black}{4.3\,(\kern-0.75pt4.2\kern-0.75pt)} & \cellcolor[rgb]{0.257,0.672,0.366}\textcolor{black}{5.2\,(\kern-0.75pt5.0\kern-0.75pt)} & \cellcolor[rgb]{0.149,0.558,0.280}\textcolor{white}{5.9\,(\kern-0.75pt5.2\kern-0.75pt)} & \cellcolor[rgb]{0.039,0.461,0.201}\textcolor{white}{6.6\,(\kern-0.75pt6.3\kern-0.75pt)} & \cellcolor[rgb]{0.000,0.373,0.150}\textcolor{white}{7.1\,(\kern-0.75pt6.3\kern-0.75pt)} \\
    18 & \cellcolor[rgb]{0.915,0.968,0.899}\textcolor{black}{1.6\,(\kern-0.75pt1.9\kern-0.75pt)} & \cellcolor[rgb]{0.728,0.892,0.702}\textcolor{black}{3.0\,(\kern-0.75pt2.9\kern-0.75pt)} & \cellcolor[rgb]{0.480,0.780,0.483}\textcolor{black}{4.2\,(\kern-0.75pt4.2\kern-0.75pt)} & \cellcolor[rgb]{0.264,0.675,0.369}\textcolor{black}{5.2\,(\kern-0.75pt4.7\kern-0.75pt)} & \cellcolor[rgb]{0.134,0.542,0.268}\textcolor{white}{6.0\,(\kern-0.75pt5.3\kern-0.75pt)} & \cellcolor[rgb]{0.031,0.454,0.194}\textcolor{white}{6.7\,(\kern-0.75pt7.1\kern-0.75pt)} & \cellcolor[rgb]{0.000,0.347,0.139}\textcolor{white}{7.3\,(\kern-0.75pt7.0\kern-0.75pt)} & \cellcolor[rgb]{0.000,0.267,0.106}\textcolor{white}{7.7\,(\kern-0.75pt7.5\kern-0.75pt)} \\
    24 & \cellcolor[rgb]{0.900,0.962,0.881}\textcolor{black}{1.8\,(\kern-0.75pt1.9\kern-0.75pt)} & \cellcolor[rgb]{0.648,0.858,0.624}\textcolor{black}{3.4\,(\kern-0.75pt3.4\kern-0.75pt)} & \cellcolor[rgb]{0.408,0.746,0.440}\textcolor{black}{4.6\,(\kern-0.75pt4.5\kern-0.75pt)} & \cellcolor[rgb]{0.212,0.625,0.330}\textcolor{black}{5.5\,(\kern-0.75pt5.7\kern-0.75pt)} & \cellcolor[rgb]{0.134,0.542,0.268}\textcolor{white}{6.0\,(\kern-0.75pt5.5\kern-0.75pt)} & \cellcolor[rgb]{0.035,0.457,0.198}\textcolor{white}{6.7\,(\kern-0.75pt6.2\kern-0.75pt)} & \cellcolor[rgb]{0.000,0.383,0.154}\textcolor{white}{7.1\,(\kern-0.75pt6.8\kern-0.75pt)} & \cellcolor[rgb]{0.000,0.292,0.116}\textcolor{white}{7.6\,(\kern-0.75pt7.3\kern-0.75pt)} \\
    30 & \cellcolor[rgb]{0.875,0.952,0.854}\textcolor{black}{2.0\,(\kern-0.75pt2.1\kern-0.75pt)} & \cellcolor[rgb]{0.634,0.852,0.611}\textcolor{black}{3.5\,(\kern-0.75pt3.1\kern-0.75pt)} & \cellcolor[rgb]{0.395,0.739,0.434}\textcolor{black}{4.6\,(\kern-0.75pt4.2\kern-0.75pt)} & \cellcolor[rgb]{0.234,0.648,0.348}\textcolor{black}{5.3\,(\kern-0.75pt5.1\kern-0.75pt)} & \cellcolor[rgb]{0.130,0.539,0.265}\textcolor{white}{6.1\,(\kern-0.75pt4.3\kern-0.75pt)} & \cellcolor[rgb]{0.044,0.465,0.204}\textcolor{white}{6.6\,(\kern-0.75pt5.0\kern-0.75pt)} & \cellcolor[rgb]{0.000,0.413,0.167}\textcolor{white}{6.9\,(\kern-0.75pt5.2\kern-0.75pt)} & \cellcolor[rgb]{0.000,0.393,0.158}\textcolor{white}{7.1\,(\kern-0.75pt5.0\kern-0.75pt)} \\
    \bottomrule
\end{tabular}
 \\[0.6em]
    \small{Intel Xeon Platinum 8360Y \quad ($p$=16, $v$=4) -- Clang 20.1.7} \\[0.2em]
    \smoll
    \begin{tabular}{r|w{c}{4.08em}w{c}{4.08em}w{c}{4.08em}w{c}{4.08em}w{c}{4.08em}w{c}{4.08em}w{c}{4.08em}w{c}{4.08em}}
    \toprule
    $M$ & $N$=32 & $N$=64 & $N$=96 & $N$=128 & $N$=160 & $N$=192 & $N$=224 & $N$=256 \\
    \midrule
    6 & \cellcolor[rgb]{0.966,0.987,0.958}\textcolor{black}{1.0\,(\kern-0.75pt1.1\kern-0.75pt)} & \cellcolor[rgb]{0.902,0.962,0.883}\textcolor{black}{2.2\,(\kern-0.75pt2.3\kern-0.75pt)} & \cellcolor[rgb]{0.827,0.933,0.803}\textcolor{black}{3.0\,(\kern-0.75pt2.8\kern-0.75pt)} & \cellcolor[rgb]{0.686,0.874,0.661}\textcolor{black}{4.4\,(\kern-0.75pt4.0\kern-0.75pt)} & \cellcolor[rgb]{0.579,0.827,0.565}\textcolor{black}{5.3\,(\kern-0.75pt4.5\kern-0.75pt)} & \cellcolor[rgb]{0.402,0.742,0.437}\textcolor{black}{6.5\,(\kern-0.75pt5.7\kern-0.75pt)} & \cellcolor[rgb]{0.257,0.672,0.366}\textcolor{black}{7.4\,(\kern-0.75pt6.4\kern-0.75pt)} & \cellcolor[rgb]{0.157,0.566,0.286}\textcolor{white}{\phantom{0}8.5\,\phantom{0}(\kern-0.75pt7.1\kern-0.75pt)} \\
    12 & \cellcolor[rgb]{0.949,0.980,0.938}\textcolor{black}{1.4\,(\kern-0.75pt1.4\kern-0.75pt)} & \cellcolor[rgb]{0.798,0.921,0.772}\textcolor{black}{3.4\,(\kern-0.75pt3.7\kern-0.75pt)} & \cellcolor[rgb]{0.624,0.847,0.602}\textcolor{black}{4.9\,(\kern-0.75pt4.6\kern-0.75pt)} & \cellcolor[rgb]{0.376,0.730,0.424}\textcolor{black}{6.6\,(\kern-0.75pt6.2\kern-0.75pt)} & \cellcolor[rgb]{0.234,0.648,0.348}\textcolor{black}{7.6\,(\kern-0.75pt7.3\kern-0.75pt)} & \cellcolor[rgb]{0.091,0.505,0.238}\textcolor{white}{9.1\,(\kern-0.75pt8.8\kern-0.75pt)} & \cellcolor[rgb]{0.013,0.439,0.182}\textcolor{white}{9.8\,(\kern-0.75pt9.4\kern-0.75pt)} & \cellcolor[rgb]{0.000,0.297,0.118}\textcolor{white}{11.0\,\phantom{0}(\kern-0.75pt9.7\kern-0.75pt)} \\
    18 & \cellcolor[rgb]{0.920,0.969,0.904}\textcolor{black}{1.9\,(\kern-0.75pt2.1\kern-0.75pt)} & \cellcolor[rgb]{0.704,0.882,0.679}\textcolor{black}{4.2\,(\kern-0.75pt4.3\kern-0.75pt)} & \cellcolor[rgb]{0.535,0.806,0.529}\textcolor{black}{5.6\,(\kern-0.75pt5.2\kern-0.75pt)} & \cellcolor[rgb]{0.257,0.672,0.366}\textcolor{black}{7.4\,(\kern-0.75pt6.7\kern-0.75pt)} & \cellcolor[rgb]{0.190,0.601,0.313}\textcolor{black}{8.1\,(\kern-0.75pt7.5\kern-0.75pt)} & \cellcolor[rgb]{0.026,0.450,0.191}\textcolor{white}{9.7\,(\kern-0.75pt10.0\kern-0.75pt)} & \cellcolor[rgb]{0.013,0.439,0.182}\textcolor{white}{9.9\,(\kern-0.75pt9.1\kern-0.75pt)} & \cellcolor[rgb]{0.000,0.267,0.106}\textcolor{white}{11.3\,(\kern-0.75pt11.0\kern-0.75pt)} \\
    24 & \cellcolor[rgb]{0.911,0.966,0.894}\textcolor{black}{2.0\,(\kern-0.75pt2.2\kern-0.75pt)} & \cellcolor[rgb]{0.690,0.876,0.665}\textcolor{black}{4.3\,(\kern-0.75pt4.6\kern-0.75pt)} & \cellcolor[rgb]{0.557,0.816,0.547}\textcolor{black}{5.4\,(\kern-0.75pt5.3\kern-0.75pt)} & \cellcolor[rgb]{0.270,0.678,0.372}\textcolor{black}{7.3\,(\kern-0.75pt7.5\kern-0.75pt)} & \cellcolor[rgb]{0.253,0.668,0.363}\textcolor{black}{7.4\,(\kern-0.75pt6.8\kern-0.75pt)} & \cellcolor[rgb]{0.112,0.524,0.253}\textcolor{white}{8.9\,(\kern-0.75pt8.6\kern-0.75pt)} & \cellcolor[rgb]{0.138,0.546,0.271}\textcolor{white}{8.7\,(\kern-0.75pt8.6\kern-0.75pt)} & \cellcolor[rgb]{0.001,0.428,0.173}\textcolor{white}{\phantom{0}10.0\,\phantom{0}(\kern-0.75pt9.6\kern-0.75pt)} \\
    30 & \cellcolor[rgb]{0.913,0.967,0.896}\textcolor{black}{2.0\,(\kern-0.75pt2.2\kern-0.75pt)} & \cellcolor[rgb]{0.704,0.882,0.679}\textcolor{black}{4.2\,(\kern-0.75pt4.0\kern-0.75pt)} & \cellcolor[rgb]{0.607,0.840,0.588}\textcolor{black}{5.0\,(\kern-0.75pt4.6\kern-0.75pt)} & \cellcolor[rgb]{0.383,0.733,0.427}\textcolor{black}{6.6\,(\kern-0.75pt6.3\kern-0.75pt)} & \cellcolor[rgb]{0.364,0.724,0.418}\textcolor{black}{6.7\,(\kern-0.75pt4.8\kern-0.75pt)} & \cellcolor[rgb]{0.190,0.601,0.313}\textcolor{black}{8.1\,(\kern-0.75pt6.2\kern-0.75pt)} & \cellcolor[rgb]{0.301,0.693,0.387}\textcolor{black}{7.1\,(\kern-0.75pt5.2\kern-0.75pt)} & \cellcolor[rgb]{0.186,0.597,0.310}\textcolor{black}{\phantom{0}8.2\,\phantom{0}(\kern-0.75pt5.8\kern-0.75pt)} \\
    \bottomrule
\end{tabular}
 \\[0.6em]
    \small{Raspberry Pi 5 \quad ($p$=4, $v$=2) -- Clang 21.1.8} \\[0.2em]
    \smoll
    \begin{tabular}{r|w{c}{4.08em}w{c}{4.08em}w{c}{4.08em}w{c}{4.08em}w{c}{4.08em}w{c}{4.08em}w{c}{4.08em}w{c}{4.08em}}
    \toprule
    $M$ & $N$=32 & $N$=64 & $N$=96 & $N$=128 & $N$=160 & $N$=192 & $N$=224 & $N$=256 \\
    \midrule
    6 & \cellcolor[rgb]{0.408,0.746,0.440}\textcolor{black}{3.0\,(\kern-0.75pt3.1\kern-0.75pt)} & \cellcolor[rgb]{0.039,0.461,0.201}\textcolor{white}{4.2\,(\kern-0.75pt4.4\kern-0.75pt)} & \cellcolor[rgb]{0.000,0.282,0.112}\textcolor{white}{4.8\,(\kern-0.75pt4.4\kern-0.75pt)} & \cellcolor[rgb]{0.000,0.267,0.106}\textcolor{white}{4.9\,(\kern-0.75pt4.6\kern-0.75pt)} & \cellcolor[rgb]{0.000,0.398,0.160}\textcolor{white}{4.5\,(\kern-0.75pt3.9\kern-0.75pt)} & \cellcolor[rgb]{0.000,0.383,0.154}\textcolor{white}{4.5\,(\kern-0.75pt3.9\kern-0.75pt)} & \cellcolor[rgb]{0.000,0.373,0.150}\textcolor{white}{4.5\,(\kern-0.75pt4.0\kern-0.75pt)} & \cellcolor[rgb]{0.000,0.408,0.164}\textcolor{white}{4.4\,(\kern-0.75pt3.8\kern-0.75pt)} \\
    12 & \cellcolor[rgb]{0.249,0.664,0.360}\textcolor{black}{3.4\,(\kern-0.75pt3.6\kern-0.75pt)} & \cellcolor[rgb]{0.414,0.749,0.443}\textcolor{black}{3.0\,(\kern-0.75pt3.2\kern-0.75pt)} & \cellcolor[rgb]{0.376,0.730,0.424}\textcolor{black}{3.1\,(\kern-0.75pt3.2\kern-0.75pt)} & \cellcolor[rgb]{0.445,0.764,0.458}\textcolor{black}{2.9\,(\kern-0.75pt2.9\kern-0.75pt)} & \cellcolor[rgb]{0.491,0.785,0.492}\textcolor{black}{2.8\,(\kern-0.75pt2.7\kern-0.75pt)} & \cellcolor[rgb]{0.546,0.811,0.538}\textcolor{black}{2.7\,(\kern-0.75pt2.5\kern-0.75pt)} & \cellcolor[rgb]{0.552,0.814,0.542}\textcolor{black}{2.7\,(\kern-0.75pt2.5\kern-0.75pt)} & \cellcolor[rgb]{0.568,0.822,0.556}\textcolor{black}{2.6\,(\kern-0.75pt2.3\kern-0.75pt)} \\
    18 & \cellcolor[rgb]{0.700,0.880,0.674}\textcolor{black}{2.2\,(\kern-0.75pt2.5\kern-0.75pt)} & \cellcolor[rgb]{0.714,0.886,0.688}\textcolor{black}{2.2\,(\kern-0.75pt2.1\kern-0.75pt)} & \cellcolor[rgb]{0.714,0.886,0.688}\textcolor{black}{2.2\,(\kern-0.75pt2.0\kern-0.75pt)} & \cellcolor[rgb]{0.723,0.890,0.697}\textcolor{black}{2.2\,(\kern-0.75pt1.9\kern-0.75pt)} & \cellcolor[rgb]{0.718,0.888,0.693}\textcolor{black}{2.2\,(\kern-0.75pt1.9\kern-0.75pt)} & \cellcolor[rgb]{0.695,0.878,0.670}\textcolor{black}{2.2\,(\kern-0.75pt2.2\kern-0.75pt)} & \cellcolor[rgb]{0.695,0.878,0.670}\textcolor{black}{2.2\,(\kern-0.75pt2.0\kern-0.75pt)} & \cellcolor[rgb]{0.690,0.876,0.665}\textcolor{black}{2.3\,(\kern-0.75pt2.1\kern-0.75pt)} \\
    24 & \cellcolor[rgb]{0.816,0.928,0.791}\textcolor{black}{1.8\,(\kern-0.75pt1.9\kern-0.75pt)} & \cellcolor[rgb]{0.775,0.911,0.747}\textcolor{black}{2.0\,(\kern-0.75pt2.1\kern-0.75pt)} & \cellcolor[rgb]{0.761,0.905,0.734}\textcolor{black}{2.0\,(\kern-0.75pt2.0\kern-0.75pt)} & \cellcolor[rgb]{0.746,0.899,0.720}\textcolor{black}{2.1\,(\kern-0.75pt2.1\kern-0.75pt)} & \cellcolor[rgb]{0.737,0.896,0.711}\textcolor{black}{2.1\,(\kern-0.75pt1.9\kern-0.75pt)} & \cellcolor[rgb]{0.718,0.888,0.693}\textcolor{black}{2.2\,(\kern-0.75pt2.0\kern-0.75pt)} & \cellcolor[rgb]{0.718,0.888,0.693}\textcolor{black}{2.2\,(\kern-0.75pt2.0\kern-0.75pt)} & \cellcolor[rgb]{0.718,0.888,0.693}\textcolor{black}{2.2\,(\kern-0.75pt2.0\kern-0.75pt)} \\
    30 & \cellcolor[rgb]{0.805,0.924,0.780}\textcolor{black}{1.9\,(\kern-0.75pt2.0\kern-0.75pt)} & \cellcolor[rgb]{0.751,0.901,0.724}\textcolor{black}{2.1\,(\kern-0.75pt1.8\kern-0.75pt)} & \cellcolor[rgb]{0.704,0.882,0.679}\textcolor{black}{2.2\,(\kern-0.75pt2.0\kern-0.75pt)} & \cellcolor[rgb]{0.700,0.880,0.674}\textcolor{black}{2.2\,(\kern-0.75pt2.0\kern-0.75pt)} & \cellcolor[rgb]{0.690,0.876,0.665}\textcolor{black}{2.3\,(\kern-0.75pt1.5\kern-0.75pt)} & \cellcolor[rgb]{0.681,0.872,0.656}\textcolor{black}{2.3\,(\kern-0.75pt1.7\kern-0.75pt)} & \cellcolor[rgb]{0.667,0.866,0.643}\textcolor{black}{2.3\,(\kern-0.75pt1.6\kern-0.75pt)} & \cellcolor[rgb]{0.662,0.864,0.638}\textcolor{black}{2.3\,(\kern-0.75pt1.6\kern-0.75pt)} \\
    \bottomrule
\end{tabular}
\end{table}
\begin{table}[htbp]
    \centering
    \caption{Average speedup of the run time and (speedup of the worst-case run time) for \cyqpalm{} (warm start) compared to \hpipm{} (warm start).}
    \smoll
    \label{tab:app:speedup-time-warm}
    \small{Intel Core i7-11700 \quad ($p$=8, $v$=4) -- Clang 20.1.8} \\[0.2em]
    \smoll
     \\[0.6em]
    \small{Intel Core Ultra 7 265 \quad ($p$=8, $v$=4) -- Clang 21.1.8} \\[0.2em]
    \smoll
    \begin{tabular}{r|w{c}{4.08em}w{c}{4.08em}w{c}{4.08em}w{c}{4.08em}w{c}{4.08em}w{c}{4.08em}w{c}{4.08em}w{c}{4.08em}}
    \toprule
    $M$ & $N$=32 & $N$=64 & $N$=96 & $N$=128 & $N$=160 & $N$=192 & $N$=224 & $N$=256 \\
    \midrule
    6 & \cellcolor[rgb]{0.918,0.968,0.901}\textcolor{black}{2.9\,(\kern-0.75pt1.8\kern-0.75pt)} & \cellcolor[rgb]{0.794,0.919,0.768}\textcolor{black}{\phantom{0}6.0\,(\kern-0.75pt4.0\kern-0.75pt)} & \cellcolor[rgb]{0.672,0.868,0.647}\textcolor{black}{\phantom{0}8.2\,(\kern-0.75pt4.5\kern-0.75pt)} & \cellcolor[rgb]{0.502,0.791,0.501}\textcolor{black}{10.8\,\phantom{0}(\kern-0.75pt5.4\kern-0.75pt)} & \cellcolor[rgb]{0.376,0.730,0.424}\textcolor{black}{12.5\,\phantom{0}(\kern-0.75pt6.2\kern-0.75pt)} & \cellcolor[rgb]{0.245,0.660,0.357}\textcolor{black}{14.3\,\phantom{0}(\kern-0.75pt7.1\kern-0.75pt)} & \cellcolor[rgb]{0.175,0.585,0.301}\textcolor{white}{15.8\,\phantom{0}(\kern-0.75pt9.3\kern-0.75pt)} & \cellcolor[rgb]{0.112,0.524,0.253}\textcolor{white}{17.1\,\phantom{0}(\kern-0.75pt9.1\kern-0.75pt)} \\
    12 & \cellcolor[rgb]{0.890,0.958,0.870}\textcolor{black}{3.8\,(\kern-0.75pt2.8\kern-0.75pt)} & \cellcolor[rgb]{0.723,0.890,0.697}\textcolor{black}{\phantom{0}7.2\,(\kern-0.75pt4.0\kern-0.75pt)} & \cellcolor[rgb]{0.530,0.804,0.524}\textcolor{black}{10.4\,(\kern-0.75pt6.8\kern-0.75pt)} & \cellcolor[rgb]{0.376,0.730,0.424}\textcolor{black}{12.5\,\phantom{0}(\kern-0.75pt7.3\kern-0.75pt)} & \cellcolor[rgb]{0.205,0.617,0.324}\textcolor{black}{15.2\,\phantom{0}(\kern-0.75pt9.2\kern-0.75pt)} & \cellcolor[rgb]{0.100,0.513,0.244}\textcolor{white}{17.4\,\phantom{0}(\kern-0.75pt9.8\kern-0.75pt)} & \cellcolor[rgb]{0.000,0.413,0.167}\textcolor{white}{19.5\,(\kern-0.75pt10.8\kern-0.75pt)} & \cellcolor[rgb]{0.000,0.297,0.118}\textcolor{white}{21.4\,(\kern-0.75pt12.9\kern-0.75pt)} \\
    18 & \cellcolor[rgb]{0.864,0.947,0.843}\textcolor{black}{4.4\,(\kern-0.75pt3.1\kern-0.75pt)} & \cellcolor[rgb]{0.618,0.845,0.597}\textcolor{black}{\phantom{0}9.1\,(\kern-0.75pt5.7\kern-0.75pt)} & \cellcolor[rgb]{0.402,0.742,0.437}\textcolor{black}{12.1\,(\kern-0.75pt9.0\kern-0.75pt)} & \cellcolor[rgb]{0.216,0.629,0.333}\textcolor{black}{14.9\,(\kern-0.75pt10.4\kern-0.75pt)} & \cellcolor[rgb]{0.130,0.539,0.265}\textcolor{white}{16.9\,(\kern-0.75pt12.4\kern-0.75pt)} & \cellcolor[rgb]{0.022,0.446,0.188}\textcolor{white}{18.9\,(\kern-0.75pt13.4\kern-0.75pt)} & \cellcolor[rgb]{0.000,0.357,0.144}\textcolor{white}{20.5\,(\kern-0.75pt15.3\kern-0.75pt)} & \cellcolor[rgb]{0.000,0.267,0.106}\textcolor{white}{22.0\,(\kern-0.75pt14.6\kern-0.75pt)} \\
    24 & \cellcolor[rgb]{0.842,0.938,0.819}\textcolor{black}{4.9\,(\kern-0.75pt3.6\kern-0.75pt)} & \cellcolor[rgb]{0.541,0.809,0.533}\textcolor{black}{10.2\,(\kern-0.75pt6.7\kern-0.75pt)} & \cellcolor[rgb]{0.320,0.702,0.397}\textcolor{black}{13.2\,(\kern-0.75pt9.2\kern-0.75pt)} & \cellcolor[rgb]{0.179,0.589,0.304}\textcolor{white}{15.8\,(\kern-0.75pt11.2\kern-0.75pt)} & \cellcolor[rgb]{0.125,0.535,0.262}\textcolor{white}{16.9\,(\kern-0.75pt13.0\kern-0.75pt)} & \cellcolor[rgb]{0.018,0.443,0.185}\textcolor{white}{18.9\,(\kern-0.75pt14.7\kern-0.75pt)} & \cellcolor[rgb]{0.000,0.357,0.144}\textcolor{white}{20.4\,(\kern-0.75pt14.6\kern-0.75pt)} & \cellcolor[rgb]{0.000,0.347,0.139}\textcolor{white}{20.6\,(\kern-0.75pt16.2\kern-0.75pt)} \\
    30 & \cellcolor[rgb]{0.820,0.930,0.795}\textcolor{black}{5.4\,(\kern-0.75pt3.3\kern-0.75pt)} & \cellcolor[rgb]{0.552,0.814,0.542}\textcolor{black}{10.0\,(\kern-0.75pt7.3\kern-0.75pt)} & \cellcolor[rgb]{0.307,0.696,0.390}\textcolor{black}{13.4\,(\kern-0.75pt9.1\kern-0.75pt)} & \cellcolor[rgb]{0.216,0.629,0.333}\textcolor{black}{15.0\,(\kern-0.75pt10.5\kern-0.75pt)} & \cellcolor[rgb]{0.146,0.554,0.277}\textcolor{white}{16.5\,(\kern-0.75pt13.5\kern-0.75pt)} & \cellcolor[rgb]{0.057,0.476,0.213}\textcolor{white}{18.2\,(\kern-0.75pt13.5\kern-0.75pt)} & \cellcolor[rgb]{0.052,0.472,0.210}\textcolor{white}{18.3\,(\kern-0.75pt13.7\kern-0.75pt)} & \cellcolor[rgb]{0.061,0.480,0.216}\textcolor{white}{18.2\,(\kern-0.75pt13.8\kern-0.75pt)} \\
    \bottomrule
\end{tabular}
 \\[0.6em]
    \small{Intel Xeon Platinum 8360Y \quad ($p$=16, $v$=4) -- Clang 20.1.7} \\[0.2em]
    \smoll
    \begin{tabular}{r|w{c}{4.08em}w{c}{4.08em}w{c}{4.08em}w{c}{4.08em}w{c}{4.08em}w{c}{4.08em}w{c}{4.08em}w{c}{4.08em}}
    \toprule
    $M$ & $N$=32 & $N$=64 & $N$=96 & $N$=128 & $N$=160 & $N$=192 & $N$=224 & $N$=256 \\
    \midrule
    6 & \cellcolor[rgb]{0.935,0.975,0.922}\textcolor{black}{2.9\,(\kern-0.75pt1.8\kern-0.75pt)} & \cellcolor[rgb]{0.850,0.941,0.827}\textcolor{black}{\phantom{0}6.7\,(\kern-0.75pt4.3\kern-0.75pt)} & \cellcolor[rgb]{0.765,0.907,0.738}\textcolor{black}{\phantom{0}9.5\,\phantom{0}(\kern-0.75pt5.6\kern-0.75pt)} & \cellcolor[rgb]{0.574,0.824,0.561}\textcolor{black}{14.3\,\phantom{0}(\kern-0.75pt7.2\kern-0.75pt)} & \cellcolor[rgb]{0.474,0.778,0.479}\textcolor{black}{16.6\,\phantom{0}(\kern-0.75pt8.4\kern-0.75pt)} & \cellcolor[rgb]{0.276,0.681,0.375}\textcolor{black}{20.6\,(\kern-0.75pt10.4\kern-0.75pt)} & \cellcolor[rgb]{0.182,0.593,0.307}\textcolor{white}{23.5\,(\kern-0.75pt13.4\kern-0.75pt)} & \cellcolor[rgb]{0.074,0.491,0.225}\textcolor{white}{26.9\,(\kern-0.75pt14.0\kern-0.75pt)} \\
    12 & \cellcolor[rgb]{0.920,0.969,0.904}\textcolor{black}{3.8\,(\kern-0.75pt2.7\kern-0.75pt)} & \cellcolor[rgb]{0.746,0.899,0.720}\textcolor{black}{\phantom{0}10.0\,(\kern-0.75pt5.7\kern-0.75pt)} & \cellcolor[rgb]{0.563,0.819,0.552}\textcolor{black}{14.6\,\phantom{0}(\kern-0.75pt9.8\kern-0.75pt)} & \cellcolor[rgb]{0.339,0.712,0.406}\textcolor{black}{19.3\,(\kern-0.75pt10.8\kern-0.75pt)} & \cellcolor[rgb]{0.223,0.637,0.339}\textcolor{black}{22.1\,(\kern-0.75pt14.0\kern-0.75pt)} & \cellcolor[rgb]{0.078,0.494,0.228}\textcolor{white}{26.8\,(\kern-0.75pt14.9\kern-0.75pt)} & \cellcolor[rgb]{0.026,0.450,0.191}\textcolor{white}{28.3\,(\kern-0.75pt16.3\kern-0.75pt)} & \cellcolor[rgb]{0.000,0.267,0.106}\textcolor{white}{33.1\,(\kern-0.75pt20.1\kern-0.75pt)} \\
    18 & \cellcolor[rgb]{0.894,0.959,0.874}\textcolor{black}{5.1\,(\kern-0.75pt3.6\kern-0.75pt)} & \cellcolor[rgb]{0.634,0.852,0.611}\textcolor{black}{13.0\,(\kern-0.75pt8.1\kern-0.75pt)} & \cellcolor[rgb]{0.513,0.796,0.511}\textcolor{black}{15.7\,(\kern-0.75pt12.0\kern-0.75pt)} & \cellcolor[rgb]{0.245,0.660,0.357}\textcolor{black}{21.4\,(\kern-0.75pt14.4\kern-0.75pt)} & \cellcolor[rgb]{0.208,0.621,0.327}\textcolor{black}{22.6\,(\kern-0.75pt17.3\kern-0.75pt)} & \cellcolor[rgb]{0.061,0.480,0.216}\textcolor{white}{27.2\,(\kern-0.75pt18.3\kern-0.75pt)} & \cellcolor[rgb]{0.044,0.465,0.204}\textcolor{white}{27.7\,(\kern-0.75pt20.5\kern-0.75pt)} & \cellcolor[rgb]{0.000,0.307,0.123}\textcolor{white}{32.0\,(\kern-0.75pt21.9\kern-0.75pt)} \\
    24 & \cellcolor[rgb]{0.887,0.956,0.866}\textcolor{black}{5.5\,(\kern-0.75pt4.0\kern-0.75pt)} & \cellcolor[rgb]{0.629,0.850,0.606}\textcolor{black}{13.1\,(\kern-0.75pt9.0\kern-0.75pt)} & \cellcolor[rgb]{0.513,0.796,0.511}\textcolor{black}{15.7\,(\kern-0.75pt11.6\kern-0.75pt)} & \cellcolor[rgb]{0.257,0.672,0.366}\textcolor{black}{21.0\,(\kern-0.75pt14.8\kern-0.75pt)} & \cellcolor[rgb]{0.282,0.684,0.378}\textcolor{black}{20.5\,(\kern-0.75pt15.0\kern-0.75pt)} & \cellcolor[rgb]{0.134,0.542,0.268}\textcolor{white}{25.2\,(\kern-0.75pt18.5\kern-0.75pt)} & \cellcolor[rgb]{0.149,0.558,0.280}\textcolor{white}{24.6\,(\kern-0.75pt17.2\kern-0.75pt)} & \cellcolor[rgb]{0.065,0.483,0.219}\textcolor{white}{27.1\,(\kern-0.75pt20.0\kern-0.75pt)} \\
    30 & \cellcolor[rgb]{0.887,0.956,0.866}\textcolor{black}{5.5\,(\kern-0.75pt3.3\kern-0.75pt)} & \cellcolor[rgb]{0.672,0.868,0.647}\textcolor{black}{12.0\,(\kern-0.75pt9.7\kern-0.75pt)} & \cellcolor[rgb]{0.557,0.816,0.547}\textcolor{black}{14.7\,(\kern-0.75pt10.8\kern-0.75pt)} & \cellcolor[rgb]{0.364,0.724,0.418}\textcolor{black}{18.9\,(\kern-0.75pt12.8\kern-0.75pt)} & \cellcolor[rgb]{0.395,0.739,0.434}\textcolor{black}{18.2\,(\kern-0.75pt13.7\kern-0.75pt)} & \cellcolor[rgb]{0.219,0.633,0.336}\textcolor{black}{22.2\,(\kern-0.75pt16.2\kern-0.75pt)} & \cellcolor[rgb]{0.320,0.702,0.397}\textcolor{black}{19.7\,(\kern-0.75pt14.4\kern-0.75pt)} & \cellcolor[rgb]{0.216,0.629,0.333}\textcolor{black}{22.3\,(\kern-0.75pt16.3\kern-0.75pt)} \\
    \bottomrule
\end{tabular}
 \\[0.6em]
    \small{Raspberry Pi 5 \quad ($p$=4, $v$=2) -- Clang 21.1.8} \\[0.2em]
    \smoll
    \begin{tabular}{r|w{c}{4.08em}w{c}{4.08em}w{c}{4.08em}w{c}{4.08em}w{c}{4.08em}w{c}{4.08em}w{c}{4.08em}w{c}{4.08em}}
    \toprule
    $M$ & $N$=32 & $N$=64 & $N$=96 & $N$=128 & $N$=160 & $N$=192 & $N$=224 & $N$=256 \\
    \midrule
    6 & \cellcolor[rgb]{0.463,0.773,0.470}\textcolor{black}{8.5\,(\kern-0.75pt4.9\kern-0.75pt)} & \cellcolor[rgb]{0.091,0.505,0.238}\textcolor{white}{12.9\,(\kern-0.75pt8.3\kern-0.75pt)} & \cellcolor[rgb]{0.000,0.322,0.129}\textcolor{white}{15.4\,(\kern-0.75pt8.6\kern-0.75pt)} & \cellcolor[rgb]{0.000,0.267,0.106}\textcolor{white}{16.1\,(\kern-0.75pt8.0\kern-0.75pt)} & \cellcolor[rgb]{0.000,0.413,0.167}\textcolor{white}{14.4\,(\kern-0.75pt6.9\kern-0.75pt)} & \cellcolor[rgb]{0.000,0.373,0.150}\textcolor{white}{14.8\,(\kern-0.75pt7.7\kern-0.75pt)} & \cellcolor[rgb]{0.000,0.368,0.148}\textcolor{white}{14.9\,(\kern-0.75pt8.3\kern-0.75pt)} & \cellcolor[rgb]{0.000,0.388,0.156}\textcolor{white}{14.7\,(\kern-0.75pt7.7\kern-0.75pt)} \\
    12 & \cellcolor[rgb]{0.351,0.718,0.412}\textcolor{black}{9.5\,(\kern-0.75pt6.9\kern-0.75pt)} & \cellcolor[rgb]{0.420,0.752,0.446}\textcolor{black}{\phantom{0}8.9\,(\kern-0.75pt5.0\kern-0.75pt)} & \cellcolor[rgb]{0.389,0.736,0.430}\textcolor{black}{\phantom{0}9.2\,(\kern-0.75pt6.3\kern-0.75pt)} & \cellcolor[rgb]{0.463,0.773,0.470}\textcolor{black}{\phantom{0}8.5\,(\kern-0.75pt4.8\kern-0.75pt)} & \cellcolor[rgb]{0.491,0.785,0.492}\textcolor{black}{\phantom{0}8.2\,(\kern-0.75pt5.0\kern-0.75pt)} & \cellcolor[rgb]{0.530,0.804,0.524}\textcolor{black}{\phantom{0}7.8\,(\kern-0.75pt4.6\kern-0.75pt)} & \cellcolor[rgb]{0.546,0.811,0.538}\textcolor{black}{\phantom{0}7.6\,(\kern-0.75pt4.5\kern-0.75pt)} & \cellcolor[rgb]{0.535,0.806,0.529}\textcolor{black}{\phantom{0}7.7\,(\kern-0.75pt4.6\kern-0.75pt)} \\
    18 & \cellcolor[rgb]{0.681,0.872,0.656}\textcolor{black}{6.0\,(\kern-0.75pt4.1\kern-0.75pt)} & \cellcolor[rgb]{0.639,0.854,0.615}\textcolor{black}{\phantom{0}6.6\,(\kern-0.75pt4.3\kern-0.75pt)} & \cellcolor[rgb]{0.672,0.868,0.647}\textcolor{black}{\phantom{0}6.1\,(\kern-0.75pt4.9\kern-0.75pt)} & \cellcolor[rgb]{0.681,0.872,0.656}\textcolor{black}{\phantom{0}6.0\,(\kern-0.75pt4.1\kern-0.75pt)} & \cellcolor[rgb]{0.700,0.880,0.674}\textcolor{black}{\phantom{0}5.8\,(\kern-0.75pt4.5\kern-0.75pt)} & \cellcolor[rgb]{0.690,0.876,0.665}\textcolor{black}{\phantom{0}5.9\,(\kern-0.75pt4.1\kern-0.75pt)} & \cellcolor[rgb]{0.690,0.876,0.665}\textcolor{black}{\phantom{0}5.9\,(\kern-0.75pt4.2\kern-0.75pt)} & \cellcolor[rgb]{0.681,0.872,0.656}\textcolor{black}{\phantom{0}6.0\,(\kern-0.75pt3.9\kern-0.75pt)} \\
    24 & \cellcolor[rgb]{0.775,0.911,0.747}\textcolor{black}{4.8\,(\kern-0.75pt3.6\kern-0.75pt)} & \cellcolor[rgb]{0.695,0.878,0.670}\textcolor{black}{\phantom{0}5.8\,(\kern-0.75pt4.1\kern-0.75pt)} & \cellcolor[rgb]{0.709,0.884,0.684}\textcolor{black}{\phantom{0}5.7\,(\kern-0.75pt4.0\kern-0.75pt)} & \cellcolor[rgb]{0.709,0.884,0.684}\textcolor{black}{\phantom{0}5.7\,(\kern-0.75pt4.1\kern-0.75pt)} & \cellcolor[rgb]{0.723,0.890,0.697}\textcolor{black}{\phantom{0}5.5\,(\kern-0.75pt4.3\kern-0.75pt)} & \cellcolor[rgb]{0.704,0.882,0.679}\textcolor{black}{\phantom{0}5.7\,(\kern-0.75pt4.0\kern-0.75pt)} & \cellcolor[rgb]{0.704,0.882,0.679}\textcolor{black}{\phantom{0}5.8\,(\kern-0.75pt3.8\kern-0.75pt)} & \cellcolor[rgb]{0.728,0.892,0.702}\textcolor{black}{\phantom{0}5.5\,(\kern-0.75pt4.3\kern-0.75pt)} \\
    30 & \cellcolor[rgb]{0.765,0.907,0.738}\textcolor{black}{5.0\,(\kern-0.75pt3.2\kern-0.75pt)} & \cellcolor[rgb]{0.714,0.886,0.688}\textcolor{black}{\phantom{0}5.7\,(\kern-0.75pt4.4\kern-0.75pt)} & \cellcolor[rgb]{0.681,0.872,0.656}\textcolor{black}{\phantom{0}6.1\,(\kern-0.75pt4.3\kern-0.75pt)} & \cellcolor[rgb]{0.690,0.876,0.665}\textcolor{black}{\phantom{0}5.9\,(\kern-0.75pt4.2\kern-0.75pt)} & \cellcolor[rgb]{0.695,0.878,0.670}\textcolor{black}{\phantom{0}5.9\,(\kern-0.75pt4.6\kern-0.75pt)} & \cellcolor[rgb]{0.681,0.872,0.656}\textcolor{black}{\phantom{0}6.1\,(\kern-0.75pt4.4\kern-0.75pt)} & \cellcolor[rgb]{0.690,0.876,0.665}\textcolor{black}{\phantom{0}5.9\,(\kern-0.75pt4.2\kern-0.75pt)} & \cellcolor[rgb]{0.695,0.878,0.670}\textcolor{black}{\phantom{0}5.9\,(\kern-0.75pt4.1\kern-0.75pt)} \\
    \bottomrule
\end{tabular}
\end{table}
\begin{table}[htbp]
    \centering
    \caption{Average speedup of the run time per iteration for \cyqpalm{} (cold start) compared to \hpipm{} (cold start).}
    \label{tab:app:speedup-time-per-iter-cold}
    \small{Intel Core i7-11700 \quad ($p$=8, $v$=4) -- Clang 20.1.8} \\[0.2em]
    \smoll
     \\[0.6em]
    \small{Intel Core Ultra 7 265 \quad ($p$=8, $v$=4) -- Clang 21.1.8} \\[0.2em]
    \smoll
    \begin{tabular}{r|w{c}{3.08em}w{c}{3.08em}w{c}{3.08em}w{c}{3.08em}w{c}{3.08em}w{c}{3.08em}w{c}{3.08em}w{c}{3.08em}}
    \toprule
    $M$ & $N$=32 & $N$=64 & $N$=96 & $N$=128 & $N$=160 & $N$=192 & $N$=224 & $N$=256 \\
    \midrule
    6 & \cellcolor[rgb]{0.942,0.978,0.930}\textcolor{black}{1.4} & \cellcolor[rgb]{0.839,0.937,0.815}\textcolor{black}{2.3} & \cellcolor[rgb]{0.737,0.896,0.711}\textcolor{black}{3.0} & \cellcolor[rgb]{0.629,0.850,0.606}\textcolor{black}{3.7} & \cellcolor[rgb]{0.519,0.798,0.515}\textcolor{black}{4.2} & \cellcolor[rgb]{0.402,0.742,0.437}\textcolor{black}{4.8} & \cellcolor[rgb]{0.326,0.706,0.400}\textcolor{black}{5.1} & \cellcolor[rgb]{0.234,0.648,0.348}\textcolor{black}{5.6} \\
    12 & \cellcolor[rgb]{0.913,0.967,0.896}\textcolor{black}{1.7} & \cellcolor[rgb]{0.770,0.909,0.743}\textcolor{black}{2.8} & \cellcolor[rgb]{0.596,0.835,0.579}\textcolor{black}{3.8} & \cellcolor[rgb]{0.427,0.755,0.449}\textcolor{black}{4.7} & \cellcolor[rgb]{0.238,0.652,0.351}\textcolor{black}{5.5} & \cellcolor[rgb]{0.138,0.546,0.271}\textcolor{white}{6.3} & \cellcolor[rgb]{0.044,0.465,0.204}\textcolor{white}{6.9} & \cellcolor[rgb]{0.000,0.383,0.154}\textcolor{white}{7.4} \\
    18 & \cellcolor[rgb]{0.887,0.956,0.866}\textcolor{black}{2.0} & \cellcolor[rgb]{0.672,0.868,0.647}\textcolor{black}{3.4} & \cellcolor[rgb]{0.408,0.746,0.440}\textcolor{black}{4.7} & \cellcolor[rgb]{0.219,0.633,0.336}\textcolor{black}{5.7} & \cellcolor[rgb]{0.112,0.524,0.253}\textcolor{white}{6.5} & \cellcolor[rgb]{0.018,0.443,0.185}\textcolor{white}{7.1} & \cellcolor[rgb]{0.000,0.337,0.135}\textcolor{white}{7.7} & \cellcolor[rgb]{0.000,0.267,0.106}\textcolor{white}{8.1} \\
    24 & \cellcolor[rgb]{0.857,0.944,0.835}\textcolor{black}{2.2} & \cellcolor[rgb]{0.579,0.827,0.565}\textcolor{black}{3.9} & \cellcolor[rgb]{0.326,0.706,0.400}\textcolor{black}{5.1} & \cellcolor[rgb]{0.182,0.593,0.307}\textcolor{white}{6.0} & \cellcolor[rgb]{0.091,0.505,0.238}\textcolor{white}{6.6} & \cellcolor[rgb]{0.018,0.443,0.185}\textcolor{white}{7.1} & \cellcolor[rgb]{0.000,0.363,0.146}\textcolor{white}{7.5} & \cellcolor[rgb]{0.000,0.282,0.112}\textcolor{white}{8.0} \\
    30 & \cellcolor[rgb]{0.827,0.933,0.803}\textcolor{black}{2.4} & \cellcolor[rgb]{0.541,0.809,0.533}\textcolor{black}{4.1} & \cellcolor[rgb]{0.314,0.699,0.394}\textcolor{black}{5.1} & \cellcolor[rgb]{0.182,0.593,0.307}\textcolor{black}{5.9} & \cellcolor[rgb]{0.082,0.498,0.231}\textcolor{white}{6.6} & \cellcolor[rgb]{0.001,0.428,0.173}\textcolor{white}{7.2} & \cellcolor[rgb]{0.000,0.378,0.152}\textcolor{white}{7.4} & \cellcolor[rgb]{0.000,0.378,0.152}\textcolor{white}{7.4} \\
    \bottomrule
\end{tabular}
 \\[0.6em]
    \small{Intel Xeon Platinum 8360Y \quad ($p$=16, $v$=4) -- Clang 20.1.7} \\[0.2em]
    \smoll
    \begin{tabular}{r|w{c}{3.08em}w{c}{3.08em}w{c}{3.08em}w{c}{3.08em}w{c}{3.08em}w{c}{3.08em}w{c}{3.08em}w{c}{3.08em}}
    \toprule
    $M$ & $N$=32 & $N$=64 & $N$=96 & $N$=128 & $N$=160 & $N$=192 & $N$=224 & $N$=256 \\
    \midrule
    6 & \cellcolor[rgb]{0.955,0.983,0.945}\textcolor{black}{1.3} & \cellcolor[rgb]{0.883,0.955,0.862}\textcolor{black}{2.5} & \cellcolor[rgb]{0.809,0.925,0.783}\textcolor{black}{3.4} & \cellcolor[rgb]{0.672,0.868,0.647}\textcolor{black}{4.7} & \cellcolor[rgb]{0.579,0.827,0.565}\textcolor{black}{5.4} & \cellcolor[rgb]{0.408,0.746,0.440}\textcolor{black}{\phantom{0}6.7} & \cellcolor[rgb]{0.301,0.693,0.387}\textcolor{black}{\phantom{0}7.4} & \cellcolor[rgb]{0.182,0.593,0.307}\textcolor{white}{\phantom{0}8.6} \\
    12 & \cellcolor[rgb]{0.933,0.974,0.919}\textcolor{black}{1.7} & \cellcolor[rgb]{0.770,0.909,0.743}\textcolor{black}{3.8} & \cellcolor[rgb]{0.591,0.832,0.574}\textcolor{black}{5.4} & \cellcolor[rgb]{0.333,0.709,0.403}\textcolor{black}{7.2} & \cellcolor[rgb]{0.223,0.637,0.339}\textcolor{black}{8.1} & \cellcolor[rgb]{0.078,0.494,0.228}\textcolor{white}{\phantom{0}9.7} & \cellcolor[rgb]{0.026,0.450,0.191}\textcolor{white}{10.2} & \cellcolor[rgb]{0.000,0.302,0.121}\textcolor{white}{11.5} \\
    18 & \cellcolor[rgb]{0.902,0.962,0.883}\textcolor{black}{2.3} & \cellcolor[rgb]{0.653,0.860,0.629}\textcolor{black}{4.9} & \cellcolor[rgb]{0.480,0.780,0.483}\textcolor{black}{6.2} & \cellcolor[rgb]{0.216,0.629,0.333}\textcolor{black}{8.2} & \cellcolor[rgb]{0.171,0.582,0.298}\textcolor{white}{8.7} & \cellcolor[rgb]{0.009,0.435,0.179}\textcolor{white}{10.3} & \cellcolor[rgb]{0.005,0.432,0.176}\textcolor{white}{10.4} & \cellcolor[rgb]{0.000,0.267,0.106}\textcolor{white}{11.8} \\
    24 & \cellcolor[rgb]{0.890,0.958,0.870}\textcolor{black}{2.5} & \cellcolor[rgb]{0.644,0.856,0.620}\textcolor{black}{4.9} & \cellcolor[rgb]{0.507,0.793,0.506}\textcolor{black}{6.0} & \cellcolor[rgb]{0.238,0.652,0.351}\textcolor{black}{7.9} & \cellcolor[rgb]{0.227,0.641,0.342}\textcolor{black}{8.1} & \cellcolor[rgb]{0.091,0.505,0.238}\textcolor{white}{\phantom{0}9.5} & \cellcolor[rgb]{0.121,0.531,0.259}\textcolor{white}{\phantom{0}9.3} & \cellcolor[rgb]{0.000,0.423,0.171}\textcolor{white}{10.5} \\
    30 & \cellcolor[rgb]{0.890,0.958,0.870}\textcolor{black}{2.5} & \cellcolor[rgb]{0.644,0.856,0.620}\textcolor{black}{4.9} & \cellcolor[rgb]{0.557,0.816,0.547}\textcolor{black}{5.6} & \cellcolor[rgb]{0.307,0.696,0.390}\textcolor{black}{7.4} & \cellcolor[rgb]{0.307,0.696,0.390}\textcolor{black}{7.4} & \cellcolor[rgb]{0.168,0.578,0.295}\textcolor{white}{\phantom{0}8.7} & \cellcolor[rgb]{0.276,0.681,0.375}\textcolor{black}{\phantom{0}7.6} & \cellcolor[rgb]{0.175,0.585,0.301}\textcolor{white}{\phantom{0}8.7} \\
    \bottomrule
\end{tabular}
 \\[0.6em]
    \small{Raspberry Pi 5 \quad ($p$=4, $v$=2) -- Clang 21.1.8} \\[0.2em]
    \smoll
    \begin{tabular}{r|w{c}{3.08em}w{c}{3.08em}w{c}{3.08em}w{c}{3.08em}w{c}{3.08em}w{c}{3.08em}w{c}{3.08em}w{c}{3.08em}}
    \toprule
    $M$ & $N$=32 & $N$=64 & $N$=96 & $N$=128 & $N$=160 & $N$=192 & $N$=224 & $N$=256 \\
    \midrule
    6 & \cellcolor[rgb]{0.249,0.664,0.360}\textcolor{black}{3.8} & \cellcolor[rgb]{0.000,0.418,0.169}\textcolor{white}{4.8} & \cellcolor[rgb]{0.000,0.267,0.106}\textcolor{white}{5.4} & \cellcolor[rgb]{0.000,0.307,0.123}\textcolor{white}{5.2} & \cellcolor[rgb]{0.044,0.465,0.204}\textcolor{white}{4.6} & \cellcolor[rgb]{0.035,0.457,0.198}\textcolor{white}{4.7} & \cellcolor[rgb]{0.065,0.483,0.219}\textcolor{white}{4.5} & \cellcolor[rgb]{0.078,0.494,0.228}\textcolor{white}{4.5} \\
    12 & \cellcolor[rgb]{0.160,0.570,0.289}\textcolor{white}{4.1} & \cellcolor[rgb]{0.376,0.730,0.424}\textcolor{black}{3.4} & \cellcolor[rgb]{0.376,0.730,0.424}\textcolor{black}{3.4} & \cellcolor[rgb]{0.445,0.764,0.458}\textcolor{black}{3.2} & \cellcolor[rgb]{0.507,0.793,0.506}\textcolor{black}{3.0} & \cellcolor[rgb]{0.557,0.816,0.547}\textcolor{black}{2.9} & \cellcolor[rgb]{0.585,0.829,0.570}\textcolor{black}{2.8} & \cellcolor[rgb]{0.602,0.837,0.583}\textcolor{black}{2.7} \\
    18 & \cellcolor[rgb]{0.613,0.842,0.592}\textcolor{black}{2.7} & \cellcolor[rgb]{0.662,0.864,0.638}\textcolor{black}{2.5} & \cellcolor[rgb]{0.690,0.876,0.665}\textcolor{black}{2.4} & \cellcolor[rgb]{0.704,0.882,0.679}\textcolor{black}{2.4} & \cellcolor[rgb]{0.718,0.888,0.693}\textcolor{black}{2.3} & \cellcolor[rgb]{0.709,0.884,0.684}\textcolor{black}{2.4} & \cellcolor[rgb]{0.704,0.882,0.679}\textcolor{black}{2.4} & \cellcolor[rgb]{0.704,0.882,0.679}\textcolor{black}{2.4} \\
    24 & \cellcolor[rgb]{0.751,0.901,0.724}\textcolor{black}{2.2} & \cellcolor[rgb]{0.732,0.894,0.706}\textcolor{black}{2.3} & \cellcolor[rgb]{0.732,0.894,0.706}\textcolor{black}{2.3} & \cellcolor[rgb]{0.737,0.896,0.711}\textcolor{black}{2.3} & \cellcolor[rgb]{0.728,0.892,0.702}\textcolor{black}{2.3} & \cellcolor[rgb]{0.718,0.888,0.693}\textcolor{black}{2.3} & \cellcolor[rgb]{0.723,0.890,0.697}\textcolor{black}{2.3} & \cellcolor[rgb]{0.732,0.894,0.706}\textcolor{black}{2.3} \\
    30 & \cellcolor[rgb]{0.742,0.897,0.715}\textcolor{black}{2.2} & \cellcolor[rgb]{0.695,0.878,0.670}\textcolor{black}{2.4} & \cellcolor[rgb]{0.681,0.872,0.656}\textcolor{black}{2.5} & \cellcolor[rgb]{0.672,0.868,0.647}\textcolor{black}{2.5} & \cellcolor[rgb]{0.676,0.870,0.652}\textcolor{black}{2.5} & \cellcolor[rgb]{0.672,0.868,0.647}\textcolor{black}{2.5} & \cellcolor[rgb]{0.672,0.868,0.647}\textcolor{black}{2.5} & \cellcolor[rgb]{0.676,0.870,0.652}\textcolor{black}{2.5} \\
    \bottomrule
\end{tabular}
\end{table}

\begin{table}[htbp]
    \centering
    \caption{Cache sizes of the different processors used in the benchmarks.}
    \label{tab:app:cache-sizes}
    \begin{tabular}{lrrrr}
        \toprule
        Processor           & L1i Cache    & L1d Cache    & L2 Cache       & L3 Cache        \\
        \midrule
        Core i7-11700       & 32\,KiB/core & 48\,KiB/core & 512\,KiB/core  & 16\,MiB shared  \\
        Core Ultra 7 265    & 64\,KiB/core & 48\,KiB/core & 3\,MiB/core    & 30\,MiB shared  \\
        Xeon Platinum 8360Y & 32\,KiB/core & 48\,KiB/core & 1280\,KiB/core & 54\,MiB shared  \\
        Raspberry Pi 5      & 64\,KiB/core & 64\,KiB/core & 512\,KiB/core  & 2\,MiB shared   \\
        \bottomrule
    \end{tabular}
\end{table}

\clearpage
\medskip
\ifusebiblatex
\printbibliography
\fi

\bigskip
\section*{Statements and Declarations}

\noindent\textbf{Conflict of interest} The authors declare that they have no conflict of interest. \\

\noindent\textbf{Acknowledgements}
This work was supported by Fonds Wetenschappelijk Onderzoek (FWO) PhD grant
11M9523N; FWO projects G081222N, G033822N; and KU Leuven internal
funding C14/24/103. \\

\noindent\textbf{Code availability} Source code for the \textsc{Cyqlone} and \textsc{CyQPALM} solvers
(including the numerical experiments in this article) is available at \url{https://github.com/kul-optec/cyqlone}. \\

\noindent\textbf{Data availability} The raw benchmark data files used to generate \Cref{fig:box-time-spring-mass-wang-boyd-2008,fig:box-time-spring-mass-wang-boyd-2008-M,fig:scaling-abs-time_per_iter-spring-mass-wang-boyd-2008-scaling}, \Cref{tab:speedup-time-cold,tab:speedup-time-warm,tab:speedup-time-per-iter-cold,tab:app:speedup-time-cold,tab:app:speedup-time-warm,tab:app:speedup-time-per-iter-cold} are available at \url{https://doi.org/10.5281/zenodo.18879648}.

\end{document}